 \documentclass[12pt,titlepage]{book}
 \usepackage{emlines}
 \usepackage{emlines2}
 \usepackage{makeidx} 

\makeindex 

\newtheorem{teo}{Teorema}[chapter]
\newtheorem{lm}{Lema}[chapter]
\newtheorem{obs}{Observa\c{c}\~{a}o}[chapter]

\def\em{\it}

\makeindex

\begin{document}

\bibliographystyle{alfa} 
\pagenumbering{roman}
\begin{titlepage}

\begin{center}
\Huge\bf
\mbox{} \\ \mbox{} \\ 
Esparsidade, Estrutura,\\
Estabilidade e Escalamento\\ 
em \'{A}lgebra Linear Computacional\\
\Large
\vspace{3.5cm}
IX Escola de Computa\c{c}\~{a}o\\ 
24 a 31 de Julho de 1994, Recife.\\ 
\vspace{3.5cm}
\large 
Julio M. Stern\\
Departamento de Ci\^{e}ncia de Computa\c{c}\~{a}o do\\
Instituto de Matem\'{a}tica e Estat\'{\i}stica da\\
Universidade de S\~{a}o Paulo\\
\normalsize 
\vfill
\end{center}
 
\end{titlepage}
 
\tableofcontents

\pagenumbering{arabic}

\setcounter{page}{5} 
%
 
\addcontentsline{toc}{chapter}{Pref\'{a}cio para a IX Escola de 
Computa\c{c}\~{a}o} 

\chapter*{Pref\'{a}cio para a IX Escola de 
Computa\c{c}\~{a}o} 

\markboth{Pref\'{a}cio}{Pref\'{a}cio}

Propomos, para a IX Escola de Computa\c{c}\~{a}o, o curso acompanhado de
livro texto: {\em \bf Esparsidade, Estrutura, Escalamento e
Estabilidade em \'{A}lgebra Linear Computacional}.  Embora 
abrangente o suficiente para a compreens\~{a}o dos
principais problemas da \'{a}rea, a t\^{o}nica do curso \'{e} a
solu\c{c}\~{a}o de sistemas lineares esparsos e-ou estruturados.  Na
primeira se\c{c}\~{a}o da introdu\c{c}\~{a}o damos uma vis\~{a}o
panor\^{a}mica sobre o tema e a organiza\c{c}\~{a}o do livro. 
Neste pref\'{a}cio ressaltamos algumas motiva\c{c}\~{o}es 
para a inclus\~{a}o do curso na IX Escola de
Computa\c{c}\~{a}o, e esclarecemos a forma e o conte\'{u}do das
palestras a serem dadas durante o mesmo. 

\section{Import\^{a}ncia da \'{A}rea}

Grande parte do processamento de dados em ci\^{e}ncia envolve
computa\c{c}\~{a}o num\'{e}rica, e a maior parte desta
computa\c{c}\~{a}o num\'{e}rica s\~{a}o rotinas b\'{a}sicas de
\'{a}lgebra linear.  V\'{a}rias aplica\c{c}\~{o}es em engenharia,
f\'{\i}sica, pesquisa operacional, administra\c{c}\~{a}o ou
matem\'{a}tica, geram problemas lineares cujas matrizes s\~{a}o esparsas
(entre tantos outros: resolu\c{c}\~{a}o de equa\c{c}\~{o}es
diferenciais, an\'{a}lise de sistemas, circuitos ou estruturas,
programa\c{c}\~{a}o linear e n\~{a}o linear, otimiza\c{c}\~{a}o de
fluxos em redes, etc...).  Ademais, observa-se que a densidade destas
matrizes decresce com a dimens\~{a}o do problema: tipicamente apenas
$n^{\frac{3}{2}}$ ou $n\log(n)$ de seus $n^2$ elementos s\~{a}o n\~{a}o
nulos.  Assim, quanto maiores e mais complexas os problemas (e
usu\'{a}rios sempre querem resolver modelos maiores), mais importante se
torna usar eficientemente a esparsidade e a estrutura das matrizes
envolvidas.  Por exemplo, na \'{a}rea de otimiza\c{c}\~{a}o, problemas
hoje considerados de grande porte envolvem milh\~{o}es de
equa\c{c}\~{o}es, sendo a maioria destes problemas altamente
estruturados, muito esparsos, e usando mais de 99\% do tempo de
resolu\c{c}\~{a}o em rotinas b\'{a}sicas de \'{a}lgebra linear!

\section{Interdisciplinaridade}

Um aspecto que, a nosso ver, torna o tema interessante para um evento
como a { Escola de Computa\c{c}\~{a}o} \'{e} sua
interdisciplinaridade, vejamos: 
 \begin{itemize}

\item A resolu\c{c}\~{a}o dos problema de \'{a}lgebra  linear envolve os
aspectos cl\'{a}ssicos de complexidade, converg\^{e}ncia e estabilidade
de an\'{a}lise num\'{e}rica. 

\item O tratamento de esparsidade e estrutura \'{e} um problema
essencialmente combinat\'{o}rio, envolvendo teoria de grafos,
hipergrafos, e heur\'{\i}sticas para solu\c{c}\~{a}o de problemas de
programa\c{c}\~{a}o matem\'{a}tica discreta. 

\item A paraleliza\c{c}\~{a}o destes algoritmos, ou o desenvolvimento de
novos m\'{e}todos, em  ambientes t\~{a}o diversos como m\'{a}quinas de
mem\'{o}ria compartilhada ou redes de esta\c{c}\~{o}es de trabalho, vem
liderando as pesquisas na \'{a}rea nos \'{u}ltimos anos.  \end{itemize}

\section{Serventia do Livro Texto} 

O curso que se segue foi recentemente montado como a disciplina MAC-795,
M\'{e}todos Computa-cio-nais da \'{A}lgebra Linear, no programa de
mestrado do Departamento de Ci\^{e}ncia da Computa\c{c}\~{a}o da
Universidade de S\~{a}o Paulo.  Cursos de m\'{e}todos computacionais da
\'{a}lgebra linear tem se popularizado nos \'{u}ltimos anos em muitos
departamentos de computa\c{c}\~{a}o, engenharia e matem\'{a}tica
aplicada.  Era minha inten\c{c}\~{a}o escrever mais um cap\'{\i}tulo
sobre m\'{e}todos iterativos mas, ao perceber que meu rascunho duplicava
o tamanho do presente livro, resolvi postergar a tarefa para outra
ocasi\~{a}o. 

\section{Plano de Aulas}

Em cinco palestras de 90 minutos \'{e} poss\'{\i}vel dar um
bom panorama da \'{a}rea, sua import\^{a}ncia, problemas, m\'{e}todos, e
perspectivas.  O material das palestras foi distribu\'{\i}do da
seguinte forma: \begin{enumerate}

\item {\bf Introdu\c{c}\~{a}o:} Vis\~{a}o geral da \'{a}rea, sua
import\^{a}ncia, origem de problemas de grande porte, e aspectos
essenciais para a sua solu\c{c}\~{a}o. 

\item {\bf Esparsidade, Elimina\c{c}\~{a}o Assim\'{e}trica:}
Minimiza\c{c}\~{a}o de preenchimento local, outras heur\'{\i}sticas,
atualiza\c{c}\~{o}es de base. 

\item {\bf Esparsidade, Elimina\c{c}\~{a}o Sim\'{e}trica:} \'{A}rvores
de elimina\c{c}\~{a}o, heur\'{\i}sticas de ordena\c{c}\~{a}o,
heur\'{\i}sticas de dissec\c{c}\~{a}o. 

\item {\bf Estrutura:} Esparsidade macrosc\'{o}pica, m\'{e}todos de
blocos, heur\'{\i}sticas de redu\c{c}\~{a}o a formas blocadas. 

\item {\bf Paralelismo:} Uso de esparsidade $\times$ uso de estrutura,
granularidade, e potencialidades de ambientes com diferentes
coeficientes entre velocidade de processamento e velocidade de
comunica\c{c}\~{a}o.  \end{enumerate}

\section{Coment\'{a}rio sobre a Bibliografia de Suporte}

Para a parte cl\'{a}ssica de an\'{a}lise num\'{e}rica h\'{a} uma farta
variedade de livros texto, como: [Stewart-73] e [Golub-83]. 
 Existem alguns livros ou colet\^{a}neas sobre matrizes esparsas em
ambiente seq\"{u}encial, sendo os principais: [Rose-72], [Tewarson-73],
[Bunch-76], [Duff-79], [George-81] e tamb\'{e}m [Pissan-84] e [Duff-86]. 
Com exce\c{c}\~{a}o dos dois \'{u}ltimos, estes livros j\'{a} est\~{a}o
bastante desatualizados e esgotados.

\'{A}lgebra Linear Computacional, sob o enfoque de Computa\c{c}\~{a}o
Paralela, \'{e} um campo de intensa pesquisa atual:
Exclusivamente para matrizes densas h\'{a} varias obras publicadas,
x1como: [Bertsekas-89] e [Dongarra-91].  Colet\^{a}neas de artigos sobre
\'{a}lgebra linear computacional em ambiente paralelo, trazendo alguns
artigos sobre matrizes esparsas, s\~{a}o muito poucos; como: [Carey-89],
[Vorst-89] e [Gallivan-90].  Nestas resenhas encontra-se tamb\'{e}m uma
extensiva bibliografia comentada da \'{a}rea.

\section{Agradecimentos} 

No DCC-IME-USP, Departamento de Ci\^{e}ncia da Computa\c{c}\~{a}o da
Instituto de Matem\'{a}tica e Estat\'{\i}stica da Universidade de
S\~{a}o Paulo, contei sempre com a ajuda e o encorajamento de muitos
colegas, como os Professores Marcos D. Gubitoso, Arnaldo Mandel, Kunio
Okuda, Siang W. Song e Routo Terada.  Sou especialmente grato ao
coordenador do grupo de Programa\c{c}\~{a}o Matem\'{a}tica,
Professor Carlos Humes Jr. 
  No CEMCAP-IME-USP, Centro de Matem\'{a}tica e Computa\c{c}\~{a}o
Aplicadas, tive sempre o apoio dos Professores Marco Antonio Raupp, Pedro
A. Morettin e Carlos A. B. Pereira.  No NOPEF-USP, N\'{u}cleo de
Otimiza\c{c}\~{a}o e Processos Estoc\'{a}sticos Aplicados \`{a} Economia
e Finan\c{c}as, contei com o companheirismo dos Professores Marcos
Eug\^{e}nio da Silva, Jos\'{e} Carlos Santos e Jo\~{a}o Carlos Prandini. 
Em v\'{a}rias oportunidades recebi a colabora\c{c}\~{a}o do 
Eng. F\'{a}bio Nakano e do Prof. Jacob Zimbarg sobrinho.

 Partes deste livro baseiam-se em trabalhos feitos com o Professor 
Stephen A. Vavasis, da Universidade de Cornell.  A
apresenta\c{c}\~{a}o de alguns t\'{o}picos foi inspirada em disciplinas
ministradas pelo Professor Thomas F. Coleman, na mesma Universidade. 
 Em 1992 montamos a disciplina MAC-795, M\'{e}todos Computacionais da
\'{A}lgebra Linear, para o programa de mestrado do DCC-IME-USP.  A
apostila que acompanhou o curso serviu de base para este livro.  Devo a
Paulo R\'{e}gis Zanj\'{a}como, presentemente na Universidade de Cornell,
e a Professora Celma O. Ribeiro, da Escola Polit\'{e}cnica da USP,
v\'{a}rias cr\'{\i}ticas, coment\'{a}rios, e sugest\~{o}es que em muito
melhoraram aquela primeira vers\~{a}o deste texto. 

  Dos promotores da IX Escola de Computa\c{c}\~{a}o, inclusive do
an\^{o}nimo, cr\'{\i}tico e bem humorado referee, e da Universidade
Federal de Pernambuco (UFP-Recife), recebi todo o necess\'{a}rio
suporte. 

A todos estes amigos, \`{a} minha esposa, Marisa, e a meus filhos, 
Rafael, Ana Carolina e Deborah, minha gratid\~{a}o. 
   
 \clearpage
 \clearpage 

\chapter{INTRODU\c{C}\~{A}O}

\section{Panorama do Livro e seu Tema}

O objetivo deste curso \'{e} o estudo de m\'{e}todos computacionais para
a resolu\c{c}\~{a}o de sistemas lineares, $Ax=b$.  Enfatizaremos quatro
aspectos da constru\c{c}\~{a}o de bons m\'{e}todos computacionais:

\begin{itemize}
\item
{\bf Esparsidade}: Como resolver eficientemente sistemas cujas matrizes 
tenham baixa densidade de elementos n\~{a}o nulos. 
\item
{\bf Estrutura}: Como resolver eficientemente sistemas esparsos cujos
elementos n\~{a}o nulos est\~{a}o dispostos com uma certa regularidade
na matriz de coeficientes. 
 \item
 {\bf Escalamento}: Como minimizar erros de arredondamento gerados ao
operarmos com n\'{u}meros de ordem de grandeza muito diferente, numa
representa\c{c}\~{a}o de precis\~{a}o limitada. 
 \item
 {\bf Estabilidade}: Como lidar com o efeito cumulativo dos erros de
arredondamento na solu\c{c}\~{a}o final gerada pelo algoritmo. 
 \end{itemize}
 \index{Esparsidade} \index{Estrutura} 
 \index{Escalamento} \index{Estabilidade} 

As ferramentas relevantes para a an\'{a}lise e implementa\c{c}\~{a}o de
algoritmos eficientes de \'{a}lgebra linear est\~{a}o espalhados entre
diversas \'{a}reas, \'{a}reas estas afins-mas-nem-tanto, como Teoria de
Grafos e Hipergrafos, \'{A}lgebra Linear, An\'{a}lise Num\'{e}rica, e
Teoria de Processamento Paralelo.  O prop\'{o}sito destas notas \'{e}
apresentar este material de forma razoavelmente coerente e did\'{a}tica. 

O Cap\'{\i}tulo 2 exp\~{o}e o m\'{e}todo de Gauss para solu\c{c}\~{a}o
de sistemas lineares, invers\~{a}o, ou fatora\c{c}\~{a}o de matrizes. 
 Os m\'{e}todos para matrizes esparsas lidam com a estrutura da
disposi\c{c}\~{a}o dos elementos n\~{a}o nulos dentro da matriz, e esta
estrutura \'{e} convenientemente descrita e manipulada em termos da
teoria de grafos.  O Cap\'{\i}tulo 3 \'{e} um resumo de conceitos 
b\'{a}sicos desta teoria. Conceitos mais avan\c{c}ados (ou menos 
usuais) como grafos cordais, separadores, e hipergrafos, s\~{a}o 
tratados em outros cap\'{\i}tulos.   

No Cap\'{\i}tulo 4 estudamos t\'{e}cnicas para tratar sistemas esparsos
assim\'{e}tricos, principalmente preenchimento local, e a
heur\'{\i}stica P3.  Desenvolvimentos posteriores deste tema
encontram-se tamb\'{e}m no cap\'{\i}tulo 7. 
 O Cap\'{\i}tulo 5 exp\~oem as Fatora\c{c}\~{o}es QR e de Cholesky, e
sua utiliza\c{c}\~{a}o na solu\c{c}\~{a}o de problemas de quadrados 
m\'{\i}nimos, programa\c{c}\~{a}o quadr\'{a}tica, e constru\c{c}\~{a}o de
projetores.  No Cap\'{\i}tulo 6 estudamos t\'{e}cnicas para tratar
sistemas esparsos sim\'{e}tricos; incluindo toda a necess\'{a}ria teoria
de grafos cordais. 

O capitulo 7 trata da estrutura de um sistema, que pode ser vista
como regularidades no padr\~{a}o de esparsidade, ou como a decomposi\c{c}\~{a}o 
do sistema em sub-sistemas acoplados. Neste cap\'{\i}tulo estudamos
brevemente a paraleliza\c{c}\~{a}o de alguns dos algoritmos, tema 
que j\'{a} aparece implicitamente em cap\'{\i}tulos anteriores. 
 O Capitulo 8 \'{e} dedicado \`{a} an\'{a}lise do efeito dos
inevit\'{a}veis erros de arredondamento nos procedimentos
computacionais.  No Capitulo 9 analisamos quanto estes erros podem
degradar a solu\c{c}\~{a}o final de um problema.  Estes dois
cap\'{\i}tulos tratam dos aspectos de An\'{a}lise Num\'{e}rica que,
embora n\~{a}o sendo a t\^{o}nica do curso, n\~{a}o podem ser
desprezados. 

O Cap\'{\i}tulo 10 trata do problema de mudan\c{c}a de base, i.\'{e}., de
atualizar a inversa de uma matriz quando nesta se substitui uma
\'{u}nica coluna.  Este problema \'{e} de fundamental import\^{a}ncia para
muitos algoritmos de Otimiza\c{c}\~{a}o. 


Parte da avalia\c{c}\~{a}o numa disciplina como a proposta deve ser
feita com exerc\'{\i}cios-programa, como os dados ao longo das notas. 
\'{E} recomend\'{a}vel o uso de uma linguagem estruturada, que inclua
entre seus tipos b\'{a}sicos, n\'{u}meros reais de precis\~{a}o simples
e dupla, inteiros, campos de bits e ponteiros, que permita a f\'{a}cil
constru\c{c}\~{a}o e manipula\c{c}\~{a}o de estruturas, e para a qual
haja compiladores em computadores dos mais diversos portes.  As
linguagens C, C++ e FORTRAN-90 s\~{a}o sugest\~{o}es naturais. 
 Ambientes iterativos para c\'{a}lculo matricial, como Matlab,
s\~{a}o um grande est\'{\i}mulo \`{a} experimenta\c{c}\~{a}o e 
compara\c{c}\~{a}o de diversos de m\'{e}todos num\'{e}ricos,  
facilitando a r\'{a}pida prototipagem e teste de algoritmos. 
Apresentamos uma introdu\c{c}\~{a}o a este tipo de ambiente como 
ap\^{e}ndice.   

\section{Nota\c{c}\~{o}es}

Utilizaremos letras mai\'{u}sculas, $A,B,\ldots$, para denotar matrizes, 
e letras min\'{u}sculas, $a,b,\ldots$, para vetores.  Numa matriz ou
vetor, os \'{\i}ndices de linha ou coluna ser\~{a}o, respectivamente,
subscritos ou superscritos \`{a} direita.  Assim: $A_i^j$ \'{e} o
elemento na $i$-\'{e}sima linha e na $j$-\'{e}sima coluna da matriz $A$;
$b_i$ \'{e} o elemento na $i$-\'{e}sima linha do vetor-coluna $b$; e
$c^j$ \'{e} o elemento na $j$-\'{e}sima coluna do vetor-linha $c$. 
 \index{Matriz!nota\c{c}\~{a}o} 
 \index{Matriz!\'{\i}ndices} 
 
\'{I}ndices ou sinais \`{a} esquerda, superscritos ou subscritos,
identificam vetores ou matrizes distintas.  Assim:$A$, $'A$, $^{j}A$,
$_{i}A$, $_{i}^{j}A$, $\ldots$ s\~{a}o matrizes distintas.  Tamb\'{e}m a
letra $\delta$ forma, \`{a} esquerda de uma letra latina o nome de um
vetor ou matriz, como $\delta a$ ou $\delta A$.  As demais letras gregas
usaremos geralmente para escalares ou fun\c{c}\~{o}es. 

Freq\"{u}entemente nos referimos aos blocos componentes de vetores ou
matrizes, como por exemplo: se $b$ e $c$ s\~{a}o vetores linha 
$a = \left[ \begin{array}{cc} b & c \end{array} \right]$ \'{e} um vetor 
linha cujos primeiros elementos s\~{a}o os elementos de $b$, e os
elementos seguintes s\~{a}o os elementos de $c$. 
 \index{Matriz!blocada} 

Analogamente, 
$$ A = 
\left[ \begin{array}{cccc}
_1^1A & _1^2A & _1^3A \\ _2^1A & _2^2A & _2^3A \\ 
_3^1A & _3^2A & _3^3A 
\end{array} \right] $$
representa uma matriz blocada, desde que as dimens\~{o}es dos blocos
sejam compat\'{\i}veis. 

Se $A$ \'{e} uma matriz, $A_i$ representa a $i$-\'{e}sima linha de $A$,
e $A^j$ sua $j$-\'{e}sima coluna, de modo que se $A$ \'{e} $mxn$,
$$
A = \left[ \begin{array}{cccc} 
A^1 & A^2 & \ldots & A^n \end{array} \right]  =
\left[ \begin{array}{c} A_1 \\ A_2 \\ \vdots \\ A_m 
\end{array} \right] 
$$

Uma {\bf matriz diagonal} ser\'{a} denotada como  
  \index{Matriz!diagonal} 
$$ diag(d) = D = \left[ \begin{array}{cccc}
d_1 & & & \\ & d_2 & & \\ & & \ddots & \\ & & & d_n 
\end{array} \right] \ , \ \   
d = \left[ \begin{array}{c} 
d_1 \\ d_2 \\ \vdots \\ d_n 
\end{array} \right]   
$$ 

O {\bf vetor unit\'{a}rio}, {\bf 1}, \'{e} o vetor, linha ou coluna, em que 
todos os elementos s\~{a}o iguais a 1, isto \'{e}, 
${\bf 1} = \left[ \begin{array}{cccc} 1 & 1 & \ldots & 1 
\end{array} \right]$.
 Assim a matriz identidade \'{e} $I = diag({\bf 1})$.
O $j$-\'{e}simo {\bf versor} de dimens\~{a}o $n$ \'{e} $I^j$, 
a $j$-\'{e}sima coluna da matriz identidade de dimens\~{a}o $n$.   
 A transposta, a inversa e a transposta da inversa de uma matriz A
s\~{a}o denotadas, respectivamente, por $A'$, $A^{-1}$ e $A^{-t}$. 
 \index{Vetor!unit\'{a}rio} 
 \index{Vetor!versor} 
 \index{Matriz!identidade} 
 \index{Matriz!transposta} 

O {\bf determinante} de uma matriz $A$, $nxn$, \'{e}
$$ det(A) = \sum_p S(p) A^1_{p(1)}A^2_{p(2)}\ldots A^n_{p(n)} $$
onde 
$p = \left[ \begin{array}{cccc} p^1 & p^2 & \ldots & p^n 
\end{array} \right]$
\'{e} uma {\bf permuta\c{c}\~{a}o} dos elementos de
$^0p = \left[ \begin{array}{cccc} 1 & 2 & \ldots & n 
\end{array} \right]$.  
O sinal de permuta\c{c}\~{a}o, $S(p)$, \'{e} +1 ou -1 conforme
o n\'{u}mero de trocas de pares de elementos que \'{e} necess\'{a}rio 
fazer em $p$ para retornar a $^0p$, seja par ou \'{\i}mpar.
 O conjunto dos elementos em cada um dos termos na somat\'{o}ria da
defini\c{c}\~{a}o do determinante \'{e} denominado uma {\bf diagonal} da
matriz.  A diagonal correspondente a pemuta\c{c}\~{a}o $^0p$ \'{e}
denominada {\bf diagonal principal}. 
 Uma matriz quadrada \'{e} singular se tiver determinante nulo.  O {\bf
posto} de uma matriz $A$, $rank(A)$, \'{e} a dimens\~{a}o de sua maior
sub-matriz quadrada n\~{a}o singular. 
 \index{Matriz!determinante} 
 \index{Permuta\c{c}\~{a}o} 
 \index{Matriz!diagonal em} 
 \index{Matriz!posto} 

Dado um sistema $Ax = b$, com matriz de coeficientes $A \ nxn$ n\~{a}o
singular, a {\bf Regra Cramer} nos diz que o sistema tem por \'{u}nica
solu\c{c}\~{a}o 
 \index{Regra de Cramer}  
 $$ x\ \mid \ x_i = det(\left[ \begin{array}{ccccccc} 
    A^1 & \ldots & A^{i-1} & b & A^{i+1} & A^n 
    \end{array} \right] ) / det(A). 
 $$

O {\bf tra\c{c}o} de uma matriz $A,\ nxn$, \'{e} 
a soma de seus elementos na diagonal principal:  
  \index{Matriz!tra\c{c}o} 
$$ tr(A) = \sum_{i=1}^n A_i^i \ \ .$$

A {\bf matriz booleana} associada \`{a} matriz $A,\ B(A)$, \'{e} a matriz, 
da mesma dimens\~{a}o de $A$, em que $B(A)_i^j = 1$ se $A_i^j\neq 0$, e
$B(A)_i^j =0$ se $A_i^j =0$.  
 O complemento de uma matriz booleana $B$, \'{e} a matriz booleana $\bar
B$, da mesma dimens\~{a}o de $B$, tal que
${\bar B}_i^j = 1 \Leftrightarrow B_i^j = 0$. 
 \index{Matriz!booleana} 

Os conjuntos $N = \{1,2,...,n\}$ e $M = \{1,2,...,m\}$ ser\~{a}o
freq\"{u}entemente usados como {\bf dom\'{\i}nios de \'{\i}ndices}. 
Assim, se $A$ \'{e} uma matriz $mxn$, faz sentido falar dos elementos
$A_i^j ,\ i\in M,\ j\in N$, 
  \index{Matriz!\'{\i}ndices} 

O n\'{u}mero de {\bf elementos n\~{a}o nulos}, ENNs, numa matriz A, 
mxn, \'{e} 
 \index{Matriz!elementos n\~{a}o nulos} 
 $$enn(A) = \sum_{i,j=1}^{m,n} B(A)_i^j = {\bf 1}'B(A){\bf 1}$$ 
 de modo que $enn(A_i)$ e $enn(A^j)$ s\~{a}o, respectivamente, o
n\'{u}mero de elementos n\~{a}o nulos na $i$-\'{e}sima linha e na
$j$-\'{e}sima coluna de $A$. 

Dada uma fun\c{c}\~{a}o $\varphi ( )$, definida num dom\'{\i}nio $D$, e
$X \subset D$, definimos seu {\bf argumento m\'{\i}nimo} 
 \index{Argmax, Argmin} 
 $V=arg\min_{x\in X} \varphi (x)$ e {\bf argumento m\'{a}ximo}
$U=arg\max_{x\in X}\varphi (x)$ como, respectivamente, os conjuntos 
$V,\ U \subset D$ que minimizam ou maximizam a fun\c{c}\~{a}o $\varphi$ 
em $X$.  Assim, se por exemplo, 
$\varphi = x^2 , \ x\in \Re$, $X=[-5,5]$ e $Y=]-5,5[$,
ent\~{a}o $arg\min_X \varphi (x) = arg\min_Y \varphi (x) = \{ 0 \}$,
$arg\max_X \varphi (x) = \{-5,5\}$, 
$arg\max_Y \varphi (x) = \emptyset$.

Para realizar experi\^{e}ncias computacionais utilizaremos por vezes os
sistemas lineares de dimens\~{a}o $n$ e solu\c{c}\~{a}o $x={\bf 1}$, com
as seguintes matrizes de coeficientes e vetores independentes:
 \index{Matriz!binomial} 
 \index{Matriz!Hilbert} 
 \index{Matriz!tridiagonal} 
\begin{itemize}
\item {\bf Binomial}:
$$ 
^nB_i^j = \left\{ \begin{array}{ccc}
\left( \begin{array}{c} i \\ j-1 \end{array} \right) & 
 \Leftrightarrow & j= 1,2,\ldots i \\ 
0 & \Leftrightarrow & j>i \end{array} \right. \ \ ,\ \  
^nb = \left[ \begin{array}{c} 2^1 \\ 2^2 \\ \vdots \\ 2^n    
\end{array} \right] 
$$
\item {\bf Hilbert}:
$$
^nH_i^j = 1/(i+j-1) \ \ ,\ \   
^nh_i = \sum_{j=1}^n 1/(i+j-1) 
$$
\item {\bf Tridiagonal}:
$$
^nT_i^j = \left\{ \begin{array}{ccc}
-2 & \Leftrightarrow & i=j \\
-1 & \Leftrightarrow & |i-j|=1 \\
0 & & \mbox{caso contrario}
\end{array} \right. \ \ ,\ \ 
^nt = \left\{ \begin{array}{ccc}
1 & \Leftrightarrow & i\in \{1,n\} \\
0 & & \mbox{caso contrario}
\end{array} \right. 
$$ 
\end{itemize}

\section{Representa\c{c}\~{a}o de Matrizes Esparsas}

Ao representarmos uma matriz esparsa gostar\'{\i}amos de utilizar o
m\'{\i}nimo de mem\'{o}ria, idealmente apenas a posi\c{c}\~{a}o e o
valor dos ENNs, mas ao mesmo tempo ter a maior flexibilidade
poss\'{\i}vel para acessar, modificar, inserir ou remover um elemento
qualquer.  Estes objetivos s\~{a}o algo conflitantes, o que motiva o uso
de diversas representa\c{c}\~{o}es: 
  \index{Matriz!representa\c{c}\~{a}o} 

Representa\c{c}\~{a}o {\bf est\'{a}tica por linhas}: Para uma matriz
$A,\ m\times n$ com $enn(A)=l$, s\~{a}o utilizados um vetor real,
$aias(k)\ k\in L$, um vetor inteiro $aijs(k),\ k\in L=\{1,\ldots l\}$, e
um segundo vetor de inteiros, $aif(i),\ i\in M$.  Os vetores aias e aijs
listam os valores e os \'{\i}ndices de coluna dos ENN's, linha por
linha.  $aif(i)$ nos d\'{a} o fim, ou a posi\c{c}\~{a}o do \'{u}ltimo
elemento da linha $i$ em $aias$ e $aijs$. 

Representa\c{c}\~{a}o {\bf est\'{a}tica por colunas}: Para uma matriz
$A,\ m\times n$ com $enn(A)=l$, s\~{a}o utilizados um vetor real,
$ajas(k)\ k\in L$, um vetor inteiro $ajis(k),\ k\in L$, e um segundo
vetor de inteiros, $ajf(j),\ j\in N$.  Os vetores $ajas$ e $ajis$ listam
os valores e os \'{\i}ndices de linha dos ENNs, coluna por coluna. 
$ajf(i)$ nos d\'{a} a posi\c{c}\~{a}o do \'{u}ltimo elemento da coluna
$j$ em $ajas$ e $ajis$.

Representa\c{c}\~{a}o de {\bf lista ligada por linhas}: Cada ENN
corresponde a uma c\'{e}lula contendo: O valor do ENN, o \'{\i}ndice de
coluna e um ponteiro para a c\'{e}lula do pr\'{o}ximo ENN na linha.  As
c\'{e}lulas dos \'{u}ltimos ENNs em cada linha cont\'{e}m o ponteiro
nulo, e um vetor de $m$ ponteiros, $ancorai(i)$, nos d\'{a} a primeira
c\'{e}lula de cada linha. 

Representa\c{c}\~{a}o de {\bf lista ligada por colunas}: Cada ENN
corresponde a uma c\'{e}lula contendo: O valor do ENN, o \'{\i}ndice de
linha e um ponteiro para a c\'{e}lula do pr\'{o}ximo ENN na coluna.  As
c\'{e}lulas dos \'{u}ltimos ENNs em cada coluna cont\'{e}m o ponteiro
nulo, e um vetor de $n$ ponteiros, $ancoraj(j)$, nos d\'{a} a primeira
c\'{e}lula de cada coluna.

Representa\c{c}\~{a}o de {\bf rede}: Cada ENN corresponde a uma
c\'{e}lula contendo o valor do ENN $A_i^j$, os \'{\i}ndices $i$ e $j$, e
ponteiros para a c\'{e}lula do pr\'{o}ximo ENN na linha e na coluna. 
Dois vetores $ancorai(i)$ e $ancoraj(j)$ cont\'{e}m ponteiros para as
primeiras c\'{e}lulas em cada linha e coluna. 

Representa\c{c}\~{a}o de {\bf rede dupla}: An\'{a}loga a
representa\c{c}\~{a}o de rede, mas cada c\'{e}lula cont\'{e}m al\'{e}m
de ponteiros para as pr\'{o}ximas c\'{e}lulas ao {\bf sul} (pr\'{o}ximo
ENN na coluna) e ao {\bf leste} (pr\'{o}ximo ENN na linha), tamb\'{e}m
ponteiros para as c\'{e}lulas ao {\bf oeste} e ao {\bf norte}, i.e. 
para os ENN anteriores na linha e na coluna.

No exemplo 1 temos as diversas 
representa\c{c}\~{o}es da matriz:

 $$A = \left[ \begin{array}{ccccc} 
 & 102 & & 104 & \\  201 & & & & \\ 301 & 304 & & & \\
 & & & & 405 \\ 501 & 502 & 503 & & \\  
 \end{array} \right] \ \ \ , \ \ l=enn(A)=9\ .   
 $$

 Representa\c{c}\~{a}o est\'{a}tica por linhas: 
 $$ aias = \left[ \begin{array}{ccccccccc} 
    102 & 104 & 201 & 301 & 304 & 405 & 501 & 502 & 503 
    \end{array} \right] \ , $$  
 $$ aijs = \left[ \begin{array}{ccccccccc} 
    2 & 4 & 1 & 1 & 4 & 5 & 1 & 2 & 3 
    \end{array} \right] \  , $$  
 $$ aif = \left[ \begin{array}{ccccc} 
    2 & 3 & 5 & 6 & 9  
    \end{array} \right] \ . $$ 
  
 Representa\c{c}\~{a}o est\'{a}tica por colunas: 
 $$ ajas = \left[ \begin{array}{ccccccccc} 
    201 & 301 & 501 & 102 & 502 & 503 & 104 & 304 & 405 
    \end{array} \right] \ , $$  
 $$ ajis = \left[ \begin{array}{ccccccccc} 
    2 & 3 & 5 & 1 & 5 & 5 & 1 & 3 & 4 
    \end{array} \right] \  , $$  
 $$ ajf = \left[ \begin{array}{ccccc} 
    3 & 5 & 6 & 8 & 9  
    \end{array} \right] \ . $$

Representa\c{c}\~{a}o de lista ligada por linhas: 

$$ \begin{array}{ccccc} 
   \mbox{ancorai(1)}\rightarrow & 
   \left( \mbox{2, 102, }\rightarrow \right) & 
   \left( \mbox{4, 104, }\dashv \right) \\  
   \mbox{ancorai(2)}\rightarrow & 
   \left( \mbox{1, 201, }\dashv \right) \\ 
   \mbox{ancorai(3)}\rightarrow & 
   \left( \mbox{1, 301, }\rightarrow \right) & 
   \left( \mbox{4, 304, }\dashv \right) \\ 
   \mbox{ancorai(4)}\rightarrow & 
   \left( \mbox{5, 405, }\dashv \right) \\  
   \mbox{ancorai(5)}\rightarrow & 
   \left( \mbox{1, 501, }\rightarrow \right) & 
   \left( \mbox{2, 502, }\rightarrow \right) & 
   \left( \mbox{3, 503, }\dashv \right)   
   \end{array} $$

Representa\c{c}\~{a}o de lista ligada por colunas: 

$$ \begin{array}{ccccc} 
   \mbox{ancoraj(1)} & \mbox{ancoraj(2)} & \mbox{ancoraj(3)} &  
   \mbox{ancoraj(4)} & \mbox{ancoraj(5)} \\ 
   \downarrow & \downarrow & \downarrow & \downarrow & \downarrow \\ 
   \left( \mbox{2, 201, }\downarrow \right) & 
   \left( \mbox{1, 102, }\downarrow \right) & 
   \left( \mbox{5, 503, }\perp \right) & 
   \left( \mbox{1, 104, }\downarrow \right) & 
   \left( \mbox{4, 405, }\perp \right) \\  
   \left( \mbox{3, 301, }\downarrow \right) & 
   \left( \mbox{5, 502, }\perp \right) &  & 
   \left( \mbox{3, 304, }\perp \right) &  \\  
   \left( \mbox{5, 501, }\perp \right) & & & & \\   
   \end{array} $$  

 Representa\c{c}\~{a}o em rede: 
 { 
 $$ \begin{array}{cccccc} 
    & \mbox{ancoraj(1)} & \mbox{ancoraj(2)} & 
   \mbox{ancoraj(3)} & \mbox{ancoraj(4)} & \mbox{ancoraj(5)} \\ 
    & \downarrow  & \downarrow  & 
   \downarrow  & \downarrow  & \downarrow \\ 
   \mbox{ancorai(1)}\rightarrow & & 
   \left( \begin{array}{ccc} 1 & 2 & 102 \\ 
     & \downarrow & \rightarrow \end{array} \right) & & 
   \left( \begin{array}{ccc} 1 & 4 & 104 \\ 
     & \downarrow & \rightarrow \end{array} \right) & \\  
   \mbox{ancorai(2)}\rightarrow & 
   \left( \begin{array}{ccc} 2 & 1 & 201 \\ 
     & \downarrow & \dashv \end{array} \right) & & & & \\ 
   \mbox{ancorai(3)}\rightarrow & 
   \left( \begin{array}{ccc} 3 & 1 & 301 \\ 
     & \downarrow & \rightarrow \end{array} \right) & & &  
   \left( \begin{array}{ccc} 3 & 4 & 304 \\ 
     & \perp & \dashv \end{array} \right) & \\  
   \mbox{ancorai(4)}\rightarrow & & & & & 
   \left( \begin{array}{ccc} 4 & 5 & 405 \\ 
     & \perp & \dashv \end{array} \right) \\  
   \mbox{ancorai(5)}\rightarrow & 
   \left( \begin{array}{ccc} 5 & 1 & 501 \\ 
     & \perp & \rightarrow \end{array} \right) &  
   \left( \begin{array}{ccc} 5 & 2 & 502 \\ 
     & \perp & \rightarrow \end{array} \right) &  
   \left( \begin{array}{ccc} 5 & 3 & 503 \\ 
     & \perp & \dashv \end{array} \right) & &  \\  
 \end{array} $$ }

\section*{Exerc\'{\i}cios}

\begin{enumerate}

\item 
Prove que 
  \index{Matriz!transposta} 
\begin{enumerate} 
\item $(AB)' = B'A'$. 
\item $(AB)^{-1} = A^{-1}B^{-1}$. 
\item $(A')^{-1} = (A^{-1})'$. 
\end{enumerate}

\item
Na representa\c{c}\~{a}o est\'{a}tica por linhas, \'{e} realmente
necess\'{a}rio termos $aif$ al\'{e}m de $aias$ e $aijs$?

\item Considere uma matriz esparsa de estrutura regular e conhecida, por
exemplo tri-diagonal, para a qual bastaria conhecermos o valor dos ENNs
numa dada seq\"{u}\^{e}ncia, por exemplo linha por linha.  Suponha que
um inteiro ou ponteiro ocupa 2 bytes e que um real ocupa 4 bytes. 
D\^{e} o coeficiente de uso de mem\'{o}ria de cada uma das outras
representa\c{c}\~{o}es em rela\c{c}\~{a}o a est\'{a}
representa\c{c}\~{a}o minimal. 

\item 
\begin{enumerate}
 \item
 Considere uma matriz esparsa $A$ de densidade, i.e.  fra\c{c}\~{a}o de
elementos n\~{a}o nulos, ${\alpha}^2$.  Suponha que $A$ est\'{a} 
representada em rede.  Queremos adicionar um novo ENN $A_i^j$ \`{a}
matriz.  Supondo que os ENN est\~{a}o aleatoriamente distribu\'{\i}dos,
e tomando acesso a uma c\'{e}lula como opera\c{c}\~{a}o elementar, qual
a complexidade esperada desta opera\c{c}\~{a}o? Explique sucintamente
como realizar a opera\c{c}\~{a}o. 
 \item 
 Nas mesmas condi\c{c}\~{o}es, queremos
substituir duas linhas $A_i$ e $A_{i+1}$, respectivamente, pelas
combina\c{c}\~{o}es lineares ${\tilde A}_i = c*A_i + s*A_{i+1}$
e ${\tilde A}_{i+1} = c*A_{i+1} - s*A_i $.  Novamente queremos um
algoritmo, descrito sum\'{a}ria mas claramente, e a complexidade
esperada.  Dicas e curiosidades:
\begin{enumerate}
\item 
Assumindo que $\alpha << 1$, o caso em que estamos interessados, a
densidade esperada das novas linhas \'{e} $2*\alpha$.  Por qu\^{e}? 
\item
O algoritmo para a segunda quest\~{a}o deve ser algo melhor que a
mera repeti\c{c}\~{a}o do algoritmo da primeira quest\~{a}o.  
\item 
Esta transforma\c{c}\~{a}o linear, tomando as constantes $c$ e $s$
como o coseno e o seno de um \^{a}ngulo $\theta$, \'{e} uma
``rota\c{c}\~{a}o de Givens'', a ser estudada no cap\'{\i}tulo 3. 
\end{enumerate}
\end{enumerate}

\item 
Dados $u$ e $w$  $n\times 1$,  $A$ e $B$ $n\times n$, e $k\in N^*$, 
prove que
 \index{Matriz!tra\c{c}o}  
\begin{enumerate} 
 \item $tr(A+B) = tr(A) + tr(B)$.
 \item $tr(AB) = tr(BA)$.
 \item $\mbox{tr}(uw') = w'u$. 
 \item $\mbox{tr}(Auw') = \mbox{tr}(uw'A) = w'Au$.
 \item $(uw'A)^k = (\mbox{tr}(uw'A))^{k-1} \ (uw'A)$. 
 \item $(Auw')^k = (\mbox{tr}(Auw'))^{k-1} \ (Auw')$. 
\end{enumerate}  

\end{enumerate}
 
 \clearpage
 \clearpage 
\setcounter{chapter}{1} 
\chapter{FATORA\c{C}\~{A}O LU }   
\begin{center}
{\Large Solu\c{c}\~{a}o de Sistemas Lineares }
\end{center}

\section{Permuta\c{c}\~{o}es e Opera\c{c}\~{o}es Elementares}

Uma {\bf matriz de permuta\c{c}\~{a}o} \'{e} uma matriz obtida pela
permuta\c{c}\~{a}o de linhas ou colunas na matriz identidade.  Realizar,
na matriz identidade, uma dada permuta\c{c}\~{a}o de linhas, nos fornece
a correspondente matriz de permuta\c{c}\~{a}o de linhas; Analogamente,  
uma permuta\c{c}\~{a}o de colunas da identidade fornece a correspondente 
matriz de permuta\c{c}\~{a}o de colunas. 
 \index{Matriz!de permuta\c{c}\~{a}o} 

Dada uma (matriz de) permuta\c{c}\~{a}o de linhas, $P$ e uma (matriz de)
permuta\c{c}\~{a}o de colunas, $Q$, o correspondente vetor de
\'{\i}ndices de linha (coluna) permutados s\~{a}o
$$ p= (P \left[ \begin{array}{c} 1\\ 2\\ \vdots \\ m 
  \end{array} \right] )'$$
$$ q= \left[ \begin{array}{cccc} 1 & 2 & \ldots & n 
  \end{array} \right] Q $$

\begin{lm} Realizar uma permuta\c{c}\~{a}o de linhas (de colunas) numa
matriz qualquer $A$, de modo a obter a matriz permutada $\tilde A$,
equivale a multiplic\'{a}-la, \`{a} esquerda (\`{a} direita), pela
correspondente matriz de permuta\c{c}\~{a}o de linhas (de colunas). 
Ademais, se $p$ ($q$) \'{e} o correspondente vetor de \'{\i}ndices de
linha (de coluna) permutados,
 $$ {\tilde A}_i^j = (P A)_i^j = A_{p(i)}^j $$
$$ {\tilde A}_i^j = (A Q)_i^j = A_i^{q(j)} \ .$$
\end{lm}

Exemplo 1:
Dadas as matrizes  
$$
A=\left[ \begin{array}{ccc} 11 & 12 & 13 \\ 21 & 22 & 23 \\
 31 & 32 & 33 \end{array} \right] \ , \  \ 
P=\left[ \begin{array}{ccc} 0 & 0 & 1 \\ 1 & 0 & 0 \\ 0 & 1 & 0 
 \end{array} \right] \ , \  \ 
Q=\left[ \begin{array}{ccc}0 & 1 & 0 \\ 0 & 0 & 1 \\ 1 & 0 & 0 
 \end{array} \right] \ ,$$ 
$$ 
p=q=\left[ \begin{array}{ccc} 3 & 1 & 2 \end{array} \right] \ , \ \ 
PA=\left[ \begin{array}{ccc} 31 & 32 & 33 \\ 11 & 12 & 13 \\ 
 21 & 22 & 23 \end{array} \right] \ ,\ \ 
AQ=\left[ \begin{array}{ccc} 13 & 11 & 12 \\ 23 & 21 & 22 \\ 
 33 & 31 & 32  \end{array} \right] \ . $$ 

Uma matriz quadrada, A, \'{e} {\bf sim\'{e}trica} sse for igual a
transposta, isto \'{e}, sse $A=A'$.  Uma {\bf permuta\c{c}\~{a}o
sim\'{e}trica} de uma matriz quadrada $A$ \'{e} uma permuta\c{c}\~{a}o da
forma $\tilde A =PAP'$, onde $P$ \'{e} uma matriz de permuta\c{c}\~{a}o. 
 Uma matriz quadrada, $A$, \'{e} {\bf ortogonal} sse sua inversa for
igual a sua transposta, isto \'{e}, sse $A^{-1}=A'$. 
  \index{Matriz!sim\'{e}trica} 
  \index{Matriz!ortogonal} 

\begin{lm}
(a) Matrizes de permuta\c{c}\~{a}o s\~{a}o ortogonais.
(b) Uma permuta\c{c}\~{a}o sim\'{e}trica de uma matriz 
sim\'{e}trica \'{e} ainda uma matriz sim\'{e}trica. 
\end{lm}

Estudaremos a seguir alguns algoritmos para solu\c{c}\~{a}o do sistema 
$Ax=b$.  Tais algoritmos s\~{a}o denominados diretos pois, a menos do erro
de arredondamento, fornecem diretamente a solu\c{c}\~{a}o do sistema. 
 A ess\^{e}ncia destes algoritmos s\~{a}o transforma\c{c}\~{o}es sobre o
sistema que deixam inalteradas sua solu\c{c}\~{a}o. 

S\~{a}o {\bf opera\c{c}\~{o}es elementares} sobre uma matriz, 
$n\times m$, qualquer:
 \index{Opera\c{c}\~{a}o Elementar} 
\begin{enumerate} 
\item Multiplicar uma linha por um escalar n\~{a}o nulo. 
\item Adicionar, a uma linha, uma outra linha da matriz.
\item Subtrair de uma linha, uma outra linha da matriz multiplicada por 
um escalar n\~{a}o nulo. 
\end{enumerate}

Uma opera\c{c}\~{a}o elementar aplicada \`{a} matriz identidade $I$, 
$n\times n$, produz a {\bf matriz elementar} correspondente, $E$. 
  \index{Matriz!elementar} 

\begin{lm}
Realizar uma opera\c{c}\~{a}o elementar sobre uma matriz qualquer, 
$A$, equivale a multiplic\'{a}-la, \`{a} esquerda, pela matriz 
elementar correspondente. 
\end{lm} 

Exemplo 3.
$$
A=\left[ \begin{array}{ccc} 11 & 12 & 13 \\ 
  21 & 22 & 23 \\ 31 & 32 & 33 \end{array} \right] \ , \  \
E=\left[ \begin{array}{ccc} 1 & 0 & 0 \\ 0 & 1 & 0 \\ 
  -31/11 & 0 & 1 \end{array} \right] \ , $$
$$  
EA=\left[ \begin{array}{ccc} 11 & 12 & 13 \\ 21 & 22 & 23 \\ 
  0 & 32-12*31/11 & 33-13*31/11 \end{array} \right] \ . 
$$

\begin{lm}
Toda matriz elementar \'{e} invers\'{\i}vel.  Ademais,  
$Ax=b \Leftrightarrow EAx=Eb$. 
\end{lm} 

A matriz aumentada correspondente ao sistema $Ax=b$, $A\ n\times n$ 
\'{e} a matriz $[A\ b]$, $n\times n+1$.
Obviamente a matriz aumentada do sistema $EAx=Eb$ \'{e} a matriz 
$E[A\ b]$. 
  \index{Matriz!aumentada} 

 Uma matriz quadrada $A$ \'{e} dita {\bf triangular} superior se todos
os elementos abaixo da diagonal principal s\~{a}o nulos, e \'{e} dita
triangular estritamente superior se tamb\'{e}m os elementos da diagonal
s\~{a}o nulos.  Analogamente, definimos matriz triangular inferior e
 estritamente inferior.  Assim, $A$ \'{e}
  \index{Matriz!triangular} 
 \begin{itemize}
 \item Triangular Superior 
  $\Leftrightarrow (i > j \Rightarrow A_i^j =0 )$.
 \item Triangular Estritamente Superior  
  $\Leftrightarrow (i \geq j \Rightarrow A_i^j =0 )$.
 \item Triangular Inferior  
  $\Leftrightarrow (i < j \Rightarrow A_i^j =0 )$. 
 \item Triangular Estritamente Inferior  
  $\Leftrightarrow (i\leq j \Rightarrow A_i^j =0 )$.
 \end{itemize}

Os m\'{e}todos diretos que estudaremos funcionam por ser f\'{a}cil
resolver um sistema Ux = b se U \'{e} triangular superior, isto \'{e}
Em um tal sistema podemos calcular por substitui\c{c}\~{a}o, 
nesta ordem, as componentes $x_n,\ x_{n-1},\ldots , x_1$, pois
\begin{eqnarray*}
x_n & = & b_n / U_n^n \\
x_{n-k} & = & (b_{n-k} - \sum_{j=n-k+1}^n U_{n-k}^j x_j )/ U_{n-k}^{n-k}
\end{eqnarray*}
 Tal m\'{e}todo funciona se $U_j^j\neq 0,\ \forall j\in N$.  
Como $det(U)= \prod_{j=1}^n U_j^j$, esta condi\c{c}\~{a}o 
equivale a termos um sistema bem determinado. 
Estudaremos a seguir como levar, atrav\'{e}s de opera\c{c}\~{o}es
elementares, um sistema qualquer \`{a} forma triangular.

\section{M\'{e}todo de Gauss} 

O m\'{e}todo de Gauss leva o sistema \`{a} forma triangular atrav\'{e}s
de opera\c{c}\~{o}es elementares do tipo 3.  Tentaremos faz\^{e}-lo
usando, nesta ordem, a primeira linha para anular os elementos abaixo da
diagonal na primeira coluna, a segunda linha para anular os elementos
abaixo da diagonal na segunda coluna, etc...  O exemplo 2 ilustra o
processo. Indicamos tamb\'{e}m os fatores pelos quais multiplicamos
cada linha antes de subtra\'{\i}-la de outra. 

Exemplo 4.

$$ 
\left[ ^0A\ ^0b\right]=\left[ \begin{array}{rrrr} 2 & 1 & 3 & 1 \\
  2 & 3 & 6 & 2 \\ 4 & 4 & 6 & 6 \end{array} \right] 
 \begin{array}{c}  \\ 1 \\ 2 \end{array}  \ \  \rightarrow \ \      
\left[ ^1A\ ^1b\right]=\left[ \begin{array}{rrrr} 2 & 1 & 3 & 1 \\ 
  0 & 2 & 3 & 1 \\ 0 & 2 & 0 & 4 \end{array} \right]   
 \begin{array}{c}  \\  \\ 1 \end{array}  \ \ \rightarrow \ \  $$ 
$$
\left[ ^2A\ ^2b\right]=\left[ \begin{array}{rrrr} 2 & 1 & 3 & 1 \\  
  0 & 2 & 3 & 1 \\ 0 & 0 & -3 & 3 \end{array} \right] \ \  \rightarrow \ \  
x = \left[ \begin{array}{r} 1\\ 2\\ -1 \end{array} \right] $$

Em cada etapa denominaremos a linha que est\'{a} sendo multiplicada e
subtra\'{\i}da \`{a}s demais de {\bf linha piv\^{o}}, seu elemento
diagonal de {\bf elemento piv\^{o}} e os fatores de
multiplica\c{c}\~{a}o de {\bf multiplicadores}. 
 Posteriormente faremos uso dos multiplicadores utilizados no processo e,
portanto, devemos armazen\'{a}-los.  Com este fim, definimos as matrizes,
$n\times n$, $^1M$, $^2M$, \ldots $^{n-1}M=M$, onde: 
Se $j\leq k$ e $i>j$, ent\~{a}o $^kM_i^j$ \'{e} o multiplicador utilizado 
na $k$-\'{e}sima etapa para anular o elemento na $i$-\'{e}sima linha e 
$j$-\'{e}sima coluna; Caso contrario, $^kM_i^j =0$. 
 \index{Multiplicadores} 
 \index{Pivoteamento} 

Observe que os elementos n\~{a}o nulos de $^kM$ correspondem a elementos
nulos em $^kA$, e vice-versa.  Podemos pois poupar mem\'{o}ria, guardando 
uma \'{u}nica matriz, $^kA\ +\ ^kM$.  
Assim, no Exemplo 4, $\left[ ^kA+ ^kM\ \mid \ ^kb\right] ,\ k=0,1,2$, 
pode ser representado como:
$$  
\left[ \begin{array}{rrrr} 2 & 1 & 3 & 1 \\ 
  2 & 3 & 6 & 2 \\ 4 & 4 & 6 & 6 \end{array} \right] \ \  \rightarrow \ \      
\left[ \begin{array}{rrrr} 2 & 1 & 3 & 1 \\  {\em 1} & 2 & 3 & 1 \\ 
  {\em 2} & 2 & 0 & 4 \end{array} \right] \ \  \rightarrow \ \
\left[ \begin{array}{rrrr} 2 & 1 & 3 & 1 \\   
  {\em 1} & 2 & 3 & 1 \\ {\em 2} & {\em 1} & -3 & 3 \end{array} \right] \ \ 
$$  

\section{Pivoteamento}

Consideremos a hip\'{o}tese do surgimento de um zero na posi\c{c}\~{a}o
do piv\^{o} da etapa seguinte do processo de triangulariza\c{c}\~{a}o, 
isto \'{e}, em $^{k-1}A$ anula-se o elemento $^{k-1}A_k^k$. 
 Se na coluna $^{k-1}A^k$ houver algum elemento n\~{a}o nulo na linha
$l$, $l>k$, podemos permutar as linhas $k$ e $l$ e continuar o processo 
de triangulariza\c{c}\~{a}o. 
Uma permuta\c{c}\~{a}o de linhas, para trocar o elemento piv\^{o} a ser
utilizado, denomina-se um {\bf pivoteamento}. 
A cada etapa queremos lembrar-nos de quais os pivoteamentos realizados
durante o processo.  Para tanto, guardamos os vetores de \'{\i}ndices de
linha permutados da permuta\c{c}\~{a}o corrente em rela\c{c}\~{a}o ao  
sistema original. Estes s\~{a}o os vetores de permuta\c{c}\~{a}o, 
$^1p$, $^2p$, \ldots $^{n-1}p=p$. 
No que tange ao armazenamento dos multiplicadores, devemos convencionar
se estes ser\~{a}o ou n\~{a}o permutados, junto com as respectivas
linhas de $A+M$, nas opera\c{c}\~{o}es de pivoteamento. 
Adotaremos, por hora, a conven\c{c}\~{a}o de SIM, permut\'{a}-los, junto
com os pivoteamentos. 

O exemplo 5 ilustra o processo para uma matriz de dimens\~{a}o $k=4$, 
apresentando a matriz $\left[ ^kA+ ^kM \right]$ juntamente 
com o vetor de permuta\c{c}\~{a}o, $^kp$ para cada etapa da 
triangulariza\c{c}\~{a}o: 

$$
\left[ \begin{array}{rrrr}
 2 & 1 & 9 & -1 \\ 1 & 3 & 7 & 7 \\ 
 2 & 8 & 4 & 2 \\  3 & 9 & 6 & 6 \\ \end{array} \right]  
\begin{array}{c}  1 \\  2 \\ 3 \\ 4  \end{array}\ \ \rightarrow \ \    
\left[ \begin{array}{rrrr}
 3 & 9 & 6 & 6 \\ 
 {\em 1/3} & 0 & 5 & 5 \\ 
 {\em 2/3} & 2 & 0 & -2 \\  
 {\em 2/3} & -5 & 5 & -5 \\ \end{array} \right]
 \begin{array}{c}  4 \\ 2 \\ 3 \\ 1  \end{array}  \ \ \rightarrow \ \ $$ 
$$
 \left[ \begin{array}{rrrr}
  3 & 9 & 6 &  6 \\ 
 {\em 2/3} & -5 & 5 & -5 \\ 
 {\em 2/3} & {\em -2/5} & 2 & -4 \\  
 {\em 1/3} & {\em 0} & 5 & 5 \\ \end{array} \right]
\begin{array}{c}  4 \\ 1 \\ 3 \\ 2  \end{array}\ \ \rightarrow \ \ 
\left[ \begin{array}{rrrr}
  3 & 9 & 6 & 6 \\ 
 {\em 2/3} & -5 & 5 & -5 \\ 
 {\em 1/3} & {\em 0} &  5 &  5 \\  
 {\em 2/3} & {\em -2/5} & {\em 2/5} & -6 \\ \end{array} \right] 
\begin{array}{c}  4 \\  1 \\ 2 \\ 3  \end{array}  $$

\begin{obs}{\rm 
 Note que se o sistema \'{e} bem determinado, um pivoteamento que
permite o prosseguimento do processo de triangulariza\c{c}\~{a}o
\'{e} sempre poss\'{\i}vel: Se na $k$-\'{e}sima etapa, para 
$l=k\ldots n$, $^{k-1}A_l^k =0$, isto \'{e} o elemento na posi\c{c}\~{a}o 
piv\^{o} e todos os elementos abaixo dele na coluna $^{k-1}A^k$ se anulam, 
ent\~{a}o as linhas $^{k-1}A_l$, $l=k\ldots n$, s\~{a}o linearmente 
dependentes e $det(\ ^{k-1}A)=0 \Rightarrow det(\ ^0A)=0$, o que 
contraria a hip\'{o}tese do sistema ser bem determinado.  
}\end{obs}

\begin{obs}{\rm 
 Assuma sabermos de antem\~{a}o um vetor de permuta\c{c}\~{a}o,
$^{n-1}p = p$, vi\'{a}vel no processo de triangulariza\c{c}\~{a}o de um
dado sistema $\left[ ^0A\ ^0b\right]$. 
Neste caso, poder\'{\i}amos permutar as linhas do sistema original como
indicado no vetor de permuta\c{c}\~{a}o, obtendo um novo sistema
$P\left[ ^0A\ ^0b\right]$, cujas equa\c{c}\~{o}es s\~{a}o as mesmas que
as do sistema original, a menos da ordem em que est\~{a}o escritas. 
Poder\'{\i}amos ent\~{a}o triangularizar este sistema sem necessidade de
nenhum pivoteamento, obtendo ao final o mesmo sistema equivalente,
$\left[ ^{n-1}A\ ^{n-1}b\right]$. 
}\end{obs} 

 Note que a triangulariza\c{c}\~{a}o da matriz dos coeficientes  
de um sistema \'{e} completamente independente do vetor dos termos
independentes. 
Suponha termos triangularizado o sistema $\left[ ^0A\ ^0b\right]$, isto
\'{e}, que temos um sistema equivalente e triangular 
$\left[ ^{n-1}A=U\ ^{n-1}b\right]$, tendo sido preservados os 
multiplicadores e as permuta\c{c}\~{o}es utilizadas, $M$ e $p$.  
Se for necess\'{a}rio resolver um segundo sistema 
$\left[ ^0A\ ^0b\right]$, que difere do primeiro apenas pelo vetor dos 
termos independentes, n\~{a}o ser\'{a} necess\'{a}rio retriangularizar a
matriz $A$; Bastar\'{a} efetuar no novo vetor de termos
independentes, $^0c$, as opera\c{c}\~{o}es indicadas pelo vetor de
permuta\c{c}\~{a}o $p$ e pela matriz de multiplicadores $M$.

 Uma pol\'{\i}tica, que mais tarde demonstraremos \'{u}til, \'{e}
realizarmos pivoteamentos para colocar como elemento piv\^{o}, em cada
etapa, o elemento de m\'{a}ximo m\'{o}dulo.  Esta estrat\'{e}gia, o
{\bf pivoteamento parcial}, garante termos todos os multiplicadores 
$\mid M_i^j \mid \ \leq 1$.
O pivoteamento parcial \'{e} de fundamental import\^{a}ncia
para a estabilidade num\'{e}rica do processo, controlando a
propaga\c{c}\~{a}o de erros de arrendondamento. 
 \index{Pivoteamento!parcial} 

\begin{obs}{\rm 
 Para realizar uma opera\c{c}\~{a}o de pivoteamento n\~{a}o \'{e}
necess\'{a}rio  efetivamente permutar linhas da matriz $A$: 
Basta, na $k$-\'{e}sima etapa da triangulariza\c{c}\~{a}o utilizar o 
\'{\i}ndice $^{k-1}p(i)$ ao inv\'{e}s do \'{\i}ndice de linha $i$.  
}\end{obs}

\begin{obs}{\rm 
 Uma maneira alternativa de guardar pivoteamentos realizados ao longo do
processo de triangulariza\c{c}\~{a}o \'{e} o {\bf vetor de
pivoteamentos}, $t$, onde $t(j)=i$ significa que ao eliminarmos a coluna
$j$ permutamos para a posi\c{c}\~{a}o de pivo (linha $j$) o elemento na
linha $i$.  O vetor de pivoteamentos do exemplo anterior \'{e}
 $\left[ \begin{array}{cccc} 4 & 1 & 
  2 & 3 \end{array} \right] '$. 
 Com o vetor de pivoteamentos adotaremos a
conven\c{c}\~{a}o de N\~{A}O permutar os multiplicadores na coluna
$^jM$, nos pivoteamentos $k>j$. 
}\end{obs} 
 \index{Pivoteamento!vetor de}

\section{Lema da Fatora\c{c}\~{a}o}

\begin{lm}[da Fatora\c{c}\~{a}o] Seja $^0A$ uma matriz invers\'{\i}vel
triangulariz\'{a}vel pelo m\'{e}todo de Gauss sem nenhum pivoteamento e
sejam, conforme a nomenclatura adotada,
 $A\ =\ ^oA,\ ^1A,\ \ldots \  ^{n-1}A$,   e  
 $ ^1M,\ ^2M,\ \ldots \  ^{n-1}M\ =\ M$,
respectivamente, a $k$-transformada da matriz $A$ e a $k$-\'{e}sima 
matriz dos multiplicadores. 
 \index{Fatora\c{c}\~{a}o!lema da} 
Fazendo $U\ =\ ^{n-1}A$ e $L\ =\ M + I$, temos que $L$ e $U$ s\~{a}o 
matrizes triangulares, respectivamente inferior e superior, e que 
$A=LU$.
\end{lm} 

Demonstra\c{c}\~{a}o.

Verifiquemos inicialmente que a transforma\c{c}\~{a}o linear que leva 
$^{k-1}A$ em $^kA$ \'{e} dada por $^kT$, i.e. $^kA= ^kT\ ^{k-1}A$, onde
$$^kT= \left[ \begin{array}{cccccc} 
1 & & & & & 0 \\
0 & \ddots & & & & \\
 & & 1 & & & \vdots \\ 
\vdots & & -M_{k+1}^k & & & \\
 & & \vdots & & \ddots & 0 \\ 
0 & & -M_n^k & & & 1 \\ 
\end{array} \right] \ ,$$ 
 de modo que
 $ U =\ ^{n-1}A =\ ^{n-1}T\ ^{n-2}T \ldots \ ^2T\ ^1T\ A = TA$.
Observando ainda que $^kT$ \'{e} invers\'{\i}vel  e que
$$^kT^{-1} =\ ^kL = \left[ \begin{array}{cccccc} 
1 & & & & & 0 \\
0 & \ddots & & & & \\
 & & 1 & & & \vdots \\ 
\vdots & & +M_{k+1}^k & & & \\
 & & \vdots & & \ddots & 0 \\ 
0 & & +M_n^k & & & 1 \\ 
\end{array} \right] \ ,$$ 
 Ademais, temos que
$T^{-1} =\ ^1T^{-1}\ ^2T^{-1} \ldots \ ^{n-1}T^{-1} =\  ^1L \ldots \ ^{n-1}L$
e \'{e} f\'{a}cil verificar que
$^1L\ ^2L \ldots \ ^{n-1}L= L$, donde $A=LU$. Q.E.D.

\begin{teo}[LU]
Seja A uma matriz invers\'{\i}vel $n\times n$.  Ent\~{a}o
existe uma (matriz de) permuta\c{c}\~{a}o de linhas $P$ de modo que
$\tilde A=PA$ \'{e} triangulariz\'{a}vel pelo m\'{e}todo de Gauss e
$\tilde A = L U$. 
\end{teo}
 \index{Fatora\c{c}\~{a}o!LU} 

Demonstra\c{c}\~{a}o.

O teorema segue trivialmente do lema da fatora\c{c}\~{a}o e das 
observa\c{c}\~{o}es 2.1 e 2.2.

 A solu\c{c}\~{a}o de um sistema linear $Ax=b$ pode ser expressa
diretamente em termos da fatora\c{c}\~{a}o LU da matriz dos
coeficientes, i.e., se $\tilde A =PA=LU$, ent\~{a}o $P'LUx=b$, e
$x=U^{-1}L^{-1}Pb$.  Assim, em muitas aplica\c{c}\~{o}es, o conhecimento
expl\'{\i}cito de $A^{-1}$ pode ser substitu\'{\i}do com vantagem pelo
conhecimento da fatora\c{c}\~{a}o da matriz. 

\begin{obs}{\rm 
 Considere os dois algoritmos seguintes, o de substitui\c{c}\~{a}o e o
de multilpica\c{c}\~{a}o pela inversa, para solu\c{c}\~{a}o do sistema
$Ux=b$: O primeiro algoritmo acessa a matriz $U$ por linha, enquanto o
segundo acessa a matriz $U$ por coluna.  }\end{obs}

\begin{verbatim}
     x = b ;
     x(n) = x(n) / U(n,n) ; 
     for i = n-1:-1:1 
       x(i) = ( x(i) - U(i,i+1:n)*x(i+1:n) ) / U(i,i) ;
     end 
\end{verbatim}

\begin{verbatim}
     x = b ;
     for j = n:-1:2 
       x(j) = x(j) / U(j,j) ; 
       x(1:j-1) = x(1:j-1) - x(j)*U(1:j-1,j) ;
     end 
     x(1) = x(1) / U(1,1) ; 
\end{verbatim}

\section{O M\'{e}todo de Doolittle}

Usando o Teorema LU podemos determinar L e U diretamente da
equa\c{c}\~{a}o de decomposi\c{c}\~{a}o, $LU=A$, tomando, nesta ordem,
para $i=1\ldots n$, as equa\c{c}\~{o}es das componentes de $A$ na
$i$-\'{e}sima linha e na $i$-\'{e}sima coluna, isto \'{e},
 \index{Fatora\c{c}\~{a}o!Doolittle} 
$$
\mbox{para}\ i=1\ldots n \left\{ \begin{array}{c}
 L_i U= A_i \\ LU^i=A^i \end{array} \right. $$
temos
$$
\mbox{para}\ i=1\ldots n \left\{ \begin{array}{cc}
 \mbox{para}\ j=i\ldots n & U_i^j = A_i^j -\sum_{k=1}^{i-1} M_i^k U_k^j \\
 \mbox{para}\ j=i+1\ldots n & M_j^i = (A_j^i -\sum_{k=1}^{i-1} M_j^k U_k^i)
 /U_i^i  \end{array} \right. $$

Note que s\'{o} escrevemos as equa\c{c}\~{o}es para os termos
inc\'{o}gnitos, isto \'{e}, ou acima da diagonal em U, e abaixo da
diagonal em L.  Tamb\'{e}m as somat\'{o}rias foram interrompidas quando
os termos remanescentes s\~{a}o todos nulos.  

O c\'{a}lculo do elemento $(M + U)_i^j$ , envolve, na ordem prescrita,
apenas elementos j\'{a} calculados de $(M+U)$.  Ademais, o c\'{a}lculo
de $(M+U)_i^j$ \'{e} a \'{u}ltima ocasi\~{a}o em que se necessita os
elementos $A_i^j$; portanto podemos armazenar $(M+U)_i^j$ no lugar de
$A_i^j$. Estas s\~{a}o as matrizes $A=\ ^0D,\ldots ,\ ^nD=M+U$.

Exemplo 6 - Usando o m\'{e}todo de Doolittle para triangularizar a
matriz do Exemplo 4, onde n\~{a}o h\'{a} necessidade de pivoteamento,
temos
 
$$ 
^0D=A=\left[ \begin{array}{rrr} 2 & 1 & 3 \\
  2 & 3 & 6 \\ 4 & 4 & 6 \end{array} \right] \ \  \rightarrow $$
$$
\left[ \begin{array}{rrr} {\em 2} & {\em 1} & {\em 3} \\ 
  2 & 3 & 6 \\ 4 & 4 & 6 \end{array} \right] \ \ \rightarrow \ \  
\left[ \begin{array}{rrr} {\em 2} & {\em 1} & {\em 3} \\ 
  {\em 1} & 3 & 6 \\ {\em 2} & 4 & 6 \end{array} \right] 
  =\ ^1D   \ \  \rightarrow $$ 
$$
\left[ \begin{array}{rrr} {\em 2} & {\em 1} & {\em 3} \\ 
  {\em 1} & {\em 2} & {\em 3} \\ {\em 2} & 4 & 6 \end{array} \right] 
  \ \ \rightarrow \ \  
\left[ \begin{array}{rrr} {\em 2} & {\em 1} & {\em 3} \\ 
  {\em 1} & {\em 2} & {\em 3} \\ {\em 2} & {\em 1} & 6 \end{array} \right] 
  =\ ^2D \ \ \rightarrow  $$  
$$ 
\left[ \begin{array}{rrr} {\em 2} & {\em 1} & {\em 3} \\ 
  {\em 1} & {\em 2} & {\em 3} \\ {\em 2} & {\em 1} & {\em -3} 
  \end{array} \right] =\ ^3D = M+U $$

Para realizar os pivoteamentos necess\'{a}rios \`{a} passagem de $^kD$ a
$^{k+1}D$ \'{e} necess\'{a}rio examinar os poss\'{\i}veis elementos
piv\^{o}s nas linhas $i=k\ldots n$, isto \'{e}, os elementos em $^kA^k$. 
 Para tanto, basta calcular, para $i=k\ldots n$, 
 $$^kv_i =\ ^0A_i^k -\sum_{l=1}^{k-1} M_i^l U_l^k\ .$$

Exemplo 7 - Triangularizando pelo m\'{e}todo de Doolittle a matriz, $A$,
com pivoteamento parcial, temos os $^kp,\ ^kD,\ ^kv$, para 
$k=0\ldots n$, como segue:

$$
\begin{array}{c}  1 \\  2 \\ 3 \\ 4  \end{array}
\left[ \begin{array}{rrrr}
 2 & 1 & 9 & -1 \\ 1 & 3 & 7 & 7 \\ 
 2 & 8 & 4 & 2 \\  3 & 9 & 6 & 6 \\ \end{array} \right]  
\begin{array}{c}  2 \\  1 \\ 2 \\ 3  \end{array}\ \ \rightarrow \ \ 
\begin{array}{c}  4 \\  2 \\ 3 \\ 1  \end{array}  
\left[ \begin{array}{rrrr}
 {\em 3} & {\em 9} & {\em 6} & {\em 6} \\ {\em 1/3} & 3 & 7 & 7 \\ 
 {\em 2/3} & 8 & 4 & 2 \\  {\em 2/3} & 1 & 9 & -1 \\ \end{array} \right]
\begin{array}{c}  * \\ 0 \\ 2 \\ -5  \end{array}  \ \ \rightarrow \ \ $$ 
$$
\begin{array}{c}  4 \\  1 \\ 3 \\ 2  \end{array}
\left[ \begin{array}{rrrr}
 {\em 3} & {\em 9} & {\em 6} & {\em 6} \\ 
 {\em 2/3} & {\em -5} & {\em 5} & {\em -5} \\ 
 {\em 2/3} & {\em -2/5} & 4 & 2 \\  
 {\em 1/3} & {\em 0} & 7 & 7 \\ \end{array} \right]
\begin{array}{c}  * \\ * \\ 2 \\ 5  \end{array}\ \ \rightarrow \ \ 
\begin{array}{c}  4 \\  1 \\ 2 \\ 3  \end{array}  
\left[ \begin{array}{rrrr}
 {\em 3} & {\em 9} & {\em 6} & {\em 6} \\ 
 {\em 2/3} & {\em -5} & {\em 5} & {\em -5} \\ 
 {\em 1/3} & {\em 0} & {\em 5} & {\em 5} \\  
 {\em 2/3} & {\em -2/5} & {\em 2/5} & 2 \\ \end{array} \right]
\ \rightarrow \ \ $$
$$ 
\begin{array}{c}  4 \\  1 \\ 2 \\ 3  \end{array}  
\left[ \begin{array}{rrrr}
 {\em 3} & {\em 9} & {\em 6} & {\em 6} \\ 
 {\em 2/3} & {\em -5} & {\em 5} & {\em -5} \\ 
 {\em 1/3} & {\em 0} & {\em 5} & {\em 5} \\  
 {\em 2/3} & {\em -2/5} & {\em 2/5} & {\em -6} \\ \end{array} \right] $$

\section{Complexidade}

Contemos agora o n\'{u}mero de opera\c{c}\~{o}es aritm\'{e}ticas
necess\'{a}rias \`{a} triangulariza\c{c}\~{a}o de uma matriz de
coeficientes e \`{a} solu\c{c}\~{a}o de um sistema linear. 
Na fase de triangulariza\c{c}\~{a}o da matriz dos coeficientes,
vemos que o c\'{a}lculo de $^kA$ a partir de $^{k-1}A$ requer  
$(n-k)$ divis\~{o}es e $(n-k)^2$ somas e produtos.  No total s\~{a}o
portanto requeridos
 \index{Fatora\c{c}\~{a}o!complexidade}  
\begin{eqnarray*}
\sum_{k=1}^{n-1} (n-k) = n(n-1)/2 & & \mbox{divis\~oes, e} \\ 
\sum_{k=1}^{n-1} (n-k)^2 = n(n^2-1)/3 + n(n-1)/2 & & 
                              \mbox{produtos e subtra\c{c}\~{o}es.}
\end{eqnarray*} 
 isto \'{e}, s\~{a}o necess\'{a}rios da ordem de
 $n^3/3\ +O(n^2)$ produtos e subtra\c{c}\~{o}es, e 
 $n^2/2\ +O(n)$ divis\~{o}es. 

Analogamente, dada a matriz dos multiplicadores e o vetor das
permuta\c{c}\~{o}es, $M$ e $p$, o tratamento de um vetor de termos
independentes requer $n(n-1)/2$ produtos e subtra\c{c}\~{o}es.
Finalmente, a solu\c{c}\~{a}o do sistema linear triangularizado requer
$n$ divis\~{o}es e $n(n-1)/2$ produtos e subtra\c{c}\~{o}es.

\section*{Exerc\'{\i}cios}
\begin{enumerate}

\item Comente como cada forma de representa\c{c}\~{a}o da matriz $A$,
por linhas ou por colunas, favorece o uso de um dos algoritmos
apresentados na observa\c{c}\~{a}o 2.5.  Escreva algoritmos similares 
para a solu\c{c}\~{a}o do sistema $Lx=b$.  

\item Programe e implemente, em linguagem C, C++, ou FORTRAN-90, 
      fun\c{c}\~{o}es para:
\begin{enumerate}
\item Fatora\c{c}\~{a}o LU com pivoteamento parcial. 
      \item Solu\c{c}\~{a}o dos sistemas $Lx=b$ e $Ux=b$, 
      \item Um argumento adicional deve indicar a forma de 
      representa\c{c}\~{a}o das matrizes, densa, est\'{a}tica por 
      linha, ou est\'{a}tica por coluna.
      \item No caso da representa\c{c}\~{a}o densa, um segundo argumento 
      deve indicar o uso do vetor de permuta\c{c}\~{o}es ou pivoteamentos. 
      No caso das representa\c{c}\~{o}es est\'{a}ticas, escreva o programa
      como lhe parecer mais conveniente.
\end{enumerate} 

\item
 Mostre que calcular explicitamente a inversa de $A$ implica saber
resolver $n$ sistemas lineares, $Ax=I^j$. Dada a fatora\c{c}\~{a}o 
$A=LU$, qual o custo de computar explicitamente a inversa $A^{-1}$?  

\end{enumerate}

 \clearpage
 \clearpage 
\setcounter{chapter}{2} 
\chapter{RESUMO DE TEORIA DOS GRAFOS}

\section{Conceitos B\'{a}sicos}

Um {\bf grafo} \'{e} um par ordenado $G=(V_G,\Gamma _G)$ onde o primeiro
elemento \'{e} um conjunto finito, o conjunto de {\bf v\'{e}rtices}, e o
segundo elemento \'{e} uma fun\c{c}\~{a}o $\Gamma _G: V_G \mapsto
P(V_G)$,no conjunto das partes de $V_G$, a {\bf fun\c{c}\~{a}o de
filia\c{c}\~{a}o}. 
 Tomaremos $V_G$ um conjunto indexado, 
$V_G = \{v_1, v_2,\ldots v_n\}$ ou ent\~{a}o tomaremos $V_G$ 
como sendo o pr\'{o}prio conjunto dos $n$ primeiros inteiros 
positivos, $N=\{1, 2,\ldots n\}$. 
 Para n\~{a}o sobrecarregar a nota\c{c}\~{a}o escreveremos, quando 
n\~{a}o houver ambig\"{u}idade,  
 $(V_G, \Gamma _G)$   como   $(V,\Gamma)$.
 \index{V\'{e}rtice} 
  \index{Grafo!fun\c{c}\~{a}o de filia\c{c}\~{a}o} 

\'{E} comum representarmos um grafo por conjunto de pontos (os
v\'{e}rtices) e um conjunto de {\bf arestas} (setas) que v\~{a}o de cada
v\'{e}rtice para os seus filhos, isto \'{e}, de cada $v\in V$ para os
elementos em $\Gamma (v)$. 
 Podemos definir o grafo atrav\'{e}s dos seus v\'{e}rtices
e de suas arestas, $G=(V_G,A_G)$, onde cada aresta \'{e} um par ordenado
de v\'{e}rtices e $(i,j)\in A_G \Leftrightarrow j\in \Gamma_G (i)$. 
 Uma terceira maneira de definir um grafo \'{e} pelo par $G=(V_G,B_G)$ 
onde $B_G$ \'{e} a {\bf matriz de adjac\^{e}ncia}, a matriz Booleana 
tal que $B_i^j=1 \Leftrightarrow j\in \Gamma (i)$.
  \index{Grafo!aresta} 
  \index{Matriz!de adjac\^{e}ncia} 

Exemplo 1:  

Considere o grafo de v\'{e}rtices $N=\{1,2,3,4,5,6\}$ e fun\c{c}\~{a}o
de filia\c{c}\~{a}o
$\Gamma (1)=\emptyset$, $\Gamma (2)= \{2\}$, $\Gamma (3)=\emptyset$, 
$\Gamma (4)=\{3,5,6\}$, $\Gamma (5)=\{3,4,5\}$ e $\Gamma (6)=\{5\}$.
Suas arestas e matriz de adjac\^{e}ncia s\~{a}o, \\ 
$A_G = \{(2,2),(4,3),(4,5),(4,6),(5,3),(5,4),(5,5),(6,5)\}$, e 
$$
 B_G = \left[ \begin{array}{cccccc} 
 0 & 0 & 0 & 0 & 0 & 0 \\ 0 & 1 & 0 & 0 & 0 & 0 \\ 
 0 & 0 & 0 & 0 & 0 & 0 \\ 0 & 0 & 1 & 0 & 1 & 1 \\ 
 0 & 0 & 1 & 1 & 1 & 0 \\ 0 & 0 & 0 & 0 & 1 & 0 
 \end{array} \right] \hspace{2cm}  
\begin{array}{rcccl}
 1 & & 3 & \leftarrow & 5 \hookleftarrow \\ 
 & & \uparrow & \nearrow \! \swarrow & \uparrow \\ 
 \hookrightarrow 2 & & 4 & \rightarrow & 6 
\end{array} $$

Dado um subconjunto de v\'{e}rtices $W\subset V$, definimos
$\Gamma (W) \equiv \cup _{w\in W} \Gamma (w)$. Definimos tamb\'{e}m
$\Gamma ^0 (i) = \{i\}$, $\Gamma ^1 (i) \equiv \Gamma (i)$ e, 
para $k>1$, $\Gamma ^k (i) \equiv \Gamma (\Gamma ^{k-1} (i) )$.
Estes s\~{a}o, para $k=1,2,3,\ldots$ os filhos, netos, bisnetos,
etc. do v\'{e}rtice $i$. Finalmente definimos os {\bf descendentes} de 
$i$ por $\bar \Gamma (i) = \cup _{k=0}^\infty \Gamma ^k (i)$, e o 
{\bf fecho transitivo} $\bar G= (V_G, \bar \Gamma _G)$.
  \index{V\'{e}rtice!descendentes} 
  \index{Grafo!fecho transitivo} 

Exemplo 2:

Damos a matriz de adjac\^{e}ncia de um grafo, e do respectivo  
fecho transitivo:
$$
B = \left[ \begin{array}{cccc} 
0 & 1 & 0 & 0 \\ 0 & 1 & 1 & 0 \\ 0 & 0 & 0 & 1 \\ 0 & 0 & 0 & 0 
\end{array} \right] \ \ , \ \ \ 
\bar B = \left[ \begin{array}{cccc} 
1 & 1 & 1 & 1 \\ 0 & 1 & 1 & 1 \\ 0 & 0 & 1 & 1 \\ 0 & 0 & 0 & 1
\end{array} \right] $$

Dado o grafo $G=(V,\Gamma)$ definimos a {\bf fun\c{c}\~{a}o de 
paternidade} $i\in \Gamma ^{-1} (j) \Leftrightarrow j\in \Gamma (i)$.
Definimos tamb\'{e}m seu {\bf grafo inverso},  
$G ^{-1} = (V , \Gamma ^{-1})$.
  \index{Grafo!inverso} 
 \index{Grafo!fun\c{c}\~{a}o de paternidade} 
 \index{V\'{e}rtice!ascendentes} 
 
\begin{lm}
Dado um grafo definido por sua matriz de adjac\^{e}ncia, $G=(N,B)$, 
seu grafo inverso \'{e} $G^{-1}=(N,B')$.
\end{lm}

Um {\bf sub-grafo} de $G=(V,A)$ \'{e} um grafo $G'=(V',A')$, onde
$V'\subset V$ e $A'$ \'{e} um subconjunto de arestas de $A$ que tem
ambos os v\'{e}rtices em $V'$.   O sub-grafo induzido por um
subconjunto de v\'{e}rtices $V'$ \'{e} $G=(V',A')$ onde $A'$ \'{e}
m\'{a}ximo, e o sub-grafo induzido por um subconjunto de arestas $A'$
\'{e} $G=(V',A')$ onde $V'$ \'{e} m\'{\i}nimo. 
 \index{Grafo!sub-grafo} 
  \index{Grafo!induzido} 

Uma aresta que parte e chega no mesmo v\'{e}rtice \'{e} dita um {\bf
loop}.  Duas arestas, $(i,j)$ e $(k,h)$, s\~{a}o ditas {\bf
cont\'{\i}guas} se a primeira chega no v\'{e}rtice de que parte a
segunda, isto \'{e} se $j=k$.  Um {\bf passeio} \'{e} uma
seq\"{u}\^{e}ncia n\~{a}o vazia e finita de arestas cont\'{\i}guas.  Um
{\bf circuito} \'{e} um passeio que parte e chega no mesmo v\'{e}rtice,
isto \'{e}, o primeiro v\'{e}rtice da primeira aresta \'{e} o segundo
v\'{e}rtice da \'{u}ltima aresta.  Uma {\bf trilha} \'{e} um passeio
onde n\~{a}o se repete nenhuma aresta e um {\bf caminho} \'{e} uma
trilha na qual n\~{a}o h\'{a} (subseq\"{u}\^{e}ncias que sejam)
circuitos.  Uma trilha na qual o \'{u}nico circuito \'{e} toda a trilha
\'{e} uma {\bf ciclo}. 
 \index{Loop}
 \index{Passeio} 
 \index{Circuito} 
 \index{Trilha} 
 \index{Caminho} 
 \index{Ciclo} 
 
Notemos que um passeio $C=(w_0,w_1), (w_1,w_2),\ldots (w_{k-1},w_k))$,
\'{e} igualmente bem determinado pela seq\"{u}\^{e}ncia de seus v\'{e}rtices
$C=(w_0,w_1,w_2,\ldots w_k)$.

Ter\'{\i}amos assim, no Exemplo 1, exemplos de:
\begin{itemize}
\item loop: $((2,2))$.
\item passeio: $((4,5),(5,5),(5,5),(5,4),(4,5)$.
\item circuito: $(5,5,4,6,5,5)$.
\item trilha: $(5,4,6,5,5,3)$.
\item caminho: $(6,5,4,3)$.
\item ciclo: $(5,4,6,5) ou (5,5)$.
\end{itemize}

Um grafo \'{e} {\bf ac\'{\i}clico} se n\~{a}o cont\'{e}m ciclos.
Uma {\bf \'{a}rvore} de {\bf raiz} $v$ \'{e} um grafo ac\'{\i}clico, 
$H=(V_H,\Gamma _H)\ ,\ v\in V_H$,  onde todos os 
v\'{e}rtices t\^{e}m no m\'{a}ximo um pai e
apenas a raiz n\~{a}o tem pai.  Os v\'{e}rtices sem filhos de uma
\'{a}rvore dizem-se folhas da \'{a}rvore. 
Uma cole\c{c}\~{a}o de \'{a}rvores \'{e} denominada {\bf floresta}.
Uma floresta cobre um grafo $G$ sse \'{e} um subgrafo de $G$ que
cont\'{e}m todos os seus v\'{e}rtices. 
 Seguem exemplos de florestas que cobrem o grafo do Exemplo 1, estando
assinaladas as ra\'{\i}zes. 
 \index{Arvore} 
 \index{Arvore!floresta} 
 \index{Arvore!raiz}  
 \index{Arvore!folha} 

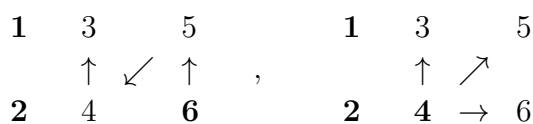
$$
\begin{array}{rcccl}
 {\bf 1} & & 3 &  & 5  \\ 
 & & \uparrow & \swarrow & \uparrow \\ 
 {\bf 2} & & 4 &  & {\bf 6} 
\end{array} 
\ \ \ \ , \ \ \ \ \ \     
\begin{array}{rcccl}
 {\bf 1} & & 3 &  & 5  \\ 
 & & \uparrow & \nearrow  &  \\ 
 {\bf 2} & & {\bf 4} & \rightarrow & 6 
\end{array} 
$$

\begin{lm}
Dada uma \'{a}rvore, $H=(V,\Gamma )$, de raiz $v$, e um v\'{e}rtice
$w\in V$, existe um \'{u}nico passeio que parte de $v$ e termina
em $w$, a saber,
 $(v, \Gamma ^{-k}(w),\ldots \Gamma ^{-2}(w), \Gamma ^{-1}(w), w)$.
Ademais, este passeio \'{e} um caminho. 
\end{lm}

\section{Rela\c{c}\~{o}es de Ordem}

Uma {\bf ordem} num conjunto $S$ \'{e} uma rela\c{c}\~{a}o, $\leq $, tal que
$\forall a,b,c \in S$, temos que
\begin{enumerate}
\item $a\leq  a$.
\item $a\leq  b\ \wedge \ b\leq   c \Rightarrow  a\leq  c$.
\end{enumerate}
i.e., uma rela\c{c}\~{a}o {\bf reflexiva} e {\bf transitiva}.
Alternativamente denotaremos $a\leq b$ por $b\geq a$.
 \index{Ordem} 
 \index{Rela\c{c}\~{a}o!reflexiva} 
 \index{Rela\c{c}\~{a}o!transitiva} 

Sendo $=$ a rela\c{c}\~{a}o identidade e $\not \leq$ a nega\c{c}\~{a}o
da rela\c{c}\~{a}o de ordem, dizemos que a ordem \'{e}
\begin{itemize}
\item {\bf total}, sse $a \not \leq b \Rightarrow  b\leq  a$.
\item {\bf parcial}, sse $(a\leq  b\ \wedge \ b\leq a) \Rightarrow  a=b$.
\item {\bf boa}, se \'{e} parcial e total.
\end{itemize}
 \index{Ordem!total} 
 \index{Ordem!parcial} 
 \index{Ordem!boa} 

\begin{lm} 
Dado um grafo, $G=(N,\Gamma )$, a fun\c{c}\~{a}o de descend\^{e}ncia,
$\bar \Gamma$, define uma ordem em seus v\'{e}rtices, a ordem natural, 
$\leq$, definida por por
$j\geq i \Leftrightarrow j\in \bar \Gamma (i)$.  
Note, por\'{e}m, que a ordem natural n\~{a}o \'{e}, em geral,
nem parcial nem total.
\end{lm}

Exemplo 3:

Nos grafos seguintes, de matrizes de adjac\^{e}ncia $B_1$,
$B_2$, $B_3$ e $B_4$,  
$$
\left[ \begin{array}{cccc}
0 & 1 & 0 & 0 \\ 1 & 0 & 1 & 0 \\ 0 & 0 & 0 & 0 \\ 0 & 0 & 0 & 0 
\end{array} \right] \ , \ \ 
\left[ \begin{array}{cccc}
0 & 1 & 0 & 0 \\ 1 & 0 & 0 & 0 \\ 0 & 0 & 0 & 1 \\ 0 & 0 & 1 & 0 
\end{array} \right] \ , \ \ 
\left[ \begin{array}{cccc}
0 & 1 & 1 & 0 \\ 0 & 0 & 0 & 1 \\ 0 & 0 & 0 & 1 \\ 0 & 0 & 0 & 0 
\end{array} \right] \ , \ \ 
\left[ \begin{array}{cccc}
0 & 1 & 0 & 0 \\ 0 & 0 & 1 & 0 \\ 0 & 0 & 0 & 1 \\ 0 & 0 & 0 & 0 
\end{array} \right] 
$$
temos ordens naturais
\begin{itemize}
\item em $G_1$: Nem parcial, pois $1\leq 2\leq 1$, 
                nem total, pois $1 \not \leq 3 \not \leq 1$.
\item em $G_2$: Total, mas n\~{a}o parcial, 
                pois $1 \leq 4 \leq 1$.
\item em $G_3$: Parcial, mas n\~{a}o total, 
                pois $2 \not \leq 3 \not \leq 2$.
\item em $G_4$: Boa.
\end{itemize} 

Uma {\bf equival\^{e}ncia} num conjunto S \'{e} uma rela\c{c}\~{a}o, 
$\sim$, tal que $\forall a,b,c \in S$:
\begin{enumerate}
\item $a \sim a$.
\item $a \sim b \ \wedge \ b \sim c \Rightarrow  a \sim c$. 
\item $a \sim b \Rightarrow b \sim a$.
\end{enumerate}
i.e., uma rela\c{c}\~{a}o reflexiva, transitiva e {\bf sim\'{e}trica}.
 \index{Rela\c{c}\~{a}o!sim\'{e}trica}
 \index{Equival\^{e}ncia} 

A {\bf classe de equival\^{e}ncia} de um elemento qualquer, $a\in S$,
\'{e} o sub-conjunto de $S$: 
$[a]=\{x\in S\mid x \sim a\}$.
  \index{Equival\^{e}ncia!classe de} 

Uma {\bf parti\c{c}\~{a}o} de um conjunto $S$ \'{e} uma cole\c{c}\~{a}o 
$P$ de sub-conjuntos n\~{a}o vazios de $S$ tal que:
  \index{Parti\c{c}\~{a}o}  
\begin{enumerate}
\item $\forall X, Y \in P,\ X\neq Y \Rightarrow X \cap Y =\emptyset$.
\item $\cup _{X\in P} X =S$.
\end{enumerate}
i.e., uma cole\c{c}\~{a}o de conjuntos disjuntos que reunidos \'{e}
igual a $S$.

\begin{lm} 
Dado $S$, um conjunto munido de uma ordem,
$\forall a,b\in S$, $a\sim b \Leftrightarrow a\leq b\ \wedge b\leq a$,
define uma rela\c{c}\~{a}o de equival\^{e}ncia.
Esta \'{e} a equival\^{e}ncia induzida pela ordem.    
Dado $S$ um conjunto munido de uma equival\^{e}ncia, 
o conjunto das classes de equival\^{e}ncia em $S$
\'{e} uma parti\c{c}\~{a}o de $S$. 
\end{lm} 

As {\bf componentes fortemente conexas}, CFC, de um grafo $G$
s\~{a}o as classes de equival\^{e}ncia induzida pela ordem natural
nos v\'{e}rtices de $G$. $G$ \'{e} dito fortemente conexo sse possui 
uma \'{u}nica CFC. 
 \index{Grafo!CFC} 

Exemplo 4:

As CFC dos grafos do Exemplo 3 s\~{a}o:
\begin{itemize}
\item em $G_1$:  $\tilde V =\{\{1,2\},\{3\},\{4\}\}$.
\item em $G_2$: $\tilde V =\{\{1,2\},\{3,4\}\}$.
\item em $G_3$ e $G_4$: $\tilde V =\{\{1\},\{2\},\{3\},\{4\}\}$.
\end{itemize}  

O {\bf grafo reduzido}, de um grafo $G=(V,\Gamma )$, \'{e} o grafo 
$G=(\tilde V,\tilde \Gamma )$, que tem por v\'{e}rtices 
as CFCs de $G$: $\tilde V= \{V_1, V_2,\ldots V_k\}$, 
e a fun\c{c}\~{a}o de filia\c{c}\~{a}o, $\tilde \Gamma$, definida por
 \index{Grafo!reduzido} 
$$ 
V_s \in \tilde \Gamma (V_r) \Leftrightarrow 
\exists i\in V_r,\ j\in V_s \mid j\in \Gamma (i),
\ (V_r \neq V_s)\ .$$

Exemplo 5:

As matrizes de adjac\^{e}ncia dos grafos reduzidos dos grafos do 
Exemplo 3, com os v\'{e}rtices correspondendo as CFCs na ordem em 
que aparecem no exemplo anterior, s\~{a}o $\tilde B_1$,
$\tilde B_2$, $\tilde B_3$ e $\tilde B_4$: 
$$
\left[ \begin{array}{ccc}
0 & 1 & 0 \\ 0 & 0 & 0 \\ 0 & 0 & 0 \end{array} \right] , \ \  
\left[ \begin{array}{cc}
0 & 0 \\ 0 & 0 \end{array} \right] , \ \  
\left[ \begin{array}{cccc}
0 & 1 & 1 & 0 \\ 0 & 0 & 0 & 1 \\ 0 & 0 & 0 & 1 \\ 0 & 0 & 0 & 0 
\end{array} \right] , \ \  
\left[ \begin{array}{cccc}
0 & 1 & 0 & 0 \\ 0 & 0 & 1 & 0 \\ 0 & 0 & 0 & 1 \\ 0 & 0 & 0 & 0 
\end{array} \right] .$$

Dado um conjunto $S$ e $\leq$ uma ordem em $S$, a {\bf ordem reduzida}
\'{e} a rela\c{c}\~{a}o de ordem no conjunto das classes
(da equival\^{e}ncia induzida pela ordem $\leq$), 
$\tilde S$, definida por
$$
\forall X,Y\in \tilde S \ , \ \ X \leq Y 
\Leftrightarrow \exists x\in X, \wedge 
\exists y\in Y \mid x\leq y\ .$$

\begin{lm}  
A ordem reduzida (nas CFCs) de um grafo $G=(V,\Gamma )$,
\'{e} a ordem natural do grafo reduzido 
$\tilde G=(\tilde V, \tilde \Gamma )$, i.e., se
$\tilde V =\{V_1,\ldots V_k\}$,
 \index{Ordem!reduzida} 
\begin{eqnarray*}
V_r \leq V_s & \Leftrightarrow & V_s \in \bar {\tilde \Gamma} (V_r) \\ 
 & \Leftrightarrow & 
 \exists i\in V_r,\ j\in V_s \mid j\in \bar \Gamma (i) \\
 & \Leftrightarrow & 
 \exists i\in V_r,\ j\in V_s \mid j \geq i \ .
\end{eqnarray*} 
Ademais, 
\begin{enumerate}
\item A ordem natural de $G$ \'{e} parcial sse,
a menos de loops, $G$ \'{e} ac\'{\i}clico.
\item Grafos reduzidos s\~{a}o ac\'{\i}clicos e ordens reduzidas 
s\~{a}o parciais. 
\end{enumerate} 
\end{lm}

 \begin{teo}[Hoffman]: 
 Um grafo $G=(V,\Gamma )$ \'{e} fortemente conexo sse dado qualquer
sub-conjunto pr\'{o}prio dos v\'{e}rtices $W \subset V$, $W\neq V$,
houver uma aresta de um v\'{e}rtice em $W$ para um v\'{e}rtice fora de
$W$. 
 \end{teo}
 \index{Teorema!Hoffman} 

Demonstra\c{c}\~{a}o:

 Se houver $W$ um subconjunto pr\'{o}prio de $V$ tal que $\Gamma
(W)=\emptyset$, ent\~{a}o n\~{a}o h\'{a} caminho de nenhum v\'{e}rtice
$w\in W$ para nenhum v\'{e}rtice $v\in V$, e $G$ n\~{a}o \'e fortemente
conexo. 
 Presupondo a condi\c{c}\~{a}o 
 $\not \exists W\neq V \mid \Gamma (W)=\emptyset$, 
 provemos que existe um caminho do v\'{e}rtice $w$ ao v\'{e}rtice $v$,
$\forall w, v \in V$.  Tomemos inicialmente $W_0=\{w\}$.  
 A condi\c{c}\~{a}o garante a exist\^{e}ncia de um caminho de
comprimento (n\'{u}mero de arestas) 1 de $w$ a algum v\'{e}rtice 
 $x_1\neq w$. 
 Se $x_1 =v$ a prova est\'{a} completa.  Caso contr\'{a}rio tomemos
$W_1=\{w,x_1\}$.  A condi\c{c}\~{a}o garante a exist\^{e}ncia de um
caminho $c_2$ de comprimento $|c_2|\leq 2$ de $w$ para algum v\'{e}rtice
$x_2\neq w, x_1$.  Repetindo o argumento at\'{e} obtermos $x_k =v$, 
$k<n$, conclu\'{\i}mos a demonstra\c{c}\~{a}o, QED. 

Dado um conjunto $S$, munido de uma ordem $\leq$, a ordem de $S$,
uma boa ordem, $<=$, em $S$ \'{e} dita uma {\bf ordem coerente} 
(com $\leq$) sse:
 \begin{enumerate}
 \item $\forall a,b\in S,\ a\leq b\ \wedge \ b\not \leq a 
      \Rightarrow  a <= b \ .$
\item $\forall a,b,c\in S,\ a <= b <= c \ \wedge \ a \sim c 
      \Rightarrow a \sim b \sim c \ .$  
\end{enumerate}

O primeiro crit\'{e}rio de coer\^{e}ncia determina que o reordenamento
se subordine \`{a} ordem parcial do grafo reduzido; O segundo
crit\'{e}rio determina que se {\bf discriminem} (n\~{a}o se misturem)
v\'{e}rtices em CFCs incompar\'{a}veis por\'{e}m distintas. 

Uma boa ordem, ou reordenamento, no conjunto $N=\{1,2,\ldots n\}$,
$q_1 <= q_2 <= \ldots <= q_n$, corresponde a um vetor de permuta\c{c}\~{a}o
$q=[q_1,q_2,\ldots q_n]= [\sigma (1), \sigma (2), \ldots \sigma(n)]$.
Um reordenamento coerente dos v\'{e}rtices de um grafo $G$ \'{e} 
um reordenamento coerente com a ordem natural do grafo.  
Um reordenamento coerente num grafo ac\'{\i}clico \'{e} tamb\'{e}m 
dito um (re){\bf ordenamento topol\'{o}gico} (dos v\'{e}rtices) 
deste grafo. 
 \index{Ordem!coerente}
 \index{Ordem!discriminante} 
 \index{Ordem!topol\'{o}gica}  

Exemplo 6:

Listamos a seguir alguns reordenamentos nos grafos do
Exemplo 3, indicando se satisfazem aos dois crit\'{e}rios de
coer\^{e}ncia. 
No grafo $G_1$: $q=[4;1,2;3]\ (1,2)$, $q=[1,3,2,4]\ (1)$,
 $q=[3,4,1,2]\ (2)$, $q=[1,3,4,2]\ ()$.
No grafo $G_2$: $q=[3,4,1,2]\ (1,2)$, $q=[1,3,2,4]\ (1)$.
No grafo $G_3$: $q=[1,3,2,4]\ (1,2)$, $q=[1,2,3,4]\ (1,2)$, e estes
s\~{a}o os \'{u}nicos reordenamentos coerentes.
No grafo $G_4$: O \'{u}nico reordenamento coerente \'{e} $q=[1,2,3,4]$.

Para reordenar coerentemente os v\'{e}rtices de um grafo precisamos 
pois determinar o grafo reduzido, e no grafo reduzido uma ordem 
topol\'{o}gica. Esta tarefa pode ser levada a cabo com o algoritmo
de Tarjan, a ser visto a seguir.

\section{Busca em Profundidade }

A busca em profundidade \'{e} um procedimento que nos permite visitar, a
partir de um determinado v\'{e}rtice $v$ de $G=(N,\Gamma )$, todos os
seus descendentes.  Cada v\'{e}rtice ser\'{a} marcado como ``j\'{a}
visitado" ou ``n\~{a}o visitado", sendo ``n\~{a}o visitado" o estado
inicial de todos os v\'{e}rtices. 
 \index{Busca!em profundidade} 

Na busca em profundidade, partindo de $v$, seguiremos um caminho, sempre
por v\'{e}rtices n\~{a}o visitados, o mais longe poss\'{\i}vel.  Ao
passar por um v\'{e}rtice, marc\'{a}-lo-emos como ``j\'{a} visitado". 
N\~{a}o sendo mais poss\'{\i}vel prosseguir de um dado v\'{e}rtice, isto
\'{e}, quando estivermos num v\'{e}rtice sem filhos n\~{a}o visitados,
retornaremos ao v\'{e}rtice de onde este foi atingido e prosseguimos na
busca em profundidade.  Quando tivermos voltado ao v\'{e}rtice inicial,
v, e j\'{a} tivermos visitado todos os seus filhos, teremos terminado a
busca em profundidade.  Os v\'{e}rtices visitados e as arestas
utilizadas para ating\'{\i}-los formam uma \'{a}rvore de raiz v, que 
\'{e} a {\bf \'{a}rvore da busca}. 

Iniciando novas buscas em profundidade, nas quais consideramos
``n\~{a}o visitados" apenas os v\'{e}rtices que n\~{a}o pertencem a
nenhuma \'{a}rvore previamente formada, teremos uma floresta que cobre o
grafo $G$, a floresta de busca. 
 \index{Busca!\'{a}rvore} 

Uma maneira de rotularmos (ou reordenarmos) os v\'{e}rtices de $G$
durante a busca em profundidade \'{e} pela {\bf ordem de
visita\c{c}\~{a}o}, $Ov( )$, onde $Ov(i)=k$ significa que o v\'{e}rtice
$i$ foi o $k$-\'{e}simo v\'{e}rtice a ser visitado na forma\c{c}\~{a}o
da floresta de busca.  Uma maneira alternativa de rotular (reordenar) os
v\'{e}rtices de uma floresta de busca \'{e} a {\bf ordem de retorno}:
$Or( )$, \'{e} a ordem em que verificamos j\'{a} termos visitado todos
os filhos de um v\'{e}rtice e retornamos ao seu pai na \'{a}rvore (ou
terminamos a \'{a}rvore se se tratar de uma raiz). 
 \index{Busca!ordem de visita\c{c}\~{a}o}  
 \index{Busca!ordem de retorno}  

Estando os v\'{e}rtices de um grafo bem ordenados por algum
crit\'{e}rio, por exemplo pelos seus \'{\i}ndices, a floresta de busca
em profundidade can\^{o}nica \'{e} aquela em que tomamos os
v\'{e}rtices, tanto para ra\'{\i}zes de novas \'{a}rvores quando para
visita\c{c}\~{a}o dentro de uma busca, na ordem estabelecida por este
crit\'{e}rio. No exemplo 7, a primeira \'{e} a floresta can\^{o}nica.

Exemplo 7:

Um grafo $G$, e v\'{a}rias das poss\'{\i}veis florestas de busca que o 
cobrem, est\~{a}o dadas na figura seguinte. Apresentamos estas florestas
com os v\'{e}rtices numerados primeiro pela ordem de visita\c{c}\~{a}o, 
depois pela ordem de retorno. Indicamos com um circunflexo as raizes 
da busca.   

$$
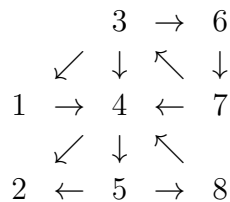
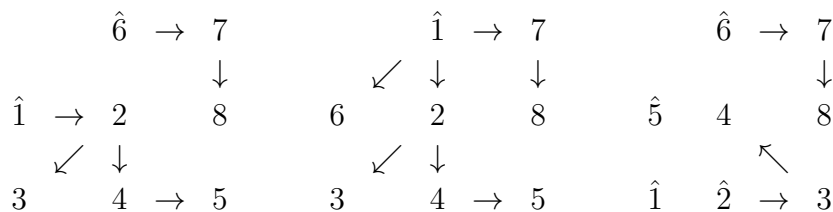
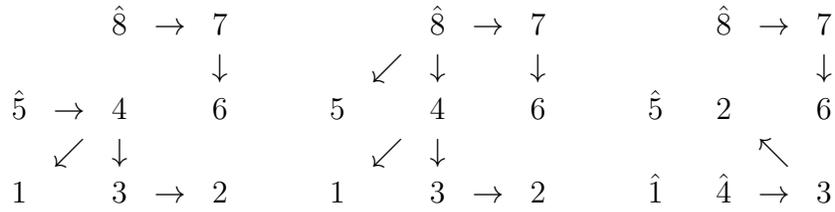
\begin{array}{ccccc}
 & & 3 & \rightarrow & 6 \\ 
 & \swarrow & \downarrow & \nwarrow & \downarrow \\ 
1 & \rightarrow & 4 & \leftarrow & 7 \\ 
 & \swarrow & \downarrow & \nwarrow & \\ 
2 & \leftarrow & 5 & \rightarrow & 8 
\end{array} $$
$$
\begin{array}{ccccc}
 & & \hat 6 & \rightarrow & 7 \\ 
 & & & & \downarrow \\ 
 \hat 1 & \rightarrow & 2 & & 8 \\ 
 & \swarrow & \downarrow &  & \\ 
3 & & 4 & \rightarrow & 5  
\end{array} \hspace{1cm}   
\begin{array}{ccccc}
 & & \hat 1 & \rightarrow & 7 \\ 
 & \swarrow & \downarrow & & \downarrow \\ 
 6 & & 2 & & 8 \\ 
 & \swarrow & \downarrow & & \\ 
 3 & & 4 & \rightarrow & 5 
\end{array} \hspace{1cm}   
\begin{array}{ccccc}
 & & \hat 6 & \rightarrow & 7 \\ 
 & & & & \downarrow \\ 
 \hat 5 & & 4 & & 8 \\ 
 & & & \nwarrow & \\ 
 \hat 1 & & \hat 2 & \rightarrow & 3  
\end{array} $$
$$
\begin{array}{ccccc}
 & & \hat 8 & \rightarrow & 7 \\ 
 & & & & \downarrow \\ 
 \hat 5 & \rightarrow & 4 & & 6 \\ 
 & \swarrow & \downarrow &  & \\ 
 1 & & 3 & \rightarrow & 2  
\end{array} \hspace{1cm}   
\begin{array}{ccccc}
 & & \hat 8 & \rightarrow & 7 \\ 
 & \swarrow & \downarrow & & \downarrow \\ 
 5 & & 4 & & 6 \\ 
 & \swarrow & \downarrow & & \\ 
 1 & & 3 & \rightarrow & 2 
\end{array} \hspace{1cm}   
\begin{array}{ccccc}
 & & \hat 8 & \rightarrow & 7 \\ 
 & & & & \downarrow \\ 
 \hat 5 & & 2 & & 6 \\ 
 & & & \nwarrow & \\ 
 \hat 1 & & \hat 4 & \rightarrow & 3  
\end{array} $$

Descrevemos agora o {\bf algoritmo de Tarjan} para determina\c{c}\~{a}o 
das componentes fortemente conexas de um grafo $G$. 
 \index{Algoritmo!Tarjan} 
 \begin{enumerate}
\item 
Considerando a boa ordem dos \'{\i}ndices, construa a floresta de
busca can\^{o}nica em $G$, marcando seus v\'{e}rtices pela ordem de
retorno;
\item
Considere $G^{-1}$ com os v\'{e}rtices reordenados, isto \'{e} 
rotulados, na ordem inversa da ordem de retorno estabelecida no passo 1;
\item
Considerando a boa ordem estabelecida no passo 2, construa
a floresta de busca can\^{o}nica em $G^{-1}$.
\end{enumerate}

\begin{teo}[Tarjan]: 
 Cada \'{a}rvore em $G^{-1}$ constru\'{\i}da no passo 3 do algoritmo de
Tarjan cobre os v\'{e}rtices de exatamente uma componente fortemente
conexa de $G$.  Mais ainda, a ordem de visita\c{c}\~{a}o (bem como a
ordem de retorno), na floresta can\^{o}nica do passo 3 do Algoritmo de
Tarjan, \'{e} um reordenamento coerente dos v\'{e}rtices de $G$. 
 \end{teo}

Demonstra\c{c}\~{a}o:

Se $v,w\in V$ est\~{a}o numa mesma componente de $G$ ent\~{a}o certamente
pertencem a uma mesma \'{a}rvore em $G^{-1}$ (bem como em $G$). 
Se $x$ \'{e} um v\'{e}rtice de $G$, denotaremos o n\'{u}mero que marca o
v\'{e}rtice $x$, pela ordem de retorno em $G$, por $Or(x)$.
Se $w$ \'{e} um v\'{e}rtice de uma \'{a}rvore de raiz $v$, em 
$G^{-1}$, ent\~{a}o:
\begin{enumerate}
\item $v$ descende de $w$, em $G$, pois $w$ descende de $v$ em $G^{-1}$.
\item $w$ descende de $v$, em $G$. 
\end{enumerate}

Para justificar a segunda afirma\c{c}\~{a}o notemos que, 
por constru\c{c}\~{a}o  $Or(v) > Or(w)$. Isto significa que
ou $w$ foi visitado durante a busca em profundidade
a partir de $v$, ou $w$ foi visitado antes de $v$.  A primeira
hip\'{o}tese implica em (2) enquanto a segunda hip\'{o}tese \'{e}
imposs\'{\i}vel pois, como por (1) $v$ descende de $w$, jamais
poder\'{\i}amos ter deixado $v$ ``n\~{a}o visitado" ap\'{o}s concluir uma
busca em profundidade que visitasse $w$.  Q.E.D.

Exemplo 8:

A figura seguinte apresenta os passos na determina\c{c}\~{a}o pelo
algoritmo de Tarjan, das CFCs do exemplo 7, e nos d\'{a} um
reordenamento coerente de seus v\'{e}rtices. 
A floresta de BEP em $G$ est\'{a} rotulada pela ordem de retorno, 
e a floresta em $G^{-1}$ pela ordem de visita\c{c}\~{a}o. 
O reordenamento coerente dos v\'{e}rtices do grafo original obtido
neste exemplo \'{e} $q=[ 3,6,7\mid 1\mid 4,5,8\mid 2 ]$.

$$
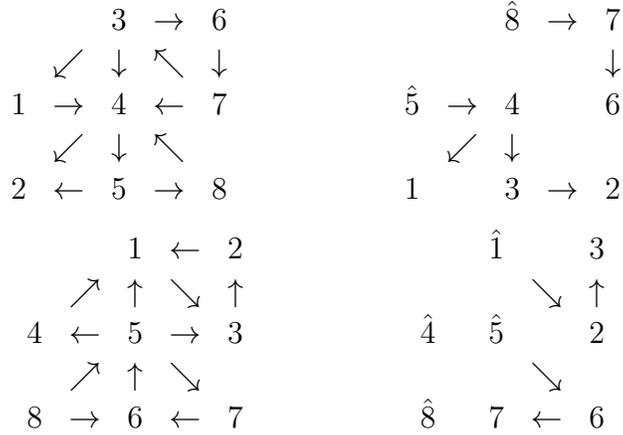
\begin{array}{ccccc}
 & & 3 & \rightarrow & 6 \\ 
 & \swarrow & \downarrow & \nwarrow & \downarrow \\ 
 1 & \rightarrow & 4 & \leftarrow & 7 \\ 
 & \swarrow & \downarrow & \nwarrow & \\ 
 2 & \leftarrow & 5 & \rightarrow & 8 
\end{array} \hspace{2cm} 
\begin{array}{ccccc}
 & & \hat 8 & \rightarrow & 7 \\ 
 & & & & \downarrow \\ 
 \hat 5 & \rightarrow & 4 & & 6 \\ 
 & \swarrow & \downarrow & & \\ 
 1 & & 3 & \rightarrow & 2  
\end{array} $$
$$
\begin{array}{ccccc}
 & & 1 & \leftarrow & 2 \\ 
 & \nearrow & \uparrow & \searrow & \uparrow \\ 
 4 & \leftarrow & 5 & \rightarrow & 3 \\ 
 & \nearrow & \uparrow & \searrow & \\ 
 8 & \rightarrow & 6 & \leftarrow & 7 
\end{array} \hspace{2cm} 
\begin{array}{ccccc}
 & & \hat 1 & & 3 \\ 
 & & & \searrow & \uparrow \\ 
 \hat 4 &  & \hat 5 & & 2 \\ 
 & & & \searrow & \\ 
 \hat 8 & & 7 & \leftarrow & 6  
\end{array} $$

\section{Grafos Sim\'{e}tricos e Casamentos}

Um grafo $G=(N,\Gamma )$ \'{e} dito {\bf sim\'{e}trico} sse 
$\forall i,j\in N,\ j\in \Gamma (i) \Rightarrow i\in \Gamma (j)$.  
\'{E} usual representarmos um grafo sim\'{e}trico substituindo cada par
de arestas, $(i,j)$ e $(j,i)$, por uma aresta sem orienta\c{c}\~{a}o ou
{\bf lado}, $\{i,j\}$.  Podemos definir um grafo sim\'{e}trico atrav\'{e}s de
seus v\'{e}rtices e dos seus lados, $G=(N,E)$, onde cada lado \'{e} um
conjunto de dois v\'{e}rtices e 
$\{i,j\}\in E \Leftrightarrow i\in \Gamma (j)$. 
 \index{Grafo!sim\'{e}trico} 
 \index{Casamento} 

Um $m$-{\bf casamento}, num grafo sim\'{e}trico, \'{e} um conjunto de $m$
lados onde s\~{a}o distintos todos os v\'{e}rtices;
$M=\{\{u_1,u_2\},\{u_3,u_4\},\ldots \{u_{2m-1},u_{2m}\}\}$.
Os v\'{e}rtices em $M$ dizem-se casados, os de mesmo lado dizem-se
companheiros, e os v\'{e}rtices fora do casamento dizem-se solteiros.
Um $m$-casamento \'{e} m\'{a}ximo em $G$ se n\~{a}o houver em $G$ um 
$m+1$ casamento e \'{e} perfeito se n\~{a}o deixar nenhum v\'{e}rtice
solteiro. 

Dado um grafo sim\'{e}trico $G=(V,E)$ e nele um $m$-casamento, $M$, 
um {\bf caminho $M$-alternado} \'{e} um caminho
$C=(\{w_0,w_1\},\{w_1,w_2\},\ldots \{w_{k-1},w_k\})$
que tem lados, alternadamente, dentro e fora do casamento.
Um caminho $M$-alternado diz-se um {\bf caminho de aumento} se
come\c{c}a e termina em v\'{e}rtices solteiros. 
 \index{Caminho!de aumento} 
 \index{Caminho!alternado}

\begin{teo}[Berge]: 
Dado um grafo sim\'{e}trico, $G=(V,E)$ e nele um $m$-casamento, 
$M=\{\{u_1,u_2\},\{u_3,u_4\},\ldots \{u_{2m},u_{2m}\}\}$, 
este casamento \'{e} m\'{a}ximo se n\~{a}o houver caminho de aumento.
\end{teo}
 \index{Teorema!Berge} 

Demonstra\c{c}\~{a}o:

Suponha que existe um caminho de aumento \\  
$C=(\{w_0,w_1\},\{w_1,w_2\},\ldots \{w_{2k-1},w_{2k+1}\})$
ent\~{a}o
\begin{eqnarray*} 
 M' = M \oplus C &=&  
   M - \{\{w_1,w_2\},\{w_3,w_4\},\ldots \{w_{2k-1},w_{2k}\}\} \\ 
 & &  + \{\{w_0,w_1\},\{w_2,w_3\},\ldots \{w_{2k},w_{2k+1}\}\}
\end{eqnarray*}
\'{e} um $m+1$-casamento.

Se, por outro lado, $M$ n\~{a}o \'{e} m\'{a}ximo, seja $M$'um
m+1-casamento e considere o grafo auxiliar $H$ de lados 
$E_H=M'\oplus M$, sendo $V_H=V_G$.  
Cada v\'{e}rtice de $H$ \'{e} isolado, ou pertence a algum lado de $E_H$
e no m\'{a}ximo a dois lados, um de $M$ e outro de $M'$. 
Cada componente de $H$ \'{e} pois, ou um v\'{e}rtice isolado, ou um
ciclo de lados alternadamente em $M$ e $M'$, ou um caminho de lados 
alternadamente em cada um dos casamentos. 
Como $H$ tem mais lados de $M$'que de $M$, h\'{a} ao menos uma 
componente de tipo caminho que come\c{c}a e termina com lados de $M'$. 
Este caminho \'{e} $M$-alternado e seus v\'{e}rtices extremos s\~{a}o,por
constru\c{c}\~{a}o, solteiros em $M$. 
Temos assim um caminho de aumento para $M$.
Q.E.D.

Exemplo 9: 

Apresentamos agora: Na primeira linha: $M$: um $m$-casamento num grafo,
$C$: um caminho de aumento, e $M'$: o $m+1$-casamento $M'=M\oplus C$. 
Na segunda linha: $M$: um $m$-casamento, $M'$: um $m+1$-casamento, e
$H$: o grafo $H=M'\oplus M$. 

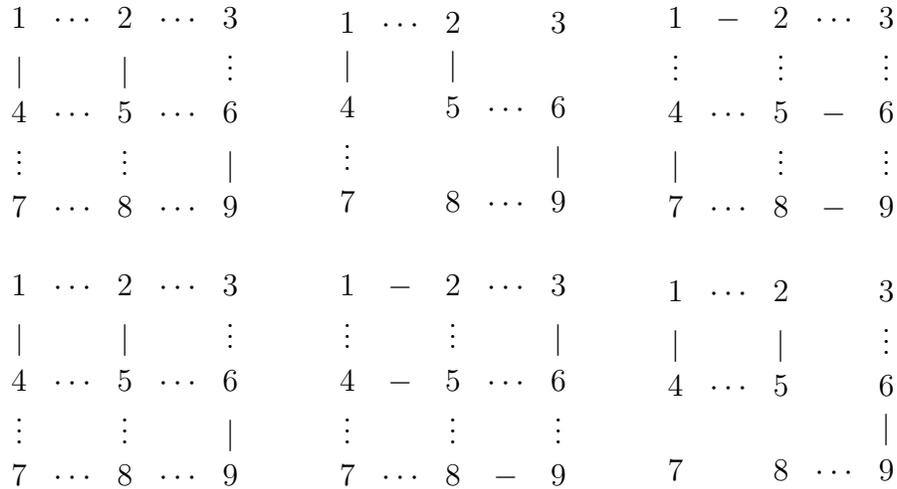
$$
\begin{array}{ccccc} 
 1 & \cdots & 2 & \cdots & 3 \\ 
 | & & | & & \vdots \\ 
 4 & \cdots & 5 & \cdots & 6 \\ 
 \vdots & & \vdots & & | \\ 
 7 & \cdots & 8 & \cdots & 9 
\end{array} \hspace{1cm} 
\begin{array}{ccccc} 
 1 & \cdots & 2 & & 3 \\ 
 | & & | & &  \\ 
 4 & & 5 & \cdots & 6 \\ 
 \vdots & & & & | \\ 
 7 & & 8 & \cdots & 9 
\end{array} \hspace{1cm} 
\begin{array}{ccccc} 
 1 & - & 2 & \cdots & 3 \\ 
 \vdots & & \vdots & & \vdots \\ 
 4 & \cdots & 5 & - & 6 \\ 
 | & & \vdots & & \vdots \\ 
 7 & \cdots & 8 & - & 9 
\end{array}$$ 

$$
\begin{array}{ccccc} 
 1 & \cdots & 2 & \cdots & 3 \\ 
 | & & | & & \vdots \\ 
 4 & \cdots & 5 & \cdots & 6 \\ 
 \vdots & & \vdots & & | \\ 
 7 & \cdots & 8 & \cdots & 9 
\end{array} \hspace{1cm} 
\begin{array}{ccccc} 
 1 & - & 2 & \cdots & 3 \\ 
 \vdots & & \vdots & & | \\ 
 4 & - & 5 & \cdots & 6 \\ 
 \vdots & & \vdots & & \vdots \\ 
 7 & \cdots & 8 & - & 9 
\end{array} \hspace{1cm} 
\begin{array}{ccccc} 
 1 & \cdots & 2 & & 3 \\ 
 | & & | & & \vdots \\ 
 4 & \cdots & 5 & & 6 \\ 
   & & & & | \\ 
 7 & & 8 & \cdots & 9 
\end{array}$$

Um grafo $G=(V,\Gamma )$ \'{e} {\bf bipartido} se houver uma 
biparti\c{c}\~{a}o  de seus v\'{e}rtices tal que todos os 
lados tenham um v\'{e}rtice em cada peda\c{c}o da 
biparti\c{c}\~{a}o, i.e. se for poss\'{\i}vel encontrar 
conjuntos  \\ 
$X,Y \mid V=X\cup Y\ \wedge \ X\cap Y=\emptyset \ \wedge \ 
\forall e\in E, e=\{x,y\}\ x\in X, y\in Y$.   
 \index{Grafo!bipartido} 

\begin{teo}[Hall]: 
Dado $G=(V,\Gamma )$ um grafo sim\'{e}trico e bipartido, 
com biparti\c{c}\~{a}o  $V=X\cup Y$, existe em $G$ um casamento que casa
todos os v\'{e}rtices de $X$ sse vale a condi\c{c}\~{a}o de Hall:
$\forall S\subseteq X,  \# \Gamma(S) \geq \# S$.
isto \'{e}, qualquer subconjunto de $X$ tem, coletivamente, ao menos
tantos filhos quanto elementos.
\end{teo} 
 \index{Teorema!Hall} 

Demonstra\c{c}\~{a}o

Se existe um casamento $M$ que casa todos os v\'{e}rtices em $X$,
ent\~{a}o os lados do casamento garantem ao menos a igualdade na
condi\c{c}\~{a}o de Hall.  Por outro lado, se houver um casamento
m\'{a}ximo, $M$, que deixe algum $x\in X$ solteiro, exibiremos um
$S\subseteq X$ que viola a condi\c{c}\~{a}o de Hall. 

Seja $Z$ o conjunto dos v\'{e}rtices conexos a $x$ por caminhos
M-alternados.  Pelo teorema de Berge sabemos que n\~{a}o h\'{a} solteiro
em $Z$, pois caso contr\'{a}rio ter\'{\i}amos um caminho de aumento e
$M$ n\~{a}o seria m\'{a}ximo.  Assim, um caminho M-alternado maximal,
isto \'{e}, que n\~{a}o pode ser continuado, que come\c{c}a no solteiro
$x\in S\subseteq X$, deve necessariamente terminar, com um n\'{u}mero
par de lados, num casado $x'\in S\subseteq X$. 

Considerando, pois, $S=(Z\cap X)+x$ e $T= Z\cap Y$, temos que 
$\# S= \# T +1$ e $\Gamma (S)=T$, donde 
$\# \Gamma (S)= \# T = \# S -1 < \# S$,
mostrando que $S$ viola a condi\c{c}\~{a}o de Hall.  Q.E.D.

\section{O Algoritmo H\'{u}ngaro}

Dado um grafo sim\'{e}trico $G=(V,E)$, bipartido em conjuntos
de mesma cardinalidade $X$ e $Y$, o algoritmo h\'{u}ngaro fornece um
casamento perfeito ou encontra um conjunto $S\subset X$ que
viola a condi\c{c}\~{a}o de Hall.
 \index{Algoritmo!H\'{u}ngaro} 

Seja $M$ um casamento que deixa $x\in X$ solteiro.
 Uma \'{a}rvore de raiz $x$, $H=(V_H,E_H)$, subgrafo de $G$,
\'{e} M-alternada se for uma \'{a}rvore onde qualquer caminho partindo
de $x$, $C=(x,w_1,w_2,\ldots )$ for M-alternado. 
 Observemos que se uma dada \'{a}rvore M-alternada de raiz $x$ tiver uma
folha solteira, ent\~{a}o o caminho $C$ que vai da raiz a esta folha
\'{e} um caminho de aumento e $M'=M\oplus C$ \'{e} um novo casamento
onde, al\'{e}m de todos os v\'{e}rtices casados em $M$, $x$ tamb\'{e}m
\'{e} casado. 

Se, todavia, encontrarmos uma \'{a}rvore M-alternada de raiz solteira,
onde todas as folhas s\~{a}o casadas e que seja m\'{a}xima, isto \'{e},
tal que n\~{a}o exista nenhum v\'{e}rtice em $V_G-V_H$, adjacente a uma
folha de $H$ teremos encontrado um $S\subset X$ que viola a 
condi\c{c}\~{a}o de Hall, $S=V_H\cap X$ e $T=V_H\cap Y = \Gamma (S)$.

O algoritmo H\'{u}ngaro faz o seguinte: dado um grafo bipartido e
sim\'{e}trico com um casamento que deixa $x\in X$ solteiro, gerar a
\'{a}rvore can\^{o}nica de busca em profundidade M-alternada, at\'{e}
encontrar um caminho de aumento, isto \'{e}, uma folha solteira.  Caso
isto n\~{a}o seja poss\'{\i}vel constatarmos a viola\c{c}\~{a}o da
condi\c{c}\~{a}o de Hall.

Exemplo 10:

No exemplo 10 temos dois exemplos de grafo e casamento, nos quais est\'{a}
assinalado um v\'{e}rtice $w$ solteiro e sem mais nenhum pretendente
dispon\'{\i}vel.  Temos as respectivas \'{a}rvores M-alternadas,
marcadas pela ordem de visita\c{c}\~{a}o, os caminhos de aumento, $C$, e
os novos casamentos, $M'=M\oplus C$, que casam o v\'{e}rtice $w$. 

O algoritmo H\'{u}ngaro, bem como qualquer algoritmo casamenteiro
conhecido, n\~{a}o \'{e} linear; Portanto vale a pena procurar
heur\'{\i}sticas mais eficientes para a procura de um casamento
perfeito.  Dado um grafo e um casamento, os pretendentes de um
v\'{e}rtice solteiro s\~{a}o os solteiros em seu conjunto de adjac\^{e}ncia,
e o peso de um solteiro, $\rho (x)$, \'{e} seu numero de pretendentes. 
A id\'{e}ia subjacente em todas estas heur\'{\i}sticas \'{e} a de casar
primeiro os v\'{e}rtices mais exigentes, i.e. de menor peso, de modo a
preservar o m\'{a}ximo de pretendentes para os v\'{e}rtices ainda
solteiros. 

Na {\bf Heur\'{\i}stica do M\'{\i}nimo Simples}, HMS, dado um m-casamento
imperfeito, procuraremos obter um m+1-casamento adicionando ao casamento
um lado que una v\'{e}rtices ainda solteiros.  Escolhemos este lado, de
modo a casar um v\'{e}rtice de m\'{\i}nimo peso, entre os ainda
solteiros, com um de m\'{\i}nimo peso dentre seus pretendentes.  Sempre
que houver um v\'{e}rtice isolado, i.e.  sem pretendentes ou de peso
zero, procuraremos cas\'{a}-lo pelo algoritmo H\'{u}ngaro. 
 Na {\bf Heur\'{\i}stica do M\'{\i}nimo Par}, HMP, escolhemos um lado a 
ser acrescentado ao casamento onde um dos v\'{e}rtices \'{e} de peso
m\'{\i}nimo, e o outro tenha peso m\'{\i}nimo dentre todos os
pretendentes a v\'{e}rtices de peso m\'{\i}nimo.  Como a HMP pode ser
bem mais custosa podemos, por exemplo, utilizar a HMS enquanto o peso
m\'{\i}nimo dos v\'{e}rtices ainda solteiros for grande, i.e.  acima de
um dado limite, e a HMP caso contr\'{a}rio.  

Exemplo 11: 

O exemplo seguinte ilustra aaplica\c{c}\~{a}o do algoritmo h\'{u}ngaro
com a heur\'{\i}stica min-min.  No exemplo indicamos o peso, ou
n\'{u}mero de pretendentes, de cada v\'{e}rtice.  Indicamos ainda, a
cada passo, o v\'{e}rtice de n\'{\i}nimo peso a ser casado, bem como, em
it\'{a}lico, seu pretendente de m\'{\i}nimo peso.

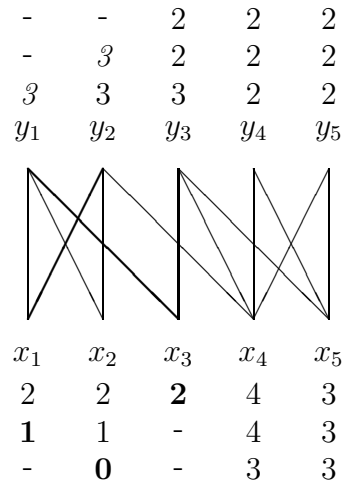
\mbox{}\\ 

\begin{center} 
\unitlength=1.00mm
\special{em:linewidth 0.4pt}
\linethickness{0.4pt}
\begin{picture}(50.00,70.00)
\put(10.00,55.00){\makebox(0,0)[cc]{$y_1$}}
\put(20.00,55.00){\makebox(0,0)[cc]{$y_2$}}
\put(30.00,55.00){\makebox(0,0)[cc]{$y_3$}}
\put(40.00,55.00){\makebox(0,0)[cc]{$y_4$}}
\put(50.00,55.00){\makebox(0,0)[cc]{$y_5$}}
\put(10.00,25.00){\makebox(0,0)[cc]{$x_1$}}
\put(30.00,25.00){\makebox(0,0)[cc]{$x_3$}}
\put(40.00,25.00){\makebox(0,0)[cc]{$x_4$}}
\put(50.00,25.00){\makebox(0,0)[cc]{$x_5$}}
\put(10.00,60.00){\makebox(0,0)[cc]{{\it 3}}}
\put(20.00,60.00){\makebox(0,0)[cc]{3}}
\put(30.00,60.00){\makebox(0,0)[cc]{3}}
\put(40.00,60.00){\makebox(0,0)[cc]{2}}
\put(50.00,60.00){\makebox(0,0)[cc]{2}}
\put(10.00,65.00){\makebox(0,0)[cc]{-}}
\put(20.00,65.00){\makebox(0,0)[cc]{{\it 3}}}
\put(30.00,65.00){\makebox(0,0)[cc]{2}}
\put(40.00,65.00){\makebox(0,0)[cc]{2}}
\put(50.00,65.00){\makebox(0,0)[cc]{2}}
\put(10.00,70.00){\makebox(0,0)[cc]{-}}
\put(20.00,70.00){\makebox(0,0)[cc]{-}}
\put(30.00,70.00){\makebox(0,0)[cc]{2}}
\put(40.00,70.00){\makebox(0,0)[cc]{2}}
\put(50.00,70.00){\makebox(0,0)[cc]{2}}
\put(10.00,20.00){\makebox(0,0)[cc]{2}}
\put(20.00,20.00){\makebox(0,0)[cc]{2}}
\put(30.00,20.00){\makebox(0,0)[cc]{{\bf 2}}}
\put(40.00,20.00){\makebox(0,0)[cc]{4}}
\put(50.00,20.00){\makebox(0,0)[cc]{3}}
\put(10.00,15.00){\makebox(0,0)[cc]{{\bf 1}}}
\put(20.00,15.00){\makebox(0,0)[cc]{1}}
\put(30.00,15.00){\makebox(0,0)[cc]{-}}
\put(40.00,15.00){\makebox(0,0)[cc]{4}}
\put(50.00,15.00){\makebox(0,0)[cc]{3}}
\put(10.00,10.00){\makebox(0,0)[cc]{-}}
\put(20.00,10.00){\makebox(0,0)[cc]{{\bf 0}}}
\put(30.00,10.00){\makebox(0,0)[cc]{-}}
\put(40.00,10.00){\makebox(0,0)[cc]{3}}
\put(50.00,10.00){\makebox(0,0)[cc]{3}}
\put(10.00,50.00){\line(0,-1){20.00}}
\put(10.00,50.00){\line(1,-2){10.00}}
\thicklines 
\put(10.00,50.00){\line(1,-1){20.00}}
\put(20.00,50.00){\line(-1,-2){10.00}}
\thinlines 
\put(20.00,50.00){\line(0,0){0.00}}
\put(20.00,50.00){\line(0,-1){20.00}}
\put(20.00,25.00){\makebox(0,0)[cc]{$x_2$}}
\put(20.00,50.00){\line(0,-1){20.00}}
\put(30.00,50.00){\line(0,-1){20.00}}
\put(30.00,50.00){\line(1,-2){10.00}}
\put(30.00,50.00){\line(1,-1){20.00}}
\put(40.00,50.00){\line(0,-1){20.00}}
\put(40.00,50.00){\line(1,-2){10.00}}
\put(50.00,50.00){\line(0,-1){20.00}}
\put(50.00,50.00){\line(-1,-2){10.00}}
\put(20.00,50.00){\line(1,-1){20.00}}
\end{picture}
\end{center}

A \'{a}rvore M-alternada com raiz em $x_2$, que ficou 
isolado, nos fornece o caminho de aumento 
 $$ c = \ \ x_2 \ -\ y_1\ {\bf -}\ x_3\ -\ y_3 $$ 
de modo que podemos recome\c{c}ar a heurr\'{\i}stica 
em $M \oplus c$:

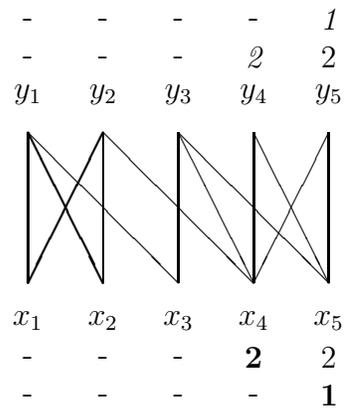
\mbox{}\\ 

\begin{center} 
\unitlength=1.00mm
\special{em:linewidth 0.4pt}
\linethickness{0.4pt}
\begin{picture}(50.00,60.00)
\put(10.00,50.00){\makebox(0,0)[cc]{$y_1$}}
\put(20.00,50.00){\makebox(0,0)[cc]{$y_2$}}
\put(30.00,50.00){\makebox(0,0)[cc]{$y_3$}}
\put(40.00,50.00){\makebox(0,0)[cc]{$y_4$}}
\put(50.00,50.00){\makebox(0,0)[cc]{$y_5$}}
\put(10.00,20.00){\makebox(0,0)[cc]{$x_1$}}
\put(30.00,20.00){\makebox(0,0)[cc]{$x_3$}}
\put(40.00,20.00){\makebox(0,0)[cc]{$x_4$}}
\put(50.00,20.00){\makebox(0,0)[cc]{$x_5$}}
\put(10.00,55.00){\makebox(0,0)[cc]{-}}
\put(20.00,55.00){\makebox(0,0)[cc]{-}}
\put(30.00,55.00){\makebox(0,0)[cc]{-}}
\put(40.00,55.00){\makebox(0,0)[cc]{{\it 2}}}
\put(50.00,55.00){\makebox(0,0)[cc]{2}}
\put(10.00,60.00){\makebox(0,0)[cc]{-}}
\put(20.00,60.00){\makebox(0,0)[cc]{-}}
\put(30.00,60.00){\makebox(0,0)[cc]{-}}
\put(40.00,60.00){\makebox(0,0)[cc]{-}}
\put(50.00,60.00){\makebox(0,0)[cc]{{\it 1}}}
\put(10.00,15.00){\makebox(0,0)[cc]{-}}
\put(20.00,15.00){\makebox(0,0)[cc]{-}}
\put(30.00,15.00){\makebox(0,0)[cc]{-}}
\put(40.00,15.00){\makebox(0,0)[cc]{{\bf 2}}}
\put(50.00,15.00){\makebox(0,0)[cc]{2}}
\put(10.00,10.00){\makebox(0,0)[cc]{-}}
\put(20.00,10.00){\makebox(0,0)[cc]{-}}
\put(30.00,10.00){\makebox(0,0)[cc]{-}}
\put(40.00,10.00){\makebox(0,0)[cc]{-}}
\put(50.00,10.00){\makebox(0,0)[cc]{{\bf 1}}}
\thicklines 
\put(10.00,45.00){\line(1,-2){10.00}}
\put(20.00,45.00){\line(-1,-2){10.00}}
\put(30.00,45.00){\line(0,-1){20.00}}
\put(40.00,45.00){\line(0,-1){20.00}}
\put(50.00,45.00){\line(0,-1){20.00}}
\thinlines  
\put(10.00,45.00){\line(0,-1){20.00}}
\put(10.00,45.00){\line(1,-1){20.00}}
\put(20.00,45.00){\line(0,0){0.00}}
\put(20.00,45.00){\line(0,-1){20.00}}
\put(20.00,20.00){\makebox(0,0)[cc]{$x_2$}}
\put(20.00,45.00){\line(0,-1){20.00}}
\put(30.00,45.00){\line(1,-2){10.00}}
\put(30.00,45.00){\line(1,-1){20.00}}
\put(40.00,45.00){\line(1,-2){10.00}}
\put(50.00,45.00){\line(-1,-2){10.00}}
\put(20.00,45.00){\line(1,-1){20.00}}
\end{picture}
\end{center}

 \clearpage
 \clearpage 
\setcounter{chapter}{3} 
\chapter{ELIMINA\c{C}\~{A}O ASSIM\'{E}TRICA}
\begin{center}
{\LARGE Esparsidade na Fatora\c{c}\~{a}o LU}
\end{center} 

A maioria dos sistemas lineares encontrados na pr\'{a}tica,
especialmente os sistemas grandes, s\~{a}o esparsos, isto \'{e}, apenas
uma pequena parte dos coeficientes do sistema n\~{a}o s\~{a}o nulos. 
Ao tratar um sistema esparso, devemos tentar manter a esparsidade do
sistema ao longo das tranforma\c{c}\~{o}es necess\'{a}rias \`{a} sua
solu\c{c}\~{a}o.  Tal medida visa, primordialmente, minimizar o
n\'{u}mero de opera\c{c}\~{o}es aritm\'{e}ticas a serem realizadas, na
pr\'{o}pria transforma\c{c}\~{a}o e na solu\c{c}\~{a}o do sistema. 
 Economia de mem\'{o}ria tamb\'{e}m \'{e} um fator importante. 
 \index{Esparsidade}

\section{Preenchimento Local}

No m\'{e}todo de Gauss as tranforma\c{c}\~{o}es 
 $^kA \rightarrow \ ^{k+1}A$ podem diminuir o n\'{u}mero de elementos
nulos.  Se tal ocorrer, dizemos que houve {\bf preenchimento} de algumas
posi\c{c}\~{o}es.  Queremos escolher o elemento piv\^{o} em $^kA$ de
modo a minimizar este preenchimento. 
 \index{Preenchimento Local} 

\begin{teo}[Tewarson] 
Seja
$$
\left[ \begin{array}{ccccc} 
 ^kA_1^1 & & \ldots & & ^kA_1^n \\
  0 & \ddots & & & \\
  \vdots & 0 & ^kA_{k+1}^{k+1} & \ldots & ^kA_{k+1}^n \\
  & \vdots & \vdots & & \vdots \\ 
  0 & 0 & ^kA_n^{k+1} & \ldots & ^kA_n^n 
\end{array} \right] $$  
seja $^kB$, $(n-k)\times (n-k)$, a matriz Booleana associada a submatriz
das $n-k$ \'{u}ltimas linhas e colunas de $^kA$, i.e. 
$$ ^kB = B(A_{k+1:n}^{k+1:n}) ,\ \ \mbox{ou} \ \ \
  ^kB_p^q = 1 \Leftrightarrow \ ^kA_{k+p}^{k+q} \neq 0 \ ,$$
e seja
 $$ ^kG = \ ^kB (\, ^k{\bar B})' \ ^kB \ .$$ 
 onde o operador complemento aplicado \`{a} matriz $B$, $\bar B$, 
troca 0's por 1's e vice-versa. 

A escolha do piv\^{o} $^kA_{k+i}^{k+j} \neq 0$
implica no prenchimento de exatamente $^kG_i^j$ posi\c{c}\~{o}es.
\end{teo}
 \index{Teorema!Tewarson} 
 \index{Matriz!booleana} 

Demonstra\c{c}\~{a}o: 

Seja $^k\tilde A$ a matriz obtida de $^kA$ por permuta\c{c}\~{a}o da
$k+1$-\'{e}sima linha e coluna com, respectivamente, a $k+i$-\'{e}sima 
linha e a $k+j$-\'{e}sima coluna. 

Os novos elementos de $^{k+1}A$ ser\~{a}o, para $p,q\ = 2\ldots (n-k)$, 
\begin{eqnarray*} 
\lefteqn{ ^{k+1}A_{k+p}^{k+q} =} \\ 
 &=& ^k{\tilde A}_{k+p}^{k+q} - \ 
     ^{k+1}M_{k+p}^{k+1} \ ^k{\tilde A}_{k+1}^{k+q} \\ 
 &=& ^k{\tilde A}_{k+p}^{k+q} - \  
     ^k{\tilde A}_{k+p}^{k+1} \ ^k{\tilde A}_{k+1}^{k+q} \ 
     / \ ^k{\tilde A}_{k+1}^{k+1} \\ 
 &=& ^kA_{k+r}^{k+s} - \  
     ^kA_{k+r}^{k+j} \ ^kA_{k+i}^{k+s} \ 
     / \ ^kA_{k+i}^{k+j} 
\end{eqnarray*}
onde
$$r=p\ \mbox{se}\ p\neq i, \ \mbox{e}\ r=1\ \mbox{se}\ p=i 
\ ;\ \ \ \mbox{e}\ \ \ 
 s=q\ \mbox{se}\ q\neq j, \ \mbox{e}\ s=1\ \mbox{se}\ q=j \ .$$

Haver\'{a} preenchimento de uma posi\c{c}\~{a}o  em $^{k+1}A$ sempre que,\\ 
para $r,s = 1\ldots (n-k)$, $r\neq i$ e $s\neq j$, 
$$ ^kB_r^s =0 \ \ \wedge \ \ 
   ^kB_r^j =1 \ \ \wedge \ \ 
   ^kB_i^s =1 $$
e portanto o total de preenchimentos, correspondente ao piv\^{o} 
$(k+i,\ k+j)$, \'{e}:
\begin{eqnarray*} 
 \lefteqn{ \sum_{r=1,\ r\neq i}^{n-k} \sum_{s=1,\ s\neq j}^{n-k}
           \ ^k{\bar B}_r^s \ ^kB_r^j \ ^kB_i^s =}\\
 &=& \sum_{r,s=1}^{n-k} \ ^kB_r^j  ((\ ^k\bar B)')_s^r \ ^kB_i^s \\ 
 &=& (\ ^kB(\ ^k\bar B)'\: ^kB)_i^j = G_i^j 
\end{eqnarray*}
Na pen\'{u}ltima passagem usamos que os termos $r=i$ ou $s=j$ 
s\~{a}o todos nulos, \\ pois $B_p^q \bar B_p^q = 0$. QED.

\begin{obs}{\rm  
Uma aproxima\c{c}\~{a}o para a matriz $^kG$ \'{e} a 
{\bf matriz de Markowitz}: 
 $$^kF_i^j = (\ ^kB\ {\bf 1}\ ^kB )_i^j \ ,$$
onde denotamos por ${\bf 1}$ a matriz quadrada de 1's. 
 A aproxima\c{c}\~{a}o $F$ \'{e} exatamente o n\'{u}mero de multiplidores
n\~{a}o nulos vezes o n\'{u}mero de elementos n\~{a}o nulos, fora o
piv\^{o}, na linha piv\^{o}.  Esta aproxima\c{c}\~{a}o \'{e} bastante
boa se $^kB$ \'{e} muito esparsa.   
}\end{obs} 
 \index{Matriz!Markowitz}

Exemplo 1:

Dada a matriz $^0A$, cujos elementos n\~{a}o nulos est\~{a}o indicados
na matriz boleana associada $B(\ ^0A)$, indique uma escolha de piv\^{o}s
que minimize preeenchimentos locais.  Assuma que ao longo do processo de
triangulariza\c{c}\~{a}o n\~{a}o h\'{a} cancelamentos, i.e., que um
elemento n\~{a}o nulo uma vez preenchido, n\~{a}o volta a se anular. 

Tomando
$$ ^0B = B(\ ^0A) =\ ^0B = \left[ \begin{array}{ccccc} 
 {\bf 1} & 1 & 1 & 0 & 0 \\ 0 & 1 & 1 & 0 & 0 \\ 
 0 & 1 & 0 & 1 & 0 \\ 0 & 0 & 0 & 1 & 1 \\ 
 1 & 1 & 1 & 0 & 1 \end{array} \right] \ \ , \ \ \ 
 ^0G = \left[ \begin{array}{ccccc} 
 {\bf 0} & 3 & 1 & 5 & 3 \\ 0 & 1 & 0 & 3 & 2 \\ 
 2 & 3 & 3 & 1 & 2 \\ 3 & 6 & 5 & 1 & 1 \\ 
 1 & 6 & 3 & 6 & 3 \end{array} \right] $$
Assim, $arg\min_{(i,j)\: \mid \: ^0B_i^j=1} \: ^0G_i^j \ = \{(1,1),(2,3)\}$.

Escolhendo $(i,j)=(1,1)$, i.e. $^0A_1^1$ como piv\^{o}, temos o 
preenchimento de $^0G_1^1=0$ posi\c{c}\~{o}es em $^1A$,
$$ ^1B = B(\ ^1A+\ ^1M) = \left[ \begin{array}{ccccc} 
 1 & 1 & 1 & 0 & 0 \\ {\em 0} & 1 & {\bf 1} & 0 & 0 \\ 
 {\em 0} & 1 & 0 & 1 & 0 \\ {\em 0} & 0 & 0 & 1 & 1 \\ 
 {\em 1} & 1 & 1 & 0 & 1 \end{array} \right] \ \ , \ \ \ 
 ^1G = \left[ \begin{array}{cccc} 
 1 & {\bf 0} & 3 & 2 \\ 2 & 2 & 1 & 2 \\ 
 4 & 3 & 1 & 1 \\ 3 & 1 & 4 & 2  \end{array} \right] $$
Assim, $arg\min_{(i,j)\: \mid \: ^1B_i^j=1} \: ^1G_i^j \ = \{(1,2)\}$.

Escolhendo $(i,j)=(1,2)$, i.e. $^1A_2^3$ como piv\^{o}, temos o 
preenchimento de $^1G_1^2=0$ posi\c{c}\~{o}es em $^2A$,
$$ ^2B = B(\ ^2A+\ ^2M) = \left[ \begin{array}{ccccc} 
 1 & 1 & 1 & 0 & 0 \\ {\em 0} & 1 & 1 & 0 & 0 \\ 
 {\em 0} & {\em 0} & {\bf 1} & 1 & 0 \\ {\em 0} & {\em 0} & 0 & 1 & 1 \\ 
 {\em 1} & {\em 1} & 1 & 0 & 1 \end{array} \right] \ \ , \ \ \ 
 ^1G = \left[ \begin{array}{ccc} 
 {\bf 1} & 1 & 2 \\ 2 & 1 & 1 \\ 1 & 2 & 1 \end{array} \right] $$
Assim, $arg\min_{(i,j)\: \mid \: ^2B_i^j=1} \: ^2G_i^j \ 
 = \{(1,1),(1,2),(2,2),(2,3),(3,1),(3,3)\}$.

Escolhendo $(i,j)=(1,1)$, i.e. $^2A_3^3$ como piv\^{o}, temos o 
preenchimento de $^2G_1^1=1$ posi\c{c}\~{o}es em $^3A$,
$$ ^3B = B(\ ^3A+\ ^3M) = \left[ \begin{array}{ccccc} 
 1 & 1 & 1 & 0 & 0 \\ {\em 0} & 1 & 1 & 0 & 0 \\ 
 {\em 0} & {\em 0} & 1 & 1 & 0 \\ {\em 0} & {\em 0} & {\em 0} & 1 & 1 \\ 
 {\em 1} & {\em 1} & {\em 1} & {\bf 1} & 1 \end{array} \right] \ \ , \ \ \ 
 ^3B = \left[ \begin{array}{cc} 1 & 1 \\ {\bf 1} & 1 
  \end{array} \right] $$

Como $^3B={\bf 1}$, \'{e} obvio que n\~{a}o haver\'{a} mais
preenchimento, quaisquer que que sejam os piv\^{o}s doravante
selecionados.  Tomando, por exemplo, o piv\^{o} $^3A_5^4$, 
temos finalmente temos a fatora\c{c}\~{a}o $\tilde A = PAQ = LU$, com
$$ B(\tilde A) = \left[ \begin{array}{ccccc} 
 1 & 1 & 1 & 0 & 0 \\ 0 & 1 & 1 & 0 & 0 \\ 
 0 & 0 & 1 & 1 & 0 \\ 1 & 1 & 1 & 0 & 1 \\ 
 0 & 0 & 0 & 1 & 1 \end{array} \right] \ \ , \ \ \ 
 B(M+U) = \left[ \begin{array}{ccccc} 
 1 & 1 & 1 & 0 & 0 \\ {\em 0} & 1 & 1 & 0 & 0 \\ 
 {\em 0} & {\em 0} & 1 & 1 & 0 \\ {\em 1} & {\em 1} & {\em 1} & 1 & 1 \\ 
 {\em 0} & {\em 0} & {\em 0} & {\em 1} & 1 \end{array} \right] $$ 
sendo os vetores de indices permutados
$^4p=[1,2,3,5,4]$ e $^4q=[1,3,2,4,5]$.
 Neste processo de triangulazira\c{c}\~{a}o, tivemos assim o preenchimento
de apenas 1 posi\c{c}\~{a}o, em $^3A_5^4$.  //

\begin{obs}{\rm 
Para n\~{a}o fazer a sele\c{c}\~{a}o de piv\^{o}s apenas pelos
crit\'{e}rios para minimiza\c{c}\~{a}o de preenchimento local,
desprezando assim completamente os apectos de estabilidade num\'{e}rica,
podemos ``vetar'' escolhas de piv\^{o}s muito pequenos, i.\'{e}., que
impliquem no uso de multiplicadores muito grandes,
$|M_i^j| \geq multmax \gg 1$.  Se nosso sistema for bem equilibrado, 
podemos usar como crit\'{e}rio de veto o tamanho do proprio piv\^{o}, 
$|U_i^i|\leq pivomin \ll 1$. 
}\end{obs} 
 \index{Pivoteamento}

\section{Pr\'{e}-Posicionamento de Piv\^{o}s}

Minimizar o preenchimento local \'{e} um processo custoso, tanto no
n\'{u}mero de opera\c{c}\~{o}es envolvidas no c\'{a}lculo de $^kG$, como
por termos de acessar toda a matriz $^kA$, o que \'{e} um grande
inconveniente para o armazenamento eficiente da matriz.  Ademais,
minimizar o preenchimento local, i.e.  aquele causado por cada piv\^{o}
individualmente, \'{e} uma estrat\'{e}gia gulosa que pode ter um fraco
desempenho global, i.e.  causar muito mais preenchimento que a
seq\"{u}\^{e}ncia \'{o}tima de $n$ piv\^{o}s.  Visando superar as
defici\^{e}ncias dos m\'{e}todos locais estudaremos a 
{\bf heur\'{\i}stica P3}, ou Procedimento de 
Pr\'{e}-determina\c{c}\~{a}o de Piv\^{o}s. 
 \index{P3}

A heur\'{\i}stica P3 procura uma permuta\c{c}\~{a}o de linhas e colunas,
$B=PAQ$, que reduza o preenchimento na fatora\c{c}\~{a}o $B=LU$.  O
procedimento que implementa a heur\'{\i}stica P3 posiciona as colunas de
$A$ em $B$ e numa matriz tempor\'{a}ria, $S$, e permuta linhas de todas
as colunas, tentando produzir uma $B$ que seja quase triangular
inferior.  As colunas de $B$ que n\~{a}o tem ENNs acima da diagonal
denominam-se {\bf colunas triangulares}.  Colunas triangulares s\~{a}o
sempre posicionadas com um ENN na diagonal, que ser\'{a} usado como
piv\^{o} da coluna na fatora\c{c}\~{a}o $B=LU$.  As colunas de $B$ que
tem ENNs acima da diagonal denominam-se {\bf espinhos}.  Na
fatora\c{c}\~{a}o LU de $B$, somente poder\'{a} ocorrer preenchimento
dentro dos espinhos.  Ademais, n\~{a}o haver\'{a} preenchimento acima do
mais alto (menor \'{\i}ndice de linha) ENN num espinho.  Portanto,
queremos minimizar o n\'{u}mero de espinhos e, como objetivo
secund\'{a}rio, minimizar a altura dos espinhos acima da diagonal. 
 \index{Espinhos} 
 \index{Coluna!triangular} 
 \index{Coluna!altura} 
 \index{P3} 

A heur\'{\i}stica P3 procede como segue: 
\begin{enumerate} 

\item Compute o n\'{u}mero de ENNs em cada linha e coluna de $A$, ou os
pesos de linha e coluna, respectivamente, $pr(i)$ e $pc(j)$. 

\item Compute $h= arg\min_i pr(i)$ e $\rho = pr(h)$. 

\item Compute a $\rho$-{\bf altura} de cada coluna, $\Theta _\rho (j)$,\\ 
 o numero de ENNs na coluna $j$ em linhas de peso $\rho$. 

\item Compute $t= arg\max_j \Theta _\rho (j)$.

\begin{enumerate}

\item Se $\rho =1$, seja 
 $h\in \{i\ \mid \ pr(i)=1 \wedge A(i,t)\neq 0 \}$;
 posicione $t$ como a primeira coluna de $A$, $h$ como a primeira
 linha, e aplique P3 recursivamente a $A$(2:m,2:n). 
 A coluna exclu\'{\i}da torna-se a \'{u}ltima coluna de $B$. 

\item Se $\rho >1$, posicione a coluna $t$ como a \'{u}ltima coluna de $A$,
 e aplique P3 recursivamente a $A$(:,1:n-1).  A coluna exclu\'{\i}da
 torna-se a primeira coluna de $S$. 

\item Se $\rho =0$, posicione $h$ como a primeira linha, reposicione a
 primeira coluna de $S$ (\'{u}ltima a ser exclu\'{\i}da de $A$) como a 
 \'{u}ltima coluna de $B$, e aplique P3 recursivamente a $A$(2:m,:). 

\end{enumerate} \end{enumerate} 

Em P3, o caso 
\begin{itemize}

\item $\rho =1$ corresponde ao posicionamento de uma coluna 
  triangular em $B$. 

\item $\rho >1$ corresponde a impossibilidade de posicionar uma
  coluna triangular.  Assim (temporariamente) eliminamos uma
  coluna de $A$ e prosseguimos com P3.  O crit\'{e}rio de
  sele\c{c}\~{a}o para a coluna a eliminar visa produzir o caso 
  $\rho =1$ nos pr\'{o}ximos passos de P3.  Em caso de empate na
  $\rho$-altura, poder\'{\i}amos usar como desempate a 
  $(\rho +1)$-altura. 

\item $\rho =0$ corresponde a n\~{a}o haver em $A$ uma coluna que
  pud\'{e}ssemos posicionar em $B$ de modo a continuar formando $B$
  n\~{a}o singular.  Ent\~{a}o reintroduzimos, como pr\'{o}xima coluna
  de $B$, a primeira coluna em $S$.  Escolhemos a primeira em $S$
  visando minimizar a altura do espinho acima da diagonal.  Existe o
  perigo de que este espinho n\~{a}o venha a prover um piv\^{o} n\~{a}o
  nulo, ou que ap\'{o}s cancelamentos, o piv\^{o} seja muito pequeno para
  garantir a estabilidade num\'{e}rica da fatora\c{c}\~{a}o.  Se $A$ \'{e}
  esparsa e tivermos apenas alguns espinhos, a melhor solu\c{c}\~{a}o
  \'{e} intercalar a fatora\c{c}\~{a}o num\'{e}rica e simb\'{o}lica, o que
  nos d\'{a} a chance de rejeitar piv\^{o}s inst\'{a}veis.  M\'{e}todos para
  lidar com piv\^{o}s inaceit\'{a}veis na fase de fatora\c{c}\~{a}o
  num\'{e}rica, ap\'{o}s termos completado uma fase independente de
  fatora\c{c}\~{a}o simb\'{o}lica, s\~{a}o discutidos posteriormente. 

\end{itemize}

Exemplo

Apliquemos o P3 \`{a} matriz booleana $A$.  
\`{A} medida que o P3 prossegue representaremos a matriz particionada 
 $\left[ \begin{array}{ccc} B, A, S \end{array} \right]$, 
 cujos ENNs ser\~{a}o denotados respectivamente $b$, $a$ e $s$.  Os
\'{\i}ndices de linha e coluna est\~{a}o a esquerda e acima da matriz. 
Os pesos, $\rho$, e as $\rho$-alturas, $\Theta _\rho$, a direita e
abaixo da matriz.  Os piv\^{o}s, ou espinhos, escolhidos a cada passo
est\~{a}o em negrito. 

$$ 
\begin{array}{c|ccccc|c} 
p\backslash q & 1 & {\bf 2} & 3 & 4 & 5 & \rho \\
\hline 
            1 & a & a       & a & a &   & 4 \\ 
            2 & a & a       &   &   &   & 2 \\ 
            3 &   &         & a & a &   & 2 \\ 
            4 &   & a       &   & a & a & 3 \\  
            5 &   & a       &   & a & a & 3 \\
\hline
\Theta _2     & 1 & 1       & 1 & 1 & 0 &   \\ 
\Theta _3     & 1 & 3       & 1 & 3 & . &   \\            
\Theta _4     & 2 & 4       & 2 & 4 & . &   
\end{array}  
\ \ , \ \ \  
\begin{array}{c|ccccc|c} 
p\backslash q & 1       & 3 & 4 & 5 & 2 & \rho \\
\hline 
            1 & a       & a & a &   & s & 3 \\ 
            2 & {\bf a} &   &   &   & s & 1 \\ 
            3 &         & a & a &   &   & 2 \\ 
            4 &         &   & a & a & s & 2 \\  
            5 &         &   & a & a & s & 2 \\
\hline
\Theta _1     & 1       & 0 & 0 & 0 & . & \\
\mbox{} & & & & & & \\               
\mbox{} & & & & & &              
\end{array} $$

$$  
\begin{array}{c|ccccc|c} 
p\backslash q & 1       & 3 & {\bf 4} & 5 & 2 & \rho \\
\hline 
            2 & {\bf b} &   &         &   & s & . \\ 
            1 & b       & a & a       &   & s & 2 \\ 
            3 &         & a & a       &   &   & 2 \\ 
            4 &         &   & a       & a & s & 2 \\  
            5 &         &   & a       & a & s & 2 \\
\hline
\Theta _2     & .       & 2 & 4       & 2 & . &  
\end{array}  
\ \ , \ \ \ 
\begin{array}{c|ccccc|c} 
p\backslash q & 1       & 3       & 5 & 4 & 2 &\rho \\
\hline 
            2 & {\bf b} &         &   &   & s & . \\ 
            1 & b       & {\bf a} &   & s & s & 1 \\ 
            3 &         & a       &   & s &   & 1 \\ 
            4 &         &         & a & s & s & 1 \\  
            5 &         &         & a & s & s & 1 \\
\hline
\Theta _1     & .       & 2       & 2 & . & . & 
\end{array}  
$$

$$
\begin{array}{c|ccccc|c} 
p\backslash q & 1       & 3       & 5 & 4       & 2 & \rho \\
\hline 
            2 & {\bf b} &         &   &         & s & . \\ 
            1 & b       & {\bf b} &   & s       & s & . \\ 
            3 &         & b       &   & {\bf s} &   & 0 \\ 
            4 &         &         & a & s       & s & 1 \\  
            5 &         &         & a & s       & s & 1 \\
\hline
\mbox{} & & & & & & 
\end{array}  
\ \ , \ \ \  
\begin{array}{c|ccccc|c} 
p\backslash q & 1       & 3       & 4       & 5       & 2 &\rho \\
\hline 
            2 & {\bf b} &         &         &         & s & . \\ 
            1 & b       & {\bf b} & b       &         & s & . \\ 
            3 &         & b       & {\bf b} &         &   & . \\ 
            4 &         &         & b       & {\bf a} & s & 1 \\  
            5 &         &         & b       & a       & s & 1 \\
\hline
\Theta _1     & .       & .       & .       & 2       & . & 
\end{array}  $$ 

$$
\begin{array}{c|ccccc} 
p\backslash q & 1       & 3       & 4       & 5       & 2 \\
\hline 
            2 & {\bf b} &         &         &         & b \\ 
            1 & b       & {\bf b} & b       &         & b \\ 
            3 &         & b       & {\bf b} &         & + \\ 
            4 &         &         & b       & {\bf b} & b \\  
            5 &         &         & b       & b       & {\bf b}  
\end{array}  $$ 
Na permuta\c{c}\~{a}o final obtida pela P3, dada pelos vetores de
permuta\c{c}\~{a}o $p$ e $q$, indicamos com um $+$ as posic{c}\~{o}es 
a serem preenchidas no processo de elimina\c{c}\~{a}o. 
Uma nota\c{c}\~{a}o mais compacta para o mesmo exemplo seria:

$$ 
\begin{array}{c|ccccc|cccccc} 
                      s &   & 1 &   & 2 &   & & &      & & \\ 
\bar p\backslash \bar q & 1 & 5 & 2 & 3 & 4 & & & \rho & & \\
\hline 
            2 & x & x & {\bf x} & x &   & 4 & 3 & 2 & 1 & . & . \\ 
            1 & {\bf x} & x &   &   &   & 2 & 1 & . & . & . & . \\ 
            3 &   &   & x & {\bf x} &   & 2 & 2 & 2 & 1 & 0 & . \\ 
            4 &   & x &   & x & {\bf x} & 3 & 2 & 2 & 1 & 1 & 1 \\  
            5 &   & {\bf x} &   & x & x & 3 & 2 & 2 & 1 & 1 & 1 \\
\hline
\Theta _2     & 1 & 1 & 1 & 1 & 0 &   \\ 
\Theta _3     & 1 & 3 & 1 & 3 & . &   \\
\hline 
\Theta _2     & . & . & 1 & 4 & 1             
\end{array}  $$

No exemplo seguinte reintroduzimos um espinho com $0$ na posi\c{c}\~{a}o
piv\^{o}.  Todavia, no processo de elimina\c{c}\~{a}o, esta
posi\c{c}\~{a}o ser\'{a} preenchida antes da sua utiliza\c{c}\~{a}o como
piv\^{o}. 

$$ 
\begin{array}{c|cccccc|ccccccc} 
                      s &   & 1 & 2 &   &   &   \\ 
\bar p\backslash \bar q & 1 & 6 & 4 & 5 & 3 & 2 \\
\hline 
 1 & {\bf x} & x & & &  &      & 2 & 1 & . & . & . & . & .\\ 
 2 & & x & x & & & {\bf x}     & 3 & 2 & 2 & 1 & . & . & . \\
 3 & & x & x & & {\bf x} &     & 3 & 2 & 2 & 1 & 1 & . & . \\ 
 5 & x & & x & {\bf x} & x &   & 4 & 4 & 3 & 2 & 2 & 1 & . \\  
 4 & x & x & {\bf 0} & & x & x & 4 & 3 & 2 & 2 & 1 & 0 & . \\ 
 6 & & {\bf x} & x & x & & x   & 4 & 3 & 3 & 2 & 2 & 1 & 0 \\ 
\hline
\Theta _2     & 1 & 1 & 0 & 0 & 0 & 0  \\ 
\Theta _3     & 1 & 3 & . & . & . & .  \\
\hline 
\Theta _2     & . & . & 2 & 0 & 2 & 2 \\             
\Theta _3     & . & . & 4 & . & 3 & 3 \\
\end{array}  $$

$$
\begin{array}{c|cccccc}
p\backslash q & 1 & 6 & 5 & 3 & 4 & 2 \\
\hline 
            1 & {\bf x} & & & & & x \\
            2 & & {\bf x} & & x & & x \\ 
            3 & & & {\bf x} & x & & x \\
            5 & x & x & x & {\bf \oplus } & & x \\ 
            4 & x & & x & x & {\bf x} & + \\ 
            6 & & x & & x & x & {\bf x} 
\end{array} $$

\begin{obs}{\rm 
Da maneira como descrevemos o P3, poder\'{\i}amos aplic\'{a}-lo
tamb\'{e}m a uma matriz retangular.  No caso de Programa\c{c}\~{a}o
linear \'{e} comum ordenarmos pelo P3 todas as colunas da matriz de
restri\c{c}\~{o}es, $A$, e depois tomarmos, a cada reinvers\~{a}o, as
colunas na base, $B$, conforme a ordem estabelecida por P3 em $A$
[Orchard-Hays68]. 
}\end{obs} 

\section*{Exerc\'{\i}cios}

\begin{enumerate}

\item Implemente a fatora\c{c}\~{a}o LU usando a heur\'{\i}stica de
Markowitz.  Antes de aceitar um piv\^{o}, assegure-se que este tem
m\'{o}dulo maior que $pivmin = 10^{-4}$, substituindo-o caso
contr\'{a}rio.  Use a representa\c{c}\~{a}o est\'{a}tica por colunas e
de rede da matriz. 

\item 
Implemente o P3 para ordenar as colunas de uma matriz retangular, 
conforme a observa\c{c}\~{a}o 6.3. Use a representa\c{c}\~{a}o 
est\'{a}tica por coluna da matriz. 

\item 
 Considere a possibilidade de termos um espinho problem\'{a}tico que, 
ap\'{o}s o processo de elimina\c{c}\~{a}o das colunas precedentes,
apresente um $0$ na posi\c{c}\~{a}o piv\^{o}. 
 Mostre que necessariamente existe um espinho, \`{a} direita do
problem\'{a}tico, cujo elemento na mesma linha do piv\^{o} do espinho
problem\'{a}tico, neste est\'{a}gio da fatora\c{c}\~{a}o da matriz,
\'{e} diferente de zero.  Descreva um procedimento para lidar com estes
casos excepcionais atrav\'{e}s da permuta\c{c}\~{a}o de colunas
espinhos.  Discuta a viabilidade de resolver o problema atrav\'{e}s de
permuta\c{c}\~{a}o de linhas. 

\end{enumerate}

 \clearpage
 \clearpage  
\setcounter{chapter}{4} 
\chapter{FATORA\c{C}\~{O}ES SIM\'{E}TRICAS e ORTOGONAIS}
\begin{center} 
 {\LARGE Projetores e Problemas Quadr\'{a}ticos}
\end{center} 

\section{Matrizes Ortogonais}

Dizemos que uma matriz quadrada e real \'{e} {\bf ortogonal} sse sua
transposta \'{e} igual a sua inversa.Dada $Q$ uma matriz ortogonal, suas
colunas formam uma base ortonormal de ${\Re}^n$, como pode ser visto da
identidade $Q'Q=I$.  O quadrado da {\bf norma quadr\'{a}tica} de um vetor 
$v$, 
 \index{Matriz!ortogonal} 
 \index{Norma!quadr\'{a}tica} 
 $$ {\| v \|}^{2} \equiv \sum_{i=1}^n (v_{i})^2 = v'v $$ 
 permanece inalterada por uma transforma\c{c}\~{a}o ortogonal, pois
 $$  {\| Qv\|}^{2} = (Qv)'(Qv) = v'Q'Qv = v'Iv = v'v = {\| v\|}^{2} .$$

\section{Fatora\c{c}\~{a}o QR}

Dada uma matriz real A $m\times n$, $m\geq n$, podemos
existe uma matriz ortogonal $Q$ tal que $A=Q \left[
\begin{array}{c} R \\ 0 \end{array} \right] $, onde $R$ \'{e} uma matriz
quadrada e triangular superior.  Esta decomposi\c{c}\~{a}o \'{e} dita
uma fatora\c{c}\~{a}o QR, ou {\bf fatora\c{c}\~{a}o ortogonal}, da matriz
$A$. Descrevemos a seguir um m\'{e}todo para fatora\c{c}\~{a}o ortogonal. 
 \index{Fatora\c{c}\~{a}o!QR}
 \index{Fatora\c{c}\~{a}o!ortogonal}

A {\bf rota\c{c}\~{a}o} de um vetor  
 $\left[ \begin{array}{c} x_1 \\ x_2 \end{array} \right] \in {\Re}^2$  
 por um \^{a}ngulo $\theta$ \'{e} dada pela transforma\c{c}\~{a}o linear 
 \index{Rota\c{c}\~{a}o!bidimensional} 
 \index{Angulo!de rota\c{c}\~{a}o} 
 $$ 
rot(\theta )x = \left[ \begin{array}{cc} 
 \cos (\theta ) & \sin (\theta ) \\ -\sin (\theta ) & \cos (\theta ) 
\end{array} \right] 
\left[ \begin{array}{c} x_1 \\ x_2 \end{array} \right] .
$$

Notemos que a rota\c{c}\~{a}o \'{e} uma
transforma\c{c}\~{a}o ortogonal, pois
$$ 
rot(\theta )' rot(\theta )
=
\left[ \begin{array}{cc} {\cos (\theta )}^2 + {\sin (\theta 
)}^2 & 0 \\ 
       0 & {\cos (\theta )}^2 + {\sin (\theta )}^2 
\end{array} \right] 
=
\left[ \begin{array}{cc} 1 & 0 \\  0 & 1 \end{array} 
\right] 
.$$

Al\'{e}m disso podemos tomar o \^{a}ngulo 
 $\theta = -\arctan (x_2 / x_1 )$ de modo que a rota\c{c}\~{a}o 
correspondente anule a segunda componente do vetor rodado 
(se $x_1=0$, tome $\theta=\pi/2$). 

 Como o produto de transforma\c{c}\~{o}es ortogonais continua ortogonal
(prove), podemos usar uma seq\"{u}\^{e}ncia de rota\c{c}\~{o}es para
levar a matriz $A$ \`{a} forma triangular superior, como veremos a
seguir. 
 
A {\bf rota\c{c}\~{a}o de Givens} \'{e} um operador linear cuja matriz
coincide com a identidade, exceto num par de linhas onde imergimos uma
matriz de rota\c{c}\~{a}o bidimencional: 
  \index{Rota\c{c}\~{a}o!de Givens} 
  \index{Matriz!de Givens} 
  $$ G(i,j,\theta ) =
\left[ \begin{array}{cccccccc} 
 1 & & & & & & & \\
   & \ddots & & & & & & \\
   & & \cos (\theta ) & & \sin (\theta ) & & & \\
   & & & \ddots & & & & \\
   & & -\sin (\theta ) & & \cos (\theta ) & & & \\
   &  & &  & &  & \ddots &  \\
   &  & &  & &  &  & 1 \\
\end{array} \right] \ .
$$
 Dizemos que a aplica\c{c}\~{a}o deste operador numa matriz $A$, $GA$,
roda as linhas $i$ e $j$ de $A$ de um \^{a}ngulo $\theta$. 

Abaixo ilustramos uma seq\"{u}\^{e}ncia de rota\c{c}\~{o}es de linhas
que leva uma matriz $5\times 3$ \`{a} forma triangular superior.  Cada
par de \'{\i}ndices, $(i,j)$, indica que rodamos estas linhas do
\^{a}ngulo apropriado para zerar a posi\c{c}\~{a}o na linha $i$, coluna
$j$.  Supomos que, inicialmente, a matriz \'{e} densa, i.e.  todos os
seus elementos s\~{a}o diferentes de zero, e ilustramos o padr\~{a}o de
esparsidade da matriz nos est\'{a}gios assinalados com um asterisco na
seq\"{u}\^{e}ncia de rota\c{c}\~{o}es. 
 $$
 (1,5) * (1,4) (1,3) (1,2) * (2,5) (2,4) (2,3) * (3,5) (3,4) *
 $$

$$
\left[ \begin{array}{ccc} 
  x & x & x \\ x & x & x \\ x & x & x \\ x & x & x \\ 0 & x & x
\end{array} \right] \ \  
\left[ \begin{array}{ccc} 
  x & x & x \\ 0 & x & x \\ 0 & x & x \\ 0 & x & x \\ 0 & x & x
\end{array} \right] \ \  
\left[ \begin{array}{ccc} 
  x & x & x \\ 0 & x & x \\ 0 & 0 & x \\ 0 & 0 & x \\ 0 & 0 & x
\end{array} \right] \ \  
\left[ \begin{array}{ccc} 
  x & x & x \\ 0 & x & x \\ 0 & 0 & x \\ 0 & 0 & 0 \\ 0 & 0 & 0
\end{array} \right]  
$$
Tomando $Q$ como a produt\'{o}ria das rota\c{c}\~{o}es de Givens, 
temos a fatora\c{c}\~{a}o $A=QR$, como procurada.

\section{Espa\c cos Vetoriais com Produto Interno}

Dados dois vetores $x,y \in {\Re}^n$, o seu {\bf produto escalar}
\'e definido como
$$ < x \mid y > \equiv x'y = \sum_{i=1}^{n} x_{i}y^{i} .$$
Com esta defini\c c\~ao, o produto escalar \'e um operador que
satisfaz as propriedades fundamentais de {\bf produto interno}, a saber: 
 \index{Produto!escalar} 
 \index{Produto!interno} 
 \index{Produto!simetria} 
 \index{Produto!linearidade} 
 \index{Produto!n\~{a}o negatividade} 
 \index{Produto!definitividade} 
\begin{enumerate}
\item $<x\mid y> = <y\mid x>$, {\bf simetria}. 
\item $<\alpha x+\beta y\mid z> = 
      \alpha <x\mid z> + \beta <y\mid z>$, {\bf linearidade}. 
\item $<x\mid x> \geq 0$ , n\~{a}o negatividade.
\item $<x\mid x>=0 \Leftrightarrow x=0$ , {\bf definitividade}.
\end{enumerate}
 
Atrav\'es do produto interno, definimos a norma:
$$ \| x \| \equiv <x\mid x>^{1/2} ;$$ 
e definimos tamb\'em o \^angulo entre dois vetores n\~{a}o nulos:
 \index{Angulo!entre vetores} 
 $$ \Theta(x,y) \equiv \arccos ( <x\mid y>  /  \| x\| \| y\| ) .$$

\section{Projetores}

Consideremos o subespa\c{c}o linear gerado pelas colunas de uma matriz
$A$, $m\times n$, $m\geq n$: $$ C(A) = \{ y=Ax, x\in {\Re}^n \}.  $$
Denominamos $C(A)$ de {\bf imagem} de $A$, e o complemento de $C(A)$,
$N(A')$, de {\bf espa\c{c}o nulo} de $A'$, 
 $$ N(A') = \{ y \mid A'y=0 \}.  $$
 O fator ortogonal $Q=[C\mid N]$ nos d\'{a} uma base ortonormal de
${\Re}^m$ onde as $n$ primeiras colunas s\~{a}o uma base ortonormal de
$C(A)$, e as $m-n$ \'{u}ltimas colunas s\~{a}o uma base de $N(A')$, 
como pode ser visto diretamente da identidade 
$Q'A=\left[ \begin{array}{c} R\\ 0 \end{array} \right]$.  
 \index{Projetor} 
 \index{Matriz!de proje\c{c}\~{a}o} 
 \index{Matriz!imagem} 
  \index{Matriz!espa\c{c}o nulo}  

Definimos a {\bf proje\c{c}\~{a}o} de um vetor $b\in {\Re}^m$ no 
espa\c{c}o das colunas de $A$, pelas rela\c{c}\~{o}es: 
$$ y = P_{C(A)}b \Leftrightarrow y\in C(A) \wedge (b-y)\perp C(A) $$ 
ou, equivalentemente, 
$$ y= P_{C(A)}b \Leftrightarrow 
   \exists x \mid y=Ax \ \wedge A'(b-y)=0. $$

No que se segue suporemos que $A$ tem posto pleno, i.e.  que suas
colunas s\~{a}o linearmente independentes.  Provemos que o projetor de
$b$ em $C(A)$ \'{e} dado pela aplica\c{c}\~{a}o linear 
$$ P_A = A(A'A)^{-1}A'.  $$ 
Se $y = A((A'A)^{-1}A'b)$, ent\~{a}o obviamente $y\in C(A)$.  
Por outro lado, $A'(b-y) = A'(I-A(A'A)^{-1}A')b = (A' - IA')b = 0$.  

\section{Quadrados M\'{\i}nimos}

Dado um sistema superdeterminado, $Ax=b$ onde a matriz $A$ $m\times n$
tem $m>n$, dizemos que $x^*$ ``resolve'' o sistema no sentido dos
{\bf quadrados m\'{\i}nimos}, ou que $x^*$ \'{e} a ``solu\c{c}\~{a}o'' 
de quadrados m\'{\i}nimos, sse $x^*$ minimiza a norma quadr\'{a}tica 
do res\'{\i}duo,
 $$ x^* = Arg \min_{x\in {\Re}^n} \| Ax - b {\|}, $$
 Dizemos tamb\'{e}m que $y=Ax^*$ \'{e} a melhor aproxima\c{c}\~{a}o, no
sentido dos quadrados m\'{\i}nimos de $b$ em $C(A)$. 
 \index{Quadrados M\'{\i}nimos}
 \index{Quadrados M\'{\i}nimos!solu\c{c}\~{a}o de} 

Como a multiplica\c{c}\~{a}o por uma matriz ortogonal deixa inalterada a
norma quadr\'{a}tica de um vetor, podemos procurar a solu\c{c}\~{a}o
deste sistema (no sentido dos quadrados m\'{\i}nimos) minimizando a
transforma\c{c}\~{a}o ortogonal do res\'{\i}duo usada na
fatora\c{c}\~{a}o QR de $A$,
$$
\| Q'(Ax-b) {\|}^2 = 
\| \left[ \begin{array}{c} R \\ 0 \end{array} \right] x - 
   \left[ \begin{array}{c} c \\ d \end{array} \right] 
{\|}^2 = 
\| Rx-c {\|}^2  +  \| 0x-d {\|}^2 .
$$
Da \'{u}ltima express\~{a}o v\^{e}-se que a solu\c{c}\~{a}o, a
aproxima\c{c}\~{a}o e o res\'{\i}duo do problema original s\~{a}o dados,
respectivamente, por
 $$
x^* = R^{-1}c \ , \ \ 
y = Ax^* \mbox{\ \ e \ \ }
z = Q \left[ \begin{array}{c} 0 \\ d \end{array} \right] .
$$
Como j\'{a} hav\'{\i}amos observado, as $m-n$ \'{u}ltimas colunas de $Q$
formam uma base ortonormal de $N(A')$, logo $z \perp C(A)$, de modo que
conclu\'{\i}mos que $y=P_{A}b$!

\section{Programa\c{c}\~{a}o Quadr\'{a}tica}

O problema de {\bf programa\c{c}\~{a}o quadr\'{a}tica} consiste em 
minimizar a fun\c{c}\~{a}o
 \index{Programa\c{c}\~{a}o Quadr\'{a}tica}
$$ f(y) \equiv (1/2)y'Wy + c'y \ , \ \ W=W'$$
sujeitos \`{a}s {\bf restri\c{c}\~{o}es} 
$$ g_{i}(y) \equiv N_{i}'y = d_i .$$
Os gradientes de $f$ e $g_i$ s\~{a}o dados, respectivamente, por 
$$ {\nabla}_{y}f = y'W + c' \ , \ \mbox{e} \  
   {\nabla}_{y}g_i = N_{i}' \ .
$$
 As {\bf condi\c{c}\~{o}es de otimalidade} de primeira ordem
(condi\c{c}\~{o}es de Lagrange) estabelecem que as restri\c{c}\~{o}es
sejam obedecidas, e que o gradiente da fun\c{c}\~{a}o sendo minimizada
seja uma combina\c{c}\~{a}o linear dos gradientes das
restri\c{c}\~{o}es.  Assim a solu\c{c}\~{a}o pode ser obtida em
fun\c{c}\~{a}o do {\bf multiplicador de Lagrange}, i.e. do  
vetor $l$ de coeficientes desta combina\c{c}\~{a}o linear, como
 \index{Condi\c{c}\~{a}o!de otimalidade} 
 \index{Condi\c{c}\~{a}o!de Lagrange} 
 \index{Multiplicadores!de Lagrange} 
 $$ N'y = d \ \wedge \ 
   y'W + c' = l' N' \ , $$ ou em forma matricial, 
 $$ \left[ \begin{array}{cc} N' & 0 \\ W & N \end{array} \right] 
 \left[ \begin{array}{c} y \\ l \end{array} \right] = 
 \left[ \begin{array}{c} d \\ c \end{array} \right] \ .  $$ 
 Este sistema de equa\c{c}\~{o}es \'{e} conhecido como o {\bf sistema
normal}.  O sistema normal tem por matriz de coeficientes uma matriz
sim\'{e}trica.  Se a forma quadr\'{a}tica $W$ for {\bf positiva definida}, 
i.e.se $\forall x\ x'Wx \geq 0 \ \wedge \ x'Wx=0 \Leftrightarrow x=0$, e
as restri\c{c}\~{o}es $N$ forem lineramente independentes, a matriz de
coeficientes do sistema normal ser\'{a} tamb\'{e}m positiva definida. 
Estudaremos a seguir como adaptar a fatora\c{c}\~{a}o de Gauss a
matrizes sim\'{e}tricas e positivas definidas. 
 \index{Sistema Normal} 
 \index{Matriz!positiva definida}

\section{Fatora\c{c}\~{a}o de Cholesky}

Ap\'{o}s a fatora\c{c}\~{a}o de Gauss de uma matriz sim\'etrica, $S=LU$,
podemos p\^{o}r em evid\^{e}ncia os elementos diagonais de $U$ obtendo
$S=LDL'$.  Se $S$ for positiva definida assim o ser\'a $D$, de modo que
podemos escrever $D=D^{1/2}D^{1/2}$, $D^{1/2}$ a matriz diagonal
contendo a raiz dos elementos em $D$.  Definindo $C=LD^{1/2}$, temos
 \index{Fatora\c{c}\~{a}o!de Cholesky}
 $S=CC'$, a {\bf fatora\c c\~ao de Cholesky} de $S$. 

Analisemos agora algumas rela\c{c}\~{o}es entre as fatora\c{c}\~{o}es de
Cholesky e ortogonal (QR), e de que maneiras a fatora\c{c}\~{a}o
ortogonal nos d\'{a} uma representa\c{c}\~{a}o da inversa. 
Primeiramente observemos que se
 $A=QR$, $$A'A = (QR)'QR = R'Q'QR = R'IR = L'L $$ i.e., o fator
triangular da fatora\c{c}\~{a}o ortogonal \'e o transposto do fator de
Cholesky da matriz sim\'{e}trica $A'A$. 

Veremos agora que, no caso de termos $A$ quadrada e de posto pleno, o
produto pela inversa, $y=A^{-1}x$, pode ser calculado sem o uso
expl\'{\i}cito do fator $Q$: 
Transpondo $AR^{-1}=Q$ obtemos $Q'=R^{-t}A'$, 
donde $$y= A^{-1}x = (QR)^{-1}x = R^{-1}Q'x = R^{-1}R^{-t}A'x \ .$$

\section*{Exerc\'{\i}cios} 
\begin{enumerate}
\item Use as propriedades fundamentais do produto interno para provar: 
  \index{Desigualdade!de Cauchy-Scwartz}
  \index{Desigualdade!triangular} 
  \index{Angulo!entre vetores}
  \begin{enumerate}
  \item A {\bf desigualdade de Cauchy-Scwartz}:
   $|<x\mid y>| \leq \| x\| \| y\|$.
   Sugest\~ao: Calcule ${\| x-\alpha y \|}^2$ para 
   $\alpha =<x\mid y>^2 / \| y\|$.
  \item A {\bf Desigualdade Triangular}:
   $\| x+y\| \leq \| x\| +\| y\| $.
  \item Em que caso temos igualdade na desigualdade de Cauchy-Schwartz?
   Relacione sua resposta com a defini\c c\~ao de \^angulo entre vetores.
  \end{enumerate}
\item  Use a defini\c c\~ao do produto interno em ${\Re}^n$ para 
 provar a Lei do Paralelogramo: 
 ${\| x+y\|}^2 + {\| x-y\|}^2 = 2{\| x\|}^2 + 2{\| y\|}^2 $.
\item Uma matriz $P$ \'{e} {\bf idempotente}, ou um projetor n\~{a}o 
ortogonal, sse $P^2 =P$.  Se $P$ \'{e} idempotente prove que:  
 \index{Matriz!idempodente} 
 \index{Matriz!sim\'{e}trica} 
 \index{Matriz!proje\c{c}\~{a}o} 
 \begin{enumerate}
 \item $R = (I-P)$ \'{e} idenpotente.
 \item ${\Re}^n = C(P) + C(R)$.
 \item Todos os autovalores de $P$ s\~{a}o $0$ ou $+1$. Sugest\~{a}o:
  Mostre que se $0$ \'{e} uma ra\'{\i}z do polin\^{o}mio 
  caracter\'{\i}stico de $P$, 
  ${\varphi}_{P}(\lambda )\equiv \det (P-\lambda I)$, ent\~{a}o
  $(1-\lambda)=1$ \'{e} ra\'{\i}z de ${\varphi}_{R}(\lambda )$.
  \end{enumerate}
\item Prove que $\forall P$ idempotente e sim\'etrico, 
  $P = P_{C(P)}$. Sugest\~ao: Mostre que $P'(I-P)=0$.
\item Prove que o operador de proje\c{c}\~{a}o num dado sub-espa\c co 
  vetorial $V$, $P_{V}$, \'e \'unico e sim\'etrico.
\item Prove o teorema de Pit\'agoras: 
  $\forall b \in {\Re}^{m} , u \in V$ temos que
  ${\| b-u \| }^2 = {\| b-P_{V}b \| }^2 + {\| P_{V}b - u \| }^2$.
\item  Formule o problema de quadrados m\'{\i}nimos como um
  problema de programa\c{c}\~{a}o quadr\'{a}tica.  
  \begin{enumerate} 
  \item Assuma dada uma base $N$ de $N(A')$.  
  \item Calcule diretamente o res\'{\i}duo, $z=b-y$, em fun\c{c}\~{a}o de $A$. 
  \end{enumerate}
\item O {\bf tra\c co} de uma matriz $A$ \'e definido por
      $ tr(A) \equiv \sum_i A_i^i$. Mostre  que
    \index{Matriz!tra\c{c}o} 
  \begin{enumerate}
  \item Se $A$, $m\times n$, tem posto pleno $\rho (A) = n$, ent\~ao
        $tr(P_A)=n$.
  \item Nas condi\c c\~oes do item anterior, definindo $R_A=(I-P_A)$,
        temos que $tr(R_A)=m-n$.
  \end{enumerate}
\item M\'{e}todo de Hauseholder: 
  A {\bf reflex\~{a}o} definida por um vetor unit\'{a}rio $u\mid u'u=1$, 
 \'{e} a transforma\c{c}\~{a}o linear $H= I-2uu'$.
   \index{Reflex\~{a}o} 
   \index{Matriz!Householder} 
  \begin{enumerate} 
   \item Interprete geometricamente a operea\c{c}\~{a}o de reflex\~{a}o. 
   \item Prove que $H=H'$, $H'=H^{-1}$, e $H^2=I$. 
   \item Dado $x\neq 0$ um vetor em $R^n$, tome 
         $$v= x \pm \| x\| 
    \left[ \begin{array}{c} 1 \\ 0 \\ \vdots \\ 0 \end{array} \right]  
    \ \ , \ \ \ u= v / \| v \| \ \ \mbox{e} \ \ \ H=I-2uu' \ \ .$$ 
    Mostre que $(Hx)_1 = \| x\|$, e que todas as demais componentes  
    de $Hx$ se anulam.   
  \item Discuta como poderiamos usar uma s\'{e}rie de reflex\~{o}es para 
        obter a fatora\c{c}\~{a}o QR de uma matriz. 
 \end{enumerate}    

\end{enumerate}

 \clearpage
 \clearpage  

\setcounter{chapter}{5} 
\chapter{ELIMINA\c{C}\~{A}O SIM\'{E}TRICA} 
\begin{center} 
{\LARGE Esparsidade na Fatora\c{c}\~{a}o de Cholesky.}
\end{center} 

\section{Grafos de Elimina\c{c}\~{a}o}

Na elimina\c{c}\~{a}o assim\'{e}trica, estudada no cap\'{\i}tulo 4,
procuramos uma permuta\c{c}\~{a}o geral, $QPAQ'$ que minimizasse o
preenchimento durante a fatora\c{c}\~{a}o $QPAQ'=LU$.  No caso
sim\'{e}trico queremos preservar a simetria da matriz, e
restringir-nos-emos a permuta\c{c}\~{o}es sim\'{e}tricas, 
$QAQ'= A_{q(i)}^{q(j)}=LL'$.  

Exemplo 1: 

Neste exemplo mostramos as posi\c{c}\~{o}es preenchidas na
fatora\c{c}\~{a}o de Cholesky de uma matriz sim\'{e}trica $A$, bem como
o preenchimento na fatora\c{c}\~{a}o de duas permuta\c{c}\~{o}es
sim\'{e}tricas da mesma matriz:
 $$ 
 \begin{array}{c} 1\\ 2\\ 3\\ 4\\ 5\\ 6 \end{array} \ \
\left[ \begin{array}{cccccc} 
    1 & x & x &   &   & x \\  
    x & 2 & x &   &   & 0 \\ 
    x & x & 3 & x &   & 0 \\ 
      &   & x & 4 &   & 0 \\ 
      &   &   &   & 5 & x \\ 
    x & 0 & 0 & 0 & x & 6 
    \end{array} \right] \ 
\left[ \begin{array}{cccccc} 
    1 & x & x & x &   &   \\  
    x & 3 & 0 & x & x &   \\ 
    x & 0 & 6 & 0 & 0 & x \\ 
    x & x & 0 & 2 & 0 & 0 \\ 
      & x & 0 & 0 & 4 & 0 \\ 
      &   & x & 0 & 0 & 5 
    \end{array} \right] \ 
\left[ \begin{array}{cccccc} 
    5 &   &   & x &   &   \\  
      & 4 &   &   & x &   \\ 
      &   & 2 &   & x & x \\ 
    x &   &   & 6 &   & x \\ 
      & x & x &   & 3 & x \\ 
      &   & x & x & x & 1 
    \end{array} \right] 
 $$ 

Assim como no cap\'{\i}tulo 4 a linguagem de teoria de grafos foi
\'{u}til para descrever o processo de elimina\c{c}\~{a}o sim\'{e}trica,
usaremos agora grafos sim\'{e}tricos para estudar a elimina\c{c}\~{a}o
sim\'{e}trica. O primeiro destes conceito a ser definido \'{e} o de 
{\bf grafos de elimina\c{c}\~{a}o}, que nada mais s\~{a}o que os grafos 
que tem por matriz de adjac\^{e}ncia as submatrizes $^kA$ 
do processo de fatora\c{c}\~{a}o (veja cap\'{\i}tulo 2):  
 Dado um grafo sim\'{e}trico $G=(N,E)$, $N=\{1,2,\ldots n\}$, e uma
ordem de elimina\c{c}\~{a}o $q=[\sigma (1),\ldots \sigma (n)]$, onde
$\sigma$ \'{e} uma permuta\c{c}\~{a}o de $[1,2,\ldots n]$, definimos o
processo de elimina\c{c}\~{a}o dos v\'{e}rtices de $G$ na ordem $q$
como a seq\"{u}\^{e}ncia de grafos de elimina\c{c}\~{a}o
$G_i=(N_i,E_i)$ onde, para $i=1\ldots n$,
$$ N_i=\{q(i), q(i+1), \ldots q(n)\}, \ \ \ 
   E_1=E, \ \ \mbox{e, para}\ i>1\ , $$
$$\{a,b\}\in E_i \Leftrightarrow \ 
  \{a,b\} \in E_{i-1} \ \ \mbox{ou} \ \ 
  \{q(i-1),a\},\ \{q(i-1),b\} \in E_{i-1} $$
 Definimos tamb\'{e}m o inverso ordem de elimina\c{c}\~{a}o, 
$\bar q(i) = k \Leftrightarrow q(k)=i$, significando que $i$ foi o
$k$-\'{e}simo v\'{e}rtice de $G$ a ser eliminado. 
 \index{Grafo!elimina\c{c}\~{a}o} 
 \index{Grafo!preenchido} 

 O {\bf grafo preenchido} \'{e} o grafo $P=(N,F)$, onde 
 $F= \cup _1^n E_i$.  
 Os lados originais e os lados preenchidos em $F$ s\~{a}o,
respectivamente, os lados em $E$ e em $F-E$. 
 %

\begin{obs}{\rm 
Ao eliminar a $j$-\'{e}sima coluna na fatora\c{c}\~{a}o de Cholesky da
matriz $QAQ'=A_{q(i)}^{q(j)}=LL'$ preenchemos exatamente as
posi\c{c}\~{o}es correspondentes aos lados preenchidos em $F$ quando 
da elimina\c{c}\~{a}o do v\'{e}rtice $q(j)$. 
}\end{obs}

Exemplo 2: 

Os grafos de elimina\c{c}\~{a}o e o grafo preenchido correspondentes
\`{a} segunda ordem de elimina\c{c}\~{a}o do exemplo 1,
$q=[1,3,6,2,4,5]$, s\~{a}o:

$$ 
 \begin{array}{cccc} 
   {\bf 1} & -      & 3 \\ 
         | & \times &  & \backslash \\ 
         2 &        & 6 & \mbox{}\; | \\ 
           &  /     &   & / \\ 
         5 &        & 4 
 \end{array} \ \ \  
 \begin{array}{cccc} 
     &        & {\bf 3} \\ 
     & /      & {\bf |} & \backslash \\ 
   2 & {\bf -} & 6 & \mbox{}\; | \\ 
     & /      &   & / \\ 
   5 &        & 4 
 \end{array} \ \ \  
 \begin{array}{ccc} 
   2 & -      & {\bf 6}  \\ 
     & \times & {\bf |}  \\ 
   5 &        & 4 
 \end{array} \ \ \  
 \begin{array}{ccc} 
   {\bf 2} &        &   \\ 
    {\bf |} & \backslash  &   \\ 
   5 & {\bf -} & 4  
 \end{array} \ \ \  
 \begin{array}{cccc} 
   5 &   -  & {\bf 4}   
 \end{array} \ \ \  
 \begin{array}{cccc} 
   1 &  -     & 3 \\ 
   | & \times & | & \backslash \\ 
   2 & -      & 6 & \mbox{}\; | \\ 
   | & \times & | & / \\ 
   5 & -      & 4 
 \end{array} 
$$

\begin{lm}[Parter] 
Se $f=\{i,j\}\in F$, ent\~{a}o ou $f$ \'{e} um lado original, i.\'{e}. 
$F\in E$, ou $f$ foi preenchido quando da elimina\c{c}\~{a}o de um
v\'{e}rtice $k \mid \bar q(k) < \min \{ \bar q(i), \bar q(j) \}$. 
\end{lm} 
 \index{Lema!Parter} 
 \index{Lema!do caminho}

\begin{lm}[do caminho]
Considere a elimina\c{c}\~{a}o dos v\'{e}rtices de $G=(N,E)$ na ordem
$q$.  Temos que $\{i,j\}\in F$ sse existe um caminho de $i$ a $j$ em
$G$, passando apenas por v\'{e}rtices eliminados antes de $i$ ou $j$,
i.e., $\exists C=(i, v_1,v_2,\ldots v_p,j)$, em $G$, com $\bar
q(1)\ldots \bar q(p) < \min \{\bar q(i), \bar q(j) \}$.
\end{lm} 

Demonstra\c{c}\~{a}o:    

$\Leftarrow$: trivial.

$\Rightarrow$: por aplica\c{c}\~{a}o recursiva do lema de Parter.

Dada a matriz $A$, $G=(N,E)=(N,B(A))$, a ordem de elimina\c{c}\~{a}o
$q$, e o respectivo grafo preenchido, consideremos o conjunto de
\'{\i}ndices de linha de ENNs na coluna $j$ do fator de Cholesky,
$L^j\mid QAQ'=LL'$:
 $$enn(L^j)=\{i\mid i>j \wedge \{q(i),q(j)\}\in F\} + \{ j\}\ .$$

Definimos a {\bf \'{a}rvore de elimina\c{c}\~{a}o}, H, por
$$\Gamma _H^{-1}(j)= \left\{ \begin{array}{c}
  j, \ \ \mbox{se}\ \ enn(L^j) = \{ j\} , \ \ \mbox{ou} \\
  \min \{i>j \mid i \in enn(L^j) \} \end{array} 
  \ , \ \ \mbox{caso contr\'{a}rio} \right. 
 $$
 Para n\~{a}o sobrecarregar a nota\c{c}\~{a}o usaremos 
$h()=\Gamma _H^{-1}(\ )$ e $g()=\Gamma _H(\ )$. 
 Assim $h(j)$, o pai de $j$ em $H$,
\'{e} simplesmente o primeiro ENN (n\~{a}o diagonal) na coluna $j$ de
$L$.  
 \index{Arvore!de elimina\c{c}\~{a}o} 

Exemplo 3:
As as \'{a}rvores de elimina\c{c}\~{a}o correspondentes ao exemplo 
1 s\~{a}o:  

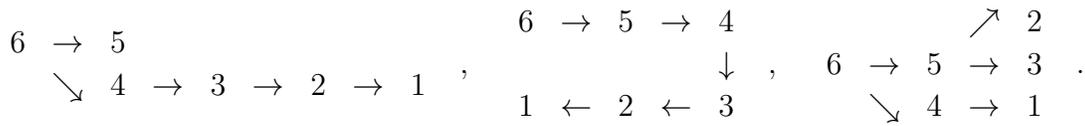
$$
\begin{array}{ccccccccc} 
6 & \rightarrow & 5 \\ 
  & \searrow  & 4 & \rightarrow & 3 & \rightarrow & 2 & \rightarrow & 1 
\end{array} \ \ , \ \ \  
\begin{array}{ccccc} 
 6 & \rightarrow & 5 & \rightarrow & 4 \\ 
   &             &   &             & \downarrow \\ 
 1 & \leftarrow & 2 & \leftarrow & 3 
\end{array} \ \ , \ \ \ 
\begin{array}{ccccc} 
   &             &   & \nearrow    & 2 \\  
 6 & \rightarrow & 5 & \rightarrow & 3 \\  
   & \searrow    & 4 & \rightarrow & 1 
\end{array} \ \ .   
 $$

\begin{teo}[da \'{a}rvore de elimina\c{c}\~{a}o]  
Dado $i>j$, $$i\in enn(L^j) \Rightarrow j\in \bar g(i) \ ,$$ i.e.,
qualquer \'{\i}ndice de linha abaixo da diagonal na coluna $j$ de $L$
\'{e} um ascendente de de $j$ na \'{a}rvore de
elimina\c{c}\~{a}o. 
\end{teo} 

Demonstra\c{c}\~{a}o: 

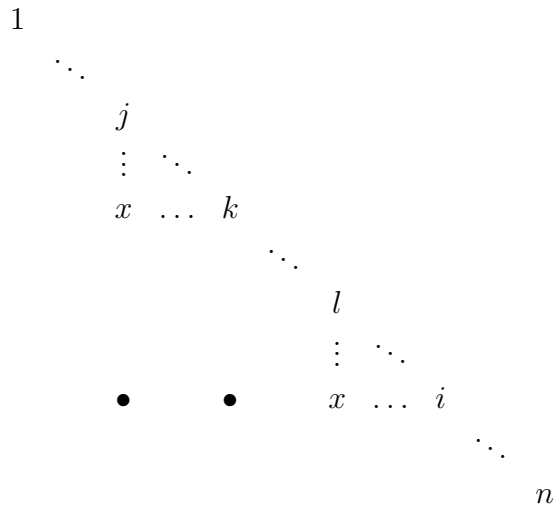
 $$ \begin{array}{ccccccccccc} 
 1 & \\ 
   & \ddots & \\ 
   &        & j \\ 
   &        & \vdots & \ddots \\ 
   &        &   x    & \ldots & k \\ 
   &        &        &        &   & \ddots \\ 
   &        &        &        &   &        &  l \\ 
   &        &        &        &   &        & \vdots & \ddots \\ 
   &        & \bullet & & \bullet &        &  x     & \ldots & i \\  
   &        &        &        &   &        &        &  &  & \ddots \\ 
   &        &        &        &   &        &        &  &  & & n \\ 
 \end{array} 
 $$

 Se $i=h(j)$ o resultado \'{e} trivial.  
Caso contr\'{a}rio (vide figura 1) seja $k=h(j)$. 
Mas $L_i^j\neq 0 \wedge L_k^j\neq 0 \Rightarrow L_i^k\neq 0$, pois
$\{q(j),q(i)\},\{q(j),q(k)\} \in G_j \Rightarrow \{q(k),q(i)\}\in G_{j+1}$.  
Agora, ou $i=h(k)$, ou aplicamos o argumento recursivamente,
reconstruindo o ramo de $H$, $(i,l,\ldots k,j)$, $i>l>\ldots >k>j$. QED.

 Pela demonstra\c{c}\~{a}o do teorema vemos que a \'{a}rvore de
elimina\c{c}\~{a}o retrata as depend\^{e}ncias entre as colunas para o
processo de fatora\c{c}\~{a}o num\'{e}rica da matriz.  Mais exatamente,
podemos eliminar a coluna $j$ de $A$ (i.e.  calcular todos os
multiplicadores na coluna $j$, denotados por $M^j$ no cap\'{\i}tulo 2, e
atualizar os elementos afetados por estes multiplicadores) sse j\'{a}
tivermos eliminado todos os descendentes de $j$ na \'{a}rvore de
elimina\c{c}\~{a}o. 

 Se pudermos realizar processamento paralelo (veja cap\'{\i}tulo 8), 
podemos eliminar simultaneamente todas as colunas em um mesmo n\'{\i}vel 
da \'{a}rvore de elimina\c{c}\~{a}o, come\c{c}ando pelas folhas, e 
terminando por eliminar a raiz.  

Exemplo 4: 
 \\ \indent   
Consideremos a elimina\c{c}\~{a}o de uma matriz com o mesmo padr\~{a}o 
de esparsidade da \'{u}ltima permuta\c{c}\~{a}o do exemplo 1. Sua 
\'{a}rvore de elimina\c{c}\~{a}o \'{e} a \'{u}ltima apresentada no 
exemplo 3. Esta \'{a}rvore tem 3 n\'{\i}veis que, das folhas para a 
raiz s\~{a}o: $\{1,3,2\}$, $\{4,5\}$, e $\{6\}$. Assim, podemos 
fatorar uma matriz com este padr\~{a}o de esparsidade em apenas 3 
etapas, como ilustrado no exemplo num\'{e}rico seguinte:  
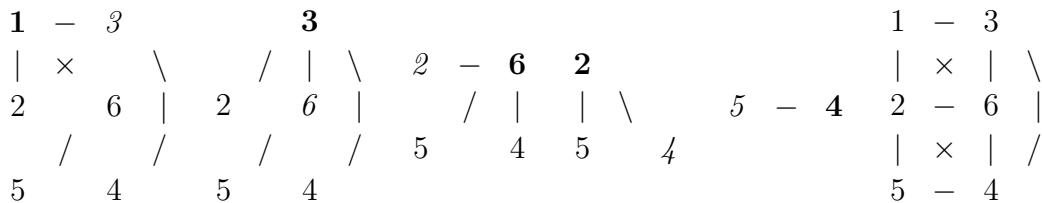
 $$ 
 \left[ \begin{array}{cccccc} 
    {\bf 1} &   &   & 7 &   &   \\  
      & {\bf 2} &   &   & 8 &   \\ 
      &   & {\bf 3} &   & 6 & 9 \\ 
    7 &   &   & 53 &   & 2 \\ 
      & 8 & 6 &   & 49 & 23 \\ 
      &   & 9 & 2 & 23 & 39 
    \end{array} \right]  
 \left[ \begin{array}{cccccc}
    1 &   &   & 7 &   &   \\
      & 2 &   &   & 8 &   \\
      &   & 3 &   & 6 & 9 \\
    {\em 7} &   &   & {\bf 4} &   & 2 \\
      & {\em 4} & {\em 2} &   & {\bf 5} & 5 \\
      &   & {\em 3} & 2 & 5 & 12 
    \end{array} \right] 
 \left[ \begin{array}{cccccc} 
    1 &   &   & 7 &   &   \\
      & 2 &   &   & 8 &   \\
      &   & 3 &   & 6 & 9 \\
    {\em 7} &   &   & 4 &   & 2 \\
      & {\em 4} & {\em 2} &   & 5 & 5 \\
      &   & {\em 3} & {\em \frac{1}{2}} & {\em 1} & {\bf 6}   
    \end{array} \right] 
$$

 O pr\'{o}ximo teorema mostra uma forma computacionalmente mais
eficiente de obter o grafo preenchido, $P=(N,F)$, e a \'{a}rvore de
elimina\c{c}\~{a}o, $H$, dado o padr\~{a}o de esparsidade da matriz
original, $G=(N,E)$, e a ordem de elimina\c{c}\~{a}o, $q$.  Nesta
vers\~{a}o simplificada dos grafos de elimina\c{c}\~{a}o, $G_j^*$, ao
eliminarmos o v\'{e}rtice $q(j)$, preenchemos apenas os lados incidentes
ao seu vizinho mais pr\'{o}ximo de ser eliminado. 
 \index{Fatora\c{c}\~{a}o!simb\'{o}lica}

\begin{teo}[da fatora\c{c}\~{a}o simb\'{o}lica]  
Definimos agora a vers\~{a}o
simplificada do processo de elimina\c{c}\~{a}o: 
\begin{eqnarray*}
 G_j^* &=& (N_j,E_j^*), \\ 
 E_1^* &=& E, \\
 h^*(j) &=& \min \{ i>j \mid \{q(j),q(i)\} \in E_j^* \} \\ 
 E_{j+1}^* &=& \{ \{a,b\} \in E_j^* \mid \bar q(a), \bar q(b) >j \} \\ 
           & &    \ \ \cup \ \ \{ \{q(h(j)),v\},\ 
               v \mid \bar q(v)\geq j \wedge \{q(j),v\}\in E_j^* \} \\ 
 F^* &=& \cup _1^n E_j^*  
\end{eqnarray*}
Afirmamos que o grafo preenchido e a \'{a}rvore de elimina\c{c}\~{a}o  
obtidos no processo simplificado coincidem com a defini\c{c}\~{a}o 
anterior, i.e., $F^*=F$ e $h^*=h$.
\end{teo}

Exemplo 5: Os grafos de elimina\c{c}\~{a}o simplificados e o grafo
preenchido referentes a segunda ordem de elimina\c{c}\~{a}o no exemplo
1, $q=[1,3,6,2,4,5]$, s\~{a}o:

$$ 
 \begin{array}{cccc} 
   {\bf 1} & -      & {\em 3} \\ 
         | & \times &  & \backslash  \\ 
         2 &        & 6 & \mbox{}\; | \\ 
           & /      &   & / \\ 
         5 &        & 4 
 \end{array} \ \ 
 \begin{array}{cccc} 
     &        & {\bf 3} \\ 
     & /      & | & {\bf \backslash } \\ 
   2 &        & {\em 6} & \mbox{}\; | \\ 
     & /      &   & / \\ 
   5 &        & 4 
 \end{array} \ \ 
 \begin{array}{cccc} 
   {\em 2} & {\bf - } & {\bf 6}  \\ 
     & /      & {\bf |} \\ 
   5 &        & 4 
 \end{array} \ \ 
 \begin{array}{cccc} 
   {\bf 2} &        &   \\ 
   {\bf | } & {\bf \backslash } &   \\ 
   5 &        & {\em 4}  
 \end{array} \ \ 
 \begin{array}{cccc} 
   {\em 5} & {\bf -}  & {\bf 4}   
 \end{array} \ \ 
 \begin{array}{cccc} 
   1 &  -     & 3 \\ 
   | & \times & | & \backslash \\ 
   2 & -      & 6 & \mbox{}\; | \\ 
   | & \times & | & / \\ 
   5 & -      & 4 
 \end{array} 
$$

O lema seguinte demonstra o teorema da fatora\c{c}\~{a}o simb\'{o}lica: 

\begin{lm} 
$$enn(L^j)=\cup_{k\in g(j)} enn(L^k) \ \cup \ enn(A^j) 
  -J +\{ j\} \ \ ,\ \ J=\{1,2,\ldots j\} \ .$$
\end{lm}

Demonstra\c{c}\~{a}o: 

$\supseteq$: 
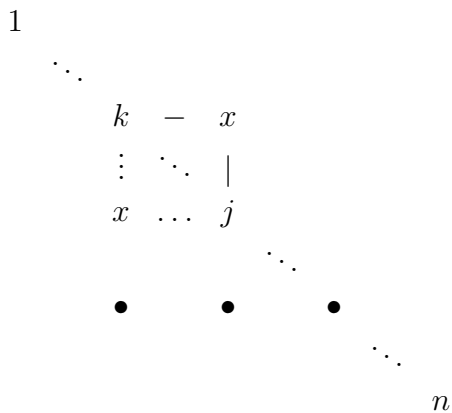
 $$ \begin{array}{ccccccccc} 
 1 & \\ 
   & \ddots & \\ 
   &        & k       &  -     &   x \\ 
   &        & \vdots  & \ddots &   | \\ 
   &        &   x     & \ldots &   j     \\ 
   &        &         &        &         & \ddots \\ 
   &        & \bullet &        & \bullet &        &  \bullet \\ 
   &        &         &        &         &        &          & \ddots \\ 
   &        &         &        &         &        &          &  & n  
 \end{array} 
 $$   
 Consideremos $k\in g(j)$ (vide figura 2). 
Se $i>j \in  enn(L^k)$, ent\~{a}o se $L_i^j$
j\'{a} n\~{a}o \'{e} um ENN, ser\'{a} preenchido ao eliminarmos $L^k$.

$\subseteq$: 
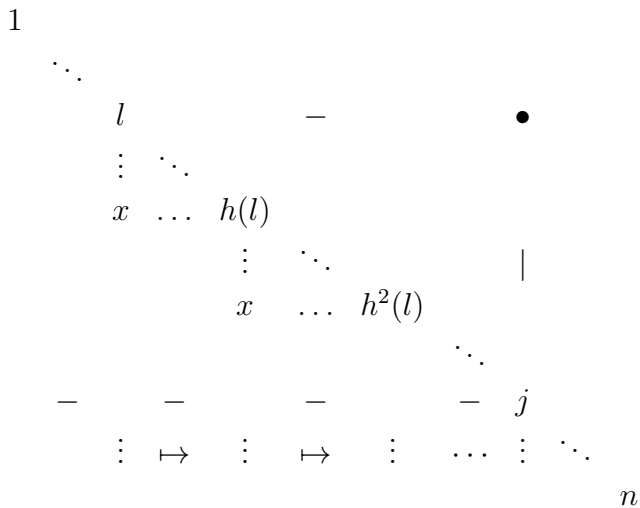
 $$ \begin{array}{ccccccccccc} 
 1 & \\ 
   & \ddots & \\ 
   &        &   l    &        &        &   -    &         &   & \bullet \\ 
   &        & \vdots & \ddots \\ 
   &        &   x    & \ldots & h(l) \\ 
   &        &        &        & \vdots & \ddots &         &   &  | \\ 
   &        &        &        &   x    & \ldots &  h^2(l) \\ 
   &        &        &        &        &        &         & \ddots \\ 
   &   -    &        &   -    &        &    -   &         &    -   & j \\  
   &        & \vdots  & \mapsto & \vdots & \mapsto & \vdots & \cdots 
         & \vdots & \ddots \\ 
   &        &        &        &   &        &        &  &  & & n \\ 
 \end{array} 
 $$   
 Seja $l<j \mid j \in enn(L^l)$, i.e., consideremos uma
coluna e cuja elimina\c{c}\~{a}o poderia causar preenchimentos na 
coluna $j$ (vide figura 3).  
Pelo teorema da \'{a}rvore de elimina\c{c}\~{a}o, $j$ \'{e} um 
ascendente de $l$, i.e., $\exists p\leq n \mid j=h^{p+1}(l)$.  
Mas pela primeira parte da prova,
 $$enn(L^l)-J \subseteq enn(L^{h(l)})-J \subseteq \ldots
  \subseteq enn(L^{h^p(j)}) \subseteq enn(L^j)- \{ j\} \ ,$$   
 de modo que qualquer v\'{e}rtice que poderia causar, durante a
elimina\c{c}\~{a}o da coluna $l$, uma preenchimento na coluna $j$, deve
necessariamente aparecer em $h^p(l) \in g(j)$.  QED.

\section{Grafos Cordais}

 Em um dado grafo sim\'{e}trico, $G=(N,E)$, definimos os seguintes
termos:
 $G=(N,E)$ \'{e} {\bf cordal} sse para qualquer ciclo
$C=(v_1,v_2,\ldots v_p,v_1),\ p\geq 3$, existe uma corda, i.e., um lado
$e\in E$ ligando dois v\'{e}rtices n\~{a}o consecutivos em $C$.  
 Um v\'{e}rtice \'{e} {\bf simplicial} sse sua elimina\c{c}\~{a}o n\~{a}o
causa preenchimento.  
 Uma ordem de elimina\c{c}\~{a}o, $q$, \'{e} {\bf perfeita} sse a
elimina\c{c}\~{a}o de nenhum dos v\'{e}rtices, na ordem $q$, causa
preenchimento.  
 Um subconjunto dos v\'{e}rtices, $S\subset N$, \'{e} um {\em
separador} entre os v\'{e}rtices $a$ e $b$, sse a remo\c{c}\~{a}o de $S$
deixa $a$ e $b$ em componentes conexas distintas. 
 Dizemos que $S\subset N$ \'{e} um {\bf separador} entre $A\subset N$ e
$B\subset N$ sse, $\forall a\in A,\ b\in B$, $S$ separa $a$ de $b$. 
 Um conjunto, $C$, satisfazendo uma propriedade, $P$, \'{e} {\em
minimal} (em rela\c{c}\~{a}o a $P$) sse nenhum subconjunto pr\'{o}prio
de $C$ satisfaz $P$.  Analogamente, um conjunto \'{e} {\bf maximal}, em
rela\c{c}\~{a}o a $P$, sse nenhum conjunto que o contenha propriamente
satisfaz $P$. 
 Um {\bf clique} \'{e} um conjunto de v\'{e}rtices $C\subset N$ onde
todos os v\'{e}rtices s\~{a}o adjacentes, i.e.,
 $\forall i,j\in C,\ \{i,j\}\in E$.  
 \index{Grafo!cordal} 
 \index{Corda} 
 \index{Ordem!de perfeita} 
 \index{Separador} 
 \index{V\'{e}rtice!simplicial} 
 \index{Maximal, Minimal} 
 index{Clique}

Exemplo 6: 

No 1o grafo do exemplo 2 vemos que 5 e 6 s\~{a}o  v\'{e}rtices 
simpliciais, $q=[5,4,2,6,3,1]$ \'{e} uma ordem de elimina\c{c}\~{a}o 
perfeita. No \'{u}ltimo grafo do exemplo 2 vemos que 
$S=\{2,3,6,4\}$ \'{e} um separador entre os v\'{e}rtices $a=1$ e $b=5$ 
(n\~{a}o minimal), $S'=S-\{3\}$ \'{e} um separador minimal, e 
$C=\{1,2,3,6\}$ \'{e} um clique.

\begin{teo}[da caracteriza\c{c}\~{a}o de grafos cordais] 
Dado $G=(N,E)$, as tr\^{e}s propriedades seguintes s\~{a}o equivalentes:
\begin{enumerate}

\item Existe em $G$ uma ordem de elimina\c{c}\~{a}o perfeita.

\item $G$ \'{e} cordal.

\item Se $S$ \'{e} um separador minimal entre $a,b \in N$, 
      ent\~{a}o $S$ \'{e} um clique.
\end{enumerate}
\end{teo} 

Demonstra\c{c}\~{a}o:

$1\Rightarrow 2$: 

Seja $q=[\sigma (1),\sigma (2),\ldots \sigma (n)]$ 
uma ordem de elimina\c{c}\~{a}o perfeita em $G$, e seja
$C=(v_1,v_2,\ldots v_p,v_1)$ um ciclo em $G$.  
Consideremos $v_k$ o primeiro v\'{e}rtice de  $C$ a ser eliminado, 
 $k=arg\min_{1\leq i\leq p}\{\bar q(v_i)\}$.  
Como $q$ \'{e} uma ordem de elimina\c{c}\~{a}o perfeita, a corda 
$\{v_{k-1},v_{k+1}\} \in E$. 

$2\Rightarrow 3$: 

 Seja $S$ um separador minimal entre $A$ e $B$, 
$a\in A$, $b\in B$, $v,w \in S$. 
 Como $S$ \'{e} minimal, existe um caminho 
$C=(a,c_1,\ldots c_p,v) \mid c_1, \ldots c_p \in A$, 
pois, caso contr\'{a}rio, $S-{v}$ continuaria separando $a$ de $b$. 
Analogamente existe um caminho $D$ conectando $a$ a $w$, com os
v\'{e}rtices intermedi\'{a}rios todos em $A$. 
 Existe pois um caminho de $v$ a $w$, cujos v\'{e}rtices
intermedi\'{a}rios est\~{a}o todos em $A$.  Seja $P$ um tal caminho de
comprimento m\'{\i}nimo, $P=(v,a_1,\ldots a_p,w)$.  Analogamente seja
$Q=(w,b_1,\ldots b_q,v)$ um caminho de $w$ a $v$ atrav\'{e}s de $B$ com
comprimento m\'{\i}nimo. 

Concatenando $P$ e $Q$ obtemos um ciclo $R=(v,a_1,\ldots
a_p,w,b_1,\ldots b_q,v)$.  Como $G$ \'{e} cordal, $R$ cont\'{e}m ao
menos uma corda, $\{x,y\}$.  Como $S$ \'{e} um separador, n\~{a}o
podemos ter $x\in A$ e $y\in B$.  Como $P$ tem comprimento m\'{\i}nimo,
n\~{a}o podemos ter $x,y \in P$ e, analogamente n\~{a}o podemos ter
$x,y\in Q$.  Assim, $\{v,w\}$ \'{e} a \'{u}nica corda poss\'{\i}vel em
$R$, e conclu\'{\i}mos que $\forall v,w\in S,\ \{v,w\}\in G$. 

Antes de $3\Rightarrow1$, provemos 2 lemas auxiliares:

\begin{lm}[hereditariedade] 
Se $G$ satisfaz a propriedade 3, ent\~{a}o qualquer subgrafo de $G$
induzido por um subconjunto de v\'{e}rtices tamb\'{e}m a satisfaz. 
\end{lm} 
 \index{Lema!hereditariedade} 

Demonstra\c{c}\~{a}o:

Seja $\tilde G=(\tilde N,\tilde E)$ o subgrafo de $G$ induzido por
$\tilde N \subset N$, e $\tilde S$ minimal em $\tilde G$ separando $a$
de $b$.  Podemos, em $G$, completar $\tilde S$, com v\'{e}rtices de
$N-\tilde N$, de modo a obter $S$, um separador minimal em $G$ entre $a$
e $b$.  Mas se $G$ satisfaz a propriedade 3, $S$ \'{e} um clique, e como
$\tilde S$ \'{e} um subgrafo de $S$, $\tilde S$ tamb\'{e}m \'{e} um 
clique. 

\begin{lm}[Gavril] 
Se $G$ satisfaz a propriedade 3 ent\~{a}o ou $G$ \'{e} completo, ou 
$G$ tem ao menos dois v\'{e}rtices simpliciais n\~{a}o adjacentes. 
\end{lm} 
 \index{Lema!Gavril} 

Demonstra\c{c}\~{a}o:

Provemos o lema por indu\c{c}\~{a}o no n\'{u}mero de v\'{e}rtices de
$G$, $n$.  Para $n=1$ o lema \'{e} trivial.  Se $G=(N,E),\ n>1$, n\~{a}o
\'{e} completo, $\exists a,b \in G \mid \{a,b\} \notin E$. 
 
Seja $S$ um separador minimal entre $a$ e $b$, partindo $G-S$ em ao
menos 2 componentes conexas distintas, $A$ e $B$, $a\in A$, $b\in B$. 
$G(S\cup A)$, o subgrafo induzido pelos v\'{e}rtices de $S\cup A$,
hereditariamente satisfaz a propriedade 3, e pela hip\'{o}tese de
indu\c{c}\~{a}o, ou $G(S\cup A)$ \'{e} um clique, ou cont\'{e}m (ao
menos) 2 v\'{e}rtices simpliciais n\~{a}o adjacentes. 

Se $G(S\cup A)$ n\~{a}o for completo, ent\~{a}o um dos seus 2
v\'{e}rtices simpliciais n\~{a}o adjacentes deve estar em $A$ (pois $S$
\'{e} um clique).  Se $G(S\cup A)$ \'{e} completo, ent\~{a}o qualquer
v\'{e}rtice de $A$ \'{e} simplicial (pois $S$ \'{e} um separador). 
Logo, existe ao menos um v\'{e}rtice simplicial, $v\in A$. 
Analogamente, existe um v\'{e}rtice simplicial $w\in B$, e $v$ n\~{a}o
\'{e} adjacente a $w$ (pois $S$ separa $v$ de $w$). 

\mbox{}\\ 
$3\Rightarrow 1$: 

Seja $G=(N,E)$ satisfazendo a propriedade 3.  Se $G$
\'{e} completo, qualquer ordem $q$ \'{e} perfeita.  Caso contr\'{a}rio,
existem 2 v\'{e}rtices simpliciais n\~{a}o adjacentes.  Podemos pois
eliminar um deles, $q(1)$, sem preenchimento, e pelo lema da
hereditariedade o subgrafo resultante continua satisfazendo a
propriedade 3. 
 Podemos ent\~{a}o usar o argumento recursivamente para obter uma ordem
de elimina\c{c}\~{a}o perfeita.  QED.

\section{Ordena\c{c}\~{o}es por Dissec\c{c}\~{a}o}

 Consideremos em $G=(N,E)$  um separador $S$ que parte 
$N-S$ com componentes conexas $N_{1},N_{2},\ldots,N_{k}$. 
 Em cada uma desta componentes podemos 
considerar um novo separador que a reparte, e assim recursivamente. 
 Podemos representar este processo pela {\bf \'{a}rvore de
dissec\c{c}\~{a}o}, $D$, onde $S$ \'{e} a raiz, uma componente separada
por $S$, ou o separador dentro dela, \'{e} filha de $S$, e as
componentes n\~{a}o repartidas s\~{a}o as folhas da \'{a}rvore. 
 \index{Ordem!dissec\c{c}\~{a}o} 
 \index{Ordem!p\'{o}s} 

 Uma {\bf p\'{o}s-ordem} dos v\'{e}rtices de uma \'{a}rvore $H$ de raiz
$r$ \'{e} uma ordem $q$, dos v\'{e}rtices de $H$, que lista os
v\'{e}rtices de cada uma das \'{a}rvores de $H-r$, recursivamente em
p\'{o}s-ordem, e, finalmente, $r$. 
  Seja $\tilde q$ uma p\'{o}s-ordem numa \'{a}rvore de dissec\c{c}\~{a}o,
$D$, de $G$.  Substituindo em $\tilde q$ cada n\'{o}, $d$, de $D$
pelos v\'{e}rtices de $G$ em $d$, temos uma {\bf ordem de 
dissec\c{c}\~{a}o}, $q$.

 \begin{lm}[Dissec\c{c}\~{a}o]
 Consideremos a elimina\c{c}\~{a}o dos v\'{e}rtices de
$G$ na ordem de dissec\c{c}\~{a}o $q$. 
 A elimina\c{c}\~{a}o de um dado v\'{e}rice, $v\in d$, s\'{o} pode
preencher lados dentro do seu n\'{o} em $D$, $d$, ou entre
(v\'{e}rtices de) $d$ e (v\'{e}rtices em) n\'{o}s ancestrais de $d$ (em
$D$), ou ainda entre (v\'{e}rtices em) ascendente de $d$ (em $D$). 
 \end{lm} 
 \index{Lema!dissec\c{c}\~{a}o} 

Demonstra\c{c}\~{a}o: Trivial, pelo lema do caminho.

Exemplo 7: 

Vejamos um exemplo de elimina\c{c}\~{a}o usando uma ordem de 
dissec\c{c}\~{a}o
 $$
 \begin{array}{ccccccccccc}
     &   &   & -  & -  & -  & 1  &            & 2  \\ 
     &   & \mbox{}\: / & &  &  / & |  &            & |  & \backslash \\ 
   3 &   & | &    & 4  &    & 5  &            & 6  &            & 7  \\ 
   | &   & | &    & |  &    & |  & \backslash & |  & \backslash \\ 
   8 & - & 9 & -  & 10 &    & 11 &            & 12 &            & 13  
 \end{array} \ \ \  
 \begin{array}{ccccc} 
           &    &   \{1,2\}  \\ 
           & /  &      |     & \backslash \\ 
  \{3,4\}  &    &   \{5,6\}  &            &  \{7\} \\ 
     |     &    &      |     & \backslash \\ 
 \{8,9,10\} &  & \{11,12\}   &            &  \{13\} 
 \end{array} 
 $$ 
 $$  
 \begin{array}{c} 
  1 \\ 2 \\ 3 \\ 4 \\ 5 \\ 6 \\ 7 \\ 8 \\ 9 \\ 10 \\ 11 \\ 12 \\ 13  
 \end{array} \ \ \ \   
  \begin{array}{ccccccccccccc} 
   9 & x & x  &     &    &    &    &    &    &     &    &  x &    \\ 
   x & 8 & 0  &  x  &    &    &    &    &    &     &    &  0 &    \\ 
   x & 0 & 10 &  0  & x  &    &    &    &    &     &    &  0 &    \\ 
     & x & 0  &  3  & 0  &    &    &    &    &     &    &  0 &    \\ 
     &   & x  &  0  & 4  &    &    &    &    &     &    &  x &    \\ 
     &   &    &     &    & 11 &    &    & x  &     &    &    &    \\ 
     &   &    &     &    &    & 12 &    & x  &  x  &    &    &    \\ 
     &   &    &     &    &    &    & 13 &    &  x  &    &    &    \\ 
     &   &    &     &    & x  & x  &    & 5  &  0  &    &  x &    \\ 
     &   &    &     &    &    & x  & x  & 0  &  6  &    &  0 &  x \\ 
     &   &    &     &    &    &    &    &    &     &  7 &    &  x \\ 
   x & 0 & 0  &  0  & x  &    &    &    & x  &  0  &    &  1 &  0 \\ 
     &   &    &     &    &    &    &    &    &  x  &  x &  0 &  2    
 \end{array}  
  $$

 O \'{u}ltimo conjunto da parti\c{c}\~{a}o, $S_k$ ou a raiz de $\tilde
G$, \'{e} em, $G$, um separador entre cada uma das componentes de
$\tilde G-k$.  Para minimizar a regi\~{a}o de poss\'{\i}vel
preenchimento (em $G$ ou $QAQ'$) gostar\'{\i}amos de ter o separador
$S_k$ ``pequeno e balanceado'', i.e. tal que 
 \index{Separador} 
\begin{itemize}
\item $\#S_k$ seja o menor poss\'{\i}vel.
\item Sub-\'{a}rvores tenham aproximadamente o mesmo n\'{u}mero de 
      v\'{e}rtices de $G$.
\end{itemize}
Recursivamente, gostar\'{\i}amos de ter como raiz de cada uma das
sub-\'{a}rvores um bom separador, i.e., pequeno e balanceado.  
Veremos a seguir v\'{a}rias heur\'{\i}sticas para obter num grafo qualquer, 
$G$, um bom separador. 

Heur\'{\i}stica de {\bf Busca em Largura}:  

 Uma busca em largura, BEL, a partir uma raiz $v\in N$, particiona
os v\'{e}rtices de $G$ em n\'{\i}veis $L_0, L_1,\ldots L_k$, definidos
por $$ L_0=\{v\}\ , \ \ L_{i+1}= adj(L_i)-L_{i-1} \ .$$
 A {\bf profundidade} do n\'{\i}vel $L_i$ \'{e} $i$, e a {\bf largura}
do n\'{\i}vel $L_i$ \'{e} $\#L_i$.  A profundidade e a largura da BEL
s\~{a}o, respectivamente, a m\'{a}xima profundidade e a m\'{a}xima
largura nos n\'{\i}veis da BEL. 
 \index{Busca!em largura}

 \begin{lm} 
 O n\'{\i}vel $L_i$
separa, em $G$, os v\'{e}rtices em n\'{\i}veis mais profundos dos em
n\'{\i}veis menos profundos que $i$.  
 \end{lm} 

 A heur\'{\i}stica de BEL procura
um separador balanceado $S\subseteq L_i$, tomando $i \approx k/2$, ou
ent\~{a}o tomando
 $i\mid \sum_1^{i-1}\#L_j <n/2 \ \wedge \ \sum_{i+1}^n\#L_j <n/2 \ .$
 
 Para obter um separador pequeno a heur\'{\i}stica procura uma raiz,
$v$, que gere uma BEL de m\'{a}xima profundidade, com o intuito de
reduzir a largura da BEL.  A {\bf dist\^{a}ncia} (em $G$) de um v\'{e}rtice
$v$ a um v\'{e}rtice $w$, $dist(v,w)$, \'{e} o comprimento, ou
n\'{u}mero de lados, do caminho mais curto entre ambos os v\'{e}rtices. 
A {\bf excentricidade} de um v\'{e}rtice $v$ \'{e} $exc(v) = \max_{w\in N}
dist(v,w)$. Um v\'{e}rtice de m\'{a}xima excentricidade se diz
{\bf perif\'{e}rico}, e sua excentricidade \'{e} o 
{\bf di\^{a}metro} de $G$. 
 \index{Grafo!di\^{a}metro} 
 \index{V\'{e}rtice!perif\'{e}rico}
 \index{V\'{e}rtice!dist\^{a}ncia}
 \index{V\'{e}rtice!excentricidade}

 Uma BEL com raiz $v$ ter\'{a} profundidade igual a excentricidade de
$v$.  Isto motiva querermos iniciar a BEL por um v\'{e}rtice
perif\'{e}rico.  Encontrar um v\'{e}rtice perif\'{e}rico \'{e} um
problema computacionalmente dif\'{\i}cil.  A {\bf heur\'{\i}stica de 
Gibbs} encontra um v\'{e}rtice {\bf quase-perif\'{e}rico} como segue:
 \begin{enumerate}
 \item Escolha como raiz um v\'{e}rtice de grau m\'{\i}nimo.
 \item Forme os n\'{\i}veis da BEL com raiz $v$, $L_1\ldots L_k$.
  Particione o n\'{\i}vel mais profundo em suas componentes conexas,
  $L_k = \cup_1^l S_j$, e tome um v\'{e}rtice de grau m\'{\i}nimo, 
  $v_j$, em cada componente.
 \item Para $j=1:l$
  \begin{itemize}
  \item Tome $v_j$ como nova raiz e encontre os n\'{\i}veis da 
   BEL, $L_1\ldots L_{k'}$
  \end{itemize}
 Ate que $k'>k$ ou $j=l$.
 \item Se o passo 3 terminou com $k'>k$, volte ao passo 2. Caso 
 contr\'{a}rio a atual raiz \'{e} um v\'{e}rtice quase-perif\'{e}rico.
 \end{enumerate}

Exemplo 8: 

Retomando o exemplo 7, a heuristica de Gibbs encontra 3 como 
v\'{e}rtice quase-perif\'{e}rico. Tomando 3 como raiz geramos a
\'{a}rvore $H$ por BEL: 
 $$ 
 \begin{array}{cccccccccccccccccccc} 
   & &  &  &  &  &  &  & 10& -& 4&  &   &  &  &  & 13 \\   
   & &  &  &  &  &  & /&   &  &  &  &   &  &  & /&  \\                       
 H=& & 3& -& 8& -& 9& -&  1& -& 5& -& 12& -& 6& -& 2& -& 7 \\ 
    \\ 
 L=& & 1&  & 2&  & 3&  &  4&  & 5&  & 6&  & 7&  &  8&  & 9  
 \end{array} 
 $$ 

Escolhendo $S_1=L_5$ como primeiro separador, e depois $S_2=L_3$ e 
$S_3=L_7$ como separadores dentro de cada uma das componetes separadas 
por $S_1$, obtemos a odem por dissec\c{c}\~{a}o 
$q=[3,8,1,10,9, 11,12,2,13,7,6, 4,5]$.    

 $$  
 \begin{array}{c} 
  1 \\ 2 \\ 3 \\ 4 \\ 5 \\ 6 \\ 7 \\ 8 \\ 9 \\ 10 \\ 11 \\ 12 \\ 13  
 \end{array} \ \ \ \   
 \begin{array}{ccccccccccccc} 
  3  & x &    &     &    &    &    &    &    &     &    &    &    \\ 
  x  & 8 &    &     &  x &    &    &    &    &     &    &    &    \\ 
     &   &  1 &     &    &    &    &    &    &     &    &  x &  x \\ 
     &   &    & 10  &  x &    &    &    &    &     &    &  x &    \\ 
     & x &    &  x  &{\bf 9}& &    &    &    &     &    &  0 &    \\ 
     &   &    &     &    & 11 &    &    &    &     &    &    &  x \\ 
     &   &    &     &    &    & 12 &    &    &     &  x &    &  x \\ 
     &   &    &     &    &    &    &  2 &    &   x &  x &    &    \\ 
     &   &    &     &    &    &    &    & 13 &     &  x &    &    \\ 
     &   &    &     &    &    &    &  x &    &   7 &  0 &    &    \\ 
     &   &    &     &    &    &  x &  x &  x &   0 &{\bf 6}& &  0 \\ 
     &   &  x &  x  &  0 &    &    &    &    &     &    &{\bf 4}& 0 \\ 
     &   &  x &     &    &  x &  x &    &    &     &  0 &  0 &{\bf 5}\\ 
 \end{array}  
  $$ 

 Note que nas linhas (colunas) correspondentes aos v\'{e}rtices do
primeiro separador, $S_1=\{1\}$, pode haver ENN's em qualquer
posi\c{c}\~{a}o.  Note tamb\'{e}m que o resto da matriz est\'{a} em
forma diagonal blocada (vide defini\c{c}\~{a}o no cap\'{\i}tulo 8), onde
cada bloco correspnde a uma das componentes separadas por $S_1$.  Esta
estrutura se repete em cada bloco, formando a estrutura ``espinha de
peixe'' caracter\'{\i}stica de ordens por dissec\c{c}\~{a}o.  Note que
esta estrutura \'{e} preservada pela fatora\c{c}\~{a}o de Cholesky. 
 \index{Matriz!espinha de peixe}

\section*{Exerc\'{\i}cios}

\begin{enumerate} 

\item 
Implemente a elimina\c{c}\~{a}o simb\'{o}lica num grafo $G$,
com a ordem $q$, computando $F$ e $H$, em tempo $O( \# enn(L))$.

\item 
Qu\~{a}o eficiente ($k$ em tempo $O(n^k)$) poder\'{\i}amos implementar o
algoritmo impl\'{\i}cito na \'{u}ltima parte do teorema de
caracteriza\c{c}\~{a}o de grafos cordais?

\item 
Considere o seguinte algoritmo para numerar os v\'{e}rtices de um
grafo na ordem $n,n-1,\ldots 2,1$:
 \index{Ordem!ORGM} 

{\bf Ordenamento Reverso por Grau M\'{a}ximo} (ORGM):
\begin{enumerate}
\item Escolha, como v\'{e}rtice $n$, um v\'{e}rtice qualquer. 
\item Escolha, como v\'{e}rtice seguinte um v\'{e}rtice ainda n\~{a}o 
 numerado adjacente a um n\'{u}mero m\'{a}ximo de v\'{e}rtices j\'{a} 
 numerados
\end{enumerate}
 Prove que ORGM define uma ordem perfeita.

\item Qu\~{a}o eficientemente podemos implementar o ORGM?

\item De um exemplo de v\'{e}rtice quase-perif\'{e}rico que 
n\~{a}o seja perif\'{e}rico. 

\end{enumerate} 
   
 \clearpage
 \clearpage 

\setcounter{chapter}{6} 
\chapter{ESTRUTURA} 
\begin{center}           
{\LARGE Acoplamento de Sub-Sistemas}
\end{center}

\section*{Estruturas Blocadas}

Consideraremos neste cap\'{\i}tulo matrizes com blocos de elementos
nulos dispostos de forma regular.  Estudaremos duas estruturas: a {\bf
triangular-blocada} superior, e a {\bf angular blocada} por colunas.  O
bloco $^r_sB$ tem dimens\~{a}o $m(r)\times n(s)$, e na estrutura
triangular blocada $m(k)=n(k)$. 
 \index{Estrutura!blocada} 
 
$$
\left[ \begin{array}{ccccc} 
 ^1_1B   & ^2_1B & ^3_1B & \ldots & ^h_1B \\
 0       & ^2_2B & ^3_2B & \ldots & ^h_2B \\ 
 0       &  0    & ^3_3B & \ldots & ^h_3B  \\
 \vdots  &       &       & \ddots & \vdots \\
 0       & 0     & 0     & \ldots & ^h_hB 
\end{array} \right] \ \ , \ \ \ 
\left[ \begin{array}{ccccc} 
 ^1_1B   & 0     & \ldots & 0 & ^h_1B \\
 0       & ^2_2B &        & 0 & ^h_2B \\ 
 \vdots  &       & \ddots &  & \vdots  \\ 
 0       & 0     &        & ^{h-1}_{h-1}B \\
 0       & 0     & \ldots & 0 & ^h_hB 
\end{array} \right] $$

Estas estruturas blocadas propiciam grandes facilidades computacionais. 
Em particular triangulariz\'{a}-las corresponde a triangularizar os
blocos diagonais, sendo que nenhum elemento n\~{a}o nulo \'{e} criado,
durante a triangulariza\c{c}\~{a}o, nos blocos nulos.

\section{Estrutura Triangular Blocada} 

Como j\'{a} visto no estudo de matrizes de permuta\c{c}\~{a}o, 
dado $G=(N,B)$, $N=\{1,\ldots n\}$, $B,\ n\times n$, sua matriz de
adjac\^{e}ncia, e um reordenamento de seus v\'{e}rtices, 
$q=[q_1,\ldots q_n]$, $q_i\in N$, ent\~{a}o a matriz de adjac\^{e}ncia
do grafo $G$ com os v\'{e}rtices reordenados (i.e., reindexados) por $q$
\'{e} $\tilde B= B_{q(i)}^{q(j)}= QBQ'$. 
 \index{Estrutura!triangular} 

 \begin{lm} 
 Considere o reordenamento coerente dos v\'{e}rtices de $G=(N,B)$ numa
ordem coerente $q=[{^1q}, {^2q} \ldots {^hq}]$.  
 Conforme a defini\c{c}\~{a}o de ordem coerente, cada bloco,
$^kq=[{^kq_1},\ldots {^kq_{n(k)}}]$, cont\'{e}m os v\'{e}rtices de uma
CFC (componente fortemente conexa) de $G$, $v_k$, e  $(v_1,\ldots
v_h)$ est\~{a}o topologicamente ordenados. 
 Neste caso a matriz de adjac\^{e}ncia do grafo reordenado, 
 $\tilde B=QBQ'$, \'{e} triangular blocada superior, de blocos 
 $^r_sB$, $n(r)\times n(s)$.  
 \end{lm}

Uma matriz que n\~{a}o possa ser, por permuta\c{c}\~{o}es de linhas e de
colunas, reduzida, isto \'{e}, posta na forma triangular-blocada, \'{e}
dita {\bf irredut\'{\i}vel}. 
  \index{Matriz!irredut\'{\i}vel}

Dada uma matriz $A$, n\~{a}o singular, procuraremos permuta\c{c}\~{o}es
de linhas e colunas, isto \'{e}, por matrizes de permuta\c{c}\~{a}o $R$
e $Q$, encontrar $\tilde A =RAQ$ que seja a ``mais fina''
parti\c{c}\~{a}o poss\'{\i}vel de $A$.  Observemos que a hip\'{o}tese de
n\~{a}o singularidade de $A$ implica na n\~{a}o singularidade dos blocos
diagonais, pois
 $$det(A)= det(\tilde A )= \prod _{k=1}^h det({^k_k{\tilde A}})\ .$$

Sejam $R$ e $Q$ matrizes de permuta\c{c}\~{a}o e seja $P=QR$, 
 $$\tilde A = R A Q = Q^{-1} Q R A Q = Q' P A Q$$ ou seja, podemos
escrever qualquer transforma\c{c}\~{a}o do tipo $RAQ$ como uma
permuta\c{c}\~{a}o de linhas, $PA$, seguida de uma permuta\c{c}\~{a}o
sim\'{e}trica $Q'(PA)Q$. 

Consideremos agora o grafo associado \`{a} matriz $PA$, $G(PA)$, que tem
por matriz de adjac\^{e}ncia o padr\~{a}o de esparsidade de $PA$,
 $B(PA)$, i.e., $$G(PA) = (N,B(PA)) = (N,\Gamma ),\ j\in \Gamma(i)
  \Leftrightarrow (PA)_i^j \neq 0\ .$$

Do \'{u}ltimo lema sabemos que a permuta\c{c}\~{a}o sim\'{e}trica que
d\'{a} a mais fina parti\c{c}\~{a}o de $PA$ \'{e} um ordenamento
coerente dos v\'{e}rtices de $G(PA)$.  Ademais, $PA$ \'{e}
irredut\'{\i}vel se $G(PA)$ \'{e} fortemente conexo.  Resta, portanto,
analisar o papel da permuta\c{c}\~{a}o de linhas, $P$, na
permuta\c{c}\~{a}o geral de linhas e colunas de $\tilde A=Q'PAQ$, onde a
permuta\c{c}\~{a}o sim\'{e}trica \'{e} dada por um ordenamento coerente
em $G(PA)$. 

Pela n\~{a}o singularidade, $A$ deve possuir ao menos uma diagonal de
elementos n\~{a}o nulos (n\~{a}o necessariamente a diagonal principal,
veja as defini\c{c}\~{o}es de determinante e diagonal).  Portanto existe
uma permuta\c{c}\~{a}o de linhas, $P$, que posiciona esta diagonal de
elementos n\~{a}o nulos na diagonal principal de $PA$.  Uma tal
permuta\c{c}\~{a}o $P$ ser\'{a} dita uma {\bf permuta\c{c}\~{a}o
pr\'{o}pria}, ao passo que uma permuta\c{c}\~{a}o que coloque algum
elemento nulo na diagonal principal ser\'{a} dita {\bf impr\'{o}pria}. 
Isto posto, mostraremos que:
 \index{Permuta\c{c}\~{a}o!pr\'{o}pria}

\pagebreak 

\begin{teo} \mbox{} 
 \begin{enumerate} 
 \item Qualquer permuta\c{c}\~{a}o pr\'{o}pria induz a mesma estrutura
de parti\c{c}\~{a}o, isto \'{e}, $h$ blocos de dimens\~{o}es $n(1)\ldots
n(h)$, com os mesmos \'{\i}ndices em cada bloco.  
 \item Qualquer permuta\c{c}\~{a}o impr\'{o}pria n\~{a}o pode induzir uma
parti\c{c}\~{a}o mais fina em $\tilde A$. 
 \end{enumerate}
 \end{teo}

Demonstra\c{c}\~{a}o: 

Pelo lema de Hoffmann $G(PA)$ \'{e} fortemente conexo, i.e., $PA$ \'{e}
irredut\'{\i}vel se qualquer conjunto de $k<n$ linhas tiver elementos
n\~{a}o nulos em ao menos $k+1$ colunas.  Como esta
caracteriza\c{c}\~{a}o independe da ordem das linhas, mostramos que a
irredutibilidade da matriz $PA$ independe da permuta\c{c}\~{a}o
pr\'{o}pria considerada. 

Da mesma forma, a irredutibilidade do bloco sudeste, $^h_h{\tilde A}$
bem como a exist\^{e}ncia da CFC $v_h$ entre as ``\'{u}ltimas'' (no 
sentido da ordem natural do grafo reduzido) componentes associadas a 
qualquer outra permuta\c{c}\~{a}o pr\'{o}pria, $P$, est\'{a} garantida, 
de modo que a invari\^{a}ncia da estrutura de $\tilde A$ segue por
indu\c{c}\~{a}o no n\'{u}mero de componente fortemente conexas (blocos). 

Mostramos agora que se $A$ tiver algum elemento diagonal nulo, ent\~{a}o
as componentes fortemente conexas de $G(A)$ s\~{a}o fus\~{o}es de
componentes fortemente conexas de $G(PA)$:
Distinguamos os zeros na diagonal de $A$ e consideremos o grafo
$G^+(PA)$ obtido de $G(PA)$ ao adicionarmos as arestas correspondentes
aos zeros distinguidos em $PA$.  As CFCs em $G(A)$ s\~{a}o exatamente as
CFCs em $G^+(PA)$, sendo que as arestas adicionais (se quando
adicionadas ao grafo reduzido de $G(PA)$ formarem ciclos) s\'{o} podem
tornar equivalentes alguns v\'{e}rtices antes em CFCs distintas. QED.  

 A permuta\c{c}\~{a}o pr\'{o}pria, $P$, pode ser vista como um casamento
perfeito entre linhas e colunas, onde casamos a coluna $j$ com a linha
$p(j)$, entre seus pretendentes
 $\Gamma ^{-1} (j) = \{i\mid B_i^j \neq 0 \}$.  
 Podemos portanto encontr\'{a}-la atrav\'{e}s do algoritmo H\'{u}ngaro,
de prefer\^{e}ncia com o emprego de alguma heur\'{\i}stica eficiente
para evitar a freq\"{u}ente gera\c{c}\~{a}o de \'{a}rvores de caminhos
de aumento. 

Exemplo 1: 

Considere as matrizes de adjac\^{e}ncia $B$, sua permuta\c{c}\~{a}o 
pr\'{o}pria $PB$, e o posterior ordenamento coerente $Q'PBQ$. 
$$
B= \left[ \begin{array}{ccc}
\mbox{\bf 0} & \mbox{\bf 1} & 0 \\ 0 & 1 & \mbox{\bf 1} \\ 
\mbox{\bf 1} & 1 & 1 \end{array} \right] \ \ , \ \ \  
PB=B^j_{p(i)}= \left[ \begin{array}{ccc} 
\mbox{\bf 1} & 1 & 1 \\ \mbox{\bf 0} & \mbox{\bf 1} & 0 \\ 
0 & 1 & \mbox{\bf 1} \end{array} \right] \ \ , $$
$$  
Q'PAQ=B^{q(j)}_{p(q(i))}= \left[ \begin{array}{ccc} 
 1 & 1 & 1 \\ 0 & 1 & 1 \\ \mbox{\bf 0} & 0 & 1 
\end{array} \right] \ \ . $$
 Observe que $G(B)$ \'{e} fortemente conexo.  Em $B$ destacamos a
diagonal a ser posicionada por $P$ na diagonal principal de $PB$, e o
zero inicialmente na diagonal principal.  Em $PB$ destacamos o mesmo
zero originalmente na diagonal de $B$.  $G(PB)$ tem 3 CFCs, que
voltariam a fundir-se numa s\'{o}, caso adicion\'{a}ssemos a aresta
correspondente ao zero destacado em $PB$.  Neste exemplo tivemos
$p=[3,1,2]$, $q=[1,3,2]$.

O {\bf procedimento P4} (Parti\c{c}\~{a}o e Pr\'{e}-Posicionamento de 
Piv\^{o}s) explora esparsidade estrutural, e posteriormente
a esparsidade local a cada bloco, como segue:
\begin{enumerate} 
\item Encontre uma permuta\c{c}\~{a}o pr\'{o}pria, $PA$, atrav\'{e}s do 
 algoritmo h\'{u}ngaro complementado com uma heur\'{\i}stica eficiente.
\item Encontre, atrav\'{e}s do algoritmo de Tarjan, um 
 ordenamento coerente $Q'PAQ$. 
\item Inverta este \'{u}ltimo ordenamento, $QPAQ'$, de modo a colocar
 a matriz na forma triangular blocada inferior. Em seguida aplique 
 a heur\'{\i}stica P3 a cada bloco diagonal.
\end{enumerate}
 \index{P4} 
 \index{Algoritmo!Tarjan} 
 \index{Algoritmo!H\'{u}ngaro} 

S\~{a}o apresentados tr\^{e}s exemplos de aplica\c{c}\~{a}o do P4.  
A disposi\c{c}\~{a}o original dos ENNs da matriz original, $A$,
corresponde aos sinais $+$, $*$, ou n\'{u}meros de 1 a 9 no corpo da
matriz.  Como prescrito no primeiro passo do P4, primeiramente
encontramos uma diagonal n\~{a}o nula, ou equivalentemente uma
permuta\c{c}\~{a}o pr\'{o}pria $PA$.  O vetor inverso de \'{\i}ndices de
linha permutados, $\bar p$, \'{e} dado \`{a} esquerda da
numera\c{c}\~{a}o original das linhas. 
Os asteriscos indicam os elementos desta diagonal. 

Como prescrito no segundo passo do P4, devemos em seguida encontrar um
reordenamento coerente, $q$, em $G(PA)$.  Para tanto aplicamos o
algoritmo de Tarjan: Primeiramente fazemos a busca em profundidade
can\^{o}nica em $G(PA)$.  Para tanto percorremos as linhas de $PA$,
gerando a primeira floresta de cada exemplo.  As ra\'{\i}zes desta busca
s\~{a}o assinaladas por um acento circunflexo, tanto na floresta como no
vetor $\bar p$.  Os n\'{u}meros no corpo da matriz correspondem a ordem
de visita\c{c}\~{a}o nesta busca.  A inversa da ordem de retorno \'{e}
dada pelo vetor $\bar b$.

O algoritmo de Tarjan exige em seguida a busca em profundidade
can\^{o}nica no grafo inverso reordenado por $b$, $G((BPAB)')$; na
verdade basta tomarmos as ra\'{\i}zes desta busca na ordem can\^{o}nica. 
Percorrendo as colunas da matriz, constru\'{\i}mos a segunda floresta em
cada exemplo, cujas ra\'{\i}zes est\~{a}o assinaladas por um acento
circunflexo, tanto na floresta como no vetor $\bar b$.  A ordem de
visita\c{c}\~{a}o nesta segunda busca em profundidade nos d\'{a} o
reordenamento coerente $QPAQ'$; apresentado explicitamente em cada
exemplo para que se reconhe\c{c}a a estrutura triangular superior. 

Finalmente, ao final dos exemplos, apresentamos o terceiro passo do P4:
a invers\~{a}o deste ordenamento coerente, $Q'PAQ$, (portanto triangular
inferior), com posterior aplica\c{c}\~{a}o do P3 a cada bloco diagonal. 
Neste ponto fica claro que o segundo e o terceiro exemplos diferem
apenas pela permuta\c{c}\~{a}o original da matriz. Neste ponto indicamos 
tamb\'{e}m os zeros preenchidos durante a fatora\c{c}\~{a}o.

\pagebreak  

Exemplo 2:
 
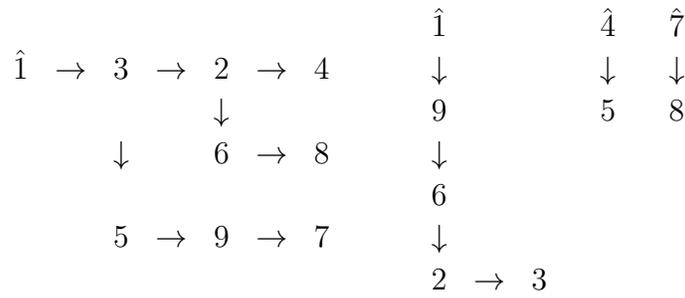
$$
\begin{array}{cccc||ccccccccc}
 & & & \bar q & 1 & 9 & 4 & 7 & 2 & 3 & 5 & 8 & 6 \\
 & & & \bar b & {\hat 1} & 9 & 3 & {\hat 7} & 4 & 2 & 6 & {\hat 8} & 5 \\
 \bar q \circ \bar p & \bar b \circ \bar p & \bar p &  
  & 1 & 2 & 3 & 4 & 5 & 6 & 7 & 8 & 9 \\
\hline \hline
 7 & 7 & 4 & 1 & & & & * & & & &  6 & \\
 2 & 4 & 5 & 2 & + & & & & * & & 4 & & 8 \\
 6 & 5 & 9 & 3 & & + & & & & & + & & * \\
 1 & 1 & \hat 1 & 4 & * & & 2 & & & 9 & + & & \\
 8 & 8 & 8 & 5 & & 7 & & & & & & * & \\
 5 & 6 & 7 & 6 & + & & + & 5 & & & * & & \\
 3 & 2 & 6 & 7 & & & + & & + & * & & & + \\
 4 & 3 & 3 & 8 & & & * & & 3 & & + & & \\
 9 & 9 & 2 & 9 & & * & & & & & & + &  \\
\end{array} $$
 
$$
\begin{array}{ccccccccccccc}
\hat 1 &\rightarrow & 3 &\rightarrow & 5 &\rightarrow & 7 &\rightarrow & 4 &\rightarrow & 8 &\rightarrow & 2 \\
   \downarrow    & &   & & \downarrow  & &   & &   & &   & &   \\ 
   6   & &   & & 9  
\end{array}$$ 

$$\begin{array}{cccccccc} 
\hat 1 & &   & & & \hat 7 & & \hat 8 \\
    \downarrow   & &   & & &        & &     \downarrow   \\ 
    4  &\rightarrow & 3 & & &        & &     9  \\ 
    \downarrow   & & \downarrow  \\ 
    2  & & 6 &\rightarrow & 5   
\end{array}$$ 
 
$$
\begin{array}{c||cccccc|c|cc}
p/q & 1 & 2 & 3 & 4 & 5 & 6 & 7 & 8 & 9 \\
\hline \hline
 1 & * & & + & + & + & & & & \\
 2 & + & * & & & + & + & & & \\
 3 & & + & * & + & & + & & & \\
 4 & & + & & * & + & & & & \\
 5 & + & & & + & * & & + & & \\
 6 & & & & & + & * & & & + \\
\hline
 7 & & & & & & & * & + & \\
\hline
 8 & & & & & & & & * & + \\
 9 & & & & & & & & + & * \\
\end{array} $$

\pagebreak 
 
Exemplo 3:
 
$$
\begin{array}{cccc||ccccccccc}
 & & & \bar q & 1 & 3 & 4 & 2 & 5 & 8 & 7 & 9 & 6 \\
 & & & \bar b & {\hat 1} & 6 & 2 & 9 & 3 & {\hat 7} & 5 & 8 & {\hat 4} \\
 \bar q \circ \bar p & \bar b \circ \bar p & \bar p &  
  & 1 & 2 & 3 & 4 & 5 & 6 & 7 & 8 & 9 \\
\hline \hline
 1 & 1 & {\hat 1} & 1 & * & & 2 & + & & & & & \\
 9 & 8 & 8 & 2 & & & & & & + & & * & \\
 6 & 4 & 9 & 3 & & & & & & & 9 & + & * \\
 3 & 6 & 2 & 4 & & * & & 4 & & 5 & & & \\
 2 & 9 & 4 & 5 & + & & + & * & & & & & \\
 7 & 5 & 7 & 6 & & & & & & & * & & + \\
 8 & 7 & 6 & 7 & & & & & & * & & 6 & \\
 5 & 3 & 5 & 8 & & & + & & * & & & & 8 \\
 4 & 2 & 3 & 9 & & 3 & * & + & 7 & & & & \\
\end{array} $$

$$
\begin{array}{ccccccc} 
 \hat  1 & \rightarrow  & 3 & \rightarrow  & 2 & \rightarrow  & 4 \\
        &   &   &   & \downarrow  &   &   \\ 
        &   & \downarrow  &   & 6 & \rightarrow  & 8 \\ 
        &   &   &   &   &   &   \\ 
        &   & 5 & \rightarrow  & 9 & \rightarrow  & 7  
\end{array} \hspace{1cm}    
\begin{array}{ccccccc}
 \hat 1 & & & & \hat 4 & & \hat 7 \\
    \downarrow    & & & &    \downarrow    & &   \downarrow      \\ 
    9   &  &  & & 5 & & 8 \\ 
    \downarrow  \\ 
    6 \\
    \downarrow  \\ 
    2 & \rightarrow  & 3 \\ 
\end{array} $$
 
$$
\begin{array}{c||ccccc|cc|cc}
p/q & 1 & 2 & 3 & 4 & 5 & 6 & 7 & 8 & 9 \\
\hline \hline
 1 & * & + & & + & & & & & \\
 2 & + & * & & + & & & & & \\
 3 & & + & * & & & & & + & \\
 4 & & + & + & * & + & & & & \\
 5 & & & & + & * & + & & & \\
\hline
 6 & & & & & & * & + & & + \\
 7 & & & & & & + & * & & \\
\hline
 8 & & & & & & & & * & + \\
 9 & & & & & & & & + & * \\
\end{array} $$
 
\pagebreak 

Exemplo 4: 

$$
\begin{array}{cccc||ccccccccc}
 & & & \bar q & 6 & 1 & 4 & 5 & 8 & 2 & 7 & 3 & 9 \\
 & & & \bar b & \hat 6 & \hat 1 & 2 & 5 & \hat 8 & 4 & 7 & 3 & 9 \\
 \bar q \circ \bar p & \bar b \circ \bar p & \bar p &  
  & 1 & 2 & 3 & 4 & 5 & 6 & 7 & 8 & 9 \\
\hline \hline
 1 & 1 & \hat 2 & 1 &   & * & 6 &   &   & + &   &   &   \\ 
 8 & 8 &      5 & 2 &   &   &   &   & * &   &   &   & 4 \\  
 7 & 7 &      7 & 3 & + &   &   &   & 3 &   & * &   &   \\  
 3 & 3 &      8 & 4 &   &   &   &   &   & + &   & * & + \\  
 2 & 4 &      6 & 5 &   & + & + &   &   & * &   &   &   \\  
 6 & 6 & \hat 1 & 6 & * &   &   &   &   &   & 2 &   &   \\  
 9 & 9 &      9 & 7 &   &   &   &   & + &   &   &   & * \\  
 5 & 5 &      4 & 8 &   &   & + & * &   &   & + &   &   \\  
 4 & 2 &      3 & 9 &   &   & * & 7 &   & 8 &   & 9 &   \\  
\end{array}$$ 

$$
\begin{array}{ccccccc} 
\hat 1 &\rightarrow & 7 &\rightarrow & 5 &\rightarrow & 9 \\ 
 \\ 
\hat 2 &\rightarrow & 3 &\rightarrow & 4 \\ 
       & & \downarrow  &\searrow &    \\ 
       & & 6 & & 8 
\end{array}  \hspace{1cm} 
\begin{array}{ccccc} 
\hat 1 & & \hat 6 & & \hat 8 \\ 
    \downarrow   & &     \downarrow   & &     \downarrow   \\ 
    4  & &     7  & &     9  \\ 
    \downarrow   \\ 
    3  \\ 
    \downarrow   \\ 
    2  &\rightarrow & 5 
\end{array}$$ 

$$
\begin{array}{c||ccccc|cc|cc}
p/q & 1 & 2 & 3 & 4 & 5 & 6 & 7 & 8 & 9 \\
\hline \hline
 1 & * & + &   & + &   &   &   &   &   \\
 2 & + & * &   & + &   &   &   &   &   \\
 3 &   & + & * &   &   &   &   &   & + \\
 4 &   & + & + & * & + &   &   &   &   \\
 5 &   &   &   & + & * &   & + &   &   \\
\hline 
 6 &   &   &   &   &   & * & + &   &   \\
 7 &   &   &   &   &   & + & * & + &   \\
\hline 
 8 &   &   &   &   &   &   &   & * & + \\
 9 &   &   &   &   &   &   &   & + & *   
\end{array}$$

\pagebreak 

Exemplos 2 a 4; Passo 3 do P4: 

$$
\begin{array}{c||cc|c|cccccc}
p/q & 1 & 2 & 3 & 4 & 5 & 6 & 7 & 8 & 9 \\
\hline \hline
 1 & * & + & & & & & & & \\
 2 & + & * & & & & & & & \\
\hline
 3 & & + & * & & & & & & \\
\hline
 4 & + & 0 & & * & & & & & + \\
 5 & & & + & & * & & + & & + \\
 6 & & & & & & * & + & & + \\
 7 & & & & + & + & + & 0 & & + \\
 8 & & & & + & + & + & 0 & * & 0 \\
 9 & & & & & + & & + & + & * \\
\end{array} $$

$$
\begin{array}{c||cc|cc|ccccc}
p/q & 1 & 2 & 3 & 4 & 5 & 6 & 7 & 8 & 9 \\
\hline \hline
 1 & * & + & & & & & & & \\
 2 & + & * & & & & & & & \\
\hline
 3 & & & * & + & & & & & \\
 4 & + & 0 & + & * & & & & & \\
\hline
 5 & & & & + & * & & & & + \\
 6 & & & & & + & * & + & & + \\
 7 & & + & & & & + & * & & 0 \\
 8 & & & & & & & + & * & + \\
 9 & & & & & & & + & + & * \\
\end{array} $$

$$
\begin{array}{c||cc|cc|ccccc}
p/q & 1 & 2 & 3 & 4 & 5 & 6 & 7 & 8 & 9 \\
\hline \hline
 1 & * & + & & & & & & & \\
 2 & + & * & & & & & & & \\
\hline
 3 & & + & * & + & & & & & \\
 4 & & & + & * & & & & & \\
\hline
 5 & & & + & 0 & * & & & & + \\
 6 & & & & & + & * & + & & + \\
 7 & + & 0 & & & & * & + & & 0 \\
 8 & & & & & & & + & * & + \\
 9 & & & & & & & + & + & * 
\end{array} $$

\pagebreak

\section{Estrutura Angular Blocada}

 Suponhamos dada uma matriz quadrada e n\~{a}o singular na forma angular
blocada, com blocos diagonais $^1B,\ldots {^hB}$, $^kB\ m(k)\times n(k)$, 
$d(k)\equiv m(k)-n(k)\geq 0$, e os correspondentes blocos nas
colunas residuais $^1C\ldots {^hC}$, $^kC\ m(k)\times n(h+1)$. 
 Podemos usar rota\c{c}\~{o}es de Givens para fatorar cada um dos blocos
diagonais,
 $^kB= {^hQ}\ \left[ \begin{array}{c} {^kV}\\ 0 \end{array} \right]$, 
 onde $^kV$ \'{e} triangular superior $n(k)\times n(k)$, e o bloco de
zeros \'{e} $d(k)\times n(k)$. 
 \index{Estrutura!angular} 

$$
\left[ \begin{array}{cccc} 
^1B & & & ^1C \\ & \ddots & & \vdots \\ & & ^hB & ^hC \\ 
\end{array} \right] \ \ , \ \ \ 
\left[ \begin{array}{cccc}
^1V & & & ^1W \\ 0 & & & ^1Z \\ 
 & \ddots & & \vdots \\ 
 & & ^hV & ^hW \\ & & 0 & ^hZ 
\end{array} \right] \ .$$  

Para completar a fatora\c{c}\~{a}o QR da matriz original, respeitando a
estrutura de blocos, permutamos os blocos $^kZ$ para as \'{u}ltimas linhas
da matriz, formando o bloco quadrado $Z$ de dimens\~{a}o 
$\sum _1^h d(k) = n(h+1)$.  Finalmente completamos a fatora\c{c}\~{a}o 
QR do bloco sudeste, $Z=QS$. 
 \index{Fatora\c{c}\~{a}o!blocada} 

$$
\left[ \begin{array}{cccc}
^1V & & & ^1W \\
 & \ddots & & \vdots \\  
 & & ^hV & ^hW \\
 & & & ^1Z \\ 
 & & & \vdots \\ 
 & & & ^hZ 
\end{array} \right] \ \ , \ \ \  
\left[ \begin{array}{cccc}
^1V & & & ^1W \\
 & \ddots & & \vdots \\  
 & & ^hV & ^hW \\
 & & & S 
\end{array} \right] \ .$$   

Como permuta\c{c}\~{o}es s\~{a}o apenas um tipo especial de
transforma\c{c}\~{o}es ortogonais, obtivemos o fator triangular da
fatora\c{c}\~{a}o QR respeitando a estrutura diagonal blocada da matriz
original, $B=QU$.  A inversa da matriz original seria dada por
$B^{-1}=U^{-1}Q'$, mas usando que $Q'=U^{-t}B'$ temos
$B^{-1}=U^{-1}U^{-t}B'$.  Isto \'{e}, obtivemos uma fatora\c{c}\~{a}o de
$B$ onde todos os fatores herdam (e s\~{a}o computados de acordo com) a
estrutura angular blocada da matriz original. 
 \index{Permuta\c{c}\~{a}o} 
 \index{Matriz!ortogonal} 

Observando que $B'B=U'Q'QU=U'U$, vemos que uma maneira alternativa de
computar o fator triangular da fatora\c{c}\~{a}o ortogonal de $B$, \'{e}
computar o fator de Cholesky  da matriz simetrizada $B'B$:
$$
\left[ \begin{array}{cccc}
 {^1B'}{^1B} & & & {^1B'}{^1C} \\
  & \ddots & & \vdots \\
 & & {^hB'}{^hB} & {^hB'}{^hC} \\
 {^1C'}{^1B} & \ldots & {^hC'}{^hB} & {^0Z} 
\end{array} \right] \ .$$
Ao eliminarmos os $h$ blocos das linhas residuais de $B'B$, formamos o
bloco sudeste $Z= {^0Z}+ \sum_1^h {^kZ}$, a ser fatorado na \'{u}ltima etapa
do processo, $Z=S'S$.  Ao final, obtemos exatamente o fator triangular
da fatora\c{c}\~{a}o QR da matriz original.

\section{Parti\c{c}\~{a}o de Hipergrafos}

Um {\bf hipergrafo} \'{e} um par ordenado $G=(V,C)$ onde o primeiro
elemento \'{e} um conjunto finito, o conjunto de v\'{e}rtices, e o
segundo elemento \'{e} uma matriz booleana, a {\bf matriz de
incid\^{e}ncia}.  Se $|V|=m$, cada coluna de $C$, $m\times n$, nos
d\'{a} os v\'{e}rtices sobre os quais incide o correspondente
(hiper)lado.  No caso particular de todos as colunas terem exatamente 2
ENN's temos na matriz de incid\^{e}ncia mais uma representa\c{c}\~{a}o
de um grafo sim\'{e}trico.  Em hipergrafos todavia um lado gen\'{e}rico
pode incidir sobre mais de 2 v\'{e}rtices. 
 \index{Hipergrafo} 
 \index{Parti\c{c}\~{a}o} 
 \index{Matriz!incid\^{e}ncia} 

O problema de permutar $PAQ$ para forma angular blocada pode ser visto
como um problema de parti\c{c}\~{a}o no hipergrafo $G=(M,B(A))$,
$M=\{1,\ldots m\}$, $B(A)$ a matriz booleana associada \`{a} matriz $A$,
$m\times n$.  No problema de parti\c{c}\~{a}o damos a cada linha $i\in
M$ de $B(A)$ uma cor $p(i)\in H=\{1,\ldots h\}$. A cor de cada lado
\'{e} definida como o conjunto de cores dos v\'{e}rtices sobre os quais
este incide, $q(j)=\{p(i)\mid A_i^j\neq 0\}$.  Lados multicoloridos
correspondem a colunas residuais na forma angular blocada, e os lados de
cor $q(j)=k$ correspondem \`{a}s colunas no bloco angular formado pelos
ENN's $A_i^j \mid p(i)=k \wedge q(j)=\{k\}$, $k\in H=\{1,\ldots h\}$. 

Para definir o problema de parti\c{c}\~{a}o falta-nos uma fun\c{c}\~{a}o
objetivo a ser minimizada.  Em vista das aplica\c{c}\~{o}es j\'{a}
apresentadas, e outras a serem apresentadas no pr\'{o}ximo capitulo,
queremos ter:
\begin{itemize}
\item
Aproximadamente o mesmo n\'{u}mero de linhas em cada bloco.
\item
Poucas colunas residuais.
\end{itemize}
Com estes prop\'{o}sitos \'{e} natural considerar a seguinte
fun\c{c}\~{a}o de custo de uma dada parti\c{c}\~{a}o de linhas em cores
$p: M\mapsto H$:
$$ f(p)= c(p) + \alpha \sum_{k=1}^h (m/h -s(k))^2 \ ,$$
$$ \mbox{onde}\ \  c(p)= |\{j\in N \mid |q(j)|\geq 2\} | \  ,   
   \mbox{e}\ \ s(k)= |\{i\in M \mid p(i)=k \} | \ .$$

Mesmo casos especiais deste problema s\~{a}o NP-dif\'{\i}ceis.  Por
exemplo: Seja $A$ a matriz de incid\^{e}ncia de um grafo, $m$ um
m\'{u}ltiplo exato de $h=2$, e faca $\alpha$ suficientemente grande para
garantir que todos os blocos tenham exatamente $m/h$ linhas.  Este \'{e}
o problema exato de 2- parti\c{c}\~{a}o em grafos, e a vers\~{a}o de
reconhecimento deste problema \'{e} NP-Completa; veja problema ND14 em
[Garey79].  Um algoritmo de anulamento simulado com perturba\c{c}\~{o}es
m\'{e}tricas para resolver este problema \'{e} apresentado em [Stern92]. 
 \index{Anulamento Simulado} 

\section{Paralelismo}

Um dos fatores mais importantes no desenvolvimento de algoritmos \'{e} a
possibilidade de realizar v\'{a}rias etapas de um procedimento em
paralelo.  No restante deste cap\'{\i}tulo adaptaremos algumas das
fatora\c{c}\~{o}es anteriormente estudadas para as estruturas blocadas,
visando paralelizar etapas independentes.  Vejamos a seguir alguns
conceitos b\'{a}sicos de computa\c{c}\~{a}o paralela. 
 \index{Paralelismo} 
 
 Existem v\'{a}rios modelos te\'{o}ricos de computador paralelo, e
in\'{u}meras inst\^{a}ncias e implementa\c{c}\~{o}es destes modelos em
m\'{a}quinas reais.  O modelo mais simples \'{e} o de {\bf mem\'{o}ria
compartilhada}.  Neste modelo v\'{a}rios processadores, t\^{e}m acesso a
uma mem\'{o}ria comum.  Neste modelo, a descri\c{c}\~{a}o de uma
algoritmo paralelo envolve basicamente dois fatores:
 \begin{itemize}
\item Como distribuir o trabalho entre os processadores.
\item Como sincronizar as diversas etapas do algoritmo.
\end{itemize}
 Este modelo \'{e} conceitualmente simples e elegante; todavia
limita\c{c}\~{o}es da nossa tecnologia inviabilizam a constru\c{c}\~{a}o
de m\'{a}quinas de mem\'{o}ria compartilhada com mais de uns poucos (da
ordem de dez) processadores. 
 \index{Mem\'{o}ria!compartilhada} 

O {\bf modelo de rede} \'{e} mais gen\'{e}rico.  Nele o computador \'{e}
visto como um grafo: cada v\'{e}rtice, {\bf n\'{o}}, ou {\bf
processador} representa: um processador propriamente dito, uma {\bf
mem\'{o}ria local}, i.e., acess\'{\i}vel somente a este processador, e
portas de comunica\c{c}\~{a}o.  Cada aresta representa uma {\bf via de
comunica\c{c}\~{a}o} inter-n\'{o}s.  Note que no modelo de mem\'{o}ria
compartilhada, a comunica\c{c}\~{a}o entre os processadores podia ser
feita de maneira trivial atrav\'{e}s da mem\'{o}ria; todavia no modelo
de rede \'{e} preciso saber os detalhes da arquitetura da m\'{a}quina
para especificar um terceiro aspecto do algoritmo:
 \begin{itemize} 
 \item A comunica\c{c}\~{a}o entre os processadores.  
 \end{itemize} 
 Estes detalhes incluem a disposi\c{c}\~{a}o das vias de
comunica\c{c}\~{a}o, ou {\bf topologia}, a velocidade de
comunica\c{c}\~{a}o em rela\c{c}\~{a}o a velocidade de processamento, a
possibilidade ou n\~{a}o de haver comunica\c{c}\~{o}es simult\^{a}neas
em vias distintas, etc.  Algumas destas topologias comumente
empregadas, s\~{a}o: {\bf Estrela}, {\bf Barra}, {\bf Anel}, {\bf Grade}, 
{\bf Toro}, {\bf Hipercubo} e {\bf Borboleta}. 
  \index{Mem\'{o}ria!distribuida} 
  \index{Comunica\c{c}\~{a}o}   

Como exemplo de algoritmo paralelo, calculemos a m\'{e}dia de $n$
n\'{u}meros numa rede com $p$ processadores.  Suponhamos que
inicialmente tenhamos $n/p$ destes n\'{u}meros em cada uma das
mem\'{o}rias locais.  Por simplicidade suponhamos que $n$ \'{e} um
m\'{u}ltiplo de $p$, e que $n>>p$.  Na primeira fase do algoritmo cada
processador, $k$, calcula a m\'{e}dia dos $n/p$ n\'{u}meros em sua
mem\'{o}ria local, $m(p)$.  Se os processadores s\~{a}o todos iguais
(rede homog\^{e}nea), cada processador completa sua tarefa em $1/p$ do
tempo necess\'{a}rio para calcular a m\'{e}dia geral num computador com
apenas um processador deste mesmo tipo.  Na segunda fase do algoritmo
reunimos as m\'{e}dias parciais para calcular a m\'{e}dia geral,
$m(0)=(1/p)\sum_{k=1}^p m(k)$.  Examinemos como calcular esta m\'{e}dia
geral em duas redes com topologia de anel, onde cada qual:
 \begin{enumerate}
 \item N\~{a}o permite comunica\c{c}\~{o}es em paralelo.
 \item Permite comunica\c{c}\~{o}es em paralelo via segmentos de arco 
      n\~{a}o superpostos.
 \end{enumerate}
Novamente por simplicidade, suporemos que $p=2^q$. 

Na primeira rede, para $k=1:p$, 
\begin{enumerate}
\item Calcule no, n\'{o} $k$,  
      $s(k)=s(k-1)+m(k)$.
\item Transmita $s(k)$ ao processador $k+1$. 
\end{enumerate}
 Neste procedimento temos as somas parciais das m\'{e}dias 
 $s(k)=\sum_{i=1}^k m(k)$, a condi\c{c}\~{a}o de inicializa\c{c}\~{a}o
\'{e} $s(0)=0$, e ao t\'{e}rmino do algoritmo podemos computar, no n\'{o}
$p$, a m\'{e}dia $m(0)=s(p)/p$.

Na segunda rede, para $i=1:q$, 
\begin{itemize}
\item  $j=2^i$;  em paralelo, para $k=j:p$, 
\begin{enumerate}
\item  Calcule, no n\'{o} $k$,  
       $r(i,k)=r(i-1,k)+r(i-1,k-j/2)$.
\item Transmita $r(i,k)$ do n\'{o} $k$ para o n\'{o} $k+j$.
\end{enumerate} 
\end{itemize} 
Neste procedimento temos as somas parciais das m\'{e}dias $r(i,k)=
\sum_{l=k-j+1}^k m(l)$, a condi\c{c}\~{a}o de inicializa\c{c}\~{a}o
\'{e} $r(0,k)=m(k)$, e ao t\'{e}rmino do algoritmo podemos computar, 
no n\'{o} $p$, a m\'{e}dia $m(0)=r(q,p)/p$. 

Em geral, a transmiss\~{a}o de dados entre n\'{o}s \'{e} muito mais
lenta que a manipula\c{c}\~{a}o destes dados localmente e, ao medir a
complexidade de um algoritmo, contamos separadamente os trabalhos de
processamento e comunica\c{c}\~{a}o.

\section{Fatora\c{c}\~{o}es Blocadas}

Das se\c{c}\~{o}es anteriores vemos que quase todo o trabalho na
fatora\c{c}\~{a}o QR de uma matriz angular blocada, $A=QU$, ou da
fatora\c{c}\~{a}o de Cholesky $A'A=U'U$, consiste na aplica\c{c}\~{a}o
repetitiva de algumas opera\c{c}\~{o}es simples sobre os blocos.  Para
tirar vantagem desta modularidade em algoritmos para fatora\c{c}\~{a}o e
atualiza\c{c}\~{a}o de matrizes com estrutura angular blocada, definimos
a seguir algumas destas opera\c{c}\~{o}es, e damos sua complexidade em
n\'{u}mero de opera\c{c}\~{o}es de ponto flutuante. 
  \index{Fatora\c{c}\~{a}o!blocada} 

\begin{enumerate}
\item Compute a {\em fatora\c{c}\~{a}o de Cholesky parcial}, 
      eliminando as primeiras $n$ colunas da matriz blocada
    \[  \left[ \begin{array}{cc} F & G \\ G^{t} & 0 \\ \end{array} \right] \]
     para obter
    \[  \left[ \begin{array}{cc} V & W \\ 0 & Z \\ 
\end{array}  \right]  \]
     onde $F=F^{t}$ \'{e} $n\times n$, e $G$ \'{e} $n\times l$.  
     Isto requer
     $(1/6)n^{3} + (1/2)n^{2}l + (1/2)nl^{2} + O(n^{2}+l^{2})$ FLOPs.

\item Compute a {\em transforma\c{c}\~{a}o inversa parcial}, i.e. $u$, em
      \[  \left[ \begin{array}{cc} V & W \\ O & I \\ \end{array}  \right]^{t} 
          \left[ \begin{array}{c} u^{1} \\ u^{2} \\ \end{array}  \right] =
          \left[ \begin{array}{c} y^{1} \\ y^{2} \\ \end{array}  \right]  \]
      onde $V$ \'{e} $n\times n$ triangular superior, $W$ \'{e} $n\times l$,
      $0$ e $I$ s\~{a}o as matrizes zero e identidade, e $u$ e $y$
      s\~{a}o vetores coluna.
      Isto requer $(1/2)n^{2} + nl + O(n+l)$ FLOPs.

\item {\em Reduza a triangular superior} uma matriz de Hessenberg, i.e.,
      aplique a seq\"{u}\^{e}ncia de rota\c{c}\~{o}es de Givens  
      $G(1,2,\theta ),G(2,3,\theta ) \ldots G(n-1,n,\theta )$ 
      \`{a} matriz blocada 
      \[  \left[ \begin{array}{cc} V & W \\ \end{array}  \right]  \]
      onde $V$ \'{e} $n\times n$ Hessenberg superior, e $W$ \'{e} $n\times l$,
      para reduzir $V$ a triangular superior.
      Isto requer $2n^{2} + 4nl + O(n^{2}+l^{2})$ FLOPs.

\item {\em Reduza a triangular superior} uma matriz blocada 
      coluna -- tri\^{a}ngulo superior, i.e., 
      aplique a seq\"{u}\^{e}ncia de rota\c{c}\~{o}es de Givens 
      $G(n-1,n,\theta ), G(n-2,n-1,\theta ), \ldots G(1,2,\theta )$ 
      \`{a} matriz blocada 
      \[  \left[ \begin{array}{cc} u & V \\ \end{array}  \right]  \] 
      onde $u$ \'{e} um vetor coluna $n\times 1$, e $V$ \'{e} $n\times n$ 
      triangular superior, de modo a reduzir $u$ \`{a} um \'{u}nico ENN na 
      primeira linha, assim transformando $V$ de triangular para 
      Hessenberg superior. 
      Isto requer $2n^{2} + O(n)$ FLOPs.
\end{enumerate}

Em ambas as fatora\c{c}\~{o}es, $A=QU$ e $A'A=U'U$, muitas das
opera\c{c}\~{o}es nos blocos podem ser feitas independentemente. 
Portanto a estrutura angular blocada n\~{a}o s\'{o} nos d\'{a} a
possibilidade de preservar esparsidade, mas tamb\'{e}m a oportunidade de
fazer v\'{a}rias opera\c{c}\~{o}es em paralelo. 

Descreveremos uma forma de paralelizar a fatora\c{c}\~{a}o de Cholesky
numa rede de $h+1$ n\'{o}s.  Para $k=1\ldots h$ alocamos blocos das
matrizes $A$ e $U$ a n\'{o}s espec\'{\i}ficos, como segue:

\begin{itemize}
\item Os blocos $D^{k}$ $E^{k}$, $V^{k}$ e $W^{k}$ s\~{a}o 
      alocados ao n\'{o} $k$.
\item Os blocos sudeste, $Z$ e $S$, s\~{a}o alocados ao n\'{o} $0$ 
     (ou $h+1$).
\end{itemize}

Expressaremos a complexidade do algoritmo em termos da soma e do
m\'{a}ximo das dimens\~{o}es dos blocos. 
\[ dbsum = \sum_{1}^{h} m(k) \]
\[ dbmax = \max \{m(1), \ldots, m(h), n(h+1)\} .  \]
 
Na an\'{a}lise de complexidade contabilizaremos o tempo de processamento,
medido em FLOPs, $pTime$, bem como a comunica\c{c}\~{a}o inter-n\'{o}s,
$INC$.  Quando $h$ opera\c{c}\~{o}es sobre blocos, 
$bop^{1} \ldots bop^{h}$, podem ser efetuadas em paralelo (em n\'{o}s
distintos), contamos seu tempo de processamento por 
$\wedge_{1}^{h} flops(bop^{k}) = \ $ 
$flops(bop^{1}) \wedge \ldots \wedge flops(bop^{h})$ \ , onde $\wedge$ 
\'{e} o operador {\em m\'{a}ximo}, e $flops(bop^{k})$ \'{e} o numero de
opera\c{c}\~{o}es de ponto flutuante necess\'{a}rio para a
opera\c{c}\~{a}o no bloco $k$, $bop(k)$.  Nas equa\c{c}\~{o}es que
seguem, $\wedge$ tem preced\^{e}ncia menor que qualquer operador
multiplicativo ou aditivo.  As express\~{o}es ``No n\'{o} {\it k=1:h}
compute'' ou ``Do n\'{o} {\it k=1:h} envie'' significam, ``Em (de) todos
os n\'{o}s $1\le k \le h$, em paralelo, compute (envie)''.  Nas
express\~{o}es de complexidade ignoraremos termos de ordem inferior. 

Damos agora uma descri\c{c}\~{a}o algor\'{\i}tmica da fatora\c{c}\~{a}o
de Cholesky blocada $bch()$:

\begin{enumerate}

\item No n\'{o} {\it k=1:h} compute os blocos 
 $(B^{k})^{t}B^{k}$, $(B^{k})^{t}C^{k}$, e $(C^{k})^{t}C^{k}$.\\ 
 $pTime=m(k)n(k)^2 + m(k)n(k)n(h+1) + m(k)n(h+1)^{2} \le 3dbmax^{3}$,\\ 
 $INC=0$. 

\item Envie $(C^{k})^{t}C^{k}$ do n\'{o} $k$ para o n\'{o} $0$, 
 onde acumulamos 
 $Z^{0}=\sum_{1}^{h} (C^{k})^{t}C^{k}$.
 $pTime = h\ n(h+1)^{2} \le h\ dbmax^{2}$ ,  
 $INC = h\ n(h+1)^{2} \le h\ dbmax^{2}$

\item No n\'{o} $k$ compute a fatora\c{c}\~{a}o de Cholesky parcial,
 eliminando as primeiras $n(k)$ colunas, da matriz blocada
   \[  \left[ \begin{array}{cc}
                  (B^{k})^{t}B^{k} & (B^{k})^{t}C^{k} \\
                  (C^{k})^{t}B^{k} & 0 \\
                 \end{array}  \right] \]
 obtendo
   \[  \left[ \begin{array}{cc}
                  V^{k} & W^{k} \\
                  0 & Z^{k} \\
                 \end{array}  \right]  \]
$pTime=(1/6)n(k)^{3}+(1/2)n(k)^{2}n(h+1)+(1/2)n(k)n(h+1)^{2}\le
 (7/6)dbmax^{3}$, \\$INC=0$.

\item Envie $Z^{k}$ do n\'{o} $k$ para o n\'{o} $0$, onde acumulamos 
  $Z = \sum_{0}^{h} Z^{k}$.\\ 
  $pTime = h\ n(h+1)^{2} \le h\ dbmax^{2}$ ,  
  $INC = h\ n(h+1)^{2} \le h\ dbmax^{2}$.

\item No n\'{o} $0$ fatore o bloco sudeste $S=chol(Z)$,
  onde $chol()$ indica a fatora\c{c}\~{a}o de Cholesky padr\~{a}o.\\
  $pTime = (1/6)n(h+1)^{3} \le (1/6)dbmax^{3}$ , $INC=0$.
\end{enumerate}

\begin{teo}
A fatora\c{c}\~{a}o de Cholesky blocada, $bch()$, requer n\~{a}o mais de 
$(4+1/3)dbmax^{3} + h\ dbmax^{2}$ tempo de processamento, e 
$h\ dbmax^{2}$ tempo de comunica\c{c}\~{a}o inter-n\'{o}s.  
\end{teo}

Nos passos 2 e 4, se a rede permite comunica\c{c}\~{o}es em paralelo, a
reuni\~{a}o da matriz acumulada pode ser feita em $log(h)$ passos, e
podemos substituir $h$ por $log(h)$ no \'{u}ltimo teorema.

 \clearpage
 \clearpage 
\setcounter{chapter}{7} 
\chapter{ESCALAMENTO} 
\begin{center}
{\LARGE Representa\c{c}\~{a}o em Ponto Flutuante}
\end{center}

\section{O Sistema de Ponto Flutuante}

 A representa\c{c}\~{a}o de um n\'{u}mero real, $\zeta \in R$, em
um computador tem, usualmente, precis\~{a}o finita. 
 A inexatid\~{a}o desta representa\c{c}\~{a}o introduz erros no resultado
final do processamento e o objetivo desta se\c{c}\~{a}o \'{e} obter
limites m\'{a}ximos para estes erros quando da aplica\c{c}\~{a}o do
m\'{e}todo de Gauss. 
 A representa\c{c}\~{a}o normalmente utilizada para n\'{u}meros reais
\'{e} o sistema de representa\c{c}\~{a}o em {\bf ponto  flutuante 
normalizado} de $t$ d\'{\i}gitos e base $b$, SPF, isto \'{e},
 \index{SPF!ponto flutuante} 
 \index{SPF!precis\~{a}o simples}  
 \begin{eqnarray*}
fl(\zeta ) & = & \pm 0.d_1d_2\ldots d_t * b^{\pm n}, \mbox{ou}\\
 & & \pm 0.d_1d_2\ldots d_t\ E \pm n 
\end{eqnarray*}
onde
\begin{eqnarray*}
& & d_k , n, b \in  N \\
& & 0 \leq d_k \leq b,\ d_1\neq 0 \\
& & 0 \leq n \leq emax,\ b\neq 0\ .
\end{eqnarray*}

Dado um real $\nu$ com representa\c{c}\~{a}o exata num dado SPF, a
fra\c{c}\~{a}o normalizada ser\'{a} denominada mantissa. A {\bf mantissa}
e o {\bf expoente} de $\nu$ ser\~{a}o denotados, respectivamente,
 \index{SPF!mantissa} 
 \index{SPF!expoente} 
\begin{eqnarray*} 
mant(\nu ) & = & \pm 0.d_1\ldots d_t ,\\
expo(\nu ) & = & \pm n 
\end{eqnarray*}

H\'{a} duas maneiras normalmente empregados para obter $fl(\zeta )$ a
partir de $\zeta$: {\bf truncamento} e {\bf arredondamento}.  Para
$trunc(\zeta)$ simplesmente truncamos a representa\c{c}\~{a}o em base
$b$ de $\zeta$ obtendo uma mantissa de $t$ digitos.  No arredondamento
tomamos a mantissa de $t$ digitos que melhor aproxima $\zeta$.  Assim,
para $\zeta =\Pi = 3.141592653\ldots$, e o SPF com $b=10$ e $t=5$,
$trunc(\zeta )= +0.31415\ E+1$ e $round(\zeta )= +0.31416\ E+1$. 
 \index{SPF!truncamento} 
 \index{SPF!arredondamento} 

Note que $\zeta =0$ n\~{a}o pode ser representado devido a
condi\c{c}\~{a}o $d_1\neq 0$.  Portanto, alguma representa\c{c}\~{a}o
inamb\'{\i}gua deve ser convencionada, por exemplo, 
 $fl(0)= +0.00\ldots 0\ E+0$.  Se $n > emax$ n\~{a}o h\'{a}
representa\c{c}\~{a}o poss\'{\i}vel no SPF, e dizemos que houve
``overflow" ou transbordamento. 

\section{Erros no Produto Escalar}

Definimos a {\bf unidade de erro} do SPF, $u$, como
$b^{1-t}$ no caso de usarmos truncamento, e $b^{1-t}/2$ no caso 
de usarmos arredondamento.
 \index{Erros!unidade no SPF} 
 \index{Erros!produto escalar} 

\begin{lm} 
Dado $\zeta \in R$ e $fl(\zeta )$ sua representa\c{c}\~{a}o, num
dado SPF de unidade de erro $u$, 
$\exists \delta \in \ [-u,u]\ \mid fl(\zeta )=\zeta (1+\delta )$.
\end{lm} 

Demonstra\c{c}\~{a}o: 

Pela defini\c{c}\~{a}o de SPF, se $expo( fl(\zeta ))=e$, v\^{e}-se que
$$ \frac{| fl(\zeta ) - \zeta |}{| \zeta | } = | \delta | \leq 
   \frac{ ub^{e-1}}{b^{e-1}} = u\ .$$ 

Se $\nu$ e $\mu$ s\~{a}o n\'{u}meros em ponto flutuante, isto \'{e}  
n\'{u}meros reais com representa\c{c}\~{a}o exata num dado SPF, \'{e}
perfeitamente poss\'{\i}vel que uma opera\c{c}\~{a}o aritm\'{e}tica
elementar entre eles resulte num n\'{u}mero sem representa\c{c}\~{a}o
exata.  Do lema anterior, por\'{e}m, sabemos que, qualquer que seja a
opera\c{c}\~{a}o aritm\'{e}tica, $\star \in \{+,-,*,/\}$,
$fl( \nu \star \mu ) = (\nu \star \mu )(1+ \delta )$, 
para algum $\delta \mid |\delta |\leq u$. QED.  

Exemplo 1: 

Consideremos um computador que armazena um n\'{u}mero real em 4
bytes, sendo 3 bytes para os d\'{\i}gitos da mantissa e 1 byte para o
sinal do n\'{u}mero, o sinal do expoente e os 6 bits restantes para
o m\'{o}dulo do expoente como um inteiro em base 2.  Supondo que todos
os c\'{a}lculos s\~{a}o feitos com arredondamento, calculamos, a unidade
de erro do SPF, $u$, e o maior real represent\'{a}vel, $rmax=b^{emax}$, 
no caso de usarmos base $b=256$, $b=16$ ou $b=2$.

\begin{center}
\begin{tabular}{|l|l|l|l|} \hline   
$b$ & $t$ & $u=b^{t-1}/2$ & $rmax=b^{64}$ \\ \hline \hline  
$2$ & 24 & $2^{-24}< 10^{-7}$ & $2^{64}>10^{19}$ \\ \hline 
$16=2^4$ & 6 & $2^{-21}< 10^{-6}$ & $2^{256}>10^{75}$ \\ \hline 
$256=2^8$ & 3 & $2^{-17}< 10^{-5}$ & $2^{512}>10^{150}$ \\ \hline 
\end{tabular}
\end{center}

Um recurso freq\"{u}entemente dispon\'{\i}vel, concomitantemente ao uso
de um SPF de base $b$ e $t$ d\'{\i}gitos, \'{e} o sistema de
representa\c{c}\~{a}o em ponto flutuante normalizado de {\bf precis\~{a}o
dupla}, SPFD, com mantissa de $2t$ d\'{\i}gitos, que denotamos 
$fld(\zeta )$. 
 \index{SPF!precis\~{a}o dupla}  
 \index{Erros!produto escalar} 

Este recurso \'{e} extremamente \'{u}til, se utilizado parcimoniosamente.
Podemos trabalhar com a maior parte dos dados no SPF, de precis\~{a}o
simples, e utilizar a precis\~{a}o dupla apenas nas passagens mais
cr\'{\i}ticas do procedimento.
Em analogia ao SPF, definimos a unidade de erro do SPFD, $ud$ como
$b^{1-2t}$ no caso de usarmos truncamento, e $b^{1-2t}/2$ no caso 
de usarmos arredondamento.

\'{E} interessante notar que se $\nu$ e $\mu$ s\~{a}o n\'{u}meros reais
com representa\c{c}\~{a}o exata no SFP, de precis\~{a}o simples, a
representa\c{c}\~{a}o do produto destes n\'{u}meros no SPFD \'{e} exata,
isto \'{e} $fld(\nu *\mu )=\nu *\mu$, pois o produto de dois inteiros de
t d\'{\i}gitos tem, no m\'{a}ximo, $2t$ d\'{\i}gitos. 

Estudamos agora o efeito cumulativo dos erros de representa\c{c}\~{a}o,
em todas as passagens intermedi\'{a}rias num produto interno.
Sejam $x$, $1\times n$ e $y$, $n\times 1$, vetores cujas componentes
t\^{e}m representa\c{c}\~{a}o exata, num dado SPF. Definimos 
$$ fl(xy)= fl( fl(x_ny^n)+ fl( fl(x_{n-1}y^{n-1}) +\ldots +
           fl( fl(x_2y^2) + fl(x_1y^1) ) \ldots ))\ .$$

\begin{lm} 
Dados $u\in [0,1[$, $n\in {\bf N}$ tq $nu\leq \alpha <1$, e 
${\delta}_i \mid |{\delta}_i| \leq u,\ i=1\ldots n$, ent\~{a}o
$$ 1-nu\leq \prod_1^n (1+{\delta}_i)\leq 1+(1+\alpha )nu\ .$$
\end{lm} 

Demonstra\c{c}\~{a}o: 

Como $$(1-nu)^n\leq \prod_1^n (1+{\delta}_i)\leq (1+nu)^n\ ,$$
basta provar que
\begin{itemize}
\item $(1-u)^n\geq 1-nu$.
\item $(1+u)^n\leq  1+ nu+ \alpha nu$.
\end{itemize}

Considerando a fun\c{c}\~{a}o $f(u)=(1-u)^n$ temos que
$\exists \theta \in ]0,1[ \ \ \mid $
\begin{eqnarray*}
f(u) & = & f(0)+uf'(u)+u^2f''(\theta u)/2 \\ 
 & = &  1-nu+n(n-1)(1-\theta u)^{n-2}u^2/2\ . 
\end{eqnarray*}
Da n\~{a}o negatividade do \'{u}ltimo termo segue a primeira 
inequa\c{c}\~{a}o.

Considerando que $\forall \beta \in [0,\alpha ]$,  
$$ 1+\beta \leq e^\beta  \leq 1+\beta +\alpha \beta \ ,$$   
temos que
$$ (1+u)^n \leq e^(nu) \leq 1+nu+\alpha nu \ .$$ QED. 

Corol\'{a}rio: 

Nas condi\c{c}\~{o}es do lema 2, $\exists \theta \in [-1,1]$ tq
$$ \prod_1^n (1+{\delta}_i ) = 1+\theta (1+\alpha )nu \ .$$

\begin{lm} 
Se $x$, $1\times n$ e $y$, $n\times 1$, s\~{a}o vetores cujas
componentes tem representa\c{c}\~{a}o exata num dado SPF, de unidade de
erro $u$, com $n u <1$, ent\~{a}o
$$ fl(xy)= \sum_{i=1}^n x_iy^i (1+{\theta }_i (1+\alpha )(n-i+2)u) \ .$$ 
\end{lm} 

Demonstra\c{c}\~{a}o: 

$\exists {\delta }_i,\ |{\delta }_i|\leq ud$, e ${\theta }_i,\  
 |{\theta }_i|\leq 1 \ \ \mid $
\begin{eqnarray*} 
fl(xy) & = & fl( fl(x_ny^n)+ fl( fl(x_{n-1}y^{n-1}) +\ldots \\ & &  
  + (x_3y^3(1+{\delta }_4)+(x_2y^2(1+{\delta }_2) +  
  x_1y^1(1+{\delta }_1))(1+{\delta }_3))(1+{\delta }_5)\ldots ))  \\ 
 & = & fl( fl(x_ny^n)+ fl( fl(x_{n-1}y^{n-1}) +\ldots + 
  x_3y^3(1+{\delta }_4)(1+{\delta }_5)  \\ & & 
  + x_2y^2(1+{\delta }_2)(1+{\delta }_3)(1+{\delta }_5) +
  x_1y^1(1+{\delta }_1)(1+{\delta }_3)(1+{\delta }_5)\ldots )) \\ 
 & = & x_ny^n(1+{\delta }_{2n-2})(1+{\delta }_{2n-1}) + \\ & & 
  + x_{n-1}y^{n-1}(1+{\delta }_{2n-4})(1+{\delta }_{2n-3})(1+{\delta }_{2n-1}) 
  + \ldots \\ & & + x_2y^2(1+{\delta }_2)(1+{\delta }_3)\ldots 
  (1+{\delta }_{2n-3})(1+{\delta }_{2n-1}) 
  \\ & & + x_1y^1(1+{\delta }_1)(1+{\delta }_3)\ldots 
  (1+{\delta }_{2n-3})(1+{\delta }_{2n-1}) \\  
 & = & x_ny^n(1+{\theta }_n(1+\alpha )2u) + 
      x_{n-1}y^{n-1}(1+{\theta }_{n-1}(1+\alpha )3u) + \ldots \\ & & + 
      x_2y^2(1+{\theta }_2(1+\alpha )nu) +
      x_1y^1(1+{\theta }_1(1+\alpha )(n+1)u) \ . 
\end{eqnarray*}

\begin{lm} 
Se $x$, $1\times n$ e $y$, $n\times 1$, s\~{a}o vetores cujas
componentes t\^{e}m representa\c{c}\~{a}o exata num SPF, em precis\~{a}o
simples, de unidade de erro $u$, sendo $n*ud\leq \gamma < 1$, 
ent\~{a}o 
\begin{eqnarray*}
fld(xy) & \equiv & fld( fld(x_ny^n)+ fld( fld(x_{n-1}y^{n-1}) +\ldots +
           fld( fld(x_2y^2) + fld(x_1y^1) ) \ldots )) \\
 & = & \sum_1^n x_iy^i(1+\theta _i(1+\gamma)(n-i+1)ud) 
\end{eqnarray*}
\end{lm} 

Demonstra\c{c}\~{a}o: 

$\exists {\epsilon }_i,\ |{\epsilon }_i|\leq ud$, e ${\theta }_i,\  
 |{\theta }_i|\leq 1 \ \ \mid $
\begin{eqnarray*} 
fld(xy)&=& fld(x_ny^n+fld(x_{n-1}y^{n-1}+\ldots \\ & &  
  + (x_3y^3+ (x_2y^2+ x_1y^1)(1+\epsilon _1))(1+\epsilon _2)\ldots )) \\ 
&=& x_ny^n(1+\epsilon _{n-1}) +
    x_{n-1}y^{n-1}(1+\epsilon _{n-2})(1+\epsilon _{n-1})+\ldots \\ & & 
    x_2y^2(1+\epsilon _1)\ldots (1+\epsilon _{n-1}) +
    x_1y^1(1+\epsilon _1)\ldots (1+\epsilon _{n-1}) \\ 
&=& x_ny^n(1+\theta _n(1+\gamma )ud) +
    x_{n-1}y^{n-1}(1+\theta _{n-1}(1+\gamma )2ud) + \\ & & 
    x_2y^2(1+\theta _2(1+\gamma )(n-1)ud) +
    x_1y^1(1+\theta _1(1+\gamma )nud) \ \ 
\end{eqnarray*}

\'{E} freq\"{u}ente o c\'{a}lculo em dupla precis\~{a}o, de produtos de
vetores armazenados em precis\~{a}o simples, e o posterior armazenamento
do resultado em precis\~{a}o simples, isto \'{e} o c\'{a}lculo 
$fl(fld(xy))$. 

Do \'{u}ltimo lema vemos que
$$ fl(fld(xy)) = (1+\delta )\sum_1^n 
   x_iy^i (1+\theta _i(1+\gamma )(n-i+1)ud) \ .$$

\begin{obs}{\rm  
Se $xy \gg \sum x_iy^i \theta _i(1+\gamma )(n-i+1)ud$, o que ocorre se
$nu\ll 1$ e n\~{a}o houver ``cancelamentos cr\'{\i}ticos'' na
somat\'{o}ria, o resultado final do produto \'{e} afetado de um erro da
ordem do erro introduzido por uma \'{u}nica opera\c{c}\~{a}o
aritm\'{e}tica, em precis\~{a}o simples. 
Assumiremos esta hip\'{o}tese no restante do cap\'{\i}tulo.
}\end{obs} 

Exemplo 2: 

Calculemos, no SPF decimal de 2 d\'{\i}gitos com arredondamento,
$fl(xy)$, $fld(xy)$ e $fl(fld(xy))$, onde 
$x=[7.5,\ 6.9,\ 1.3]$ e $y=[0.38,\ -0.41,\ 0.011]'$.

$$fld(xy)=fld(2.85 - 2.829 + 0.0143)= fld(0.021 + 0.0143)= 0.0353\ .$$
$$fl(fld(xy))= fl(0.0353)= 0.035 \ .$$
$$fl(xy)= fl(2.9 - 2.8 - 0.014)= fl(0.1 + 0.014)= 0.11 \ .$$

\section{Escalamento de Matrizes}

Da an\'{a}lise de erros no produto interno, e considerando que 
a maioria das opera\c{c}\~{o}es nos algoritmos de triangulariza\c{c}\~{a}o 
podem ser agrupadas na forma de produtos internos, fica
evidente que \'{e} conveniente termos todos os elementos da matriz da
mesma ordem de grandeza, i.e., termos a matriz bem {\bf equilibrada}. 
Evitamos assim ter multiplicadores muito pequenos e soma de termos de
ordem de grandeza muito diferentes.  Se tal n\~{a}o ocorre para a 
matriz de um dado sistema, $Ax=b$, podemos procurar $x$ atrav\'{e}s 
da solu\c{c}\~{a}o de um outro sistema melhor equilibrado. 
 \index{Escalamento} 
 \index{Matriz!equilibrada} 
       
Se $E= diag(e_1,e_2,\ldots e_n)$ e $D= diag(d_1,d_2,\ldots d_n)$, onde
$e_i, d_i \neq 0$, o sistema $EADy=Eb$ \'{e} obtido multiplicando-se a
$i$-\'{e}sima equa\c{c}\~{a}o do sistema por $e_i$, e efetuando-se a
substitui\c{c}\~{a}o de vari\'{a}veis $x=Dy$, o que equivale a multiplicar a
j-\'{e}sima coluna de $A$ por $d_j$.  Uma transforama\c{c}\~{a}o $\tilde
A=EAD$, $A\ m\times n$, $E$ e $D$ diagonais com $e_i, d_i \neq 0$, \'{e}
um {\bf escalamento} da matriz $A$. 

Exemplo 3: 
 
$$ EAD = 
   \left[ \begin{array}{ccc} 
    1 & 0 & 0 \\ 0 & 2 & 0 \\ 0 & 0 & 3 
    \end{array} \right] 
   \left[ \begin{array}{ccc} 
    1 & 1 & 1 \\ 1 & 1 & 1 \\ 1 & 1 & 1 
    \end{array} \right] 
   \left[ \begin{array}{ccc} 
    11 & 0 & 0 \\ 0 & 12 & 0 \\ 0 & 0 & 13 
    \end{array} \right] 
   = 
   \left[ \begin{array}{ccc} 
    11 & 12 & 13 \\ 22 & 24 & 26 \\ 31 & 36 & 39 
    \end{array} \right] 
 $$ 

Estudaremos o problema de escolher um escalamento que
``melhor equilibre'' uma dada matriz.
 Em primeiro lugar, vale notar que, estando num SPF de base $b$, a
escolha de $E$ e $D$ da forma $E_i^i=b^{e_i}$.  $D_j^j=b^{d_j}$, onde
$e$ e $d$ s\~{a}o vetores de elementos inteiros, \'{e} muito
conveniente, pois $\tilde A=EAD$ e $A$ ter\~{a}o elementos de mesma
mantissa, sendo o efeito do escalamento apenas o de alterar os expoentes
dos elementos da matriz,
 $$ mant(\tilde A_i^j)= mant(A_i^j) \ \ , \ \ \ 
   expo(\tilde A_i^j)= expo(A_i^j)+e_i+d^j \ .$$

\begin{obs}{\rm  
Para implementar eficientemente a fun\c{c}\~{a}o $expo( )$ e a
opera\c{c}\~{a}o de soma de um inteiro ao expoente de um n\'{u}mero
real, devemos conhecer detalhadamente o SPF usado e manipular
diretamente os campos de bits envolvidos.  
}\end{obs} 

Uma maneira de medir o grau de desequilibrio de uma matriz,
$A$ \'{e} atravez da m\'{e}dia e da vari\^{a}ncia dos expoentes de seus 
elementos: 
 $$mex(A) = \sum_{i,j \mid A_i^j\neq 0} expo(A_i^j) / enn(A) \ ,$$
 $$vex(A) = \sum_{i,j \mid A_i^j\neq 0} (expo(A_i^j)- mex )^2 / enn(A) \ .$$
 Nas somat\'{o}rias que definem  $mex$ e $vex$ excluimos 
os elementos nulos da matriz, os elementos nulos da matriz podem ser 
eliminados das opera\c{c}\~{o}es de produto interno.

 Para n\~{a}o sobrecarregar a nota\c{c}\~{a}o, doravante escrevemos
 $$\sum_i{}' = \sum_{i=1\ \mid A_i^j \neq 0}^m \  , \ \ \ 
  \sum_j {}' = \sum_{j=1\ \mid A_i^j \neq 0}^n \  , \ \ \ 
  \sum_{i,j} {} ' = \sum_{i,j=1\ \mid A_i^j \neq 0}^{m,n} \ .$$

O primeiro m\'{e}todo de escalamento que estudaremos \'{e} justamente
o {\bf m\'{e}todo da redu\c{c}\~{a}o de vari\^{a}ncia} em que procuramos
minimizar a vari\^{a}ncia dos expoentes de $EAD$.

Tomemos as matrizes de escalamento esquerda e direita como, 
respectivamente, $E=diag(round(e_i^*))$ e $D=diag(round(d_j^*))$, 
sendo os vetores $e^*$ e $d^*$ argumentos que minimiz\~{a}o 
a vari\^{a}ncia dos expoentes da matriz escalada,   
$$vex(EAD) = \sum_{i,j}{}' (expo(A_i^j) +e_i +d_j -mex(A) )^2 \ .$$

Um ponto de m\'{\i}nimo deve obedecer ao sistema
\begin{eqnarray*}
\frac{\partial vex(EAD)}{\partial e_i^*} &= 0 =&
 2\sum_j{}' (expo(A_i^j) +e_i^* +d_j^* -mex(A) ) \\
\frac{\partial vex(EAD)}{\partial d_j^*} &= 0 =&
 2\sum_i{}' (expo(A_i^j) +e_i^* +d_j^* -mex(A) ) 
\end{eqnarray*}
ou, fazendo a substituicao $e_i^+= e_i^* -mex(A)/2$,
$d_j^+= d_j^* -mex(A)/2$
\begin{eqnarray*}
\sum_j{}' e_i^+ +d_j^+ &=& -\sum_j{}' expo(A_i^j) \\
\sum_i{}' e_i^+ +d_j^+ &=& -\sum_i{}' expo(A_i^j)
\end{eqnarray*}


Se a matriz A n\~{a}o tiver elementos nulos, ent\~{a}o a solu\c{c}\~{a}o
$e^*$, $d^*$ do sistema \'{e} imediata:
\begin{eqnarray*}
e_i^* &=& (1/n) \sum_j (mex(A) -expo(A_i^j)) \\
d_i^* &=& (1/m) \sum_i (mex(A) -expo(A_i^j)) 
\end{eqnarray*}

Esta equa\c{c}\~{a}o pode nos dar uma boa aproxima\c{c}\~{a}o da
solu\c{c}\~{a}o em matrizes densas, i.\'{e}, com poucos elementos nulos,
como se consider\'{a}ssemos o ``expoente dos elementos nulos'' como
sendo expoente m\'{e}dio de A.  Este \'{e} o {\bf m\'{e}todo aproximado da
redu\c{c}\~{a}o de vari\^{a}ncia}.  Procuremos agora m\'{e}todos
heur\'{\i}sticos para determina\c{c}\~{a}o de escalamentos, que sejam
menos trabalhosos que o m\'{e}todo da redu\c{c}\~{a}o de vari\^{a}ncia. 

O {\bf m\'{e}todo da m\'{e}dia geom\'{e}trica} consiste em tomar
\begin{eqnarray*}
d_j &=& -int( (\sum_i{}' expo(A_i^j) )/ enn(A^j) ) \\
e_i &=& -int( (\sum_j{}'( expo(A_i^j) +d_j ))/ enn(A^i) ) 
\end{eqnarray*}

Assim os fatores de quilibramento s\~{a}o uma aproxima\c{c}\~{a}o
do inverso da m\'{e}dia geom\'{e}trica dos elementos n\~{a}o nulos das
colunas, ou linhas, i.e.,
\begin{eqnarray*}
d_j &=& -int( (\sum_i{}' \log(|A_i^j|) )/ enn(A^j) ) \\
e_i &=& -int( (\sum_j{}'( \log(|A_i^j|) +d_j ))/ enn(A^i) ) 
\end{eqnarray*}

Uma variante do m\'{e}todo da m\'{e}dia geom\'{e}trica \'{e} o
{\bf m\'{e}todo da m\'{e}dia max-min}, no qual tomamos
\begin{eqnarray*}
d_j &=& -int(( \max_i{}'expo(A_i^j) +\min_i{}'expo(A_i^j))/2) \\
e_i &=& -int(( \max_j{}'(expo(A_i^j)+d_j) 
               +\min_j{}'(expo(A_i^j)+d_j) )/2)
\end{eqnarray*}

Finalmente, o {\bf m\'{e}todo da norma infinito} consiste em tomar
\begin{eqnarray*}
e_i &=& -\max_j{}' expo(A_i^j) \\
d_j &=& -\max_i{}' (expo(A_i^j) +e_i)
\end{eqnarray*}

A escolha do particular m\'{e}todo a ser empregado depende bastante do
tipo de matriz a ser equilibrada e da exig\^{e}ncia que temos sobre
$vex(A)$.  Em pacotes comerciais de otimiza\c{c}\~{a}o \'{e} comum a
aplica\c{c}\~{a}o do m\'{e}todo max-min um n\'{u}mero
pr\'{e}-determinado de vezes, ou at\'{e} que vari\^{a}ncia dos expoentes
se reduza a um valor limite aceit\'{a}vel.  Este limite deve ser tomado
em fun\c{c}\~{a}o das condi\c{c}\~{o}es do problema, como por exemplo o
n\'{u}mero de condi\c{c}\~{a}o da matriz, ser definido no cap\'{\i}tulo
9, e da unidade de erro do SPF. 

Exemplo 4: 

Equilibremos a matriz $A$, dada abaixo, num SPF de base 10,
\begin{enumerate}
\item pelo m\'{e}todo aproximado de redu\c{c}\~{a}o de vari\^{a}ncia,
\item pelo m\'{e}todo da m\'{e}dia geom\'{e}trica,
\item pelo m\'{e}todo max-min,
\item pelo m\'{e}todo da norma $\infty$. 
\end{enumerate}
Para $A$ e para cada um dos escalamentos, calculemos $mex$, $vex$,
e o di\^{a}metro da matriz, definido como a diferen\c{c}a entre o maior 
e o menor expoente dos elementos de $A$.
 
$$
A= \left[ \begin{array}{ccccc} 
1 & 0.7E-4 & 0.5E-1 & 0.9E0 & -0.2E2 \\ 
0 & 0 & 0.3E-1 & 0 & 0.3E4 \\
1 & -0.8E-3 & 0 & 0 & 0.3E6 \\
0 & 0.3E-1 & -0.8E0 & 0.1E3 & 0.7E9 
\end{array} \right] $$
$$ 
expo(A)= \left[ \begin{array}{ccccc} 
0 & -4 & -1 & 0 & 2 \\ x & x & -1 & x & 4 \\ 
0 & -3 & x & x & 6 \\ x & -1 & 0 & 3 & 9 \end{array} \right]  
\ \ \ \begin{array}{c} +1 \\ 160 \\ 13 \end{array} $$

Pelo m\'{e}todo aproximado de redu\c{c}\~{a}o de vari\^{a}ncia, temos 
$e$, $d$, $expo(EAD)$ e $[mex, vex, diam ]$, respectivamente
$$
\left[ \begin{array}{c} 
int(8/5)=2 \\ int(2/5)=0 \\ int(2/5)=0 \\ int(-6/5)=-1 \\ 
\end{array} \right] \ \ 
\left[ \begin{array}{c} 
int(4/4)=1 \\ int(12/4)=3 \\ int(6/4)=2 \\ int(1/4)=0 \\ int(-17/4)=-4 
\end{array} \right] \ \ 
\left[ \begin{array}{ccccc}
3 & 1 & 3 & 2 & 0 \\ x & x & 1 & x & 0 \\ 0 & 0 & x & x & 2 \\ 
x & 3 & 1 & 2 & 4 \end{array} \right] \ \ 
\begin{array}{c}  1.57 \\ 23.4 \\ 4 \end{array} $$ 

Analogamente, pelo m\'{e}todo da m\'{e}dia geom\'{e}trica, temos
$$
\left[ \begin{array}{c} 
-int(4/5)=1 \\ -int(-1/2)=1 \\ -int(1/3)=0 \\ -int(8/4)=-2 
\end{array} \right] \ \ 
\left[ \begin{array}{c} 
-int(0/2)=0\\ -int(-8/3)=3\\ -int(-2/4)=1\\ -int(3/2)=-2\\ -int(21/4)=-5 
\end{array} \right] \ \ 
\left[ \begin{array}{ccccc}
1 & 0 & 1 & 1 & -2 \\ x & x & 1 & x & 0 \\ 0 & 0 & x & x & 1 \\ 
x & 0 & -1 & -1 & 2 \end{array} \right] \ \ 
\begin{array}{c}  0.214 \\ 14.4 \\ 4 \end{array} $$ 

Analogamente, pelo m\'{e}todo max-min, temos
$$
\left[ \begin{array}{c} 
-int(-4/2)=2 \\ -int(-2/2)=1 \\ -int(0/2)=0\\ -int(4/2)=-2  
\end{array} \right] \ \ 
\left[ \begin{array}{c} 
-int(0/2)=0\\ -int(-5/2)=3\\ -int(-1/2)=1\\ -int(3/2)=-2\\ -int(11/2)=-6 
\end{array} \right] \ \ 
\left[ \begin{array}{ccccc}
2 & 1 & 2 &  & -2 \\ x & x & 1 & x & -1 \\ 0 & 0 & x & x & 0 \\ 
x & 0 & -1 &  & 1 \end{array} \right] \ \ 
\begin{array}{c}  0.143 \\ 14.3 \\ 4 \end{array} $$ 

Analogamente, pelo m\'{e}todo da norma-$\infty$, temos
$$
\left[ \begin{array}{c} 
 0 \\ 1 \\ 0 \\ 0  \end{array} \right] \ \ 
\left[ \begin{array}{c} 
 0 \\ 1 \\ 0 \\ -3 \\ -9 \end{array} \right] \ \ 
\left[ \begin{array}{ccccc}
0 & -3 & -1 & -3 & -7 \\ x & x & 0 & x & -4 \\ 0 & -2 & x & x & -3 \\ 
x & 0 & 0 & 0 & 0 \end{array} \right] \ \ 
\begin{array}{c}  -1.64 \\ 59.2 \\ 7 \end{array} $$

 Escalamentos visando equilibrar as matrizes do problema s\~{a}o uma
etapa importante na solu\c{c}\~{a}o de sistemas lineares de grande
porte.  Note que opera\c{c}\~{o}es de escalamento n\~{a}o afetam os
elementos nulos de uma matriz, e portanto n\~{a}o alteram sua estrutura
e esparsidade.  O desempenho das heur\'{\i}sticas estudadas variam
conforme a \'{a}rea de aplica\c{c}\~{a}o; vale pois testar
experimentalmente as heur\'{\i}stcas escalamento e suas
varia\c{c}\~{o}es.

 \clearpage
 \clearpage 
\setcounter{chapter}{8} 
\chapter{ESTABILIDADE}
\begin{center}
{\LARGE Erros e Perturba\c{c}\~{o}es da Solu\c{c}\~{a}o }
\end{center} 

\section{Normas e Condi\c{c}\~{a}o} 

Uma {\bf norma}, num dado espa\c{c}o vetorial $E$, \'{e} uma 
fun\c{c}\~{a}o $$ ||\ .\ ||\ :\ E \Rightarrow {\bf R}  \mid 
   \forall x, y \in E, \alpha \in {\bf R} \ :$$
\begin{enumerate}
\item $ || x || \geq 0,\ \mbox{e}\ || x || = 0 \Leftrightarrow x = 0 $. 
\item $ || \alpha x || = | \alpha | \ || x || $.
\item $ || x + y || \leq || x || + || y || $, a desigualdade triangular. 
\end{enumerate}
 
S\~{a}o de grande interesse em ${\bf R}^n$ as {\bf {\it p}-normas}
 $$|| x ||_p = ( \sum_1^n |x_i|^p )^{1/p} \ ,$$ 
e de particular interesse a norma 1, a norma 2
(ou norma quadr\'{a}tica, ou Euclidiana) e a norma $p=+\infty$.
No caso da norma infinito, devemos tomar o limite da defini\c{c}\~{a}o 
de $p$-norma para $p\rightarrow +\infty$, ou seja, 
 $$|| x ||_p = \max_{i=1}^n |x_i| \ .$$ 
 \index{Norma!p-norma} 
 \index{Norma!induzida} 
 \index{Norma!de matriz}  

Dado um espa\c{c}o vetorial normado $(E,\ ||\ ||)$ definimos a 
{\bf norma induzida} sobre as transforma\c{c}\~{o}es lineares limitadas,
$ T:\ E \rightarrow E \ \mbox{tq}\ \exists \alpha \in {\bf R} \mid $ 
$ \forall x \in E,\  ||T(x)|| \leq \alpha ||x|| $
como sendo
$$ ||T|| \equiv \max_{x\neq 0} ( \ ||T(x)|| \ /\ ||x|| \ ) $$
ou equivalentemente, por linearidade
$$ ||T|| \equiv \max_{x\mid \ ||x||=1} ||T(x)|| \ .$$

Em $( {\bf R}^n,\ ||\ || )$ falamos da norma induzida sobre as matrizes
$n\times n$ como sendo a norma da transforma\c{c}\~{a}o associada, isto
\'{e} $||A|| = ||T||$, onde $T(x) = Ax$,

\begin{lm} 
A norma induzida sobre as matrizes em $( {\bf R}^n,\ ||\ || )$, 
goza das propriedades, para 
$\forall A, B\ n\times n,\ \alpha \in {\bf R},$
\begin{enumerate}
\item $||A|| \geq 0 \ \mbox{e}\ ||A|| = 0 \Leftrightarrow A=0 $
\item $ || A+B || \leq ||A|| + ||B|| $
\item $ || AB || \leq ||A|| \ ||B|| $ 
\end{enumerate} 
\end{lm}

\begin{lm} 
Para a norma 1 e para norma $\infty$ temos as seguintes express\~{o}es
expl\'{\i}citas da norma induzida sobre as transforma\c{c}\~{o}es, ou
matrizes,
\begin{eqnarray*}
||A||_1 &=& \max_{j=1}^n \sum_{i=1}^n  | A_i^j |  \\ 
||A||_\infty &=& \max_{i=1}^n \sum_{j=1}^n | A_i^j | 
\end{eqnarray*} 
\end{lm} 

Demonstra\c{c}\~{a}o:

Para verificar a validade das express\~{o}es observe que
\begin{eqnarray*} 
||Ax||_1 &=& \sum_{i=1}^n \mid \sum_{j=1}^n A_i^j x_j \mid 
       \leq \sum_{i=1}^n \sum_{j=1}^n |A_i^j | \ |x_j | \\ 
 &\leq & \sum_{j=1}^n |x_j | \max_{j=1}^n \sum_{i=1}^n |A_i^j | 
       = ||A||_1 \ ||x||_1 
\end{eqnarray*} 
\begin{eqnarray*} 
||Ax||_{\infty } &=& \max_{i=1}^n \mid \sum_{j=1}^n A_i^j x_j \mid  
       \leq \max_{i=1}^n \sum_{j=1}^n |A_i^j | \ |x_j | \\ 
 &\leq & \max_{j=1}^n |x_j | \max_{i=1}^n \sum_{j=1}^n |A_i^j | 
       = ||x||_\infty \ ||A||_\infty  
\end{eqnarray*} 

e que, se $k$ \'{e} o \'{\i}ndice que realiza o m\'{a}ximo na
defini\c{c}\~{a}o da norma, ent\~{a}o as igualdades s\~{a}o realizadas
pelos vetores $x = I^k$, para a norma 1, e
$x \mid  x_j = sig(A_i^j )$, para a norma $\infty$.

Definimos o {\bf n\'{u}mero de condi\c{c}\~{a}o} de uma matriz como
$$cond(A) = ||A|| \ ||A^{-1}|| \ .$$ 
 \index{Condi\c{c}\~{a}o}
 \index{Norma!condi\c{c}\~{a}o} 
 \index{Condi\c{c}\~{a}o!n\'{u}mero de} 

Exemplo 1: 

Calcule a norma da matriz de binomial de dimens\~{a}o 3, e de sua
inversa, nas normas 1 e $\infty$.  Para cada uma das normas calculadas
exiba um vetor $x$ que, multiplicado pela matriz, seja ``esticado'' de
um fator igual \`{a} pr\'{o}pria norma, i.e. 
 $x\mid \ \| Ax\| = \| A\| \| x\|$. 
 $$ ^3B = 
\left[ \begin{array}{ccc} 
\left( \begin{array}{c} 1 \\ 0 \end{array} \right) & 
\left( \begin{array}{c} 1 \\ 1 \end{array} \right) & 0 \\ 
\left( \begin{array}{c} 2 \\ 0 \end{array} \right) & 
\left( \begin{array}{c} 2 \\ 1 \end{array} \right) & 
\left( \begin{array}{c} 2 \\ 2 \end{array} \right)  \\  
\left( \begin{array}{c} 3 \\ 0 \end{array} \right) & 
\left( \begin{array}{c} 3 \\ 1 \end{array} \right) & 
\left( \begin{array}{c} 3 \\ 2 \end{array} \right) 
\end{array} \right] $$ 
$$
^3B = \left[ \begin{array}{ccc} 1 & 1 & 0 \\ 1 & 2 & 1 \\ 
 1 & 3 & 3 \end{array} \right] \ \ 
^3B^{-1} = \left[ \begin{array}{rrr} 3 & -3 & 1 \\ -2 & 3 & -1 \\ 
 1 & -2 & 1 \end{array} \right] $$ 
Para $^3B$ a soma da norma dos elementos $|A_i^j|$ por linha e por 
coluna s\~{a}o, respectivamente, 
$\left[ \begin{array}{ccc} 7 & 6 & 4 \end{array} \right]$ e 
$\left[ \begin{array}{ccc} 6 & 8 & 3 \end{array} \right]$.
Para $^3B^{-1}$, analogamente, temos 
$\left[ \begin{array}{ccc} 7 & 6 & 4 \end{array} \right]$ e 
$\left[ \begin{array}{ccc} 7 & 6 & 4 \end{array} \right]$. 
Assim, para a norma 1, $|| ^3B ||$, $|| ^3B^{-1} ||$, e $cond( ^3B)$ 
s\~{a}o respectivamente 6, 8, e 48. 
Analogamente, para norma $\infty$, temos 7, 7, e 49. 
Finalmente,  
$\left[ \begin{array}{ccc} 0 & 1 & 0 \end{array} \right] '$ e 
$\left[ \begin{array}{ccc} 1 & -1 & 1 \end{array} \right] '$ 
s\~{a}o vetores como os procurados.

\section{Perturba\c{c}\~{o}es }

Estudaremos agora uma maneira de limitar o efeito de uma
pertuba\c{c}\~{a}o nos dados do sistema linear, $Ax=b$, sobre a sua
solu\c{c}\~{a}o.

\begin{teo}[da pertuba\c{c}\~{a}o] 
Se os dados do sistema $Ax=b$ forem pertubados, isto \'{e}, alterados de
uma matriz $\delta A$ de um vetor $\delta b$, a resposta ser\'{a} 
pertubada por um $\delta x$, isto \'{e}
$(A + \delta A)(x + \delta x) = (b + \delta b)$, tal que, se  
$||\delta A|| \ ||A^{-1}|| < 1$, ent\~{a}o
$$ \frac{||\delta x||}{||x||} \leq 
   \frac{cond(A)}{1- cond(A) \ ||\delta A|| / ||A||} \left(
  \frac{||\delta b||}{||b||} +\frac{||\delta A||}{||A||} \right) \ .$$
\end{teo} 
 \index{Teorema!da Perturba\c{c}\~{a}o} 

Demonstra\c{c}\~{a}o.

$(A +\delta A)(x +\delta x) = 
Ax +A\delta x +\delta Ax +\delta A\delta x = b +\delta b$,
portanto $A\delta x = \delta b -\delta Ax -\delta A\delta x$, e 
$\delta x = A^{-1}( \delta b  -\delta Ax -\delta A\delta x )$.
Assim, 
$$||\delta x|| \leq ||A^{-1}|| \left( ||\delta b|| + 
 ||\delta A|| ||x|| + ||\delta A|| ||\delta x|| \right) \ .$$ 
 
Como tamb\'{e}m $||b|| \leq ||A|| ||x||$,
temos que 
\begin{eqnarray*}  
 \left( 1- ||A^{-1}|| ||\delta A|| \right) \frac{||\delta x||}{||x||} 
 &\leq & ||A^{-1}|| \left( \frac{||\delta b||}{||x||} 
         + ||\delta A|| \right) \\  & &  
 \leq ||A^{-1}|| \left( \frac{||A|| ||\delta b||}{||b||} 
      + ||\delta A|| \right) \ .\end{eqnarray*} 
Usando a hip\'{o}tese  $||\delta A|| ||A^{-1}|| \leq 1$, 
temos o teorema, QED.

Exemplo 2:

 Considere o sistema $(B + \delta A )x = (b + \delta b)$, onde a matriz
de coeficientes, e o vetor de termos independentes correspondem \`{a}
matriz binomial de dimens\~{a}o 3,
 $^3B$, e $^3b\mid \ ^3B{\bf 1}=\ ^3b$.  
 Os elementos da matriz e do vetor de pertuba\c{c}\~{a}o,
 $\delta A_i^j$ e $\delta b_i$, s\~{a}o aleat\'{o}rios.  
 A distribui\c{c}\~{a}o de cada um destes elementos de
perturba\c{c}\~{a}o \'{e} independente e \'{e} uma fun\c{c}\~{a}o de
dois parametros $(\alpha ,p)$: Tomemos cada elemento da
perturba\c{c}\~{a}o \'{e} no cojunto $\{0, -\alpha ,\alpha \}$ com
probabilidade, respectivamente, $[1-p, p/2, p/2]$.

Com os dados do Exemplo 1, determine um limite m\'{a}ximo para $0 \leq
\alpha \leq alfamax$ que garanta a hip\'{o}tese do teorema da
pertuba\c{c}\~{a}o.  Fa\c{c}a uma experi\^{e}ncia comparando limites de
erro e pertuba\c{c}\~{a}o da solu\c{c}\~{a}o do sistema proposto. 

Na norma 1, $||\delta A|| \leq 3\alpha$, e do Exemplo 1 sabemos que
$||B^{-1}|| \leq 7$.  Assim, a condi\c{c}\~{a}o 
$||\delta A|| ||A^{-1}|| \leq 1$ est\'{a} garantida para 
$\alpha \leq 0.04 < 1/21 $. 

Tomando como experimento de perturba\c{c}\~{a}o,  
$$ \delta A = \left[ \begin{array}{rrr} 0 & 0.01 & -0.01 \\ 
   0 & 0.01 & 0 \\ 0 & 0.01 & -0.01 \end{array} \right] \ \ e \ \ 
   \delta b = \left[ \begin{array}{r} 
   0 \\ -0.01 \\ 0 \end{array} \right] $$ 
a solu\c{c}\~{a}o do sistema pertubado \'{e} 
$x = [1.0408, 0.9388, 1.0622 ]'$, i.e. 
$\delta x = [ 0.04, -0.06, 0.06 ]'$ e $||\delta x|| = 0.16$.

Por outro lado, o Teorema da pertuba\c{c}\~{a}o nos d\'{a} o limite
$$ \frac{||\delta x||}{||x||} \leq 
   \frac{49}{1- 49*0.03/7}(\frac{0.01}{7}+\frac{0.03}{7}) = 0.35 $$
um limite superior de acordo com o resultado obtido no experimento.

\section{Erro na Fatora\c{c}\~{a}o LU}

Consideremos a resolu\c{c}\~{a}o de um sistema, $Ax=b$,
pelo m\'{e}todo de Doolittle, isto \'{e},
a decomposi\c{c}\~{a}o $A=LU$,
e a solu\c{c}\~{a}o do sistema   $Ly=b$ e $Ux=b$.

\begin{teo}[Wilkinson] %
 A solu\c{c}\~{a}o, afetada pelos erros de
arredondamento, pode ser escrita como a solu\c{c}\~{a}o exata, $x$, mais
um termo de pertuba\c{c}\~{a}o $\delta x$. 
 Assumiremos que um produto escalar em dupla precis\~{a}o \'{e} afetado
de um erro da ordem do erro de passagem para precis\~{a}o simples, 
conforme a observa\c{c}\~{a}o 8.1. 
 Nestas condi\c{c}\~{o}es,, a solu\c{c}\~{a}o calculada,
 $(x +\delta x)$, \'{e} solu\c{c}\~{a}o exata de um sistema.
 $(A +\delta A)(x +\delta x) = b$, tal que 
 $||\delta A||_{\infty} \leq (2n+1)gu$,   
onde $g = \max_{i,j} |U_i^j|$.
\end{teo} 
 \index{Teorema!Wilkinson}

Demonstra\c{c}\~{a}o.

Consideremos as linhas da matriz $A$ j\'{a} ordenadas de modo a n\~{a}o
haver necessidade de pivoteamentos.  A decomposi\c{c}\~{a}o da matriz $A$
ser\'{a} dada por, para $i=1\ldots n$ 
\begin{eqnarray*}
L_i U &=& A_i \\ 
LU^i &=& A^i 
\end{eqnarray*} 
ou, para $i=1\ldots n$, para  
\begin{eqnarray*}
j=i\ldots n 
 & & 
 U_i^j = A_i^j -\sum_{k=1}^{i-1} M_i^k U_k^j \\
j=i+1\ldots n 
 & & 
 M_j^i = (A_j^i -\sum_{k=1}^{i-1} M_j^k U_k^i) /U_i^i  
\end{eqnarray*} 

Na realidade, obteremos as matrizes afetadas de erro, 
para $i=1\ldots n$, para  
\begin{eqnarray*}
j=i\ldots n 
 & & 
 \tilde U_i^j = fl(fld( A_i^j -
  \sum_{k=1}^{i-1} \tilde M_i^k \tilde U_k^j )) \\
j=i+1\ldots n 
 & & 
 \tilde M_j^i = fl(fld( (A_j^i 
   -\sum_{k=1}^{i-1} \tilde M_j^k \tilde U_k^i) / \tilde U_i^i ))) 
\end{eqnarray*} 

Como supomos desprez\'{\i}veis os termos de $O(ud)$ (vide
observa\c{c}\~{a}o 8.1) temos, para $i=1\ldots n$, para 
\begin{eqnarray*}
j=i\ldots n 
 & & 
 \tilde U_i^j = (1+\delta _i^j)( A_i^j 
   -\sum_{k=1}^{i-1} \tilde M_i^k \tilde U_k^j ) \\
j=i+1\ldots n 
 & & 
 \tilde M_j^i = (1+\delta _j^i )(A_j^i 
   -\sum_{k=1}^{i-1} \tilde M_j^k \tilde U_k^i) /\tilde U_i^i  
\end{eqnarray*} 

Portanto, para $i=1\ldots n$, para   
\begin{eqnarray*}
j=i\ldots n 
 & & 
\sum_{k=1}^{i-1} {\tilde M}_i^k {\tilde U}_k^j 
    + 1U_i^i = {\tilde L}_i {\tilde U}^j 
    = A_i^j + {\tilde U}_i^j \delta _i^j / (1+\delta _i^i) \\
j=i+1\ldots n 
 & & \sum_{k=1}^{i-1} 
  {\tilde M}_j^k {\tilde U}_k^i  + {\tilde M}_j^i {\tilde U}_i^i 
  = {\tilde L}_j {\tilde U}^i 
  = A_j^i + {\tilde M}_j^i {\tilde U}_i^i \delta _j^i / (1+\delta _j^i )  
\end{eqnarray*} 

Notemos agora que: $|\delta _i^j|/(1+\delta _i^j) \sim 
|\delta _i^j| \leq u$, que $|M_i^j|\leq 1$, pelo pivoteamento parcial, 
e definindo $g=\max_{i,j} |U_i^j|$, temos que $\tilde L \tilde U = A+E$, 
onde $|E_i^j| \leq gu$, donde $||E||_\infty \leq ngu$.

Na solu\c{c}\~{a}o do sistema $\tilde Ly=b$ e $\tilde Ux=y$ 
calculamos, para $i=1\ldots n$,
$$
{\tilde y}_i = fl(fld((b_i -\sum_{j=1}^{i-1} 
 {\tilde L}_i^j {\tilde y}_j )/{\tilde L}_i^i )) $$ 
e novamente, para $i=1\ldots n$,
$$
{\tilde x}_i = fl(fld((\tilde y_i -\sum_{j=i+1}^{n} 
 {\tilde U}_i^j {\tilde x}_j )/{\tilde U}_i^i )) \ .$$ 

Supondo desprez\'{\i}veis os termos de $0(ud)$ 
\begin{eqnarray*} 
 {\tilde y}_i &=& (1+\delta _i) (b_i -\sum_{j=1}^{i-1} 
  {\tilde L}_i^j {\tilde y}_j )/{\tilde L}_i^i  \\
 {\tilde x}_i &=& (1+\delta _i ') (y_i -\sum_{j=i+1}^{n} 
 {\tilde U}_i^j {\tilde x}_j )/{\tilde U}_i^i  
\end{eqnarray*} 
ou
\begin{eqnarray*}
\sum_{j=1}^{i} {\tilde L}_i^j {\tilde y}_j 
  &=& {\tilde L}_i {\tilde y}
  = b_i + \tilde L_i^j {\tilde y}_i \delta _i / (1+\delta _i) \\
\sum_{j=i}^{n} {\tilde U}_i^j {\tilde x}_j 
  &=& {\tilde U}_i {\tilde x}
  = {\tilde y}_i + U_i^j {\tilde x}_i {\delta _i}' / (1+\delta _i ') \\
\end{eqnarray*}

isto \'{e}
\begin{eqnarray*}
({\tilde L} + \delta L ){\tilde y} &=& b \\
({\tilde U} + \delta U ){\tilde x} &=& \tilde y 
\end{eqnarray*} 

onde para as matrizes diagonais $\delta L$ e $\delta U$, temos
$|\delta L_i^i| \leq u$ e $|\delta U_i^i| \leq gu$.

Em suma,
\begin{eqnarray*}
b &=& (\tilde L + \delta L) \tilde y 
      = ( \tilde L + \delta L)( \tilde U + \delta U) \tilde x \\
 &=& ( \tilde L \tilde U + \tilde L \delta U + \delta L \tilde U 
       + \delta L \delta U) \tilde x \\
 &=& ( A + E + \tilde L \delta U + \delta L \tilde U 
       + \delta L \delta U ) \tilde x \\
 &=& ( A + \delta A) \tilde x
\end{eqnarray*}

Desprezando o termo $\delta L \delta U$, de $O(ud)$, 
$\delta A = E + \tilde L \delta U + \delta L \tilde U$, donde
$$ ||\delta A||_\infty = ||E|| + 
   ||\tilde L \delta U + \delta LU \tilde U|| \leq ngu + (n+1)gu $$
pois
\begin{eqnarray*}
\lefteqn{ ||\tilde L \delta U + \delta L \tilde U||_\infty \leq } \\
 & & 
\leq \left| \left| 
\left[ \begin{array}{ccc}
1 & & 0 \\ & \ddots & \\ 1 & & 1 \end{array} \right] 
\left[ \begin{array}{ccc}
gu & & 0 \\ & \ddots & \\ 0 & & gu \end{array} \right] +  
\left[ \begin{array}{ccc}
u & & 0 \\ & \ddots & \\ 0 & & u \end{array} \right] 
\left[ \begin{array}{ccc}
g & & g \\ & \ddots & \\ 0 & & g \end{array} \right] 
\right| \right| _\infty \\ 
 & & 
\leq \left| \left| 
\left[ \begin{array}{ccc}
gu & & 0 \\ & \ddots & \\ gu & & gu \end{array} \right] + 
\left[ \begin{array}{ccc}
gu & & gu \\ & \ddots & \\ 0 & & gu \end{array} \right] 
\right| \right|  \\  
 & & 
\leq \left| \left| 
\left[ \begin{array}{ccc}
2gu & & gu \\ & \ddots & \\ gu & & 2gu \end{array} \right] 
\right| \right|  = (n+1)gu  \ \ \ QED. 
\end{eqnarray*}

Para frisar a import\^{a}ncia do uso conjugado do m\'{e}todo de
Doolittle e dupla precis\~{a}o daremos, sem demonstra\c{c}\~{a}o, a
vers\~{a}o do teorema ora demonstrado para o m\'{e}todo de
Gauss, que equivale ao m\'{e}todo de Doolittle sem o recurso da dupla
precis\~{a}o:  
 Definindo $h =\max_{i,j,k} |\ ^kA_i^j|$, ter\'{\i}amos
 $||E||_\infty \leq O(n^2hu)$ e 
 $||\delta A||_\infty \leq O(n^3hu)$.
 Estes resultados tornam evidente a obrigatoriedade de calcularmos
os produtos escalares envolvidos no processo, sempre em dupla 
precis\~{a}o, a menos que tratemos de sistemas de pequena dimens\~{a}o. 
 \index{Fatora\c{c}\~{a}o!Doolittle} 
 \index{Precis\~{a}o!dupla}

Al\'{e}m do pivoteamento parcial, isto \'{e}, da escolha como elemento
piv\^{o} do maior elemento em m\'{o}dulo na coluna, poder\'{\i}amos usar
o {\bf pivoteamento total}, isto \'{e}, encolher como elemento piv\^{o}
o maior elemento em m\'{o}dulo em qualquer linha ou coluna
utiliz\'{a}vel (isto \'{e}, poderiamos tomar por piv\^{o} da
transforma\c{c}\~{a}o $^{k-1}A \rightarrow \ ^kA$, um elemento em
 $arg\max_{k\leq i, j\leq n} |\ ^{k-1}A_i^j$.                 
 Usando pivoteamento total, podemos demonstrar a exist\^{e}ncia de um
limite superior para a constante $h$, supondo que a matriz original $A$
\'{e} tal que $|A_i^j|\leq 1$, da ordem de $O(n^{(1/4)\ln (n)})$, 
contra limites de $O(2^n)$ para pivoteamento parcial. 
 O uso do pivoteamento total \'{e} todavia, extremamente custoso, pois
exige a cada passo da $O(n^2)$ compara\c{c}\~{o}es, contra $O(n)$ no
pivoteamento parcial.  Ademais, estes limites de $h$ tendem a ser
extremamente pessimistas, principalmente se $A$ for bem equilibrada. 
 \index{Pivoteamento!total}

\section*{Exerc\'{\i}cios}
 
\begin{enumerate}
\item Estabilidade de Fatora\c{c}\~{o}es Sim\'{e}tricas. 
 \index{Fatora\c{c}\~{a}o!sim\'{e}trica} 
 \index{Fatora\c{c}\~{a}o!de Cholesky}  
 \begin{enumerate}  
 \item Adapte os resultados de estabilidade da faora\c{c}\~{a}o LU para
o caso particular da fatora\c{c}\~{a}o de Cholesky. 
  \item Qual o n\'{u}mero de condi\c{c}\~{a}o de uma matriz ortogonal?
Usando a rela\c{c}\~{a}o entre o fator triangular da fatora\c{c}\~{a}o
QR e a fatora\c{c}\~{a}o de Cholesky, discuta a estabilidade da
fatora\c{c}\~{a}o QR.  
 \end{enumerate} 
\item Prove que $||\ ||_1$, $||\ ||_2$ e $||\ ||_\infty$ s\~{a}o
 efetivamente normas.
\item Desenhe em ${\bf R}^2$ a regiao $||x||\leq 1$, 
 para as normas 1, 2 e $\infty$.
\item Prove o Lema 1.
\item Equival\^{e}ncia das normas 1, 2 e $\infty$: 
 Prove em ${\bf R}^n$ que 
\begin{enumerate} 
\item $||x||_\infty \leq ||x||_1 \leq n||x||_\infty $
\item $||x||_\infty \leq ||x||_2 \leq n^{1/2}||x||_\infty $
\item $||A||_2 = \lambda$, 
      onde $\lambda ^2$ \'{e} o maior autovalor de $A'A$.
\item $||A'A||_2 = ||A||_2^2$
\item $||A||_2 \leq (\ ||A||_1 ||A||_\infty )^{1/2}$
\end{enumerate}

\item Calcule computacionalmente $cond(A)$, para as matrizes de teste
$^nT$, $A=\ ^nB$ e $A=\ ^nH$, para $n= 2,4,8,16$.  Estime graficamente o
crescimento do n\'{u}mero de condi\c{c}\~{a}o das matrizes teste em
fun\c{c}\~{a}o da dimens\~{a}o. 

\item
Resolva computacionalmente os sistemas teste 
$^8Tx=\ ^8t$, $^8Bx=\ ^8b$ e $^8Hx=\ ^8h$. 
\begin{enumerate}
\item Pelo m\'{e}todo de Gauss com pivoteamento parcial;
\item Pelo m\'{e}todo de Gauss com pivoteamento total;
\item Pelo m\'{e}todo de Doolittle com pivoteamento parcial 
      e dupla precis\~{a}o; 
\item Pelo m\'{e}todo de Doolittle com pivoteamento total e 
      dupla precis\~{a}o. 
\end{enumerate}
No item c, verifique se o erro final est\'{a} dentro do esperado, em
fun\c{c}\~{a}o da unidade de erro do SPF utilizado.

\end{enumerate}

 \clearpage
 \clearpage 

\setcounter{chapter}{9} 
\chapter{MUDAN\c{C}A de BASE}  
\begin{center} 
{\LARGE Atualiza\c{c}\~{o}es de Posto 1}
\end{center} 

Em Otimiza\c{c}\~{a}o, principalmente em Programa\c{c}\~{a}o Linear, o
termo base significa uma matriz quadrada de posto pleno.  Estudaremos
agora o problema de {\bf mudan\c{c}a de base}, aqui posto na seguinte
forma: Seja $\hat A$ a nova base, obtida da base original, $A$,
substituindo-se a coluna $A^s$ por uma nova coluna, $a$, isto \'{e}
 $\hat A = [A^1,\ldots A^{s-1},a, A^{s+1},\ldots A^n]$. 
 Se j\'{a} dispusermos da inversa de $A$ (ou na pr\'{a}tica mais
comumente de uma fatora\c{c}\~{a}o de $A$), como atualizar a inversa (ou
a fatora\c{c}\~{a}o)? Isto \'{e}, como, a partir da inversa de $A$ obter
a inversa de $\hat A$ com um trabalho muito menor que o necess\'{a}rio
para {\bf reinverter} $\hat A$, i.e., computar a inversa de novo, sem
aproveitar a informa\c{c}\~{a}o contida em $A^{-1}$ ou, na
fatora\c{c}\~{a}o $A=LU$ ou $A=QR$ ?
 \index{Atualiza\c{c}\~{a}o!base}
 \index{Atualiza\c{c}\~{a}o!posto 1}

\section{F\'{o}rmulas de Modifica\c{c}\~{a}o}

\begin{teo}[f\'{o}rmula geral de modifica\c{c}\~{a}o] 
 Dada $A$, $n\times n$ e invers\'{\i}vel, $V=A^{-1}$, $\delta A$,
$n\times n$, e $\alpha \in R$ suficientemente pequeno, podemos fazer a
expans\~{a}o
$$ (A+ \alpha \delta A)^{-1} = V+ \sum_{k=1}^\infty 
   (-\alpha)^k \; (V \; \delta A)^k \; V \ .$$
\end{teo} 

Demonstra\c{c}\~{a}o: 

Em virtude da regra de Cramer podemos, para 
$\alpha\in [0, \mbox{alphamax}]$, escrever a s\'{e}rie de Taylor para 
cada elemento de $(A+ \alpha \delta A)^{-1}$,
 $$(A+ \delta A )^{-1} = V+ \sum_{k=1}^\infty 
  \frac{{\alpha}^k}{k!} \ {\delta ^k V}$$
 Para determinar a matriz que d\'{a} a perturba\c{c}\~{a}o de ordem $k$
da inversa $\delta ^k V$ observemos que
 $$ (A+ \alpha \delta A)(V+ \alpha \ {\delta V}+ 
   \frac{{\alpha}^2}{2} \ {\delta ^2 V}+  
   \frac{{\alpha}^3}{3!} \ {\delta ^3 V}+ 
   \ldots ) = I $$
 donde
 $$ \alpha ( {\delta A}\ V+ A\ {\delta V} )+
  \alpha ^2 ( \frac{1}{2}A\ {\delta ^2 V}+ {\delta A}\ {\delta V} )+
  \alpha ^3 ( \frac{1}{6}A\ {\delta ^3 V}+ 
              \frac{1}{2} {\delta A}\ {\delta ^2 V} )+ 
  \ldots =0 $$
 de modo que
 $$ {\delta V} = -V\ {\delta A}\ V \ , \ \ 
   {\delta ^2 V} = 2V\ {\delta A}\ V\ {\delta A}\ V \ , \ \ 
 {\delta ^3 V} = ( V\ {\delta A} )^3\ V \ , \ldots$$
$$ {\delta ^k V} = (-1)^k\ k!\ ( V\ {\delta A} )^k \ V \ .$$  

Substituindo esta f\'{o}rmula geral de ${\delta ^k V}$  na s\'{e}rie de 
Taylor de $(A+ \alpha \; {\delta A} )$, segue diretamente o teorema.
QED.

\begin{teo}[f\'{o}rmula de Sherman e Morrison]  
Dada $A$, $n\times n$ e invers\'{\i}vel, $u$ e $w$ $n\times 1$, 
e  $\alpha \in R^*$, ent\~{a}o a inversa de 
$\hat A = A + \alpha uw'$ \'{e} dada por
$$ \hat A^{-1} = A^{-1} + \beta A^{-1} u w' A^{-1} \ , $$
$$ \beta = -( {\alpha}^{-1} + w' A^{-1} u )^{-1} \ . $$
\end{teo} 
 \index{Teorema!Sherman e Morrison}

Demonstra\c{c}\~{a}o: 

Da f\'{o}rmula geral de modifica\c{c}\~{a}o, supondo que
$\alpha < \mbox{alfamax}$, e das propriedades do operador
tra\c{c}o (exerc\'{\i}cio 1.5) temos: 
 \begin{eqnarray*}
 \lefteqn{ ( A+ \alpha u w' )^{-1} =}\\  
 & & = V + \sum_{k=1}^\infty (-\alpha )^k (V u w' ) V \\
 & & = V -\alpha V u w' V \sum_k   ( -\alpha \mbox{tr}( V u w' ) )^k \\
 & & = V - \frac{\alpha V u w' V}{1+ \alpha \mbox{tr}(Vuw') }\\
 & & = V - \frac{Vuw'V}{ {\alpha}^{-1}+ \mbox{tr}(Vuw') } \\
 & & = V - ( {\alpha}^{-1}+ w'Vu )^{-1}\; Vuw'V \ .
 \end{eqnarray*}
o que prova o teorema para $0< \alpha < \mbox{alphamax}$.

 A f\'{o}rmula de Sherman e Morrison \'{e} todavia uma identidade que
podemos provar diretamente para qualquer   
$\alpha > 0 \mid ( {\alpha}^{-1}+ w'A^{-1}u )\neq 0$:
Usando a formula de Sherman e Morisson para desenvolver a identidade
$\hat A^{-1} \hat A = I$, obtemos,

\begin{eqnarray*}
(V+ \beta Vuw'V)(A+ \alpha uw')  &=&  I 
 \ \ \Leftrightarrow  \\  
 \beta Vuw' + \alpha \beta (Vuw')^2  &=& -\alpha Vuw' 
 \ \ \Leftrightarrow  \\
 ( {\alpha}^{-1}+ w'Vu ) \alpha Vuw' &=& 
 -{\beta }^{-1} \alpha Vuw' 
\end{eqnarray*}
sendo esta \'{u}ltima identidade trivialmente verdadeira.
Q.E.D.

Exemplo 1: 
  
Tomando
$$A=\left[ \begin{array}{cc} 1& 0\\ 2& 1\end{array} \right]
\ , \ \  
u=\left[ \begin{array}{c} 1\\ 2\end{array} \right] \ , \ \  
w=\left[ \begin{array}{c} 3\\ 4\end{array} \right] 
\ , \ \  \alpha =1 \ ,  
$$
podemos calcular a inversa de
$$
\hat A= \left[ \begin{array}{cc} 1& 0\\ 2& 1\end{array} \right]
  + \left[ \begin{array}{cc} 3& 4\\ 6& 8\end{array} \right]
  = \left[ \begin{array}{cc} 4& 4\\ 8& 9\end{array} \right]
$$ 
calculando 
$$ \beta = -(1+3)^{-1} = -1/4 \ ,$$
$$
\hat A^{-1}= \left[ \begin{array}{cc} 1& 0\\ -2& 1\end{array} \right]
 -(1/4) \left[ \begin{array}{cc} -5& 4\\ 0& 0\end{array} \right]
 = \left[ \begin{array}{cc} 9/4& -1\\ -2& 1\end{array} \right]
$$

\begin{obs}{\rm 
Note que uma mudan\c{c}a de base poderia ser feita pela f\'{o}rmula de
Sherman e Morrison, pois se
$\hat A = [A^1,\ldots A^{j-1}, a, A^{j+1},\ldots A^n]$,
ent\~{a}o, $\hat A = A + (a - A^j)\; I_j$, 
de modo que a inversa da nova base seria, tomando $V=A^{-1}$,
\begin{eqnarray*}
 \hat A^{-1} &=& V + \beta V ( a- A^j ) I_j V \\
 &=& V + V \beta ( a- A^j ) V_j 
\end{eqnarray*}
$$ \mbox{onde} \ \ 
   \beta = -(1- I_j V (a- A^j) )^{-1} = -(V_j a)^{-1} \ .$$
}\end{obs}  

Exemplo 2:  

Se 
$$
A= \left[ \begin{array}{cc} 1& 0\\ 2& 1\end{array} \right]
\ , \ \  
a= \left[ \begin{array}{c} 1\\ 3 \end{array} \right]
\ \mbox{e} \ \ j=2 \ , $$
ent\~{a}o a inversa da nova base, $\hat A^{-1}$, pose ser calculada como segue
$$
\beta = -( \left[ \begin{array}{cc} -2 & 1 \end{array} \right]
 \left[ \begin{array}{cc} 1\\ 3\end{array} \right] )^{-1} 
 = -1 $$
\begin{eqnarray*}
\lefteqn{ \hat A^{-1} = 
 ( \left[ \begin{array}{cc} 1& 1\\ 2& 3\end{array} \right] )^{-1} =}\\
 &=& \left[ \begin{array}{cc} 1& 0\\ -2& 1 \end{array} \right] + 
 \left[ \begin{array}{cc} 1& 0\\ -2& 1 \end{array} \right] 
 (-1) ( \left[ \begin{array}{c} 1\\ 3 \end{array} \right]
    - \left[ \begin{array}{c} 0\\ 1 \end{array} \right] )
 \left[ \begin{array}{cc} -2& 1 \end{array} \right] \\
 &=& \left[ \begin{array}{cc} 1& 0\\ -2& 1 \end{array} \right] + 
 \left[ \begin{array}{c} -1\\ 0 \end{array} \right]
 \left[ \begin{array}{cc} -2& 1 \end{array} \right] 
 =  \left[ \begin{array}{cc} 3& -1\\ -2& 1 \end{array} \right] 
\end{eqnarray*}

\section{Atualiza\c{c}\~{o}es Est\'{a}veis}

O {\bf Algoritmo de Bartels e Golub}, que apresentaremos a seguir, nos
d\'{a} uma atualiza\c{c}\~{a}o est\'{a}vel da fatora\c{c}\~{a}o LU de
uma matriz.  Tomando $T=L^{-1}$, $A=LU$,
 \index{Algoritmo!Bartels e Golub}
$$T \hat A = [U^1,\ldots U^{s-1}, Ta, U^{s+1},\ldots U^n] \ .$$

Consideremos agora a permuta\c{c}\~{a}o de colunas que comuta a 
$s$-\'{e}sima e a $n$-\'{e}sima colunas, $\hat Q$
Esta permuta\c{c}\~{a}o leva $\hat A$ em $\bar A = \hat A \hat Q$,  
e $T \hat A$ em
$$
T \hat A \hat Q = 
\left[ \begin{array}{ccccccc} 
 U^1_1 & \ldots & U^{s-1}_1 & T_1a & U^{s+1}_1 & \ldots & U^n_1 \\
 \vdots & \ddots & \vdots & \vdots & & & \vdots \\
 0 &  & U^{s-1}_{s-1} & T_{s-1}a & U^{s+1}_{s-1} & & U^n_{s-1} \\
 0 & & 0 & T_s a & U^{s+1}_s & & U^n_s \\
 0 &  & 0 & T_{s+1}a & U^{s+1}_{s+1} & & U^n_{s+1} \\
 \vdots & & &\vdots & & \ddots & \\ 
 0 &\ldots  & 0 & T_na & 0 & \ldots & U^n_n 
\end{array} \right]  Q = $$
$$
H= \left[ \begin{array}{cccccccc} 
 U^1_1 & \ldots & U^{s-1}_1 & U^{s+1}_1 & \ldots &U^{n-1}_1 & U^n_1 & T_1a \\
 0 & \ddots &  & & & & & \vdots \\
 \vdots &\ddots & U^{s-1}_{s-1} & U^{s+1}_{s-1} & & 
   U^{n-1}_{s-1}& U^n_{s-1} & T_{s-1}a \\
 0 & & 0 & U^{s+1}_s & & U^{n-1}_s & U^n_s & T_sa \\
 0 & & 0 & U^{s+1}_{s+1} & \ddots & U^{n-1}_{s+1} & U^n_{s+1} & T_{s+1}a \\
 \vdots & & & &\ddots & \ddots & & \vdots \\ 
 0 &  & 0 & 0 & & U^{n-1}_{n-1} & U^n_{n-1} &  T_{n-1}a \\
 0 &\ldots  & 0 & 0 & \ldots & 0 &  U^n_n & T_na 
\end{array} \right] $$

$H$ \'{e} uma matriz de {\bf Hessenberg} superior, isto \'{e}, apenas os
elementos paralelos \`{a} diagonal principal podem ser n\~{a}o nulos no
tri\^{a}ngulo inferior. O m\'{e}todo de Bartels e Golub consiste na
aplica\c{c}\~{a}o de m\'{e}todo de Gauss, com pivotamento parcial, \`{a}
matriz $H$. 
 \index{Matriz!Hessenberg}

Observemos que a aplica\c{c}\~{a}o do m\'{e}todo de Gauss \`{a} matriz 
$H$ \'{e} extremamente simples, pois
\begin{enumerate}

\item As primeiras $s-1$ etapas de transforma\c{c}\~{a}o n\~{a}o
s\~{a}o necess\'{a}rias, pois nas colunas $1,\ldots s-1$, $H$ 
j\'{a} \'{e} triangular superior. 

\item As etapas $j= s,\ldots n-1$, que transformam 

$$^0H=\; ^{s-1}H \rightarrow \; ^sH \rightarrow \ldots 
   \; ^jH \rightarrow \ldots \; ^{n-1}H = \; \hat U $$
resumem-se em
\begin{enumerate}

\item Permutar as linhas $j$ e $j+1$ se
      $ |\; ^{j-1}H^j_{j+1} \; | > |\; ^{j-1}H^j_j \; | $, 
obtendo uma nova matriz $^{j-1}{\tilde H}$.

\item Calcular um \'{u}nico multiplicador 
 $^jN^j_{j+1}= \; ^{j-1}{\tilde H}^j_{j+1} / \; ^{j-1}{\tilde H}^j_j $.

\item Atualizar uma \'{u}nica linha
$ ^jH_{j+1} = \; ^{j-1}{\tilde H}_{j+1} 
  - \; ^jN^j_{j+1} \; ^{j-1}{\tilde H}_j $.

\end{enumerate}
\end{enumerate}

A aplica\c{c}\~{a}o do metodo de Gauss nos d\'{a} a fatora\c{c}\~{a}o 
$\tilde H =  \hat L  \hat U$, onde $\tilde H =  \hat R H$ \'{e} a matriz
obtida de $H$ pelas permuta\c{c}\~{o}es de
linhas realizadas durante a triangulariza\c{c}\~{a}o. 
Assim,  
\begin{eqnarray*}
\tilde H &=&  \hat RH =  \hat RT \bar A = \hat L \hat U \ ,\ \ \mbox{donde}\\ 
\bar A &=& L \hat R^t \hat L \hat U \ , \ \ \mbox{ou}\\
\hat V= \bar A^{-1} &=&  \hat U^{-1} \hat T \hat R T 
\end{eqnarray*}

Ap\'{o}s uma seq\"{u}\^{e}ncia de mudan\c{c}as de base, nossa 
representa\c{c}\~{a}o
da inversa teria a forma 
$$ ^kV= \; ^kU^{-1}\; ^kT\; ^kR \ldots ^1T\; ^1R\; T \ .$$

\section{Preservando Esparsidade}

O {\bf algoritmo de Saunders}, que agora examinamos, pode ser visto como
uma adapta\c{c}\~{a}o do algoritmo de Bartels e Golub visando aproveitar
a pr\'{e}via estrutura\c{c}\~{a}o da base pelo algoritmo P4. 
 \index{Algoritmo!Saunders} 

Suponhamos termos a fatora\c{c}\~{a}o $A=LU$ obtida pela
aplica\c{c}\~{a}o do m\'{e}todo de Gauss (sem pivotamento de linhas mas
com poss\'{\i}veis permuta\c{c}\~{o}es de colunas) \`{a} matriz $A$
previamente estruturada pelo P4, por exemplo

Exemplo 3: 

$$
\begin{array}{c||c|ccc|c|ccccc}
   & 1 & 2 & 3 & 4 & 5 & 6 & 7 & 8 & 9 & 10 \\ 
\hline \hline 
 1 & x &   &   &   &   &   &   &   &   &   \\
\hline 
 2 &  & x &   & x &   &   &   &   &   &   \\ 
 3 & x & x & x & 0 &   &   &   &   &   &   \\
 4 &  &   & x & x &   &   &   &   &   &   \\ 
\hline 
 5 & x &   &   & x & x &   &   &   &   &   \\ 
\hline 
 6 & x &   & x & 0 &   & x &   &   &   & x \\ 
 7 &  & x &   & x & x &   & x & x &   & x \\ 
 8 &  & x &   & 0 &   &   & x & x &   & x \\ 
 9 &  & x &   & x &   & x &   & x & x & x \\ 
 10 & x &   &   &   & x &   &   & x & x & 0    
\end{array}  $$ 
 No Exemplo 3, ``x'' indica as posi\c{c}\~{o}es originalmente n\~{a}o
nulas de A e ``0'' indica as posi\c{c}\~{o}es preenchidas durante a
triangulariza\c{c}\~{a}o.  

Seja $\tilde U$ a matriz obtida de $U$ pela permuta\c{c}\~{a}o
sim\'{e}trica, $Q'UQ$, que leva os espinhos, preservando sua ordem de
posicionamento, para as \'{u}ltimas colunas \`{a} direita. 
No Exemplo 3, ter\'{\i}amos, 
$$
\tilde U= Q'UQ= \begin{array}{c||ccccccc|ccc} 
 & 1 & 2 & 3 & 5 & 6 & 7 & 9 & 4 & 8 & 10 \\
\hline \hline 
 1 & x & & & & & & & & & \\ 
 2 & & x & & & & & & x & & \\ 
 3 & & & x & & & & & x & & \\ 
 5 & & & & x & & & & & & \\ 
 6 & & & & & x & & & & & \\ 
 7 & & & & & & x & & & & x \\ 
 9 & & & & & & & x & & x & x \\ 
\hline 
 4 & & & & & & & & x & & \\ 
 8 & & & & & & & & & x & x \\ 
10 & & & & & & & & & & x \\ 
\end{array} $$ 
Esta matriz tem estrutura
$\tilde U = \left[ \begin{array}{cc} D & E\\ 0 & F \end{array} \right]$, 
onde $D$ \'{e} diagonal e $F$ \'{e} triangular superior.

O algoritmo de Saunders atualiza a decomposi\c{c}\~{a}o  
$A=LQ\tilde UQ'$ para a nova base
\begin{eqnarray*}
\hat A &=& \left[ \begin{array}{ccccccc} 
       A^1 & \ldots & A^{s-1} & a & A^{s+1} & \ldots & A^n 
       \end{array} \right] \ \ \mbox{ou} \\
\hat A \hat Q &=& \left[ \begin{array}{ccccccc} 
       A^1 & \ldots & A^{s-1} & A^{s+1} & \ldots & A^n & a  
       \end{array} \right]
\end{eqnarray*}       
como segue:
\begin{enumerate}

\item 
Forma, pela permuta\c{c}\~{a}o $\tilde U \hat Q$, e substitui\c{c}\~{a}o da
\'{u}ltima coluna, 
$$
  W= T \hat A \hat Q = Q'T AQ \hat Q = \left[ \begin{array}{ccccccc} 
  {\tilde U}^1 & \ldots & {\tilde U}^{s-1} & {\tilde U}^{s+1} & \ldots & 
  {\tilde U}^n & Q'Ta  \end{array} \right] \ .$$

\item
Forma, pela permuta\c{c}\~{a}o sim\'{e}trica, 
\begin{eqnarray*}
 \hat W =  \hat Q' T \hat A \hat Q = 
 = \hat Q' Q' T \hat A Q  \hat Q =  \hat W = \\ 
 \left[ \begin{array}{ccccccc} 
 (W_1)' & \ldots & (W_{s-1})' & (W_{s+1})' & \ldots & (W_n)' & (W_s)' 
\end{array} \right] ' \ .
\end{eqnarray*}

Note que $\hat W$ tem uma estrutura do tipo
$$ \hat W= \left[ \begin{array}{ccc} 
   \hat D & \hat E & \hat W^n_{1:n-c-1} \\ 
   0 & \hat F & \hat W^n_{n-c:n-1} \\ 
   0 & W_s & \hat W^n_n \end{array} \right] $$ 
onde $\hat D$ \'{e} diagonal e $W_s$ s\'{o} pode ter elementos n\~{a}o
nulos nas \'{u}ltimas posi\c{c}\~{o}es \`{a} direita 
(sob $\hat F$ e $\hat W^n$ ). 

\item 
Triangulariza $\hat W$, isto \'{e}, $\hat N$, pelo m\'{e}todo de Gauss
com pivotamento parcial.  A matriz $\hat N$, constitu\'{\i}da de $\hat F$ 
e dos \'{u}ltimos elementos da linha $W_s$ e da coluna $\hat W^n$, 
\'{e} denominada {\bf n\'{u}cleo} de $\hat W$. Nas rotinas num\'{e}ricas 
pode ser conveniente guardar apenas o n\'{u}cleo, $c\times c,\ c<<n$ 
como uma matriz densa na mem\'{o}ria principal, enquanto o resto da 
matriz, usada apenas para leitura, pode ser guardada em uma 
representa\c{c}\~{a}o esparsa, e na mem\'{o}ria secund\'{a}ria.   

\end{enumerate}

No Exemplo 3 ter\'{\i}amos, ao mudar a coluna 3 da base A,
$$
\hat A = \begin{array}{c||cccccccccc} 
   & 1 & 2 & 3 & 4 & 5 & 6 & 7 & 8 & 9 & 10 \\ 
\hline \hline 
 1 & x &   & y &   &   &   &   &   &   &   \\ 
 2 &   & x &   & x &   &   &   &   &   &   \\ 
 3 & x & x &   &  &   &   &   &   &   &   \\
 4 &   &   & y & x &   &   &   &   &   &   \\ 
 5 & x &   &   & x & x &   &   &   &   &   \\ 
 6 & x &   &   &  &   & x &   &   &   & x \\ 
 7 &   & x & y & x & x &   & x & x &   & x \\ 
 8 &   & x &   &  &   &   & x & x &   & x \\ 
 9 &   & x &   & x &   & x &   & x & x & x \\ 
10 & x &   & y &   & x &   &   & x & x &     
\end{array} $$ 

$$
 W =  \begin{array}{c||cccccc|ccc|c}
   & 1 & 2 & 5 & 6 & 7 & 9 & 4 & 8 & 10 & 3 \\
\hline \hline
 1 & x &   &   &   &   &   &   &   &   & y \\ 
 2 &  & x &   &   &   &   & x &   &   &   \\ 
 3 &  &   &   &   &   &   & x &   &   &   \\
 5 &  &   & x &   &   &   &   &   &   &   \\ 
 6 &  &   &   & x &   &   &   &   & x &   \\ 
 7 &  &   &   &   & x &   &   & x & x & y \\ 
 9 &  &   &   &   &   & x &   &   & x &   \\ 
 4 &  &   &   &   &   &   & x &   &   & y \\ 
 8 &  &   &   &   &   &   &   & x & x &   \\ 
 10 &  &   &   &   &   &   &   &   & x & y   
\end{array} $$ 

$$
 \hat W =  \begin{array}{c||cccccc|ccc|c}
   & 1 & 2 & 5 & 6 & 7 & 9 & 4 & 8 & 10 & 3 \\
\hline \hline
 1 & x &   &   &   &   &   &   &   &   & y \\ 
 2 &  & x &   &   &   &   & x &   &   &   \\ 
 5 &  &   & x &   &   &   &   &   &   &   \\ 
 6 &  &   &   & x &   &   &   &   & x &   \\ 
 7 &  &   &   &   & x &   &   & x & x & y \\ 
 9 &  &   &   &   &   & x &   &   & x &   \\ 
\hline 
 4 &  &   &   &   &   &   & x &   &   & y \\ 
 8 &  &   &   &   &   &   &   & x & x &   \\ 
 10 &  &   &   &   &   &   &   &   & x & y \\
\hline  
 3 &  &   &   &   &   &   & x &   &   &  0 \\
\end{array} $$

Ap\'{o}s a triangulariza\c{c}\~{a}o de $\hat W$ teremos, sendo 
$\hat R$ as permuta\c{c}\~{o}es de linhas realizadas ao triangularizar 
$\hat W$, 

\begin{eqnarray*}
\hat R \hat W &=& \hat L \hat U = \hat R \hat Q' Q' T \hat A Q \hat Q 
\ \ \mbox{portanto} \\
\hat A &=& L Q \hat Q \hat R' \hat L \hat U \hat Q' Q' 
\ \ \mbox{e} \\ 
{\hat A}^{-1} &=& Q \hat Q {\hat U}^{-1} \hat T \hat R \hat Q' Q' T \ . 
\end{eqnarray*}

Em geral, ap\'{o}s uma seq\"{u}\^{e}ncia de mudan\c{c}as de base, 
nossa representa\c{c}\~{a}o da inversa ter\'{a} a forma 
$$
^kA^{-1} = Q \; ^1Q \ldots \; ^kQ \; ^kU^{-1} 
  \; ^kT \; ^kR \; ^kQ' \ldots \; ^1T \; ^1R \; ^1Q \; Q \; T \ .$$

\begin{obs}{\rm  
\mbox{} \\   
 \begin{enumerate} 
 \item 
A regi\~{a}o dos preenchimentos em  $^kW$
est\'{a} restrita ao tri\^{a}ngulo superior do n\'{u}cleo e sua
\'{u}ltima linha. 
 \item 
A matriz $^kL$  s\'{o} ter\'{a} elementos nulos na \'{u}ltima
linha do n\'{u}cleo.  
 \item 
A cada nova mudan\c{c}a de base, a dimens\~{a}o do n\'{u}cleo:
   \begin{enumerate} 
   \item 
Aumenta de 1, se a coluna que sai da base \'{e} uma coluna triangular
da base original.  
   \item 
Permanece constante, se a coluna que sai da base \'{e} um espinho, ou
uma coluna j\'{a} anteriormente substitu\'{\i}da. 
   \end{enumerate} 
 \item 
Cada mudan\c{c}a de base aumenta um termo na fatora\c{c}\~{a}o da
base, o que leva a uma gradativa perda de efici\^{e}ncia e ac\'{u}mulo
de erros.  Depois de um n\'{u}mero pr\'{e}-determinado de mudan\c{c}as
sobre uma base original, ou quando o ac\'{u}mulo de erros assim o
exigir, devemos fazer uma reivers\~{a}o, i.\'{e}., reiniciar o processo
aplicando a P4 e triangularizando a pr\'{o}xima base desejada.  
 \end{enumerate}
 }\end{obs} 
 \index{N\'{u}cleo} 
 \index{Espinhos}

\section{Preservando Estrutura}

Consideraremos agora o problema de atualizar a fatora\c{c}\~{a}o $A=QU$
de uma base na forma angular blocada.  Conforme argumentamos no final do
cap\'{\i}tulo 7, estamos apenas interessados no fator $U$, sendo as
transforma\c{c}\~{o}es ortogonais descartadas logo ap\'{o}s o seu uso. 
No que segue, a coluna a sair da base ser\'{a} a coluna $outj$ do bloco
$outk$, $1\leq outk\leq h+1$.  Em seu lugar ser\'{a} introduzida na base
uma coluna com a estrutura do bloco $ink$, $1\leq ink\leq h+1$. 
 \index{Estrutura!angular} 
 \index{Estrutura!blocada} 
 \index{Fatora\c{c}\~{a}o!blocada} 
 \index{Atualiza\c{c}\~{a}o!blocada}

Apresentamos agora um procedimento de atualiza\c{c}\~{a}o de $U$ por
blocos, $bup()$, que usa explicitamente a estrutura angular blocada da
base $A$ e do fator $U$.  Este procedimento ser\'{a} descrito em termos
das opera\c{c}\~{o}es simples em cada bloco apresentadas na
se\c{c}\~{a}o 7.5. 

Em $bup()$ consideraremos cinco casos distintos:
\begin{description}
\item[Caso I] $ink \ne outk , \ \ ink \ne h+1 , \ \ outk \ne h+1$. 
\item[Caso II] $ink = outk ,\ \ ink \ne h+1$. 
\item[Caso III] $ink\ne h+1, \ \ outk = h+1$.  
\item[Caso IV] $ink = h+1 , \ \ outk\ne h+1$.  
\item[Caso V] $ink = outk, \ \ ink = h+1$. 
\end{description}

Vejamos em detalhe o caso I, quando $ink$ and $outk$ s\~{a}o blocos
diagonais distintos, como mostrado na figura 1.  Neste caso os
\'{u}nicos ENN's na coluna saindo est\~{a}o no bloco
$^{\mbox{outk}}B_{\bullet}^{\mbox{outj}}$.  Analogamente os \'{u}nicos
ENN's na coluna entrando na base, $a$, est\~{a}o em $^{\mbox{ink}}a$. 

Definimos $y\equiv A'a$, e $u\equiv Q'a = U^{-t}A'a = U^{-t}y$.  Notamos
que os vetores $y$ e $u$ preservam a estrutura blocada da base.  Assim,
os ENN's em $u$ est\~{a}o nos blocos $^{ink}u$ e
$^{h+1}u$,
   \[  \left[ \begin{array}{c} ^{ink}u \\ ^{h+1}u \\ \end{array}  \right] =
       \left[ \begin{array}{cc} ^{ink}V & ^{ink}W \\
                                0 & S \\ \end{array}  \right]^{-t}
       \left[ \begin{array}{c} ^{ink}y \\ ^{h+1}y \\ \end{array}  \right] \]

Para atualizar $U$ removemos a coluna $outj$ do bloco $outk$, e
inserimos $u$ como a \'{u}ltima coluna de $^{ink}U$.  Depois apenas temos
que reduzir $U$ a uma matriz triangular superior atrav\'{e}s de
transforma\c{c}\~{o}es ortogonais. 

\begin{enumerate}

\item  Leve   
   $ \left[ \begin{array}{cc} ^{outk}V & ^{outk}W \\ \end{array} \right] $ 
   de Hessenberg a triangular superior.

\item Leve
   $ \left[ \begin{array}{cc} ^{h+1}u & S \\ \end{array} \right] $
   a triangular.  

\item Insira a primeira linha de $^{h+1}U$ , como a \'{u}ltima linha de 
      $^{ink}U$. 

\item Insira a \'{u}ltima linha de  $^{outk}U$
      como a primeira de $^{h+1}U$.

\item Leve $S$ de Hessenberg a triangular superior. 

\end{enumerate}

Os outros casos s\~{a}o bastante similares.  Os esquemas nas figuras 2,
3 e 4 s\~{a}o os an\'{a}logos do esquema na figura 1, e d\~{a}o uma
descri\c{c}\~{a}o sum\'{a}ria de cada um dos casos.

\noindent
Estes s\~{a}o os passos para o caso I :  
 
\begin{enumerate}

\item No n\'{o} $ink$, compute $y^{ink}=(B^{ink})^{t}a^{ink}$ e 
 $y^{b+1}=(C^{ink})^{t}a^{ink}$.\\
 $pTime= m(ink)n(ink) + m(ink)n(b+1) \le 2dbmax^{2}$ , 
 $INC=0$ .

\item No n\'{o} $ink$, compute a transforma\c{c}\~{a}o inversa parcial
   \[  \left[ \begin{array}{c} u^{ink} \\ z \\ \end{array}  \right] = 
       \left[ \begin{array}{cc} V^{ink} & W^{ink} \\ 
                                0 & I \\ \end{array}  \right]^{-t} 
       \left[ \begin{array}{c} y^{ink} \\ y^{b+1} \\ \end{array}  \right] \]
 Ent\~{a}o insira $u^{ink}$ como a \'{u}ltima coluna de $V^{ink}$.\\
 $pTime = (1/2)n(ink)^{2} + n(ink)n(b+1) \le (3/2)dbmax^{2}$ , $INC=0$ .

\item Do n\'{o} $ink$ envie $z$ ao n\'{o} $0$.\\
 $pTime = 0$ , $INC = n(b+1) \le dbmax$ . 

\item No n\'{o} $0$ compute $u^{b+1}=S^{-t}z$.\\    
 $pTime = (1/2)n(b+1)^{2} \le (1/2)dbmax^{2}$ , $INC = 0$.

\item 
 \begin{enumerate}
  \item No n\'{o} $outk$, remova a coluna $V^{outk}_{\bullet,outj}$ de 
   $V^{outk}$. Ent\~{a}o reduza
   $ \left[ \begin{array}{cc} V^{outk} & W^{outk} \\ \end{array} \right] $ 
   de Hessenberg para triangular superior.
  \item No n\'{o} $0$, reduza 
   $ \left[ \begin{array}{cc} u^{b+1} & S \\ \end{array} \right] $
   para triangular superior.  
 \end{enumerate}
 Observe que as opera\c{c}\~{o}es nos passos $5a$ e $5b$ s\~{a}o 
 independentes, portanto\\  
 $pTime = 2n(ink)^{2} + 4n(ink)n(b+1) \wedge 2n(b+1)^{2} \le 6dbmax^{2}$ ,
 $INC = 0$.

 \item Do n\'{o} $0$ envie o vetor $S_{1,\bullet}$ ao n\'{o} $ink$, onde
  ele \'{e} inserido como a \'{u}ltima coluna de $W^{ink}$. 
  Do n\'{o} $0$ envie o elemento $u^{b+1}_{1}$ ao n\'{o} $ink$, onde ele 
  \'{e} inserido como $U^{ink}_{n(ink)+1,n(ink)+1}$.    
  Do n\'{o} $outk$ envie vetor $W^{outk}_{n(outk),\bullet}$ ao n\'{o}
  $0$, onde ele \'{e} inserido como a primeira linha de $S$.\\
  $pTime = 0$ , $INC = 2n(b+1) + n(outk) \le 3dbmax$.

\item No n\'{o} $0$, reduza  $S$ de Hessenberg para triangular superior.\\ 
 $pTime = 2n(b+1)^{2} \le 2dbmax^{2}$ , $INC = 0$. 

\end{enumerate}

\noindent
Estes s\~{a}o os passos para o caso II : 

Os Passos 1---5 s\~{a}o exatamente como no caso I.  
\begin{enumerate}
\setcounter{enumi}{5}

\item 
 \begin{enumerate}
   \item Do n\'{o} $ink$ envie 
    $V^{ink}_{n(ink),n(ink)}$ e $W^{ink}_{n(ink),\bullet}$ ao n\'{o} $0$.

    \item No n\'{o} $0$ reduza para triangular superior a matriz 
     $2\times n(b+1)+1$ 
     \[ \left[ \begin{array}{cc} 
        V^{ink}_{n(ink),n(ink)} & W^{ink}_{n(ink),\bullet} \\ 
        u^{b+1}_{1} & S_{1,\bullet} \\  \end{array} \right] \] 
 \end{enumerate}
 $pTime = 4n(b+1) \le 4dbmax$ , $INC = n(b+1) \le dbmax$. 
   
\item 
 \begin{enumerate}
    \item Do n\'{o} $0$ envie o vetor modificado  
     $  \left[ \begin{array}{cc} 
               V^{ink}_{n(ink),n(ink)} & W^{ink}_{n(ink),\bullet} \\
                                                       
\end{array} \right] $
     de volta ao n\'{o} $ink$.
  
    \item No n\'{o} $0$, reduza  $S$ de Hessenberg para triangular superior.
 \end{enumerate}
 $pTime = 2n(b+1)^{2} \le 2dbmax^{2}$ , $INC = n(b+1) \le dbmax$.  

\end{enumerate}

\noindent
Estes s\~{a}o os passos do caso III : 

Os passos 1---4 s\~{a}o exatamente como no caso I.  
\begin{enumerate}
\setcounter{enumi}{4}

\item 
 \begin{enumerate}
   \item No n\'{o} $k=1:b$ remova a coluna $W^{k}_{\bullet,outj}$ de
    $W^{k}$.  
   \item No n\'{o} $0$ reduza
    $\left[ \begin{array}{cc} u^{b+1} & S \\ \end{array} \right] $
    para triangular superior. Remova $S_{\bullet,outj}$ de $S$.
 \end{enumerate}
 $pTime = 2n(b+1)^{2} \le 2dbmax^{2}$ , $INC = 0$.

\item Do n\'{o} $0$ envie ao n\'{o} $ink$, $u^{b+1}_{1}$ para ser inserido
      em $V^{ink}$ como $V^{ink}_{n(ink)+1,n(ink)+1}$, e 
      $S_{1,\bullet}$ para ser inserido como a \'{u}ltima coluna de $W^{ink}$.\\
 $pTime = 0$ , $INC = n(b+1) \le dbmax$.

\item No n\'{o} $0$, reduza  $S$ de Hessenberg para triangular superior.\\
 $pTime = 2n(b+1)^{2} \le 2dbmax^{2}$ , $INC = 0$. 

\end{enumerate}

\noindent 
Estes s\~{a}o os passos do caso IV :  
 
\begin{enumerate}

\item No n\'{o} {\it k=1:b}, compute $y^{k}=(B^{k})^{t}a^{k}$ e 
 $x^{k}=(C^{k})^{t}a^{k}$.\\
 $pTime = \wedge_{1}^{b} \ m(k)n(ink) + m(k)n(b+1) \le 2dbmax^{2}$ ,
 $INC = 0$.

\item No n\'{o} {\it k=1:b}, compute a transforma\c{c}\~{a}o inversa parcial
   \[  \left[ \begin{array}{c} u^{k} \\ z^{k} \\ \end{array}  \right] = 
       \left[ \begin{array}{cc} V^{k} & W^{k} \\ 
                                0 & I \\ \end{array}  \right]^{-t} 
       \left[ \begin{array}{c} y^{ink} \\ x^{k} \\ \end{array}  \right] \]
 e insira $u^{k}$ como a \'{u}ltima coluna de $W^{k}$.\\
 $pTime = \wedge_{1}^{b} \ (1/2)n(k)^{2} + n(k)n(b+1) \le (3/2)dbmax^{2}$ ,
 $INC = 0$.

\item Do n\'{o} {\it k=1:b} envie $z^{k}$ ao n\'{o} $0$, onde acumulamos
 $z = \sum_{1}^{b} z^{k}$.\\ 
 $pTime = b\ n(b+1) \le b\ dbmax$ ,  $INC = b\ n(b+1) \le b\ dbmax$. 

\item No n\'{o} $0$ compute $u^{b+1}=S^{-t}z$,
 e insira $u^{b+1}$ como a \'{u}ltima coluna de $S$.\\    
 $pTime = (1/2)n(b+1)^{2} \le (1/2)dbmax^{2}$ , $INC = 0$.

\item remova a coluna $V^{outk}_{\bullet,outj}$ de $V^{outk}$ , 
  e reduza 
  $\left[ \begin{array}{cc} V^{outk} & W^{outk} \\ \end{array}  \right]$ 
  para triangular superior. \\ 
  $pTime = 2n(outk)^{2} + 4n(outk)n(b+1) \le 6dbmax^{2}$ , $INC = 0$.

\item Envie o vetor $W^{outk}_{n(ouk),\bullet}$ do n\'{o} $outk$ ao
n\'{o} $0$, onde o inserimos como a primeira linha de $S$, e reduza $S$
a triangular.\\
 $pTime = 2n(b+1)^{2} \le 2dbmax^{2}$ , $INC = n(b+1) \le dbmax$. 
 
\end{enumerate}

\noindent 
Estes s\~{a}o os passos do caso V :  

Os passos 1---4 s\~{a}o exatamente como no caso IV.  
\begin{enumerate}
\setcounter{enumi}{4}

\item No n\'{o} {\it k=1:b}, remova a coluna 
 $W^{k}_{\bullet,outj}$ de $W^{outk}$ , 
 e insira  $u^{k}$ como a \'{u}ltima coluna de $W^{k}$.  
 No n\'{o} $0$ remova a coluna $S_{\bullet,outj}$ de $S$ , e insira  
 $u^{b+1}$ como a \'{u}ltima coluna de $S$.\\
 $pTime = 0$ , $INC = 0$.

\item No n\'{o} $0$, reduza  $S$ de Hessenberg para triangular superior. \\
 $pTime = 2n(b+1)^{2} \le 2dbmax^{2}$ , $INC = 0$. 

\end{enumerate}

A complexidade do procedimento de atualiza\c{c}\~{a}o de $U$ por blocos,
$bup()$, \'{e} dada pelo seguinte teorema: Como na se\c{c}\~{a}o 7.4, 
$dbmax = \max_{k=1}^h n(k)$. 

\begin{teo}
A complexidade de $bup()$, desconsideramos termos de ordem inferior, tem
os seguintes limitantes superiores para tempo de processamento e
comunica\c{c}\~{a}o:
\begin{center}
\begin{tabular}{|l|l|l|} \hline
 Caso & Processamento & Comunica\c{c}\~{a}o \\ \hline
 I   & $12dbmax^{2}$ & $3dbmax$ \\ \hline
 II  & $12dbmax^{2}$ & $3dbmax$ \\ \hline
 III & $8dbmax^{2}$ & $2dbmax$  \\ \hline
 IV  & $12dbmax^{2}+b\ dbmax$ & $h\ dbmax$  \\ \hline
 V   & $6dbmax^{2}+b\ dbmax$  & $h\ dbmax$  \\ \hline
\end{tabular}
\end{center}
\end{teo}

Caso o ambiente permita comunica\c{c}\~{o}es em paralelo, podemos
substituir $h$ por $log(h)$ nas express\~{o}es de complexidade. 

\pagebreak 

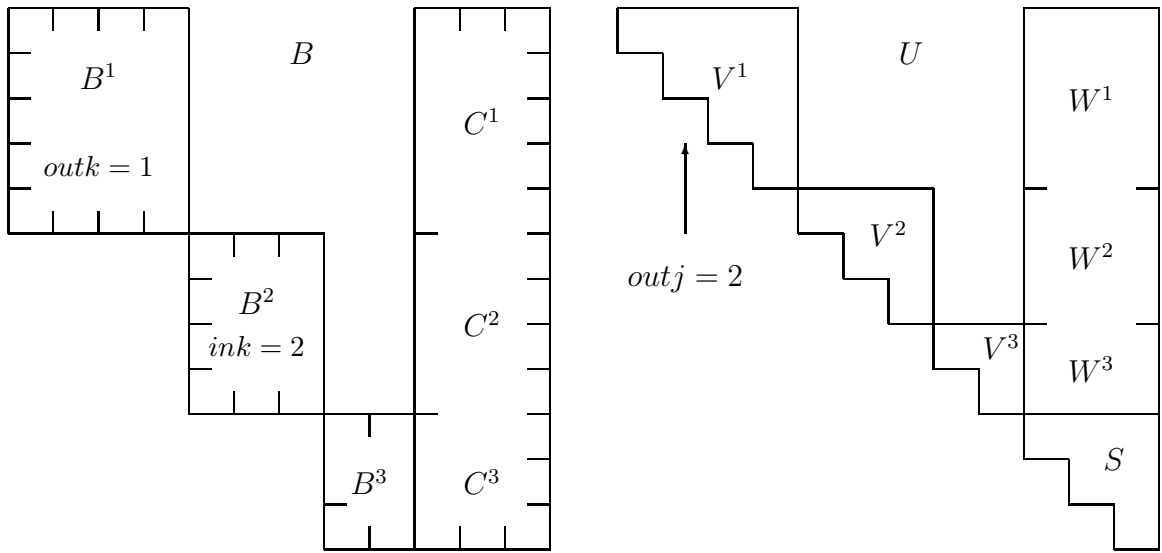
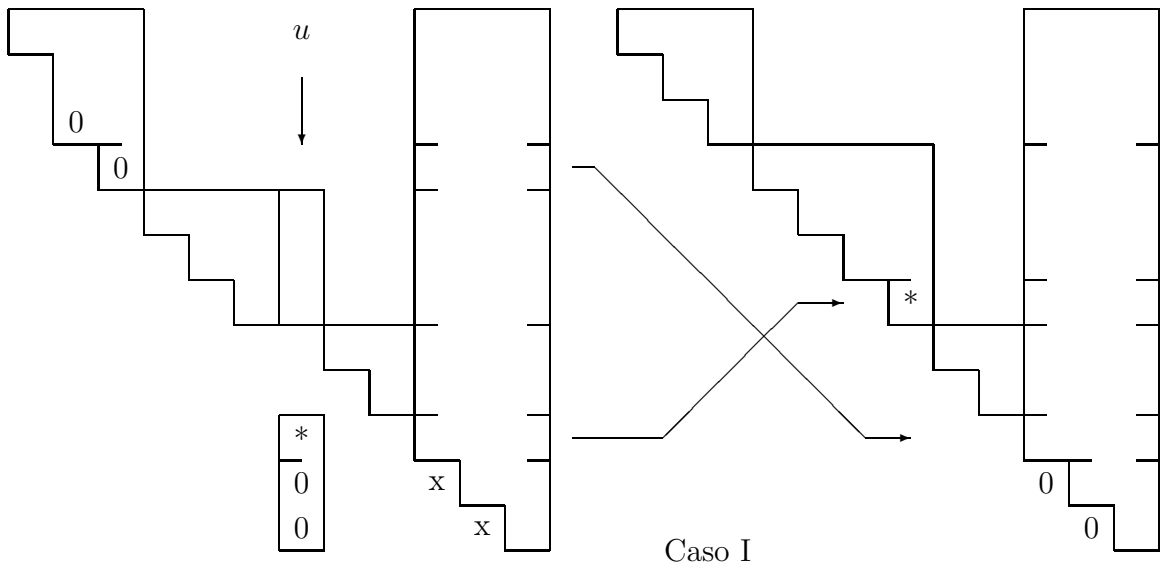
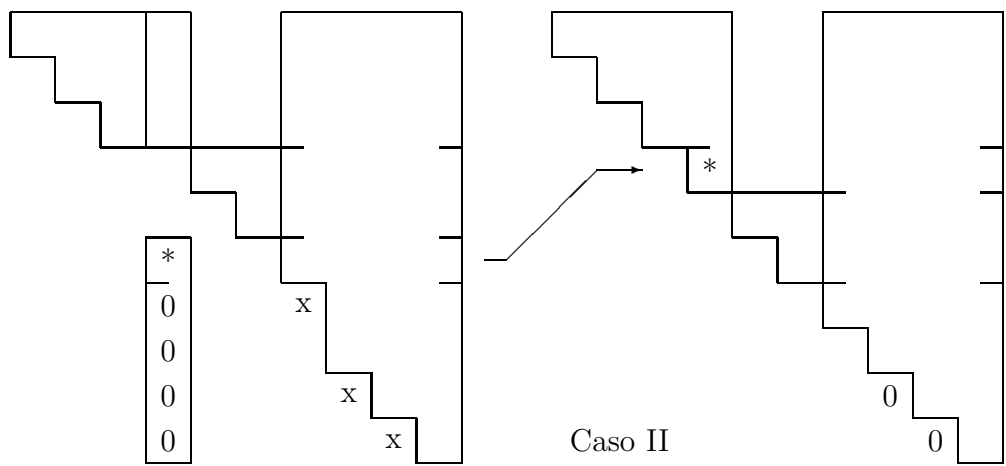
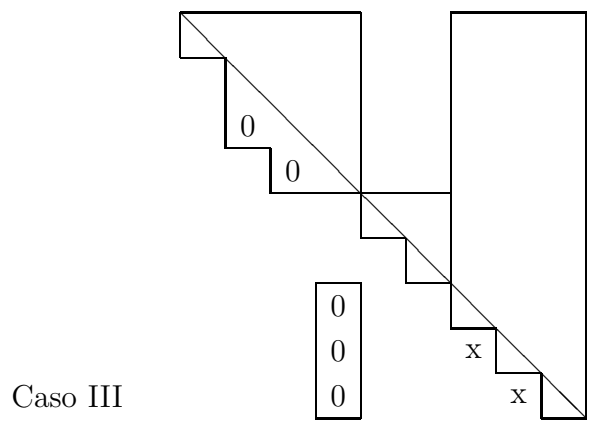
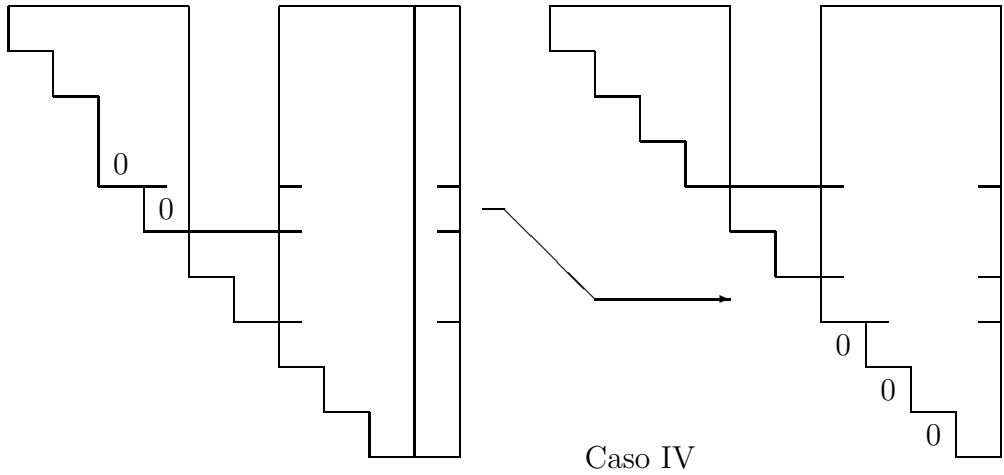
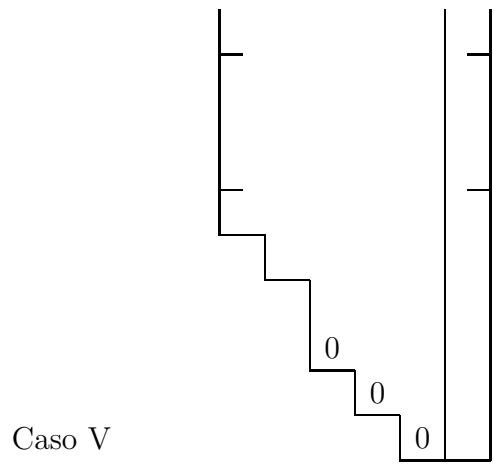
\mbox{}\\ 

\unitlength=1.00mm
\special{em:linewidth 0.4pt}
\linethickness{0.4pt}
\begin{picture}(153.00,72.00)
\put(6.00,72.00){\line(0,-1){3.00}}
\put(12.00,72.00){\line(0,-1){3.00}}
\put(18.00,72.00){\line(0,-1){3.00}}
\put(0.00,66.00){\line(1,0){3.00}}
\put(0.00,60.00){\line(1,0){3.00}}
\put(0.00,54.00){\line(1,0){3.00}}
\put(0.00,48.00){\line(1,0){3.00}}
\put(6.00,42.00){\line(0,1){3.00}}
\put(12.00,42.00){\line(0,1){3.00}}
\put(18.00,45.00){\line(0,-1){3.00}}
\put(30.00,42.00){\line(0,-1){3.00}}
\put(36.00,42.00){\line(0,-1){3.00}}
\put(24.00,36.00){\line(1,0){3.00}}
\put(24.00,30.00){\line(1,0){3.00}}
\put(24.00,24.00){\line(1,0){3.00}}
\put(30.00,18.00){\line(0,1){3.00}}
\put(36.00,18.00){\line(0,1){3.00}}
\put(48.00,18.00){\line(0,-1){3.00}}
\put(42.00,12.00){\rule{3.00\unitlength}{0.00\unitlength}}
\put(42.00,6.00){\line(1,0){3.00}}
\put(48.00,0.00){\line(0,1){3.00}}
\put(60.00,0.00){\line(0,1){3.00}}
\put(66.00,0.00){\line(0,1){3.00}}
\put(72.00,6.00){\line(-1,0){3.00}}
\put(72.00,12.00){\line(-1,0){3.00}}
\put(72.00,18.00){\line(-1,0){3.00}}
\put(72.00,24.00){\line(-1,0){3.00}}
\put(72.00,30.00){\line(-1,0){3.00}}
\put(72.00,36.00){\line(-1,0){3.00}}
\put(72.00,42.00){\line(-1,0){3.00}}
\put(72.00,48.00){\line(-1,0){3.00}}
\put(72.00,54.00){\line(-1,0){3.00}}
\put(72.00,60.00){\line(-1,0){3.00}}
\put(72.00,66.00){\line(-1,0){3.00}}
\put(66.00,72.00){\line(0,-1){3.00}}
\put(60.00,72.00){\line(0,-1){3.00}}
\put(54.00,42.00){\line(1,0){3.00}}
\put(54.00,18.00){\line(1,0){3.00}}
\put(81.00,72.00){\line(1,0){24.00}}
\put(105.00,72.00){\line(0,-1){24.00}}
\put(105.00,48.00){\line(-1,0){6.00}}
\put(99.00,48.00){\line(0,1){6.00}}
\put(99.00,54.00){\line(-1,0){6.00}}
\put(93.00,54.00){\line(0,1){6.00}}
\put(93.00,60.00){\line(-1,0){6.00}}
\put(87.00,60.00){\line(0,1){6.00}}
\put(87.00,66.00){\line(-1,0){6.00}}
\put(81.00,66.00){\line(0,1){6.00}}
\put(105.00,48.00){\line(1,0){18.00}}
\put(123.00,48.00){\line(0,-1){18.00}}
\put(123.00,30.00){\line(-1,0){6.00}}
\put(117.00,30.00){\line(0,1){6.00}}
\put(117.00,36.00){\line(-1,0){6.00}}
\put(111.00,36.00){\line(0,1){6.00}}
\put(111.00,42.00){\line(-1,0){6.00}}
\put(105.00,42.00){\line(0,1){6.00}}
\put(135.00,18.00){\line(1,0){18.00}}
\put(153.00,18.00){\line(0,-1){18.00}}
\put(153.00,0.00){\line(-1,0){6.00}}
\put(147.00,0.00){\line(0,1){6.00}}
\put(147.00,6.00){\line(-1,0){6.00}}
\put(141.00,6.00){\line(0,1){6.00}}
\put(141.00,12.00){\line(-1,0){6.00}}
\put(135.00,12.00){\line(0,1){6.00}}
\put(123.00,30.00){\line(1,0){12.00}}
\put(135.00,30.00){\line(0,-1){12.00}}
\put(135.00,18.00){\line(-1,0){6.00}}
\put(129.00,18.00){\line(0,1){6.00}}
\put(129.00,24.00){\line(-1,0){6.00}}
\put(123.00,24.00){\line(0,1){6.00}}
\put(135.00,30.00){\line(0,1){42.00}}
\put(135.00,72.00){\line(1,0){18.00}}
\put(153.00,72.00){\line(0,-1){54.00}}
\put(153.00,30.00){\line(-1,0){3.00}}
\put(153.00,48.00){\line(-1,0){3.00}}
\put(0.00,72.00){\line(1,0){24.00}}
\put(24.00,72.00){\line(0,-1){30.00}}
\put(24.00,42.00){\line(-1,0){24.00}}
\put(0.00,42.00){\line(0,1){30.00}}
\put(24.00,42.00){\line(1,0){18.00}}
\put(42.00,42.00){\line(0,-1){24.00}}
\put(42.00,18.00){\line(-1,0){18.00}}
\put(24.00,18.00){\line(0,1){24.00}}
\put(42.00,18.00){\line(1,0){12.00}}
\put(54.00,0.00){\line(-1,0){12.00}}
\put(42.00,0.00){\line(0,1){18.00}}
\put(54.00,0.00){\line(0,1){72.00}}
\put(54.00,72.00){\line(1,0){18.00}}
\put(72.00,72.00){\line(0,-1){72.00}}
\put(72.00,0.00){\line(-1,0){18.00}}
\put(12.00,63.00){\makebox(0,0)[cc]{$B^1$}}
\put(96.00,63.00){\makebox(0,0)[cc]{$V^1$}}
\put(147.00,12.00){\makebox(0,0)[cc]{$S$}}
\put(117.00,42.00){\makebox(0,0)[cc]{$V^2$}}
\put(132.00,27.00){\makebox(0,0)[cc]{$V^3$}}
\put(12.00,51.00){\makebox(0,0)[cc]{ {\small $outk=1$} }}
\put(33.00,33.00){\makebox(0,0)[cc]{$B^2$}}
\put(33.00,27.00){\makebox(0,0)[cc]{ {\small $ink=2$} }}
\put(48.00,9.00){\makebox(0,0)[cc]{$B^3$}}
\put(63.00,30.00){\makebox(0,0)[cc]{$C^2$}}
\put(63.00,9.00){\makebox(0,0)[cc]{$C^3$}}
\put(63.00,57.00){\makebox(0,0)[cc]{$C^1$}}
\put(135.00,30.00){\line(1,0){3.00}}
\put(135.00,48.00){\line(1,0){3.00}}
\put(144.00,60.00){\makebox(0,0)[cc]{$W^1$}}
\put(144.00,39.00){\makebox(0,0)[cc]{$W^2$}}
\put(144.00,24.00){\makebox(0,0)[cc]{$W^3$}}
\put(39.00,66.00){\makebox(0,0)[cc]{$B$}}
\put(120.00,66.00){\makebox(0,0)[cc]{$U$}}
\put(90.00,42.00){\vector(0,1){12.00}}
\put(90.00,36.00){\makebox(0,0)[cc]{$outj=2$}}
\end{picture} 

\mbox{}\\   
\mbox{}\\   
 
\unitlength=1.00mm
\special{em:linewidth 0.4pt}
\linethickness{0.4pt}
\begin{picture}(153.00,73.00)
\put(0.00,73.00){\line(1,0){18.00}}
\put(18.00,73.00){\line(0,-1){24.00}}
\put(18.00,49.00){\line(-1,0){6.00}}
\put(12.00,49.00){\line(0,1){6.00}}
\put(12.00,55.00){\line(-1,0){6.00}}
\put(6.00,55.00){\line(0,1){12.00}}
\put(6.00,67.00){\line(-1,0){6.00}}
\put(0.00,67.00){\line(0,1){6.00}}
\put(12.00,55.00){\line(1,0){3.00}}
\put(18.00,49.00){\line(1,0){18.00}}
\put(36.00,49.00){\line(0,-1){18.00}}
\put(36.00,31.00){\line(-1,0){6.00}}
\put(30.00,31.00){\line(0,1){6.00}}
\put(30.00,37.00){\line(-1,0){6.00}}
\put(24.00,37.00){\line(0,1){6.00}}
\put(24.00,43.00){\line(-1,0){6.00}}
\put(18.00,43.00){\line(0,1){6.00}}
\put(36.00,49.00){\line(1,0){6.00}}
\put(42.00,49.00){\line(0,-1){18.00}}
\put(42.00,31.00){\line(-1,0){6.00}}
\put(42.00,31.00){\line(1,0){12.00}}
\put(54.00,31.00){\line(0,-1){12.00}}
\put(54.00,19.00){\line(-1,0){6.00}}
\put(48.00,19.00){\line(0,1){6.00}}
\put(48.00,25.00){\line(-1,0){6.00}}
\put(42.00,25.00){\line(0,1){6.00}}
\put(54.00,31.00){\line(0,1){42.00}}
\put(54.00,73.00){\line(1,0){18.00}}
\put(72.00,73.00){\line(0,-1){54.00}}
\put(72.00,19.00){\line(0,-1){18.00}}
\put(72.00,1.00){\line(-1,0){6.00}}
\put(66.00,1.00){\line(0,1){6.00}}
\put(66.00,7.00){\line(-1,0){6.00}}
\put(60.00,7.00){\line(0,1){6.00}}
\put(60.00,13.00){\line(-1,0){6.00}}
\put(54.00,13.00){\line(0,1){6.00}}
\put(54.00,19.00){\line(1,0){3.00}}
\put(54.00,31.00){\line(1,0){3.00}}
\put(54.00,49.00){\line(1,0){3.00}}
\put(54.00,55.00){\line(1,0){3.00}}
\put(72.00,55.00){\line(-1,0){3.00}}
\put(72.00,49.00){\line(-1,0){3.00}}
\put(72.00,31.00){\line(-1,0){3.00}}
\put(72.00,19.00){\line(-1,0){3.00}}
\put(72.00,13.00){\line(-1,0){3.00}}
\put(36.00,19.00){\line(0,-1){18.00}}
\put(36.00,1.00){\line(1,0){6.00}}
\put(42.00,1.00){\line(0,1){18.00}}
\put(42.00,19.00){\line(-1,0){6.00}}
\put(36.00,13.00){\line(1,0){3.00}}
\put(81.00,73.00){\line(1,0){18.00}}
\put(99.00,73.00){\line(0,-1){18.00}}
\put(99.00,55.00){\line(-1,0){6.00}}
\put(93.00,55.00){\line(0,1){6.00}}
\put(93.00,61.00){\line(-1,0){6.00}}
\put(87.00,61.00){\line(0,1){6.00}}
\put(87.00,67.00){\line(-1,0){6.00}}
\put(81.00,67.00){\line(0,1){6.00}}
\put(99.00,55.00){\line(1,0){24.00}}
\put(123.00,55.00){\line(0,-1){24.00}}
\put(123.00,31.00){\line(-1,0){6.00}}
\put(117.00,31.00){\line(0,1){6.00}}
\put(117.00,37.00){\line(-1,0){6.00}}
\put(111.00,37.00){\line(0,1){6.00}}
\put(111.00,43.00){\line(-1,0){6.00}}
\put(105.00,43.00){\line(0,1){6.00}}
\put(105.00,49.00){\line(-1,0){6.00}}
\put(99.00,49.00){\line(0,1){6.00}}
\put(117.00,37.00){\line(1,0){3.00}}
\put(123.00,31.00){\line(1,0){12.00}}
\put(135.00,31.00){\line(0,-1){12.00}}
\put(135.00,19.00){\line(-1,0){6.00}}
\put(129.00,19.00){\line(0,1){6.00}}
\put(129.00,25.00){\line(-1,0){6.00}}
\put(123.00,25.00){\line(0,1){6.00}}
\put(135.00,31.00){\line(0,1){42.00}}
\put(135.00,73.00){\line(1,0){18.00}}
\put(153.00,73.00){\line(0,-1){54.00}}
\put(153.00,19.00){\line(0,-1){18.00}}
\put(153.00,1.00){\line(-1,0){6.00}}
\put(147.00,1.00){\line(0,1){6.00}}
\put(147.00,7.00){\line(-1,0){6.00}}
\put(141.00,7.00){\line(0,1){6.00}}
\put(141.00,13.00){\line(-1,0){6.00}}
\put(135.00,13.00){\line(0,1){6.00}}
\put(135.00,19.00){\line(1,0){3.00}}
\put(135.00,31.00){\line(1,0){3.00}}
\put(135.00,55.00){\line(1,0){3.00}}
\put(153.00,55.00){\line(-1,0){3.00}}
\put(153.00,31.00){\line(-1,0){3.00}}
\put(153.00,19.00){\line(-1,0){3.00}}
\put(153.00,13.00){\line(-1,0){3.00}}
\put(141.00,13.00){\line(1,0){3.00}}
\put(39.00,64.00){\vector(0,-1){9.00}}
\put(9.00,58.00){\makebox(0,0)[cc]{0}}
\put(15.00,52.00){\makebox(0,0)[cc]{0}}
\put(39.00,10.00){\makebox(0,0)[cc]{0}}
\put(39.00,4.00){\makebox(0,0)[cc]{0}}
\put(57.00,10.00){\makebox(0,0)[cc]{x}}
\put(63.00,4.00){\makebox(0,0)[cc]{x}}
\put(39.00,70.00){\makebox(0,0)[cc]{$u$}}
\put(39.00,16.00){\makebox(0,0)[cc]{*}}
\put(120.00,34.00){\makebox(0,0)[cc]{*}}
\put(138.00,10.00){\makebox(0,0)[cc]{0}}
\put(144.00,4.00){\makebox(0,0)[cc]{0}}
\put(75.00,52.00){\line(1,0){3.00}}
\put(78.00,52.00){\line(1,-1){36.00}}
\put(114.00,16.00){\vector(1,0){6.00}}
\put(105.00,34.00){\vector(1,0){6.00}}
\put(105.00,34.00){\line(-1,-1){18.00}}
\put(87.00,16.00){\line(-1,0){12.00}}
\put(93.00,1.00){\makebox(0,0)[cc]{Caso I}}
\put(135.00,37.00){\line(1,0){3.00}}
\put(153.00,37.00){\line(-1,0){3.00}}
\end{picture} 

\pagebreak 

\mbox{}\\   

\unitlength=1.00mm
\special{em:linewidth 0.4pt}
\linethickness{0.4pt}
\begin{picture}(132.00,144.00)
\put(0.00,144.00){\line(1,0){18.00}}
\put(18.00,144.00){\line(0,-1){18.00}}
\put(18.00,126.00){\line(-1,0){6.00}}
\put(12.00,126.00){\line(0,1){6.00}}
\put(12.00,132.00){\line(-1,0){6.00}}
\put(6.00,132.00){\line(0,1){6.00}}
\put(6.00,138.00){\line(-1,0){6.00}}
\put(0.00,138.00){\line(0,1){6.00}}
\put(18.00,144.00){\line(1,0){6.00}}
\put(24.00,144.00){\line(0,-1){18.00}}
\put(24.00,126.00){\line(-1,0){6.00}}
\put(24.00,126.00){\line(1,0){12.00}}
\put(36.00,126.00){\line(0,-1){12.00}}
\put(36.00,114.00){\line(-1,0){6.00}}
\put(30.00,114.00){\line(0,1){6.00}}
\put(30.00,120.00){\line(-1,0){6.00}}
\put(24.00,120.00){\line(0,1){6.00}}
\put(36.00,114.00){\line(0,-1){6.00}}
\put(36.00,108.00){\line(1,0){6.00}}
\put(42.00,108.00){\line(0,-1){12.00}}
\put(42.00,96.00){\line(1,0){6.00}}
\put(48.00,96.00){\line(0,-1){6.00}}
\put(48.00,90.00){\line(1,0){6.00}}
\put(54.00,90.00){\line(0,-1){6.00}}
\put(54.00,84.00){\line(1,0){6.00}}
\put(60.00,84.00){\line(0,1){30.00}}
\put(60.00,114.00){\line(0,1){12.00}}
\put(60.00,126.00){\line(0,1){18.00}}
\put(60.00,144.00){\line(-1,0){24.00}}
\put(36.00,144.00){\line(0,-1){18.00}}
\put(18.00,114.00){\line(0,-1){30.00}}
\put(18.00,84.00){\line(1,0){6.00}}
\put(24.00,84.00){\line(0,1){30.00}}
\put(24.00,114.00){\line(-1,0){6.00}}
\put(18.00,108.00){\line(1,0){3.00}}
\put(36.00,126.00){\line(1,0){3.00}}
\put(36.00,114.00){\line(1,0){3.00}}
\put(60.00,114.00){\line(-1,0){3.00}}
\put(60.00,108.00){\line(-1,0){3.00}}
\put(60.00,126.00){\line(-1,0){3.00}}
\put(72.00,144.00){\line(1,0){24.00}}
\put(96.00,144.00){\line(0,-1){24.00}}
\put(96.00,120.00){\line(-1,0){6.00}}
\put(90.00,120.00){\line(0,1){6.00}}
\put(90.00,126.00){\line(-1,0){6.00}}
\put(84.00,126.00){\line(0,1){6.00}}
\put(84.00,132.00){\line(-1,0){6.00}}
\put(78.00,132.00){\line(0,1){6.00}}
\put(78.00,138.00){\line(-1,0){6.00}}
\put(72.00,138.00){\line(0,1){6.00}}
\put(96.00,120.00){\line(1,0){12.00}}
\put(108.00,120.00){\line(0,-1){12.00}}
\put(108.00,108.00){\line(-1,0){6.00}}
\put(102.00,108.00){\line(0,1){6.00}}
\put(102.00,114.00){\line(-1,0){6.00}}
\put(96.00,114.00){\line(0,1){6.00}}
\put(108.00,108.00){\line(0,-1){6.00}}
\put(108.00,102.00){\line(1,0){6.00}}
\put(114.00,102.00){\line(0,-1){6.00}}
\put(114.00,96.00){\line(1,0){6.00}}
\put(120.00,96.00){\line(0,-1){6.00}}
\put(120.00,90.00){\line(1,0){6.00}}
\put(126.00,90.00){\line(0,-1){6.00}}
\put(126.00,84.00){\line(1,0){6.00}}
\put(132.00,84.00){\line(0,1){60.00}}
\put(132.00,144.00){\line(-1,0){24.00}}
\put(108.00,144.00){\line(0,-1){24.00}}
\put(108.00,108.00){\line(1,0){3.00}}
\put(108.00,120.00){\line(1,0){3.00}}
\put(132.00,108.00){\line(-1,0){3.00}}
\put(132.00,120.00){\line(-1,0){3.00}}
\put(132.00,126.00){\line(-1,0){3.00}}
\put(90.00,126.00){\line(1,0){3.00}}
\put(63.00,111.00){\line(1,0){3.00}}
\put(66.00,111.00){\line(1,1){12.00}}
\put(78.00,123.00){\vector(1,0){6.00}}
\put(21.00,111.00){\makebox(0,0)[cc]{*}}
\put(21.00,105.00){\makebox(0,0)[cc]{0}}
\put(21.00,99.00){\makebox(0,0)[cc]{0}}
\put(21.00,93.00){\makebox(0,0)[cc]{0}}
\put(21.00,87.00){\makebox(0,0)[cc]{0}}
\put(39.00,105.00){\makebox(0,0)[cc]{x}}
\put(45.00,93.00){\makebox(0,0)[cc]{x}}
\put(51.00,87.00){\makebox(0,0)[cc]{x}}
\put(93.00,123.00){\makebox(0,0)[cc]{*}}
\put(117.00,93.00){\makebox(0,0)[cc]{0}}
\put(123.00,87.00){\makebox(0,0)[cc]{0}}
\put(42.00,54.00){\line(1,0){24.00}}
\put(66.00,54.00){\line(0,-1){24.00}}
\put(66.00,30.00){\line(-1,0){12.00}}
\put(54.00,30.00){\line(0,1){6.00}}
\put(54.00,36.00){\line(-1,0){6.00}}
\put(48.00,36.00){\line(0,1){12.00}}
\put(48.00,48.00){\line(-1,0){6.00}}
\put(42.00,48.00){\line(0,1){6.00}}
\put(66.00,30.00){\line(1,0){12.00}}
\put(78.00,30.00){\line(0,-1){12.00}}
\put(78.00,18.00){\line(-1,0){6.00}}
\put(72.00,18.00){\line(0,1){6.00}}
\put(72.00,24.00){\line(-1,0){6.00}}
\put(66.00,24.00){\line(0,1){6.00}}
\put(78.00,18.00){\line(0,-1){6.00}}
\put(78.00,12.00){\line(1,0){6.00}}
\put(84.00,12.00){\line(0,-1){6.00}}
\put(84.00,6.00){\line(1,0){6.00}}
\put(90.00,6.00){\line(0,-1){6.00}}
\put(90.00,0.00){\line(1,0){6.00}}
\put(96.00,0.00){\line(0,1){54.00}}
\put(96.00,54.00){\line(-1,0){18.00}}
\put(78.00,54.00){\line(0,-1){24.00}}
\put(42.00,54.00){\line(1,-1){54.00}}
\put(51.00,39.00){\makebox(0,0)[cc]{0}}
\put(57.00,33.00){\makebox(0,0)[cc]{0}}
\put(81.00,9.00){\makebox(0,0)[cc]{x}}
\put(87.00,3.00){\makebox(0,0)[cc]{x}}
\put(60.00,18.00){\line(1,0){6.00}}
\put(66.00,18.00){\line(0,-1){18.00}}
\put(66.00,0.00){\line(-1,0){6.00}}
\put(60.00,0.00){\line(0,1){18.00}}
\put(63.00,15.00){\makebox(0,0)[cc]{0}}
\put(63.00,9.00){\makebox(0,0)[cc]{0}}
\put(63.00,3.00){\makebox(0,0)[cc]{0}}
\put(81.00,87.00){\makebox(0,0)[cc]{Caso II}}
\put(27.00,3.00){\makebox(0,0)[cc]{Caso III}}
\end{picture} 

\pagebreak 

\mbox{}\\  

\unitlength=1.00mm
\special{em:linewidth 0.4pt}
\linethickness{0.4pt}
\begin{picture}(132.00,144.00)
\put(0.00,144.00){\line(1,0){24.00}}
\put(24.00,144.00){\line(0,-1){30.00}}
\put(24.00,114.00){\line(-1,0){6.00}}
\put(18.00,114.00){\line(0,1){6.00}}
\put(18.00,120.00){\line(-1,0){6.00}}
\put(12.00,120.00){\line(0,1){12.00}}
\put(12.00,132.00){\line(-1,0){6.00}}
\put(6.00,132.00){\line(0,1){6.00}}
\put(6.00,138.00){\line(-1,0){6.00}}
\put(0.00,138.00){\line(0,1){6.00}}
\put(24.00,114.00){\line(1,0){12.00}}
\put(36.00,114.00){\line(0,-1){12.00}}
\put(36.00,102.00){\line(-1,0){6.00}}
\put(30.00,102.00){\line(0,1){6.00}}
\put(30.00,108.00){\line(-1,0){6.00}}
\put(24.00,108.00){\line(0,1){6.00}}
\put(36.00,102.00){\line(0,-1){6.00}}
\put(36.00,96.00){\line(1,0){6.00}}
\put(42.00,96.00){\line(0,-1){6.00}}
\put(42.00,90.00){\line(1,0){6.00}}
\put(48.00,90.00){\line(0,-1){6.00}}
\put(48.00,84.00){\line(1,0){6.00}}
\put(54.00,84.00){\line(0,1){60.00}}
\put(54.00,144.00){\line(-1,0){18.00}}
\put(36.00,144.00){\line(0,-1){30.00}}
\put(54.00,144.00){\line(1,0){6.00}}
\put(60.00,144.00){\line(0,-1){60.00}}
\put(60.00,84.00){\line(-1,0){6.00}}
\put(18.00,120.00){\line(1,0){3.00}}
\put(60.00,120.00){\line(-1,0){3.00}}
\put(60.00,114.00){\line(-1,0){3.00}}
\put(60.00,102.00){\line(-1,0){3.00}}
\put(72.00,144.00){\line(1,0){24.00}}
\put(96.00,144.00){\line(0,-1){24.00}}
\put(96.00,120.00){\line(-1,0){6.00}}
\put(90.00,120.00){\line(0,1){6.00}}
\put(90.00,126.00){\line(-1,0){6.00}}
\put(84.00,126.00){\line(0,1){6.00}}
\put(84.00,132.00){\line(-1,0){6.00}}
\put(78.00,132.00){\line(0,1){6.00}}
\put(78.00,138.00){\line(-1,0){6.00}}
\put(72.00,138.00){\line(0,1){6.00}}
\put(96.00,120.00){\line(1,0){12.00}}
\put(108.00,120.00){\line(0,-1){12.00}}
\put(108.00,108.00){\line(-1,0){6.00}}
\put(102.00,108.00){\line(0,1){6.00}}
\put(102.00,114.00){\line(-1,0){6.00}}
\put(96.00,114.00){\line(0,1){6.00}}
\put(108.00,108.00){\line(0,-1){6.00}}
\put(108.00,102.00){\line(1,0){6.00}}
\put(114.00,102.00){\line(0,-1){6.00}}
\put(114.00,96.00){\line(1,0){6.00}}
\put(120.00,96.00){\line(0,-1){6.00}}
\put(120.00,90.00){\line(1,0){6.00}}
\put(126.00,90.00){\line(0,-1){6.00}}
\put(126.00,84.00){\line(1,0){6.00}}
\put(132.00,84.00){\line(0,1){60.00}}
\put(132.00,144.00){\line(-1,0){24.00}}
\put(108.00,144.00){\line(0,-1){24.00}}
\put(132.00,120.00){\line(-1,0){3.00}}
\put(132.00,108.00){\line(-1,0){3.00}}
\put(114.00,102.00){\line(1,0){3.00}}
\put(132.00,102.00){\line(-1,0){3.00}}
\put(108.00,120.00){\line(1,0){3.00}}
\put(108.00,108.00){\line(1,0){3.00}}
\put(36.00,120.00){\line(1,0){3.00}}
\put(36.00,114.00){\line(1,0){3.00}}
\put(36.00,102.00){\line(1,0){3.00}}
\put(63.00,117.00){\line(1,0){3.00}}
\put(66.00,117.00){\line(1,-1){12.00}}
\put(78.00,105.00){\vector(1,0){18.00}}
\put(15.00,123.00){\makebox(0,0)[cc]{0}}
\put(21.00,117.00){\makebox(0,0)[cc]{0}}
\put(111.00,99.00){\makebox(0,0)[cc]{0}}
\put(117.00,93.00){\makebox(0,0)[cc]{0}}
\put(123.00,87.00){\makebox(0,0)[cc]{0}}
\put(54.00,60.00){\line(0,-1){6.00}}
\put(54.00,54.00){\line(0,-1){18.00}}
\put(54.00,36.00){\line(0,-1){6.00}}
\put(54.00,30.00){\line(1,0){6.00}}
\put(60.00,30.00){\line(0,-1){6.00}}
\put(60.00,24.00){\line(1,0){6.00}}
\put(66.00,24.00){\line(0,-1){12.00}}
\put(66.00,12.00){\line(1,0){6.00}}
\put(72.00,12.00){\line(0,-1){6.00}}
\put(72.00,6.00){\line(1,0){6.00}}
\put(78.00,6.00){\line(0,-1){6.00}}
\put(78.00,0.00){\line(1,0){6.00}}
\put(84.00,0.00){\line(0,1){60.00}}
\put(84.00,0.00){\line(1,0){6.00}}
\put(90.00,0.00){\line(0,1){60.00}}
\put(54.00,36.00){\line(1,0){3.00}}
\put(54.00,54.00){\line(1,0){3.00}}
\put(90.00,54.00){\line(-1,0){3.00}}
\put(90.00,36.00){\line(-1,0){3.00}}
\put(69.00,15.00){\makebox(0,0)[cc]{0}}
\put(75.00,9.00){\makebox(0,0)[cc]{0}}
\put(81.00,3.00){\makebox(0,0)[cc]{0}}
\put(84.00,84.00){\makebox(0,0)[cc]{Caso IV}}
\put(33.00,3.00){\makebox(0,0)[cc]{Caso V}}
\end{picture}

 \clearpage
 \clearpage 
\chapter*{Matlab}  
\addcontentsline{toc}{chapter}{Matlab} 
\markboth{Matlab}{Matlab} 

\section*{Hist\'{o}rico} 

Matlab, de Matrix Laboratory, \'{e} um ambiente interativo para
computa\c{c}\~{a}o envolvendo matrizes.  Matlab foi desenvolvido no
inicio da d\'{e}cada de 80 por Cleve Moler, no Departamento de
Ci\^{e}ncia da Computa\c{c}\~{a}o da Universidade do Novo M\'{e}xico,
EUA. \index{Matlab} 
 
Matlab coloca \`{a} disposi\c{c}\~{a}o do usu\'{a}rio, num ambiente
interativo, as bibliotecas desenvolvidas nos projetos LINPACK e EISPACK. 
Estes projetos elaboraram bibliotecas de dom\'{\i}nio p\'{u}blico para
\'{A}lgebra Linear. LINPACK tem rotinas para solu\c{c}\~{a}o de
sistemas de equa\c{c}\~{o}es lineares, e EISPACK tem rotinas para
c\'{a}lculo de autovalores.  Os manuais destes projetos s\~{a}o portanto
documenta\c{c}\~{a}o complementar \`{a} documenta\c{c}\~{a}o do Matlab. 

Vers\~{o}es posteriores de Matlab, atualmente na vers\~{a}o 4.0, foram
desenvolvidas na firma comercial MathWorks Inc., que
det\^{e}m os direitos autorais destas implementa\c{c}\~{o}es.  
As vers\~{o}es recentes do produto Matlab melhoram significativamente 
o ambiente interativo, incluindo facilidades gr\'{a}ficas de
visualiza\c{c}\~{a}o e impress\~{a}o; todavia a ``Linguagem Matlab''
manteve-se quase inalterada.  Existem v\'{a}rios interpretadores da
linguagem Matlab em dom\'{\i}nio publico, como Matlab 1.0, Octave e
rlab.  Existem tamb\'{e}m outros interpretadores comerciais de Matlab,
como CLAM.  Existem ainda v\'{a}rias Tool Boxes, bibliotecas vendidas
pela MathWorks e por terceiros, com rotinas em Matlab para \'{a}reas
espec\'{\i}ficas. 

Usaremos a grafia de nome pr\'{o}prio, {\bf Matlab}, como refer\^{e}ncia
a linguagem, o nome em mai\'{u}sculas, {\bf MATLAB}, como refer\^{e}ncia
ao produto comercial da MathWorks, e a grafia em min\'{u}sculas, 
{\bf matlab}, como refer\^{e}ncia a um interpretador gen\'{e}rico 
da linguagem Matlab.

\section*{O Ambiente} 

Para entrar no ambiente matlab, simplesmente digite ``matlab''.  O
prompt do matlab \'{e} \verb#>># , que espera por comandos.  Para sair
use o comando \verb#quit#.  Dentro do ambiente matlab, um comando
precedido do bang, \verb#!#, \'{e} executado pelo sistema operacional,
assim: usando \verb#!dir# ou \verb#!ls# ficaremos sabendo os arquivos no
diret\'{o}rio corrente, e usando \verb#!edit# , \verb#!vi# ou
\verb#!emacs# , seremos capazes de editar um arquivo.  Normalmente
Matlab distingue mai\'{u}sculas de min\'{u}sculas. 

O comando \verb#help# exibe todos os comandos e s\'{\i}mbolos
sint\'{a}ticos dispon\'{\i}veis.  O comando \verb#help nomecom# fornece
informa\c{c}\~{o}es sobre o comando de nome {\em nomecom}.  O comando
\verb#diary nomearq# abre um arquivo de nome {\em nomearq}, e ecoa tudo
que aparece na tela para este arquivo.  Repetindo o comando \verb#diary#
fechamos este arquivo. O formato dos n\'{u}meros na tela pode ser 
alterado com o comando \verb#format#. 

Os comandos \verb#who# e \verb#whos# listam as vari\'{a}veis em
exist\^{e}ncia no espa\c{c}o de trabalho.  O comando \verb#clear# limpa o
espa\c{c}o de trabalho, extinguindo todas as vari\'{a}veis.  O comando
\verb#save nomearq# guarda o espa\c{c}o de trabalho no arquivo 
{\em nomearq}, e o comando \verb#load nomearq# carrega um espa\c{c}o de
trabalho previamente guardado com o comando \verb#save#. 

Em Matlab h\'{a} dois terminadores de comando: a v\'{\i}rgula, \verb#,#
, e o ponto-e-v\'{\i}rgula, \verb#;# .  O resultado de um comando
terminado por v\'{\i}rgula \'{e} ecoado para a tela.  O terminador
ponto-e-v\'{\i}rgula n\~{a}o causa eco.  Resultados ecoados na tela
s\~{a}o atribu\'{\i}dos a vari\'{a}vel do sistema \verb#ans# (de answer,
ou resposta).  O terminador v\'{\i}rgula pode ser suprimido no
\'{u}ltimo comando da linha.  Para continuar um comando na linha
seguinte devemos terminar a linha corrente com tr\^{e}s pontos,
\verb#...# , o s\'{\i}mbolo de continua\c{c}\~{a}o.  O sinal de
porcento, \verb#%# , indica que o resto da linha \'{e} coment\'{a}rio.

\section*{Matrizes} 

O tipo b\'{a}sico do Matlab \'{e} uma matriz de n\'{u}meros complexos. 
Uma matriz real \'{e} uma matriz que tem a parte imagin\'{a}ria de todos
os seus elementos nula.  Um vetor linha \'{e} uma matriz $1\times n$,
um vetor coluna uma matriz $n\times 1$, e um escalar uma matriz
$1\times 1$. 

As atribui\c{c}\~{o}es \\  
 \verb#A = [1, 2, 3; 4, 5, 6; 7, 8, 9];#  ou \\  
 \verb#A = [1 2 3; 4 5 6; 7 8 9];# , \\  
s\~{a}o equivalentes, e atribuem \`{a} vari\'{a}vel $A$ o valor
 $$ A = \left[ \begin{array}{ccc} 
    1 & 2 & 3 \\ 4 & 5 & 6 \\ 7 & 8 & 9  \end{array} \right] 
 $$

Matrizes s\~{a}o delimitadas por colchetes, elementos de uma mesma linha
s\~{a}o separados por v\'{\i}rgulas (ou apenas por espa\c{c}os em
branco), e linhas s\~{a}o separadas por ponto-e-v\'{\i}rgula (ou pelo
caracter nova-linha).  O ap\'{o}strofe, \verb#'# , transp\~{o}em uma
matriz.  \'{E} f\'{a}cil em Matlab compor matrizes blocadas, desde que
os blocos tenham dimens\~{o}es compat\'{\i}veis! Por exemplo:
 
 $$
 A = \left[ \begin{array}{cc} 1 & 2 \\ 3 & 4 \end{array} \right] 
 \ , \ \ 
 B = \left[ \begin{array}{c} 5 \\ 6 \end{array} \right] 
 \ , \ \ 
 C = \left[ \begin{array}{c} 7 \\ 8 \\ 9  \end{array} \right] 
 \ , \ \ 
 \verb#D=[A,B;C'];# 
 \ \ D = \left[ \begin{array}{ccc} 1 & 2 & 5 \\ 3 & 4 & 6 \\ 
                 7 & 8 & 9 \end{array} \right]
 $$ 
Cuidado ao concatenar matrizes com os espacos em branco, pois estes 
s\~{a}o equivalentes a v\'{\i}rgulas, separando elementos. Assim: 
\verb#[1,2+3]==[1 5]# mas \verb#[1,2 +3]==[1,2,3]#. 
 \index{Matriz!blocada} 

H\'{a} v\'{a}rias fun\c{c}\~{o}es para gerar matrizes e vetores
especiais: \verb#zeros(m,n)#, \verb#ones(m,n)# e \verb#rand(m,n)#
s\~{a}o matrizes de dimens\~{a}o $m\times n$ com, respectivamente,
zeros, uns, e n\'{u}meros aleat\'{o}rios em [0,1].  O vetor \verb#i:k:j#
\'{e} o vetor linha $\left[i,i+k,i+2k,\ldots i+nk\right]$, onde 
 $n=\max m \mid i+mk\leq j$.  Podemos suprimir o ``passo'' $k=1$,
escrevendo o vetor \verb#i:1:j# simplesmente como \verb#i:j#.  Se $v$
\'{e} um vetor, \verb#diag(v)# \'{e} a matriz diagonal com diagonal $v$,
se $A$ \'{e} uma matriz quadrada, \verb#diag(A)# \'{e} o vetor com os
elementos da diagonal principal de $A$.  A matriz identidade de ordem
$n$ \'{e} \verb#eye(n)#, o mesmo que \verb#diag(ones(1,n))#. 
  \index{Matriz!\'{\i}ndices}

\verb#A(i,j)# \'{e} o elemento na $i$-\'{e}sima linha, $j$-\'{e}sima
coluna de $A$, $m\times n$.  Se $vi$ e $vj$ s\~{a}o vetores de
\'{\i}ndices em $A$, i.e.  vetores linha com elementos inteiros
positivos, em $vi$ n\~{a}o excedendo $m$, e em $vj$ n\~{a}o excedendo
$n$, \verb#A(vi,vj)# \'{e} a sub-matriz do elementos de $A$ com
\'{\i}ndice de linha em $vi$ e \'{\i}ndice de coluna em $vj$.  Em
particular \verb#A(1:m,j)#, ou \verb#A(:,j)# \'{e} a $j$-\'{e}sima coluna
de de $A$, e \verb#A(i,1:j)#, ou \verb#A(i,:)#, \'{e} a $i$-\'{e}sima
linha de $A$. Exemplo: 
 $$ 
 A = \left[ \begin{array}{ccc} 
     11 & 12 & 13 \\ 21 & 22 & 23 \\ 31 & 32 & 33  
     \end{array} \right] \ \ 
 vi = \left[ \begin{array}{c} 3 \\ 1 \end{array} \right] \ \  
 vj = \left[ \begin{array}{c} 2 \\ 3 \end{array} \right] \ \  
 A(vi,vj) = 
 \left[ \begin{array}{cc} 32 & 33 \\ 12 & 13 \end{array} \right]
 $$ 

As opera\c{c}\~{o}es de adi\c{c}\~{a}o, subtra\c{c}\~{a}o, produto,
divis\~{a}o e pot\^{e}ncia, \verb#+ - * / ^#, s\~{a}o as usuais da
\'{a}lgebra linear.  
 O operador de divis\~{a}o \`{a} esquerda, \verb#\#, fornece em 
 \verb#x = A\b;# uma solu\c{c}\~{a}o $x \mid A*x=b$.  
 O operador de divis\~{a}o a direita, \verb#/# , fornece em 
 \verb#a = b/A;# uma solu\c{c}\~{a}o $x \mid x*A=b$. Quando o sistema
\'{e} bem determinado isto \'{e} o mesmo que, respectivamente, 
 \verb#inv(A)*b# ou \verb#b*inv(A)#.  
 (Quando o sistema \'{e} super-determinado x \'{e} uma solu\c{c}\~{a}o no
sentido dos m\'{\i}nimos quadrados.)  
 Os operadores aritm\'{e}ticos de produto, pot\^{e}ncia e divis\~{a}o
tem a vers\~{a}o pontual, \verb#.* .^ ./ .\#  , que s\~{a}o executados
elemento a elemento.  Exemplo: 
 $$ x=\left[ \begin{array}{ccc} 5 & 5 & 5 \end{array} \right] \ , \ \ 
    y=\left[ \begin{array}{ccc} -1 & 0 & 1 \end{array} \right] \ , \ \ 
 x\, .*\, y =\left[ \begin{array}{ccc} -5 & 0 & 5 \end{array} \right]  
 $$ 

Matlab \'{e} uma linguagem declarativa (em oposi\c{c}\~{a}o a
imperativa) onde as vari\'{a}veis s\~{a}o declaradas, dimensionadas e
redimensionadas dinamicamente pelo interpretador.  Assim, se $A$
presentemente n\~{a}o existe, \verb#A=11;# declara e inicializa uma
matriz real $1\times 1$.  Em seguida, o comando \verb#A(2,3)=23;# 
redimensionaria $A$ como a matriz $2\times 3$ \verb#[11, 0, 0; 0, 0, 23]#.  
A nome da matriz nula \'{e} \verb#[ ]#.  A atribui\c{c}\~{a}o 
 \verb#A(:,2)=[];# anularia a segunda coluna de $A$, tornado-a a matriz 
 $2\times 2$ \verb#[11, 0; 0, 23]#.

\section*{Controle de Fluxo}

Os operadores relacionais da linguagem s\~{a}o \verb#< <= > >= == ~=#, 
que retornam valor $0$ ou $1$ conforme a condi\c{c}\~{a}o seja 
verdadeira ou falsa. Os operadores l\'{o}gicos, {\em n\~{a}o e ou}, 
s\~{a}o, respectivamente, \verb#~ & | #.  

Matlab possui os comandos de fluxo {\em for -- end}, {\em while -- end}
e {\em if -- elseif -- else -- end}, que tem a mesma sintaxe do Fortran,
exemplificada a seguir.  Lembre-se de {\bf n\~{a}o} escrever a palavra
\verb#elseif# como duas palavras separadas. 

\begin{verbatim}                               
for i=1:n                 if(x<0)          fatn=1;  
  for j=1:m                 signx=-1;      while(n>1) 
    H(i,j)=1/(i+j-1);     elseif(x>0)        fatn=fatn*n; 
  end                       signx=1;         n=n-1; 
end                       else             end 
                            signx=0; 
                          end                        
\end{verbatim} 

Uma considera\c{c}\~{a}o sobre efici\^{e}ncia: Matlab \'{e} uma
linguagem interpretada, e tempo de interpreta\c{c}\~{a}o de um comando
simples pode ser bem maior que seu tempo de execu\c{c}\~{a}o.  Para
tornar o programa r\'{a}pido, tente operar sobre matrizes ou blocos,
evitando loops expl\'{\i}citos que operem numa matriz elemento por
elemento.  Em outras palavras, tente evitar loops que repetem um comando
que atribui valores a elementos, por atribui\c{c}\~{o}es a vetores ou
matrizes.  As facilidades no uso de vetores de \'{\i}ndices, os
operadores aritm\'{e}ticos pontuais, e fun\c{c}\~{o}es como \verb#max# ,
\verb#min# , \verb#sort# , etc., tornam esta tarefa f\'{a}cil, e os
programas bem mais curtos e leg\'{\i}veis. 

\section*{Scripts e Fun\c{c}\~{o}es}

O fluxo do programa tamb\'{e}m \'{e} desviado pela invoca\c{c}\~{a}o de
subrotinas.  A subrotina de nome {\em nomsubr} deve estar guardada no
arquivo \verb#nomsubr.m#  ; por isto subrotinas s\~{a}o tamb\'{e}m
denominadas M-arquivos (M-files).  H\'{a} dois tipos de subrotinas:
Scripts e Fun\c{c}\~{o}es. 

Um script \'{e} simplesmente uma seq\"{u}\^{e}ncia de comandos, que
ser\~{a}o executados como se fossem digitados ao prompt do matlab. 
Subrotinas podem invocar outras subrotinas, inclusive recursivamente. 

Um M-arquivo de fun\c{c}\~{a}o em Matlab come\c{c}a com a
declara\c{c}\~{a}o da forma \\ 
 \verb#[ps1, ps2,... psm] = nomefunc( pe1, pe2,... pen )#  \\  
 A lista entre par\^{e}nteses \'{e} a lista de par\^{a}metros de entrada
da fun\c{c}\~{a}o, e a lista entre colchetes s\~{a}o os par\^{a}metros
de sa\'{\i}da.  Par\^{a}metros s\~{a}o vari\'{a}veis locais, assim como
s\~{a}o vari\'{a}veis locais todas as vari\'{a}veis no corpo da
fun\c{c}\~{a}o. 

Ao invocarmos a fun\c{c}\~{a}o {\em nomefunc} com o comando \\ 
  \verb#[as1,... asm] = nomefunc(ae1,... aen);# \\  
 Os argumentos de entrada, $ae1,\ldots aen$, s\~{a}o passados por valor
aos (i.e.  copiados nos) par\^{a}metros de entrada, e ao fim do
M-arquivo, ou ao encontrar o comando \verb#return# , os par\^{a}metros
de sa\'{\i}da s\~{a}o passados aos argumentos de sa\'{\i}da.  Exemplos:

\begin{verbatim} 
 function [mabel, mabin] = vmax(v) 
 % procura o maior elemento, em valor absoluto, 
 % dentro de um vetor, linha ou coluna.  
 %Obs: Esta funcao NAO segue os conselhos 
 %para operar sobre blocos, e nao elementos! 
                                             
 [m,n]=size(v);                             
 if( m ~= 1  &  n ~= 1 ) 
   erro; 
 else 
   K=max([m,n]); 
   mabel= abs(v(1)); mabin= 1;  
   for k=2:K 
     if( abs(v(k)) > mabel ) 
       mabel= abs(v(k)); mabin= k; 
     end%if 
   end%for  
 end%else 
\end{verbatim} 

Para referir-mo-nos, dentro de uma fun\c{c}\~{a}o, a uma vari\'{a}vel
externa, esta deve ter sido declarada uma vari\'{a}vel global com o
comando \verb#global#.  A forma de especificar uma vari\'{a}vel como global
muda um pouco de interpretador para interpretador, e mesmo de vers\~{a}o
para vers\~{a}o: Diga \verb#help global# para saber os detalhes de como
funciona a sua implementa\c{c}\~{a}o.

\section*{Bibliografia}

[Colema-88] T.F.Coleman and C.F.van Loan.  {\em A Matrix Computation Handbook}. 
SIAM Publications, Philadelphia. 

\noindent 
 [Moler-81] C.B.Moler.  {\em Matlab Manual}.  Department of Computer
Science, University of New Mexico. 

\noindent 
 [Dongar-79] J.J.Dongarra, J.R.Bunch, C.B.Moler, G.W.Stewart.  {\em
LINPACK Users' Guide}.  Society for Industrial and Applied Mathematics,
Philadelphia. 

\noindent 
 [Smith-76] B.T.Smith, J.M.Boyle, J.J.Dongarra, B.S.Garbow, Y.Ikebe,
V.C.Klema, C.B.Moler.  {\em Matrix Eigensystem Routines: EISPACK Guide}. 
Lecture Notes in Computer Science, volume 6, second edition,
Springer-Verlag. 
 
\noindent 
 [Garbow-77] B.S.Garbow, J.M.Boyle, J.J.Dongarra, C.B.Moler.  {\em Matrix
Eigensystem Routines: EISPACK Guide Extension}.  Lecture Notes in
Computer Science, volume 51, Springer-Verlag. 
           
\noindent 
 [Moler-92] C.B.Moler, J.N.Litte and S.Bangert.  {\em PC-MATLAB User's
Guide}.  The MathWorks Inc.  Sherborn, Massachusetts. 

\noindent 
 [Eaton-92] J.W.Eaton.  {\em Octave Manual}.  Chemical Engineering
Department, University of Texas at Austin.  Austin, Texas.

 \clearpage
 \clearpage 
\chapter*{Bibliografia}
\addcontentsline{toc}{chapter}{Bibliografia} 
\markboth{Bibliografia}{Bibliografia} 

Bibliografia Sum\'{a}ria de suporte para o curso: 

\noindent 
[Bertsekas-89] D.P.Bertsekas J.N.Tsitsiklis. 
{\em Parallel and Distributed Computation, Numerical Methods}. 
Prentice Hall, 1989.

\noindent 
[Bunch-76] J.R.Bunch D.J.Rose. 
{\em Sparse Matrix Computations}.
Academic Press, 1976. 

\noindent 
[Carey-89] G.F.Carey. 
{\em Parallel Supercomputing}. 
John Wiley, 1989. 

\noindent 
[Dongarra-91] J.J.Dongarra et all. 
{\em Solving Linear Systems on Vector and Shared Memory Computers}. 
SIAM, 1991. 

\noindent 
[Duff-79] I.S.Duff G.W.Stewart. 
{\em Sparse Matrix Proceedings 1978}. 
SIAM, 1979.

\noindent 
[Duff-86] I.S.Duff A.M.Erisman J.K.Reid. 
{\em Direct Methods for Sparse Matrices}.
Oxford, 1986 

\noindent 
[Gallivan-90] K.A.Gallivan et all.  
{\em Parallel Algorithms for Matrix Computations}. 
SIAM, 1990. 

\noindent 
[George-81] A.George J.W.H.Liu. 
{\em Computer Solution of Large Sparse Positive-Definite Systems}. 
Prentice Hall, 1981. 

\noindent 
[Golub-83] G.H.Golub C.F.van Loan. 
{\em Matrix Computations}. 
John Hopkins, 1983. 
 
\noindent 
[Pissan-84] S.Pissanetzky. 
{\em Sparse Matrix Technology}. 
Academic Press 1984. 

\noindent 
[Rose-72] D.J.Rose R.A.Willoughby. 
{\em Sparse Matrices}. 
Plenum Press, 1972. 

\noindent 
[Stern-92] J.M.Stern. 
{\em Simulated Annealing with a Temperature Dependent 
Penalty Function}. ORSA Journal on Computing, V-4, N-3, p-311-319, 1992.

\noindent 
[Stern-93] J.M.Stern S.A.Vavasis. 
{\em Nested Dissection for Sparse Nullspace Bases}. 
SIAM Journal on Matrix Analysis and Applications, 
V-14, N-3, p-766-775, 1993. 

\noindent 
[Stern-94] J.M.Stern S.A.Vavasis. 
{\em Active Set Algorithms for Problems in Block Angular Form}. 
Computational and Applied Mathematics, V-13, 1994.   

\noindent
[Stewart-73]  G.W.Stewart. 
{\em Introduction to Matrix Computations}. 
Academic Press, 1973. 

\noindent 
[Tewarson-73] R.P.Tewarson. 
{\em Sparse Matrices}. 
Academic Press, 1973. 

\noindent 
[Vorst-89] H.A.van der Vorst P.van Dooren. 
{\em Parallel Algorithms for Numerical Linear Algebra}. 
North Holland, 1989.

 \clearpage
 \clearpage \printindex \clearpage   
%

\end{document}